\documentclass[twoside,openright]{report}[11pt]
\pdfoutput=1
\usepackage[utf8]{inputenc}
\usepackage{environ}
\usepackage{etoolbox}
\title{Thesis}
\author{Callum Reader}
\date{January 2022}

\usepackage{xargs}
\usepackage[margin=1.5in]{geometry}
\usepackage{mathrsfs}   
\usepackage[dvipsnames]{xcolor}
\usepackage{tikz-cd}    
\usepackage{tikz}       
\usepackage{xspace}     
\usepackage{graphicx}
\usepackage{booktabs}
\usepackage{ebproof}
\usepackage{bbm}
\usepackage{leftidx}
\usepackage{listofitems}
\usepackage{caption}
\usepackage{subfig}
\usepackage[outline]{contour}
\usepackage{float}
\usepackage{graphbox}
\newcommand{\altchi}{\mathbb{X}}

\newcommand{\myinput}[1]{\fbox{\includegraphics[align=c]{#1}}}

\graphicspath{{Figures/}}   
\usepackage{comment}        
\usepackage{xstring}

\usepackage{afterpage}

\newcommand\blankpage{%
    \null
    \thispagestyle{empty}%
    \addtocounter{page}{-1}%
    \newpage}

\DeclareFontFamily{U}{mathb}{}
\DeclareFontShape{U}{mathb}{m}{n}{
  <5> <6> <7> <8> <9> <10>
  <10.95> <12> <14.4> <17.28> <20.74> <24.88>
  mathb10
  }{}
\DeclareSymbolFont{mathb}{U}{mathb}{m}{n}
\DeclareMathSymbol{\mminus}{\mathbin}{mathb}{"01}   

\usepackage[all,knot,poly]{xy}
\usepackage{standalone}
\usepackage{amsmath,amssymb,amsthm,amsfonts}
\usepackage{enumitem}
\usepackage{multirow}
\usepackage{pdfpages}
\usepackage{mathtools}
\usepackage[hidelinks]{hyperref}
\usepackage{nicefrac}
\usepackage{accents}
\usepackage{todonotes}
\usepackage[noabbrev,capitalize]{cleveref}
\usepackage{caption}
\usepackage{mathtools}
\usepackage{thmtools}
\usepackage{thm-restate}
\usepackage{tabularx}
\usepackage{scalerel}
\usepackage{geometry}
\usepackage{pgffor}
\usepackage{circuitikz}
\usepackage[nomessages]{fp}
\usepackage{xstring} 
\usepackage{xargs}
\usepackage{ifthen}
\usepackage{cmll}
\usepackage{quiver}

\RequirePackage{snapshot}

\makeatletter
\DeclareFontFamily{U}{MnSymbolA}{}
\DeclareFontShape{U}{MnSymbolA}{m}{n}{
    <-6>  MnSymbolA5
   <6-7>  MnSymbolA6
   <7-8>  MnSymbolA7
   <8-9>  MnSymbolA8
   <9-10> MnSymbolA9
  <10-12> MnSymbolA10
  <12->   MnSymbolA12}{}
\DeclareFontShape{U}{MnSymbolA}{b}{n}{
    <-6>  MnSymbolA-Bold5
   <6-7>  MnSymbolA-Bold6
   <7-8>  MnSymbolA-Bold7
   <8-9>  MnSymbolA-Bold8
   <9-10> MnSymbolA-Bold9
  <10-12> MnSymbolA-Bold10
  <12->   MnSymbolA-Bold12}{}
\DeclareSymbolFont{MnSyA}{U}{MnSymbolA}{m}{n}
\SetSymbolFont{MnSyA}{bold}{U}{MnSymbolA}{b}{n}

\DeclareRobustCommand{\overleftharpoon}{\mathpalette{\overarrow@\leftharpoonfill@}}
\DeclareRobustCommand{\overrightharpoon}{\mathpalette{\overarrow@\rightharpoonfill@}}
\def\leftharpoonfill@{\arrowfill@\leftharpoondown\mn@relbar\mn@relbar}
\def\rightharpoonfill@{\arrowfill@\mn@relbar\mn@relbar\rightharpoonup}

\DeclareMathSymbol{\leftharpoondown}{\mathrel}{MnSyA}{'112}
\DeclareMathSymbol{\rightharpoonup}{\mathrel}{MnSyA}{'100}
\DeclareMathSymbol{\mn@relbar}{\mathrel}{MnSyA}{'320}
\makeatother

\usetikzlibrary{patterns}
\usetikzlibrary{through}
\usetikzlibrary{knots}
\usetikzlibrary{decorations.markings}
\usetikzlibrary{fit}
\usetikzlibrary{external}
\usetikzlibrary{shadows,shadings,shapes.symbols}
\usetikzlibrary{math}

\usetikzlibrary{decorations.pathmorphing}

\usepackage[american, british]{babel}
\usepackage[english=british]{csquotes}

\usepackage[style=ieee, backref=true, backend=biber, hyperref=true, giveninits=true, sorting=nyt, language=british]{biblatex}
\DeclareSourcemap{
  \maps[datatype=bibtex]{
    \map{
      \step[fieldset=issn, null]
    }
  }
}
\DefineBibliographyStrings{english}{%
  backrefpage = {cited on p.~},
  backrefpages = {cited on pp.~},
}
\DefineBibliographyExtras{british}{%
  \stdpunctuation
}
\usepackage{xpatch}
    \xpatchfieldformat{eprint:arxiv}{\space}{}{}{}
    \xpatchfieldformat{doi}{\space}{}{}{}
   
\renewbibmacro{in:}{}
\DeclareFieldFormat{url}{Available\addcolon\space\url{#1}}

\tikzset{%
    symbol/.style={%
        draw=none,
        every to/.append style={%
            edge node={node [sloped, allow upside down, auto=false]{$#1$}}}
    }
}

\newcommand{\textdef}[1]{\textbf{#1}}
\newcommand{\textbook}[1]{\textit{#1}}
\newcommand{\morphindx}[4]{(#1\colon #2\rightarrow #3)\in #4}

\newcommand{\outin}[2]{
out=\IfEqCase{#1}{
{ne}{-45}
{se}{45}
{sw}{135}
{nw}{-135}
}[0],
in=\IfEqCase{#2}{
{ne}{-45}
{se}{45}
{sw}{135}
{nw}{-135}
}
}

\newcommand{\compassToDegrees}[1]{
\IfEqCase{#1}{
{ne}{-45}
{se}{45}
{sw}{135}
{nw}{-135}
}[#1]
}

\newcommand{\inedge}[1]{
to [in=#1]
}
\newcommand{\inedgetest}[1]{\expandafter\inedge\expandafter{\compassToDegrees{#1}}}

\theoremstyle{plain}
\newtheorem{theorem}{Theorem}[section]
\newtheorem*{theorem*}{Theorem}

\newtheorem{lemma}[theorem]{Lemma}
\newtheorem*{lemma*}{Lemma}

\newtheorem{proposition}[theorem]{Proposition}

\newtheorem{corollary}[theorem]{Corollary}

\renewcommand{\autoref}{\cref}          

\crefname{equation}{Equation}{Equations}
\crefname{theorem}{Theorem}{Theorems}
\crefname{lemma}{Lemma}{Lemmata}
\crefname{example}{Example}{Examples}

\theoremstyle{definition}
\newtheorem{definition}[theorem]{Definition}
\newtheorem{example}[theorem]{Example}

\newtheorem*{remark}{Remark}
\newtheorem*{notation}{Notation}

\newenvironment{lst}{\begin{enumerate}                  
[label=\emph{\roman*.}]}{\end{enumerate}}

\newcommand{\mhyphen}{\text{-}}

\newcommand{\cat}[1]{#1\mhyphen\mathrm{Cat}}         
              
\newcommand{\rep}[1]{#1\mhyphen\mathrm{Rep}}
\newcommand{\clrep}[1]{#1\mhyphen\mathrm{Rep}_{\mathrm{cl}}}         
\newcommand{\Bicat}{\mathrm{Bicat}}
\renewcommand{\phi}{\varphi}                	    

\newcommand{\blank}{\phantom{ }}                    

\newcommand{\tensor}{\otimes}                       

\newcommand{\close}[1]{\langle #1\rangle}

\newcommand{\N}{\mathbb{N}}

\renewcommand{\hom}[1]{\left[#1\right]}

\DeclareSymbolFont{bbold}{U}{bbold}{m}{n}
\DeclareSymbolFontAlphabet{\mathbbold}{bbold}

\DeclareMathSymbol{\mminus}{\mathbin}{mathb}{"01}

\newcommand{\id}{\mathrm{id}}
\newcommand{\Id}{\mathrm{Id}}
\DeclareMathOperator{\Set}{Set}
\DeclareMathOperator{\obj}{obj}

\newcommand{\Hom}{\mathrm{{Hom}}}
\newcommand{\RHom}{\mathrm{{RHom}}}

\DeclareMathOperator{\colim}{colim}

\newcommand{\ev}{\mathrm{ev}}

\DeclareMathOperator{\wlim}{lim}

\DeclareMathOperator{\Cat}{Cat}
\DeclareMathOperator{\Nat}{Nat}

\DeclareFontFamily{U}{min}{}
\DeclareFontShape{U}{min}{m}{n}{<-> udmj30}{}
\DeclareMathOperator{\Prof}{Prof}

\newcommand{\op}{{\mathrm{op}}}

\newsavebox{\pullbackdrawing}
\sbox\pullbackdrawing{%
\begin{tikzpicture}%
\draw (0,0) -- (1ex,0ex);%
\draw (1ex,0ex) -- (1ex,1ex);%
\end{tikzpicture}}

\newsavebox{\pushforwarddrawing}
\sbox\pushforwarddrawing{%
\begin{tikzpicture}[rotate=180, transform shape]%
\draw (0,0) -- (1ex,0ex);%
\draw (1ex,0ex) -- (1ex,1ex);%
\end{tikzpicture}}

\newcommand*\cocolon{%
        \nobreak
        \mskip6mu plus1mu
        \mathpunct{}%
        \nonscript
        \mkern-\thinmuskip
        {:}%
        \mskip2mu
        \relax
}

\renewcommand*\cocolon{%
        \nobreak
        \mskip6mu plus1mu
        \mathpunct{}%
        \nonscript
        \mkern-\thinmuskip
        {:}%
        \mskip2mu
        \relax
}

\makeatletter
\def\slashedarrowfill@#1#2#3#4#5{%
  $\m@th\thickmuskip0mu\medmuskip\thickmuskip\thinmuskip\thickmuskip
  \relax#5#1\mkern-7mu%
  \cleaders\hbox{$#5\mkern-2mu#2\mkern-2mu$}\hfill
  \mathclap{#3}\mathclap{#2}%
  \cleaders\hbox{$#5\mkern-2mu#2\mkern-2mu$}\hfill
  \mkern-7mu#4$%
}
\def\rightslashedarrowfill@{%
  \slashedarrowfill@\relbar\relbar\mapstochar\rightarrow}
\newcommand\xprofto[2][]{%
  \ext@arrow 0055{\rightslashedarrowfill@}{#1}{#2}}
\makeatother

\def\profto{\xprofto{}}
\tikzset{commutative diagrams/.cd, 
    prof/.style = "\shortmid"{marking, sloped},
    }

\DeclareMathSymbol\simplex  \mathord{bbold}{"01}

\DeclareFontFamily{U}{mathx}{\hyphenchar\font45}
\DeclareFontShape{U}{mathx}{m}{n}{<->mathx10}{}
\DeclareFontSubstitution{U}{mathx}{m}{n}
\DeclareSymbolFont{mathx}{U}{mathx}{m}{n}

\DeclareMathSymbol{\endint}{\mathop}{mathx}{"B3}

\DeclareFontFamily{U}{mathx}{\hyphenchar\font45}
\DeclareFontShape{U}{mathx}{m}{n}{
      <5> <6> <7> <8> <9> <10>
      <10.95> <12> <14.4> <17.28> <20.74> <24.88>
      mathx10
      }{}
\DeclareSymbolFont{mathx}{U}{mathx}{m}{n}
\DeclareFontSubstitution{U}{mathx}{m}{n}
\DeclareMathAccent{\widecheck}{0}{mathx}{"71}

\newcommand{\extend}{\mathrel{\multimapinv}}
\newcommand{\lift}{\mathrel{\multimap}}
\tikzset{shorten <>/.style={shorten >=#1,shorten <=#1}}

\DeclareMathOperator{\Span}{Span}

\newcommand{\ecirc}{\diamond}
\newcommand{\elift}{
\mathbin{
    \vcenter{
      \hbox{${-}$}
    }
      \kern-1.25pt
    {
      \vcenter{\hbox{${\scriptstyle\diamond}$}}
    }
  }
}
\newcommand{\comp}[1]{\mathring{#1}}
\newcommand{\ecomp}[1]{\protect\accentset{\ecirc}{#1}}

\newcommand{\Dim}{\mathrm{Dim}}
\newcommand{\Tr}{\mathrm{Tr}}
\newcommand{\coev}{\mathrm{coev}}

\newcommand{\Bim}{\mathrm{Bim}}
\newcommand{\DGBim}{\mathrm{DBim}}
\newcommand{\Rel}{\mathrm{Rel}}
\newcommand{\Rep}{\mathrm{Rep}}
\newcommand{\Repp}[1]{#1\mhyphen\mathrm{Rep}}
\newcommand{\Ext}{\mathrm{Ext}}
\newcommand{\Tor}{\mathrm{Tor}}

\newcommand{\MonCat}{\mathrm{MonCat}}

\newcommand{\dTr}{\mathrm{\overline{\Tr}}^{\circ}}
\newcommand{\rTr}{\mathrm{\Tr}^\circ}

\newcommand{\dDim}{\mathrm{\overline{\Dim}^{\circ}}}
\newcommand{\rDim}{\mathrm{{\Dim}^{\circ}}}

\DeclareMathSymbol{\boxbackslash}{2}{mathb}{"6E}

\newcommand{\invbullet}{\ominus}

    \usepackage[a]{esvect}
    \usepackage{array}

    \usetikzlibrary{positioning,decorations.markings,arrows.meta,calc,fit,quotes,cd,math,arrows,backgrounds,shapes.geometric}
    
    \makeatletter
\DeclareRobustCommand{\rvdots}{%
  \vbox{
    \baselineskip4\p@\lineskiplimit\z@
    \kern-\p@
    \hbox{.}\hbox{.}\hbox{.}
  }}
\makeatother
\makeatletter
\def\instring#1#2{TT\fi\begingroup
  \edef\x{\endgroup\noexpand\in@{#1}{#2}}\x\ifin@}
    \let\nvec\Vec
    \def\Vec#1{\nvec{\vphantom b\smash{#1}}}
    \renewcommand{\Vec}[1]{\overrightarrow{\vphantom a\smash{#1}}}
    \renewcommand{\Vec}[1]{
    \if\instring{#1}{ABCDEFGHIJKLMNOPQRSTUVWXYZ}
    \vv{\vphantom{A}\smash{#1}}
    \else
    \vv{\vphantom{\tiny{a}}\smash{#1}}
    \fi}

    \newcommand{\str}{\mathrm{str}}

    \newcommand{\EnCat}{\mathrm{EnCat}}

    \renewcommand{\AA}{\mathscr{A}}
    \newcommand{\BB}{\mathscr{B}}
    \newcommand{\BBB}{\underline{\BB}}
    \newcommand{\CC}{\mathscr{C}}
    \newcommand{\DD}{\mathscr{D}}
    
    \newcommand{\MM}{\mathscr{M}}

    \newcommand{\UU}{\mathscr{U}}
    
    \newcommand{\VV}{\mathscr{V}}
    \newcommand{\WW}{\mathscr{W}}
    
    \newcommand{\act}{\odot}
    
    \newcommand{\Act}{\operatorname{Spr}^\circ}
    \newcommand{\Spr}{\Act}
    \newcommand{\rSpr}{\overline{\mathrm{Spr}}^\circ}
    \newcommand{\dSpr}{\mathrm{Spr}^\circ}

    \newcommand{\tr}{\rTr}
    \newcommand{\ctr}{\dTr}
    \newcommand{\cspr}{\rSpr}
    \newcommand{\spr}{\dSpr}
    \newcommand{\Path}{\operatorname{Path}}

  \newcommand{\name}[1]{\overrightharpoon{\vphantom{h}#1}}
  \newcommand{\unname}[1]{\overleftharpoon{\vphantom{h}#1}}
  
  \makeatletter
  \newcommand{\superimpose}[2]{%
  {\ooalign{$#1\@firstoftwo#2$\cr\hfil$#1\@secondoftwo#2$\hfil\cr}}}
  \makeatother

    \newcommand{\forgetdagger}[1]{\ifthenelse{\equal{#1}{(-)}}{^{\dagger}}{#1^{\dagger}}}
    
    \newcommand{\freedagger}[1]{#1^{\dagger}}
    \let\f\freedagger
    \renewcommand{\freedagger}[1]{\ifthenelse{\equal{#1}{}}{\f{(-)}}{\f{#1}}}
    
    \newcommand{\freemulti}[1]{#1^{\triangleright}}
    \let\f\freemulti
    \renewcommand{\freemulti}[1]{\ifthenelse{\equal{#1}{}}{\f{(-)}}{\f{#1}}}

    \newcommand{\forgetmulti}[1]{#1_{\triangleright}}
    \let\f\forgetmulti
    \renewcommand{\forgetmulti}[1]{\ifthenelse{\equal{#1}{}}{\f{(-)}}{\f{#1}}}

    \let\f\freemonoidal
    \renewcommand{\forgetmulti}[1]{\ifthenelse{\equal{#1}{}}{\f{(-)}}{\f{#1}}}

    \let\f\forgetmonoidal
    \renewcommand{\forgetmulti}[1]{\ifthenelse{\equal{#1}{}}{\f{(-)}}{\f{#1}}}
    
    \renewcommand{\forgetdagger}[1]{L_{\dagger}#1}
    \renewcommand{\freedagger}[1]{R_{\dagger}#1}
    \renewcommand{\freemulti}[1]{L_{\triangleright}#1}
    \renewcommand{\forgetmulti}[1]{R_{\triangleright}#1}

    \newcommand{\FDHilb}{\mathrm{FDHilb}}

\definecolor{cb-orange}{cmyk}{0,0.5,1,0}
\definecolor{cb-light-blue}{cmyk}{0.8,0,0,0}
\definecolor{cb-green}{cmyk}{0.97,0,0.75,0}
\definecolor{cb-yellow}{cmyk}{0.1,0.05,0.9,0}
\definecolor{cb-blue}{cmyk}{1,0.5,0,0}
\definecolor{cb-red}{cmyk}{0,0.8,1,0}
\definecolor{cb-pink}{cmyk}{0.1,0.7,0,0}

\definecolor{cb-purple}{RGB}{62,35,115}

\definecolor{alt-plum}{RGB}{51,34,136}
\definecolor{alt-green}{RGB}{17,119,51}
\definecolor{alt-turq}{RGB}{68,170,153}
\definecolor{alt-blue}{RGB}{136,204,238}
\definecolor{alt-yellow}{RGB}{221,204,119}
\definecolor{alt-pink}{RGB}{204,102,119}
\definecolor{alt-purple}{RGB}{170,68,153}
\definecolor{alt-red}{RGB}{136,34,85}

\colorlet{c1}{cb-orange}
\colorlet{c2}{cb-light-blue}
\colorlet{c3}{cb-green}
\colorlet{c4}{cb-pink}
\colorlet{c5}{cb-blue}
\colorlet{c6}{cb-red}
\colorlet{c7}{cb-yellow}
\colorlet{c8}{alt-red}

\tikzset{
  blank/.style={
    inner sep = 0cm, 
    outer sep = 0cm, 
    align=center, 
    minimum size=0pt
  },
  label/.style={
    minimum size=5mm, 
    inner sep=0.1cm,
    fill=white
  },
  midarrow/.style={
    currarrow,
    sloped,
    scale=0.5,
    allow upside down
  },
  shorten <>/.style={
    shorten >=#1,
    shorten <=#1
  },
  bullet/.style={
    rectangle,
    draw,
    fill=black,
    minimum size=3pt,
    inner sep=0pt,
  },
  colourBullet/.style={
    circle,
    draw=#1,
    fill=#1,
    minimum size=3pt,
    inner sep=0pt,
  },
  double colour/.style 2 args={fill=#1, 
  path picture={
    \fill[#2, sharp corners] 
    (path picture bounding box.south) -|
    (path picture bounding box.north east) --
    (path picture bounding box.north) --
    cycle;
    }
  },
  modification/.style 2 args={
    rounded corners=2,
    double colour={#1}{#2},
    minimum size=5mm,
    inner sep = 0.1cm, 
    outer sep = 0cm, 
    align=center,
    text=white
  },
  functorNode/.style={
    modification={#1}{#1}
  },
  twocellarr/.style={
    shorten <= 10pt, 
    shorten >= 10pt
  },
  twocell/.style={
    functorNode={black}
  },
  twocellX/.style ={
    rectangle,
    thick,
    fill=white,
    rounded corners=0.5,
    inner sep = 0.1cm, 
    outer sep = 0cm, 
    align=center,
    minimum width={#1 pt},
    draw=black,
    line width=\lineWidth
    },
  mod/.style ={
    rectangle,
    thick,
    fill=white,
    rounded corners=0.5,
    inner sep = 0.1cm, 
    outer sep = 0cm, 
    align=center,
    minimum width={#1 pt},
  },
  mod2/.style ={
    rectangle,
    dashed,
    thick,
    fill=white,
    rounded corners=0.5,
    inner sep = 0.1cm, 
    outer sep = 0cm, 
    align=center,
    minimum width={#1 pt},
  },
  colChange/.style ={
    rectangle,
    thick,
    fill=black,
    rounded corners=0.5,
    inner sep = 0.1cm, 
    outer sep = 0cm, 
    align=center,
    minimum width={#1 pt},
  },
  only in/.style={
    to path={
      let
      \p1 = ($ (\tikztotarget) - (\tikztostart) $)
      in
      .. controls 
      ($(\tikztostart)!{veclen(\x1,\y1)/3}!($(\tikztotarget)+(#1:{veclen(\x1,\y1)/3})$)$)
      and
      ($(\tikztotarget)+(#1:{veclen(\x1,\y1)/3})$)
      .. (\tikztotarget)
    }
  },
  only out/.style={
    \path draw (\tikztostart) --node[](temp1){} (\tikztotarget-|)
    to path={
      let
      \p1 = ($ (\tikztotarget) - (\tikztostart) $)
      in
      .. controls 
      ($(\tikztostart)+(#1:{veclen(\x1,\y1)/3})$)
      and
      ($(\tikztotarget)!{veclen(\x1,\y1)/3}!($(\tikztostart)+(#1:{veclen(\x1,\y1)/3})$)$)
      .. (\tikztotarget)
    }
  },
  braidTop/.style={
    to path={
      let
      \p1=($(\tikztotarget)-(\tikztostart)$)
      in
      ..controls +(0,\y1/2).. (\tikztotarget)
    }
  },
  braidBottom/.style={
    to path={
      let
      \p1=($(\tikztotarget)-(\tikztostart)$)
      in
      ..controls +(\x1/2,0).. (\tikztotarget)
    }
  },
  NTDash/.style={
    dash pattern= on 0.20cm off 0.05cm,
  },
}
\makeatletter
\newcommand{\gettikzxy}[3]{%
    \tikz@scan@one@point\pgfutil@firstofone#1\relax
    \edef#2{\the\pgf@x}%
    \edef#3{\the\pgf@y}%
}
\makeatother

\newcommandx{\drawStraightDownX}[2][1=1]{
  \foreach \i in {#2}{
    \draw[thick,black] (\i) to[out=-90, in=90] ++(0,-#1) node[blank] (\i){};
    }
}

\newcommandx{\drawId}[2][1=1]{
  \drawStraightDownX[#1]{#2}
}

\newcommandx{\drawBead}[5][1=0.5,4=0.5]{
    \draw[thick, black] (#2) to[out=-90, in=90] ++(0,-#1) node [twocell] {#3} to[out=-90, in=90] ++(0,-#4) node [blank] (#5){};
}

\newcommandx{\drawBeadBottom}[4][1=0.5]{
    \draw[thick, black] (#2) to[out=-90, in=90] ++(0,-#1) node [twocell](#4) {#3};
}

\newcommandx{\drawBeadTop}[4][3=0.5]{
    \drawBead[0]{#1}{#2}[#3]{#4}
}

\newcommandx{\drawCap}[7][2=1,3=1,6=,7=]{
    \draw (#1)++(#2,0) node[blank] (#4){};
    \draw (#4)++(0.5*#3,0.5) node[blank] (capTop){};
    \draw (#4)++(#3,0) node[blank] (#5){};

    \draw[thick,black] (#4) to[out=90, in=180] node[label,pos=0.1, anchor=east] {#6} (capTop);

    \draw[thick,black] (capTop) to[out=0, in=90] node[label,pos=0.9, anchor=west] {#7} (#5);
}

\newcommandx{\drawCup}[5][3=0.5,4=,5=,usedefault]{
    \path (#1) to node[pos=0.5, blank] (temp1){} (#2);
    \draw (temp1)++(0,-#3) node[blank] (temp2){};

    \strand[thick,black] (#1) to[out=-90, in=180] 
    node[label,pos=0.1, anchor=east] {#4} (temp2);

    \strand[thick,black] (temp2) to[out=0, in=-90] node[label,pos=0.9, anchor=west] {#5} (#2);
}

\newcommandx{\drawTwoToOne}[6][1=0.5,5=0.5]{
    \path (#2) to node[blank, pos=0.5] (temp1){} (#3);
    \draw (temp1)++(0,-#1) node[twocell](temp2){#4};

    \draw[thick,black] (#2) to [out=-90, in=135] (temp2);

    \draw[thick,black] (#3) to [out=-90, in=45] (temp2);

    \draw[thick,black] (temp2) to [out=-90, in=90] ++(0,-#5) node[blank](#6){};
}
\newcommandx{\drawOneToTwo}[7][1=0.5,4=1,5=0.5,usedefault]{
    \draw (#2)++(0,-#1) node[twocell] (temp1){#3};

    \draw[thick,black] (#2) to [out=-90, in=90] (temp1);

    \draw (temp1)++(-0.5*#4,-#5) node [blank] (#6){};
    \draw (temp1)++(0.5*#4,-#5) node [blank] (#7){};

    \draw[thick, black] (temp1) to [out=-135, in=90] (#6);

    \draw[thick,black]
      (temp1) to [out=-45, in=90] (#7);
}

\newcommandx{\drawTwoToTwo}[8][1=0.5,5=1,6=0.5]{
    \path (#2) to node[blank, pos=0.5] (temp1){} (#3);
    \draw (temp1)++(0,-#1) node[twocell](temp2){#4};

    \draw[thick,black] (#2) to [out=-90, in=135] (temp2);

    \draw[thick,black] (#3) to [out=-90, in=45] (temp2);

    \draw (temp2)++(-0.5*#5,-#6) node [blank] (#7){};
    \draw (temp2)++(0.5*#5,-#6) node [blank] (#8){};

    \draw[thick,black] (temp2) to [out=-135, in=90] (#7); 

    \draw[thick,black] (temp2) to [out=-45, in=90] (#8); 
}

    \newcommandx{\drawThreeToTwo}[9][1=0.5,6=1,7=0.5]{
      \path (#2) to node[blank, pos=0.5] (temp1){} (#4);
      \draw (temp1)++(0,-#1) node[twocell](temp2){#5};

      \draw[thick,black]
      (#2) to [out=-90, in=135] (temp2);

      \draw[thick,black]
      (#3) to [out=-90, in=90] (temp2);
      
      \draw[thick,black]
      (#4) to [out=-90, in=45] (temp2);
      
      \draw (temp2)++(-0.5*#6,-#7) node [blank] (#8){};
      \draw (temp2)++(0.5*#6,-#7) node [blank] (#9){};

      \draw[thick,black]
      (temp2) to [out=-135, in=90] (#8); 

      \draw[thick,black]
      (temp2) to [out=-45, in=90] (#9); 
    }
\newcommandx{\drawTwoCell}[6][1=0.5,5=0,6=0.5,usedefault]{
    \readlist*\inputs{#2}
    \xdef\xCoordinateList{}
    \xdef\yCoordinateList{}
    \foreach \index in {1,...,\inputslen}{
        \gettikzxy{(\inputs[\index])}{\ax}{\ay}
        \xdef\xCoordinateList{\ax,\xCoordinateList}
        \xdef\yCoordinateList{\ay,\yCoordinateList}
    }
    \pgfmathsetmacro\minX{min(\xCoordinateList)}
    \pgfmathsetmacro\maxX{max(\xCoordinateList)}
    \pgfmathsetmacro\xCoordinate{(\maxX+\minX)/2}
    \pgfmathsetmacro\yCoordinate{min(\yCoordinateList)}
    
    \coordinate (temp1) at (\xCoordinate pt, \yCoordinate pt);
    \coordinate (temp1) at ($(temp1)+(0,-#1)$);
    \draw (temp1) node[twocell](temp1){#3};
    
    \pgfmathsetmacro\startAngle{180*(1-(1/(2*\inputslen)))}
    \pgfmathsetmacro\increaseAngle{180/(\inputslen)}
    
    \foreach \index in {1,...,\inputslen}{
        \pgfmathsetmacro\angle{\startAngle+((1-\index)*\increaseAngle)}
        \draw[thick,black] (\inputs[\index]) to[out=-90, in=\angle] (temp1);
    }
    \draw (temp1) node[twocell](temp1){#3};

    \ifthenelse{\equal{#4}{}}{
    }
    {
    \readlist*\outputs{#4}
    \draw (temp1)++(-0.5*#5,-#6) node [blank] (temp2){};
    
    \pgfmathsetmacro\startAngle{180*((1/(2*\outputslen))-1)}
    \pgfmathsetmacro\increaseAngle{180/(\outputslen)}

    \ifthenelse{\outputslen>1}{
    \def\outputspace{
        {#5/(\outputslen-1)}}
    }{
        \def\outputspace{0}
    }

    \foreach \index in {1,...,\outputslen}{
        \pgfmathsetmacro\angle{\startAngle+((\index-1)*\increaseAngle)}

        \draw (temp2) node[blank] (\outputs[\index]){};
        \draw[thick,black] (temp1) to[out=\angle, in=90] (\outputs[\index]);
        
        \draw (temp2)++({\outputspace},0) node[blank] (temp2){};
    }
    }
}
    \newcommandx{\drawThreeToNone}[5][1=0.5]{
      \path (#2) to node[blank, pos=0.5] (temp1){} (#4);
      \draw (temp1)++(0,-#1) node[twocell](temp2){#5};

      \draw[thick,black]
      (#2) to [out=-90, in=135] (temp2);

      \draw[thick,black]
      (#3) to [out=-90, in=90] (temp2);
      
      \draw[thick,black]
      (#4) to [out=-90, in=45] (temp2);
    }
    
    \newcommandx{\drawCapX}[9][2=1,3=0.5,4=1,5=0.5,8=,9=, usedefault]{
      \newdimen \y
      \y = #2 pt
      \ifthenelse{\lengthtest{\y > 0pt}}{
      \draw (#1)++(#2,{#3-#5}) node[blank] (#6){};
      \draw (#1)++({#2+(0.5*#4)},#3) node[blank] (temp1){};
      \draw (#1)++({#2+#4},{#3-#5}) node[blank] (#7){};
      }
      {
        \draw (#1)++(#2,{#3-#5}) node[blank] (#7){};
        \draw (#1)++({#2-(0.5*#4)},#3) node[blank] (temp1){};
        \draw (#1)++({#2-#4},{#3-#5}) node[blank] (#6){};
      }
      \strand[thick,black] 
      (#6) to[out=90, in=180] 
      node[label,pos=0.1, anchor=east] {#8} (temp1);

      \strand[thick,black]
      (temp1) to[out=0, in=90] node[label,pos=0.9, anchor=west] {#9} (#7);
    }
    \newcommandx{\drawPositiveBraid}[6][3=1,4=-1,5=temp1,6=temp1,usedefault]{
      \path (#1) to node [pos=0.5, blank](temp1){} (#2);

      \draw (temp1)++(-0.5*#4,-#3) node[blank](temp2){};
      \draw (temp1)++(0.5*#4,-#3) node[blank](temp3){};

      \strand [thick,black]
      (#1) to[out=-90, in=90] (temp3);
      \strand [thick,black]
      (#2) to[out=-90, in=90] (temp2);

      \ifthenelse{\equal{#5}{temp1}}{
      \draw (temp2) node[blank] (#2){};
      \draw (temp3) node[blank] (#1){};
      }
      {
      \draw (temp2) node[blank] (#5){};
      \draw (temp3) node[blank] (#6){};
      }
    }
    \newcommandx{\drawNegativeBraid}[6][3=1,4=-1,5=temp1,6=temp1,usedefault]{
      \draw (#1)++(0,-#3) node[blank](temp2){};
      \draw (#2)++(0,-#3) node[blank](temp3){};

      \strand [thick,black]
      (#1) to[out=-90, in=90] (temp3);
      \strand [thick,black]
      (#2) to[out=-90, in=90] (temp2);

      \ifthenelse{\equal{#5}{temp1}}{
      \draw (temp2) node[blank] (#2){};
      \draw (temp3) node[blank] (#1){};
      }
      {
      \draw (temp2) node[blank] (#5){};
      \draw (temp3) node[blank] (#6){};
      }
    }
    \newcommandx{\drawBraidX}[6][3=1,4=-1,5=temp1,6=temp1,usedefault]{
      \newdimen \y
      \y = #4 pt
      \ifthenelse{\lengthtest{\y < 0pt}}
      {
        \drawNegativeBraid{#1}{#2}[#3][#4][#5][#6]
      }
      {
        \drawPositiveBraid{#1}{#2}[#3][#4][#5][#6]
      }
    }
    \newcommandx{\drawAcrossBraid}[2]{
      \draw (#1)++(0.2,0) node[blank](temp1){};
      \draw (#2)++(-0.2,0) node[blank](temp2){};

      \draw[thick, black](temp1) to [out=180, in=0] (temp2);
    }
    \newcommand{\drawStartNodes}[1]{
      \xdef\xCoord{0}
      \foreach \i/\j/\k in {#1}
      {
        \coordinate (\i) at (\xCoord,0);
        \draw (\i) node[label,anchor=south] {\j};
        \xdef\xCoord{\xCoord+\k}
      }
    }
    \newcommand{\drawEndNodes}[1]{
      \foreach \i/\j in {#1}
      {
        \draw(\i) node[label,anchor=north]{\j};
      }
    }
\newcommandx{\drawCapXX}[5]{
  \readlist*\feet{#4}
  \readlist*\feetlabels{#5}
  \coordinate (\feet[1]) at ($(#1)+(-0.5*#2,-#3)$);
  \coordinate (\feet[2]) at ($(#1)+(0.5*#2,-#3)$);

  \strand[thick,black]
    (\feet[1]) to[out=90, in=180] 
    node[pos=0.1, label] {\feetlabels[1]} (#1) 
    to[out=0, in=90] 
    node[pos=0.9, label]{\feetlabels[2]} (\feet[2]);
}
\newcommandx{\drawBraidLeft}[4]{
  \readlist*\inputDBL{#1}
  \readlist*\colours{#2}
  \readlist*\output{#4}
  \gettikzxy{(\inputDBL[1])}{\ax}{\ay}
  \gettikzxy{(\inputDBL[2])}{\bx}{\by}
  
  \pgfmathsetmacro\xCoord{(\ax+\bx)/2}
  \pgfmathsetmacro\yCoord{min(\ay,\by)}

  \coordinate (temp1) at ($(\xCoord pt,\yCoord pt)-(0,0.5*#3)$);

  \draw[thick, natTrans={green}{blue}]
  (\inputDBL[1]) to[out=-90, in=90] ($(\bx,\yCoord)-(0,#3)$) coordinate (\output[2]);

  \draw [shorten >= 0.2*#3,thick] 
  (\inputDBL[2]) to[out=-90, in=45] (temp1);

  \draw [shorten <= 0.2*#3, thick] 
  (temp1) to[only in=90] ($(\ax,\yCoord)-(0,#3)$) coordinate (\output[1]);
}
\newcommandx{\TwoTone}[6]{
  \xdef\tempPath{
    (#1) to[out=#2, in=#3] (#4)
  }
  \draw[very thick, color=#5] \tempPath;
  \begin{scope}[overlay]
    \clip  \tempPath --(20,0)--cycle;
    \draw[color=#6,very thick] \tempPath;
  \end{scope}
}

\newcommandx{\stringTwoTone}[6][6=]{
    \ifthenelse{\equal{#6}{}}{
        \tikzset{
            optionZero/.style={
                line width=2pt, color=white
            },
            optionOne/.style={
                line width=2pt, color=#2
            },
            optionTwo/.style={
                line width=2pt, color=#3
            }
        }
    }
    {
        \tikzset{
            optionZero/.style={
                very thick, draw=white, double=white, double distance=2pt
            },
            optionOne/.style={
                very thick, draw=white, double=#2, double distance=2pt
            },
            optionTwo/.style={
                very thick, draw=white, double=#3, double distance=2pt
            }
        }
    }
    \begin{scope}
        \draw[optionZero] #1;
    \end{scope}
    \begin{scope}[overlay]
        \clip #1 #4;
        \draw[optionOne] #1;
    \end{scope}
    \begin{scope}[overlay]
        \clip #1 #5;
        \draw[optionTwo] #1;
    \end{scope}
    \begin{scope}[]
        \draw[thin, white] #1;
    \end{scope}
}
\newcommandx{\stringLeftBraid}[5][5=1]{
\setsepchar{,}
\readlist*\inputSLB{#1}
\readlist*\output{#2}
\readlist*\colours{#3}

\gettikzxy{(\inputSLB[1])}{\ax}{\ay}
\gettikzxy{(\inputSLB[2])}{\bx}{\by}
  
\pgfmathsetmacro\xCoord{(\ax+\bx)/2}
\pgfmathsetmacro\yCoord{min(\ay,\by)}

\coordinate (temp1) at ($(\xCoord pt,\yCoord pt)-(0,0.5*#5)$);

\pgfmathsetmacro\xCoord{min(\ax,\bx)}

\coordinate (temp2) at ($(\xCoord pt,\yCoord pt)-(0,#5)$);

\pgfmathsetmacro\xCoord{max(\ax,\bx)}

\coordinate (temp3) at ($(\xCoord pt,\yCoord pt)-(0,#5)$);

\draw[thick,color=#4] (\inputSLB[2]) to[only out=-90] (temp1);

\draw[thick,color=#3] (temp1) to[only in=90] (temp2);

\stringTwoTone{(\inputSLB[1]) to [out=-90,in=90] (temp3)} {#4}{#3}{
    --++(0,-2)
    --++(2,0)
    --($(\inputSLB[2])+(2,2)$)
    --($(\inputSLB[1])+(0,2)$)
    --++(0,-2)
}{
    --++(0,-2)
    --($(temp2)+(-2,-2)$)
    --($(\inputSLB[1])+(-2,2)$)
    --++(2,0)
    --++(0,-2)
}[
    braid
]
\coordinate (\output[1]) at (temp2);
\coordinate (\output[2]) at (temp3);
}
\newcommandx{\stringRightBraid}[5][5=1]{
    \setsepchar{,}
    \readlist*\inputSRB{#1}
    \readlist*\output{#2}
    \readlist*\colours{#3}
    
    \gettikzxy{(\inputSRB[1])}{\ax}{\ay}
    \gettikzxy{(\inputSRB[2])}{\bx}{\by}
      
    \pgfmathsetmacro\xCoord{(\ax+\bx)/2}
    \pgfmathsetmacro\yCoord{min(\ay,\by)}
    
    \coordinate (temp1) at ($(\xCoord pt,\yCoord pt)-(0,0.5*#5)$);

    \pgfmathsetmacro\xCoord{min(\ax,\bx)}

    \coordinate (temp2) at ($(\xCoord pt,\yCoord pt)-(0,#5)$);

    \pgfmathsetmacro\xCoord{max(\ax,\bx)}

    \coordinate (temp3) at ($(\xCoord pt,\yCoord pt)-(0,#5)$);
    
    \draw[thick,color=#3] (\inputSRB[1]) to[only out=-90] (temp1);
    
    \draw[thick,color=#4] (temp1) to[only in=90] (temp3);
    
    \stringTwoTone{
        (\inputSRB[2]) to[out=-90, in=90] (temp2)
        }{#4}{#3}{
        --++(0,-2)
        --($(temp3)+(2,-2)$)
        --($(\inputSRB[2])+(2,2)$)
        --++(-2,0)
        --++(0,-2)
    }{
        --++(0,-2)
        --++(-2,0)
        --($(\inputSRB[1])+(-2,2)$)
        --($(\inputSRB[2])+(0,2)$)
        --++(0,-2)
    }[
        braid
    ]
    \coordinate (\output[1]) at (temp2);
    \coordinate (\output[2]) at (temp3);
}

\newcommandx{\stringCapCC}[6][]{
    \readlist*\outputnodes{#2}
    \coordinate (temp1) at ($(#1)+(-0.5*#3,-#4)$);
    \coordinate (temp2) at ($(#1)+(0.5*#3,-#4)$);
    \stringTwoTone{
        (temp1) to[out=90, in=180] (#1)
        to [out=0, in=90] (temp2)
    }{#5}{#6}
    {
        --++(0,-1)
        --++(4,0)
        --++(0,#4)
        --++(0,3)
        --++(-#3,0)
        --++(-8,0)
        --++(0,-3)
        --++(0,-#4)
        --++(4,0)
        --++(0,1)
    }
    {
        --++(0,-1)
        --++(#3,0)
        --++(0,1)
    }
    \coordinate (\outputnodes[1]) at (temp1);
    \coordinate (\outputnodes[2]) at (temp2);
}
\newcommandx{\stringCapC}[5]{
    \readlist*\outputnodes{#2}
    \coordinate (temp1) at ($(#1)+(-0.5*#3,-#4)$);
    \coordinate (temp2) at ($(#1)+(0.5*#3,-#4)$);
    \draw[thick,#5] (temp1) to[out=90, in=180] (#1)
        to [out=0, in=90] (temp2);
    \coordinate (\outputnodes[1]) at (temp1);
    \coordinate (\outputnodes[2]) at (temp2);
}

\newcommandx{\stringCap}[5][3=1,4=0.5,5=black,usedefault]{
    \setsepchar{,}
    \readlist*\colours{#5}
    \ifthenelse{\colourslen>1}{
        \xdef\colourOne{\colours[1]}
        \xdef\colourTwo{\colours[2]}
        \stringCapCC{#1}{#2}{#3}{#4}{\colourOne}{\colourTwo}
    }{
        \xdef\colourOne{\colours[1]}
        \stringCapC{#1}{#2}{#3}{#4}{\colourOne}
    }
}
\newcommandx{\stringCup}[3][2=0.5,3=black,usedefault]{
    \readlist*\inputStringCup{#1}
    \coordinate (temp1) at (\inputStringCup[1]);
    \coordinate (temp3) at (\inputStringCup[2]);

    \gettikzxy{(\inputStringCup[1])}{\ax}{\ay}
    \gettikzxy{(\inputStringCup[2])}{\bx}{\by}
      
    \pgfmathsetmacro\xCoord{(\ax+\bx)/2}
    \pgfmathsetmacro\yCoord{min(\ay,\by)}

    \coordinate (temp2) at ($(\xCoord pt,\yCoord pt)+(0,-#2)$);

    \readlist*\colours{#3}
    \ifthenelse{\colourslen>1}{
        \xdef\colourOne{\colours[1]}
        \xdef\colourTwo{\colours[2]}
    \stringTwoTone{
        (temp1) to[out=-90, in=180] (temp2)
        to[out=0, in=-90] (temp3)
    }{\colourOne}{\colourTwo}
    {
        --++(0,2)
        --++(2,0)
        --($(temp2)+(\xCoord pt,-2)+(2,0)$)
        --($(temp2)+(-\xCoord pt,-2)-(2,0)$)
        --($(temp1)+(-2,2)$)
        --++(2,0)
        --++(0,-2)
    }
    {
        --++(0,2)
        --($(temp1)+(0,2)$)
        --++(0,-2)
    }
    }{
        \xdef\colourOne{\colours[1]}
        \draw[thick, \colourOne] (temp1) to[out=-90, in=180] (temp2) to[out=0, in=-90] (temp3);
    }
}

\newcommandx{\stringIdCC}[4][4=1]{
    \coordinate (temp1) at ($(#1)+(0,-#4)$);
    \stringTwoTone{
        (#1)to[out=-90,in=90](temp1)
        }{#2}{#3}{
        --++(0,-1)
        --++(-10,0)
        --++(0,2)
        --++(0,#4)
        --++(10,0)
    }
    {
        --++(0,-1)
        --++(10,0)
        --++(0,2)
        --++(0,#4)
        --++(-10,0)
    }[]
    \coordinate(#1) at (temp1);
}
\newcommandx{\stringIdC}[3][3=1]{
    \coordinate (temp1) at ($(#1)+(0,-#3)$);
    \draw [color=#2, thick] (#1) to (temp1);
    \coordinate (#1) at (temp1);
}

\newcommandx{\stringId}[3][2=1,3=black,usedefault]{
    \setsepchar{,}
    \readlist*\colours{#3}
    \ifthenelse{
        \colourslen>1
    }{
        \xdef\colourOne{\colours[1]}
        \xdef\colourTwo{\colours[2]}
        \stringIdCC{#1}{\colourOne}{\colourTwo}[#2]
    }{
        \xdef\colourOne{\colours[1]}
        \stringIdC{#1}{\colourOne}[#2]
    }
}

\newcommandx{\colourChange}[5][4=0.5,5=0.5]{
    \stringId{#1}[#4][#2]
    \coordinate (temp11) at (#1);
    \stringId{#1}[#5][#3]
    \draw (temp11) node[bullet]{};
}
\newcommandx{\stringLabel}[5][4=0.5,5=0.5]{
    \stringId{#1}[#4][#2]
    \coordinate (temp11) at (#1);
    \stringId{#1}[#5][#2]
    \draw (temp11) node[label]{#3};
}

\newcommandx{\clipPathS}[2]{
    \gettikzxy{(#1)}{\xTop}{\yTop}
    \gettikzxy{(#2)}{\xBottom}{\yBottom}

    \newdimen \xTopDim
    \xTopDim = \xTop pt
    \newdimen \xBottomDim
    \xBottomDim =\xBottom pt

    \ifthenelse{
        \lengthtest{\xTopDim <\xBottomDim}
    }{
        \xdef\leftClipPath{
            --++(0,-2)
            --($(\xTop,\yBottom)+(-2,-2)$)
            --($(#1)+(-2,2)$)
            --++(2,0)
            --++(0,-2)
        }
        \xdef\rightClipPath{
            --++(0,-2)
            --++(2,0)
            --($(\xBottom, \yTop)+(2,2)$)
            --($(#1)+(0,2)$)
            --++(0,-2)
        }
    }{
        \xdef\leftClipPath{
            --++(0,-2)
            --++(-2,0)
            --($(\xBottom, \yTop)+(-2,2)$)
            --($(#1)+(0,2)$)
            --++(0,-2)
        }  
        \xdef\rightClipPath{
            --++(0,-2)
            --($(\xTop,\yBottom)+(2,-2)$)
            --($(#1)+(2,2)$)
            --++(-2,0)
            --++(0,-2)
        }
    }
}
\newcommandx{\stringTwoCell}[7][1=0.5,5=0,6=0.5,7=twocell,usedefault]{
    \setsepchar[.]{,./}
    \readlist*\inputs{#2}
    \xdef\xCoordinateList{}
    \xdef\yCoordinateList{}
    \foreach \index in {1,...,\inputslen}{
        \gettikzxy{(\inputs[\index,1])}{\ax}{\ay}
        \xdef\xCoordinateList{\ax,\xCoordinateList}
        \xdef\yCoordinateList{\ay,\yCoordinateList}
    }
    \pgfmathsetmacro\minX{min(\xCoordinateList)}
    \pgfmathsetmacro\maxX{max(\xCoordinateList)}
    \pgfmathsetmacro\xCoordinate{(\maxX+\minX)/2}
    \pgfmathsetmacro\yCoordinate{min(\yCoordinateList)}
    
    \coordinate (temp1) at (\xCoordinate pt, \yCoordinate pt);
    \coordinate (temp1) at ($(temp1)+(0,-#1)$);
    
    \pgfmathsetmacro\startAngle{180*(1-(1/(2*\inputslen)))}
    \pgfmathsetmacro\increaseAngle{180/(\inputslen)}
    \foreach \index in {1,...,\inputslen}{
        \xdef\currentPath{
            (\inputs[\index,1]) to[out=-90, in=\angle] (temp1)
        }
        \pgfmathsetmacro\angle{\startAngle+((1-\index)*\increaseAngle)}
        \ifthenelse{\listlen\inputs[\index]=1}{
            \draw[thick,black] \currentPath;
        }{}
        \ifthenelse{\listlen\inputs[\index]=2}{
            \xdef\currentColour{\inputs[\index,2]}
            \tikzset{
                workAround/.style={
                    thick,
                    color=\currentColour
                }
            }
            \draw[workAround] \currentPath;
        }{}
        \ifthenelse{\listlen\inputs[\index]>2}{
            \xdef\colourOne{\inputs[\index,2]}
            \xdef\colourTwo{\inputs[\index,3]}
            
            \clipPathS{\inputs[\index,1]}{temp1}

            \stringTwoTone{\currentPath}{\colourOne}{\colourTwo}{
                \leftClipPath;
            }
            {
                \rightClipPath;
            }
        }{}
    }
    \ifthenelse{\equal{#4}{}}{
    }
    {
    \readlist*\outputs{#4}
    \draw (temp1)++(-0.5*#5,-#6) node [blank] (temp2){};
    
    \pgfmathsetmacro\startAngle{180*((1/(2*\outputslen))-1)}
    \pgfmathsetmacro\increaseAngle{180/(\outputslen)}

    \ifthenelse{\outputslen>1}{
    \def\outputspace{
        {#5/(\outputslen-1)}}
    }{
        \def\outputspace{0}
    }

    \foreach \index in {1,...,\outputslen}{
        \pgfmathsetmacro\angle{\startAngle+((\index-1)*\increaseAngle)}

        \draw (temp2) node[blank] (\outputs[\index,1]){};
        \xdef\currentPath{
            (temp1) to[out=\angle, in=90] (\outputs[\index,1])
        }
        \ifthenelse{\listlen\outputs[\index]=1}
        {
            \draw[thick,black] \currentPath;
        }{}
        \ifthenelse{\listlen\outputs[\index]=2}{
            \xdef\currentColour{\outputs[\index,2]}
            \tikzset{
                workAround/.style={
                    thick,
                    color=\currentColour
                }
            }
            \draw[workAround] \currentPath;
        }{}
        \ifthenelse{\listlen\outputs[\index]>2}{
            \xdef\colourOne{\outputs[\index,2]}
            \xdef\colourTwo{\outputs[\index,3]}
            
            \clipPathS{temp1}{\outputs[\index,1]}

            \stringTwoTone{\currentPath}{\colourOne}{\colourTwo}{
                \leftClipPath;
            }
            {
                \rightClipPath;
            }
        }{}
        \draw (temp2)++({\outputspace},0) node[blank] (temp2){};
    }
    }
    \draw (temp1) node[#7](temp1){#3};
    \setsepchar[.]{,./}
}

    \setsepchar[.]{;./}

\xdef\gapWidth{1.2}
\xdef\hlWidth{1.2}
\xdef\lineWidth{0.4}
\xdef\outerLineWidth{0.2}

\tikzmath{
    function angle(\a){
        return Mod(\a,360);
    };
    function bisect(\i,\o){
        if angle(\i)<angle(\o) then{
            return angle((angle(\i)+angle(\o))/2);
        }
        else{
            return angle(((angle(\i)+angle(\o))/2)-180);
        };
    };
    function shiftToFactor(\s,\x1,\x2){
        if \x1==\x2 then{
            return 0;
        } else{
            return (1+(\s/(\x1-\x2)))/(28.45274);
        };
    };
    function vecAngle(\x,\y){
        if \y<0 then {
            \s=-1;
        } else {
            \s=1;
        };
        return 180+\s*acos(\x/(veclen(\x,\y)));
    };
    function hlDblSep(\n){
        return \lineWidth+2*\gapWidth+2*(\n-1)*(\hlWidth+\gapWidth);
    };
    function hlDblSep2(\n){
        return 2*\outerLineWidth+4*\gapWidth+\lineWidth+2*\gapWidth+2*(\n-1)*(\hlWidth+\gapWidth);
    };
    function TCLShft(\n){
        return ((\lineWidth-hlDblSep(\n))/2)-\hlWidth;
    };
    function TCelShft(\l,\r){
        return (1/2)*(\r-\l)*(\hlWidth+\gapWidth);
    };
    function TCWidth(\n,\l,\r){
        return \n+2*\hlWidth+hlDblSep(\l)+hlDblSep(\r)+10;
    };
    function TCWidthX(\n,\l,\r){
        return \n+(\l+\r)*(\hlWidth+\gapWidth)+10;
    };
    function symWidth(\n,\l,\r){
        return (max(\l*2*(\hlWidth+\gapWidth),\r*2*(\hlWidth+\gapWidth)))+\n;
    };
    function testing(\l,\r){
        return (max(\l*2*(\hlWidth+\gapWidth),\r*2*(\hlWidth+\gapWidth)));
    };
    function twoCellWidth(\n,\l,\r,\n2,\l2,\r2){
        return 10+max(symWidth(\n,\l,\r), symWidth(\n2,\l2,\r2));
    };
    function outerLineSep(\n){
        return \lineWidth+2*\gapWidth+2*\n*(\hlWidth+\gapWidth);
    };
    function joinAngle(\i,\n){
        return 180*(1-(1/(2*\n)))+(1-\i)*(180/\n);
    };
    function spread(\i,\n,\w){
        if \n==1 then{
            return 0;
        }else{
            return (\i-1)*(\w/(\n-1))-(\w/2);
        };
    };
    function idRad(\i){
        return (1/2)*\lineWidth+(1/2)*\hlWidth+\gapWidth+(\i-1)*(\hlWidth+\gapWidth);
    };
}

\newcommand{\ifempty}[2]{
    \ifthenelse{\equal{}{#1}}{
        #2
    }{}
}

\newcommand{\ifnonempty}[2]{
    \ifthenelse{\equal{}{#1}}{}{
        #2
    }
}

\newcommandx{\place}[4][3=0]{
    \coordinate[shift={(#2,#3)}] (tempPlace) at (#4); 
    \coordinate (#1) at (tempPlace); 
}

\newcommand{\listLength}[2]{
    \edef\inputLL{#1}
    \xdef#2{0}
    \foreach \i in \inputLL{
        \tikzmath{
            integer \currentLength;
            \currentLength = #2+1;
        }
        \xdef#2{\currentLength}
    }
}

\newcommand{\listLengthX}[2]{
    \edef\inputLLX{#1}
    \xdef#2{0}
    \foreach \i in \inputLLX{
        \ifthenelse{\equal{\i}{}}{}{
            \tikzmath{
                integer \currentLength;
                \currentLength = #2+1;
            }
            \xdef#2{\currentLength}
        }
    }
}

\newcommand{\reverseList}[2]{
    \edef\inputRL{#1}
    \xdef#2{}
    \foreach \x[count=\i] in \inputRL{
        \ifthenelse{\equal{#2}{}}{
            \xdef#2{\x}
        }{
            \xdef#2{\x,#2}
        }
    }
}

\newcommandx{\joinTarget}[5][4=\nodesWidth,5=-1,usedefault]{
    \readlist*\inputJT{#1}
    \xdef\xCoordinateList{}
    \xdef\yCoordinateList{}

    \foreachitem \n \in \inputJT{
        \gettikzxy{(\n)}{\ax}{\ay}
        \xdef\xCoordinateList{\ax,\xCoordinateList}
        \xdef\yCoordinateList{\ay,\yCoordinateList}
    }
    \tikzmath{
        \minX=min(\xCoordinateList);
        \maxX=max(\xCoordinateList);
        \xCoordinate=(\maxX+\minX)/2;
        \minY=min(\yCoordinateList);
        \maxY=max(\yCoordinateList);
        \width=\maxX-\minX;
    }
    \ifdim #5 pt < 0 pt
        \xdef\yCoordinate{\minY}
    \else
        \xdef\yCoordinate{\maxY}
    \fi{}
    \xdef#4{\width}
    \coordinate[shift={(0,#2)}] (#3) at (\xCoordinate pt, \yCoordinate pt);
}

\newcommandx{\splitTargets}[4]{
    \readlist*\splitInput{#2}
    \foreach \i in {1,...,\splitInputlen}{
        \tikzmath{
            \p=spread(\i,\splitInputlen,#3);
        }
        \coordinate[shift={(\p,-#4)}] (\splitInput[\i]) at (#1);
    }
}

\newcommandx{\halfPath}[3][1=]{
    \begin{scope}
        \begin{pgfinterruptboundingbox}
        \clip #2 #3;
        \end{pgfinterruptboundingbox}
        \draw[#1] #2;
    \end{scope}
}

\newcommandx{\removeEmpty}[2]{
    \xdef\removeEmptyInput{#1}
    \xdef\output{}
    \foreach \x[count=\i] in \removeEmptyInput{
        \ifnonempty{\x}{
            \ifnonempty{\output}{
                \xdef\output{\output,\x}
            }
            \ifempty{\output}{
                \xdef\output{\x}
            }
        }
    }
    \xdef#2{\output}
}

\newcommandx{\highlightSide}[4]{
    \removeEmpty{#1}{\coloursI}        
    \listLength{\coloursI}{\colourslenI}
    \removeEmpty{#2}{\coloursII}
    \listLength{\coloursII}{\colourslenII}

    \foreach \c [count=\k] in \coloursI{
        \ifthenelse{\equal{\c}{}}{}{
        \halfPath[
            double,
            double=white,
            draw=\c, 
            double distance=hlDblSep(\colourslenI+\colourslenII+1-\k), 
            line width=\hlWidth pt,
            rounded corners=0.5
        ]{#3}{#4}
        }
    }
    \foreach \c [count=\k] in \coloursII{
        \ifthenelse{\equal{\c}{}}{}{
        \halfPath[
            double,
            double=white,
            draw=\c, 
            double distance=hlDblSep(\colourslenII+1-\k), 
            line width=\hlWidth pt,
            rounded corners=0.5,
            NTDash
            ]{#3}{#4}
        }
    }
}

\newcommandx{\highlight}[7][7=]{
    \xdef\coloursTempI{#2}
    \xdef\coloursTempII{#3}
    \reverseList{#4}{\coloursTempIII}    
    \highlightSide{\coloursTempI}{\coloursTempII}{#1}{#5}
    \highlightSide{\coloursTempI}{\coloursTempIII}{#1}{#6}
    
    \listLengthX{#3}{\ntleftlen}
    \listLengthX{#4}{\ntrightlen}
    \tikzmath{
        integer \nt;
        \nt=\ntleftlen+\ntrightlen;
    }
    \ifnum \nt>0
        \draw[line width=\lineWidth,NTDash,#7] #1;
    \else
        \draw[line width=\lineWidth,#7] #1;
    \fi
}

\newcommandx{\clipS}[4][3=\leftClip,4=\rightClip,usedefault]{
    \gettikzxy{(#1)}{\ix}{\iy}
    \gettikzxy{(#2)}{\ox}{\oy}

    \ifdim \iy>\oy
        \xdef\s{-1}
    \else
        \xdef\s{1}
    \fi

    \ifdim \ix <\ox
        \xdef#3{
            --++(0,\s*3)
            -|($(#1)+(-3,-\s*3)$)
            -|(#1)
        }
        \xdef#4{
            |-++(3,\s*3)
            |-($(#1)+(0,-\s*3)$)
            --(#1)
        }
    \else
        \xdef#3{
            |-++(-3,\s*3)
            |-($(#1)+(0,-\s*3)$)
            --(#1)
        }
        \xdef#4{
            --++(0,\s*3)
            -|($(#1)+(3,-\s*3)$)
            -|(#1)
        }
    \fi{}
}

\newcommandx{\drawS}[6][3=,4=,5=,6=,usedefault]{
    \clipS{#1}{#2}
    \highlight{(#1) to[out=-90, in=90] (#2)}{#3}{#4}{#5}{\leftClip}{\rightClip}[#6]
}

\newcommandx{\strS}[5][3=,4=,5=]{
    \drawS{#1}{#2}[#3][#4][#5]
    \coordinate (#1) at (#2);
}

\newcommandx{\strSX}[5][3=,4=,5=]{
    \coordinate[shift={(#2)}] (temp) at (#1);
    \drawS{#1}{temp}[#3][#4][#5]
    \coordinate (#1) at (temp);
}

\newcommandx{\drawIdTo}[6][3=,4=,5=,6=,usedefault]{
    \clipS{#1}{#2}
    \highlight{(#1)--(#1|-#2)}{#3}{#4}{#5}{\leftClip}{\rightClip}[#6]
}

\newcommandx{\strIdTo}[6][3=,4=,5=,6=,usedefault]{
    \drawIdTo{#1}{#2}[#3][#4][#5][#6]
    \coordinate (#1) at (#1|-#2);
}

\newcommandx{\strId}[6][2=,3=,4=,5=1,6=,usedefault]{
    \coordinate[shift={(0,-#5)}] (temp) at (#1);
    \strIdTo{#1}{temp}[#2][#3][#4][#6]
}

\newcommandx{\strBraid}[9][3=,4=,5=,6=,7=,8=1,9=;,usedefault]{
    \joinTarget{#1;#2}{0}{topLevel}
    \readlist\options{#9}
    \xdef\optionI{\options[1,1]}
    \xdef\optionII{\options[2,1]}

    \coordinate[shift={(0,-#8)}] (rightBottom) at (#2|-topLevel);
    \coordinate[shift={(0,-#8)}] (leftBottom) at (#1|-topLevel);

    \xdef\overPath{
        (#1) to[out=-90, in=90] (rightBottom)
    }

    \clipS{#1}{rightBottom}
    \begin{scope}
        \begin{pgfinterruptboundingbox}
            \clip \overPath\rightClip;
        \end{pgfinterruptboundingbox}
        \drawS{#2}{leftBottom}[#3,#5][#6][#7][\optionII]
    \end{scope}
    \begin{scope}
        \begin{pgfinterruptboundingbox}
            \clip \overPath\leftClip;
        \end{pgfinterruptboundingbox}
        \drawS{#2}{leftBottom}[#3,#4][#6][#7][\optionII]
    \end{scope}
    \drawS{#1}{rightBottom}[white,#3][#4][#5][\optionI]
    \coordinate (#1) at (rightBottom);
    \coordinate (#2) at (leftBottom);
}

\newcommandx{\strBraidX}[3][2=1,3=;,usedefault]{
    \readlist*\braidInput{#1}

    \xdef\topN{\braidInput[1,1]}
    \xdef\topFCol{\braidInput[1,2]}
    \xdef\topNTLCol{\braidInput[1,3]}
    \xdef\topNTRCol{\braidInput[1,4]}
    \xdef\botN{\braidInput[2,1]}
    \xdef\botFCol{\braidInput[2,2]}
    \xdef\botNTLCol{\braidInput[2,3]}
    \xdef\botNTRCol{\braidInput[2,4]}

    \joinTarget{\topN;\botN}{0}{topLevel}
    \readlist*\options{#3}
    \xdef\optionI{\options[1,1]}
    \xdef\optionII{\options[2,1]}

    \coordinate[shift={(0,-#2)}] (rightBottom) at (\botN|-topLevel);
    \coordinate[shift={(0,-#2)}] (leftBottom) at (\topN|-topLevel);

    \xdef\overPath{
        (\topN) to[out=-90, in=90] (rightBottom)
    }

    \clipS{\topN}{rightBottom}
    \begin{scope}
        \begin{pgfinterruptboundingbox}
            \clip \overPath\rightClip;
        \end{pgfinterruptboundingbox}

        \drawS{\botN}{leftBottom}[
            \topFCol,
            \topNTRCol,
            \botFCol
        ][\botNTLCol][\botNTRCol][\optionII]
    \end{scope}
    \begin{scope}
        \begin{pgfinterruptboundingbox}
            \clip \overPath\leftClip;
        \end{pgfinterruptboundingbox}
        \drawS{\botN}{leftBottom}[
            \topFCol,
            \topNTLCol,
            \botFCol
        ][\botNTLCol][\botNTRCol][\optionII]
    \end{scope}
    \drawS{\topN}{rightBottom}[
        white,\topFCol
    ][
        \topNTLCol
    ][
        \topNTRCol][\optionI]
    \coordinate (\topN) at (rightBottom);
    \coordinate (\botN) at (leftBottom);
}

\newcommandx{\drawSCirc}[6][3=,4=,5=,6=-0.5]{
    \joinTarget{#1;#2}{#6}{joinPoint}[][#6]
    \tikzmath{
        \a=sign(#6)*90;
        \shft=(3*(sign(#6)));
    }
    \highlight{
        (#1) to[out=\a, in=180] (joinPoint) to[out=0, in=\a] (#2)
    }{#3}{#4}{#5}{
        |-++(3,-\shft)
        |-($(joinPoint)+(0,\shft)$)
        -|($(#1)+(-3,-\shft)$)
        -|(#1)
    }{
        --++(0,-\shft)
        -|(#1)
    }
}

\newcommandx{\strCup}[6][3=,4=,5=,6=0.5,usedefault]{
    \drawSCirc{#1}{#2}[#3][#4][#5][-#6]
}

\newcommandx{\strCap}[6][3=,4=,5=,6=0.5,usedefault]{
    \drawSCirc{#1}{#2}[#3][#4][#5][#6]
}

\newcommandx{\getNodes}[2]{
    \readlist*\inputGN{#1}
    \xdef\output{}
    \foreach \i in {1,...,\inputGNlen}{
        \ifnum \i=1
            \xdef\output{\inputGN[\i,1]}
        \else
            \xdef\output{\output;\inputGN[\i,1]}
        \fi
    }
    \xdef#2{\output}
}

\newcommandx{\drawFork}[6][1=,2=,4=-0.5,5=-0.5,6=]{
    \getNodes{#3}{\nodes}
    \joinTarget{\nodes}{#4}{join}[][#4]
    \coordinate[shift={(0,#5)}] (exit) at (join);
    
    \readlist*\prong{#3}

    \xdef\natTrans{0}
    \ifnonempty{#2}{
        \xdef\natTrans{1}
    }
    \foreach \i in {1,...,\pronglen}{
        \ifnonempty{prong[\i,1]}{
            \xdef\natTrans{1}
        }
    }
    \ifnum \natTrans>0
        \tikzset{
            temp/.style={
                NTDash,#6,line width=\lineWidth
            }
        }
    \else
        \tikset{
            \temp./style={
                line width=\lineWidth,#6
            }
        }
    \fi

    \tikzmath{
        \s=sign(#4);
        \shft=(-3)*(sign(#4));
    }

    \pgfmathsetmacro\angle{
        -\s*joinAngle(1,\pronglen)
    }

    \clipS{\prong[1,1]}{exit}
    
    \highlightSide{#1}{#2}{
        (\prong[1,1])
        to[out=\s*90,in=\angle] (join) 
        to[out=\s*90, in=-\s*90] (exit)
    }{
        \leftClip
    }

    \draw[temp]      
    (\prong[1,1])
    to[out=\s*90,in=\angle] (join) 
    to[out=\s*90, in=-\s*90] (exit);

    \ifnum \pronglen>1
        \tikzmath{\finish=\pronglen-1;}
        \foreach \i in {1,...,\finish}{

            \xdef\currentNode{\prong[\i,1]}
            \xdef\nextNode{\prong[\i+1,1]}
            \reverseList{\prong[\i,2]}{\currentColours}
            \pgfmathsetmacro\inAngle{
                -\s*joinAngle(\i,\pronglen)
            }
            \pgfmathsetmacro\outAngle{
                -\s*joinAngle(\i+1,\pronglen)
            }
            \highlightSide{#1}{\currentColours}{
                (\currentNode)
                to[out=\s*90,in=\inAngle] (join) 
                to[out=\outAngle, in=\s*90] (\nextNode)
            }{
                --++(0,\shft)
                -|(\currentNode)
            }

            \draw[temp]      
            (\currentNode)
            to[out=\s*90,in=\inAngle] (join) 
            to[out=\outAngle, in=\s*90] (\nextNode);
            
        }
    \fi

    \reverseList{\prong[-1,2]}{\currentColours}

    \pgfmathsetmacro\angle{
        -\s*joinAngle(\pronglen,\pronglen)
    }
    
    \clipS{\prong[-1,1]}{exit}

    \highlightSide{#1}{\currentColours}{
        (\prong[-1,1])
        to[out=\s*90,in=\angle] (join) 
        to[out=\s*90, in=-\s*90] (exit)
    }{\rightClip}

    \draw[temp]         
    (\prong[-1,1])
    to[out=\s*90,in=\angle] (join) 
    to[out=\s*90, in=-\s*90] (exit);
}

\newcommandx{\strJoin}[7][1=,2=,4=0.5,5=0.5,7=,usedefault]{
    \drawFork[#1][#2]{#3}[-#4][-#5][#7]
    \coordinate (#6) at (exit);
}

\newcommandx{\strSplit}[8][1=,2=,4=0.5,5=0.5,7=1,8=]{
    \getNodes{#6}{\splitNodes}
    \tikzmath{
        \s=#4+#5;
    }
    \splitTargets{#3}{\splitNodes}{#7}{\s}
    \drawFork[#1][#2]{#6}[#5][#4][#8]
}

\newcommandx{\strMod}[8][1=,2=,3=,5=0.5,7=0.5,usedefault]{
    \coordinate[shift={(0,-#5)}] (joint) at (#4);

    \listLengthX{#1,#2}{\lflen}
    \listLengthX{#1,#3}{\rtlen}

    \tikzmath{
        \TCWidth=10+symWidth(0,\lflen,\rtlen);    
    }
    \draw (joint) node[mod={\TCWidth}](joint){#6};
    
    \ifdim #5 pt > 0 pt
        \coordinate (entry) at (#4);
        \xdef\modPathL{
            (entry)--(entry|-joint.north)-|(joint.south west)
        }
        \xdef\modPathR{
            (entry)--(entry|-joint.north)-|(joint.south east)
        }
    \else
        \coordinate (entry) at (joint.north);
        \xdef\modPathL{
            (joint.north)-|(joint.south west)
        }
        \xdef\modPathR{
            (joint.north)-|(joint.south east)
        }
    \fi

    \ifdim #7 pt > 0 pt
        \coordinate[shift={(0,-#7)}] (exit) at (joint);
        \coordinate (#8) at (exit);
        \xdef\modPathL{\modPathL-|(exit)}
        \xdef\modPathR{\modPathR-|(exit)}
    \else
        \coordinate (exit) at (joint.south);
        \xdef\modPathL{\modPathL-|(exit)}
        \xdef\modPathR{\modPathR-|(exit)}
    \fi

    \highlightSide{#1}{#2}{
        \modPathL
    }{
        --++(0,-3)
        -|($(joint.west)+(-3,0)$)
        |-($(entry)+(0,3)$)
        --(entry)
    }
    \draw[NTDash, rounded corners=0.5,line width=\lineWidth] \modPathL;

    \reverseList{#3}{\rNTCol}
    \highlightSide{#1}{\rNTCol}{
        \modPathR
    }{
        --++(0,-3)
        -|($(joint.east)+(3,0)$)
        |-($(entry)+(0,3)$)
        --(entry)
    }
    \draw[NTDash, rounded corners=0.5, line width=\lineWidth] \modPathR;
}

\newcommandx{\strTwoCell}[8][1=,3=0.5,5=0,6=0.5,8=,usedefault]{

    \readlist*\lst{#2}
    \readlist*\outlst{#7}
    \getNodes{#2}{\nodes}
    \joinTarget{\nodes}{-#3}{joint}[\nodesWidth]

    \listLengthX{\lst[1,2],\lst[1,3]}{\lflen}
    \listLengthX{\lst[-1,2],\lst[-1,4]}{\rtlen}
    \listLengthX{\outlst[1,2],\lst[1,3]}{\lflenII}
    \listLengthX{\outlst[-1,2],\lst[-1,4]}{\rtlenII}

    \tikzmath{
        \lshft=TCLShft(\lflen);
        \rshft=-TCLShft(\rtlen);
        \TCWidth=10+max(symWidth(\nodesWidth,\lflen,\rtlen),symWidth(((#5)*28.45276),\lflenII, \rtlenII));
    }
    \draw (joint) node[twocellX={\TCWidth}](joint){#4};

    \foreach \i in {1,...,\lstlen}{
        \ifnum \i=1 then
            \coordinate[xshift=\lshft] (hlcursor) at (\lst[1,1]);
            \clipS{hlcursor}{joint.west}
            \highlightSide{#1}{}{
                (hlcursor)
                --(hlcursor|-joint.north)
                -|(joint.west)
            }{\leftClip}
        \fi

        \ifnum \i=\lstlen then
            \coordinate[xshift=\rshft] (hlcursor) at (\lst[-1,1]);
            \clipS{hlcursor}{joint.east}
            \highlightSide{#1}{}{
                (hlcursor)
                --(hlcursor|-joint.north)
                -|(joint.east)
            }{\rightClip}

            \edef\fCol{\lst[-1,2]}
            \edef\lntCol{\lst[-1,3]}
            \edef\rntCol{\lst[-1,4]}

            \strIdTo{\lst[-1,1]}{joint}[\fCol][\lntCol][\rntCol]
        \else
            \listLengthX{\lst[\i,2],\lst[\i,4]}{\rtlen}
            \listLengthX{\lst[\i+1,2],\lst[\i+1,3]}{\lflen}

            \tikzmath{
                \rshft=-TCLShft(\rtlen);
                \lshft=TCLShft(\lflen);
            }

            \coordinate[xshift=\rshft] (hlcursor) at (\lst[\i,1]);
            \coordinate[xshift=\lshft] (hlcursor2) at (\lst[\i+1,1]);

            \highlightSide{#1}{}{
                (hlcursor)
                --(hlcursor|-joint.north)
                -|(hlcursor2)
            }{--++(0,3)
            -| (hlcursor)}

            \edef\fCol{\lst[\i,2]}
            \edef\lntCol{\lst[\i,3]}
            \edef\rntCol{\lst[\i,4]}

            \strIdTo{\lst[\i,1]}{joint}[\fCol][\lntCol][\rntCol]
        \fi
    }    

    \getNodes{#7}{\outNodes}
    \splitTargets{joint}{\outNodes}{#5}{#6}

    \foreach \i in {1,...,\outlstlen}{
        \ifnum \i=1 then
            \listLengthX{\outlst[1,2],\outlst[1,3]}{\lflen}

            \tikzmath{
                \lshft=TCLShft(\lflen);
            }
            \coordinate[xshift=\lshft] (hlcursor) at (\outlst[1,1]);
            \clipS{hlcursor}{joint.west}
            \highlightSide{#1}{}{
                (hlcursor)
                --(hlcursor|-joint.south)
                -|(joint.west)
            }{\leftClip}
        \fi

        \ifnum \i=\outlstlen then
            \listLengthX{\outlst[-1,2],\outlst[-1,4]}{\rtlen}

            \tikzmath{
                \rshft=-TCLShft(\rtlen);
            }
            \coordinate[xshift=\rshft] (hlcursor) at (\outlst[-1,1]);
            \clipS{hlcursor}{joint.east}
            \highlightSide{#1}{}{
                (hlcursor)
                --(hlcursor|-joint.south)
                -|(joint.east)
            }{\rightClip}

            \edef\fCol{\outlst[-1,2]}
            \edef\lntCol{\outlst[-1,3]}
            \edef\rntCol{\outlst[-1,4]}

            \drawS{joint-|\outlst[\i,1]}{\outlst[\i,1]}[\fCol][\lntCol][\rntCol]
        \else
            \listLengthX{\outlst[\i,2],\outlst[\i,4]}{\rtlen}
            \listLengthX{\outlst[\i+1,2],\outlst[\i+1,3]}{\lflen}

            \tikzmath{
                \rshft=-TCLShft(\rtlen);
                \lshft=TCLShft(\lflen);
            }

            \coordinate[xshift=\rshft] (hlcursor) at (\outlst[\i,1]);
            \coordinate[xshift=\lshft] (hlcursor2) at (\outlst[\i+1,1]);

            \highlightSide{#1}{}{
                (hlcursor)
                --(hlcursor|-joint.south)
                -|(hlcursor2)
            }{--++(0,-3)
                -| (hlcursor)
            }

            \edef\fCol{\outlst[\i,2]}
            \edef\lntCol{\outlst[\i,3]}
            \edef\rntCol{\outlst[\i,4]}

            \drawS{joint-|\outlst[\i,1]}{\outlst[\i,1]}[\fCol][\lntCol][\rntCol]
        \fi
    }
    \draw (joint) node[twocellX={\TCWidth}](joint){#4};
}

\newcommandx{\idCup}[2]{
    \xdef\colours{#2}        
    \listLength{\colours}{\colourslen}
    \foreach \c [count=\k] in \colours{
        \tikzmath{
            real \rad;
            \rad=idRad(\colourslen-\k+1);
        }
        \coordinate[xshift=-\rad pt
        ] (temp) at (#1);
        \draw[
            line width=\hlWidth pt,
            color=\c
        ]
        (temp) arc(180:360:\rad pt);
    }
}

\newcommandx{\idCap}[2]{
    \xdef\colours{#2}        
    \listLength{\colours}{\colourslen}
    \foreach \c [count=\k] in \colours{
        \tikzmath{
            real \rad;
            \rad=idRad(\colourslen-\k+1);
        }
        \coordinate[xshift=\rad pt
        ] (temp) at (#1);
        \draw[
            line width=\hlWidth pt,
            color=\c
        ]
        (temp) arc(0:180:\rad pt);
    }
}

\newcommandx{\strKey}[1]{
    \readlist*\keyIn{#1}
    \coordinate (currentNode) at (current bounding box.north);
    \begin{scope}[local bounding box=outer box]
        \foreach \i in {\keyInlen,...,1}{
            \edef\nodeName{\keyIn[\i,1]}
            \edef\nodeColour{\keyIn[\i,2]}
        
            \coordinate (temp) at (currentNode);
            \coordinate[shift={(0,0.5)}] (currentNode) at (temp);

            \draw (currentNode) node[label, anchor=east](temp){\nodeName};
            \draw[color=\nodeColour, line width=\hlWidth] (currentNode) --++(0.3,0);
        }   
    \end{scope}
    \coordinate[shift={(0.1,0)}] (topRight) at (outer box.north east);
    \draw[black] (outer box.south west) rectangle (topRight);
}

\newcommandx{\strKeyX}[1]{
    \readlist*\keyIn{#1}
    \tikzset{
        key/.pic={    
            \coordinate (currentNode) at (0,0);
            \begin{scope}[local bounding box=outer box]
                \foreach \i in {\keyInlen,...,1}{
                    \edef\nodeName{\keyIn[\i,1]}
                    \edef\nodeColour{\keyIn[\i,2]}
        
                    \coordinate (temp) at (currentNode);
                    \coordinate[shift={(0,0.5)}] (currentNode) at (temp);

                    \draw (currentNode) node[label, anchor=east](temp){\nodeName};
                    \draw[color=\nodeColour, line width=\hlWidth] (currentNode) --++(0.3,0);
                }  
            \end{scope}
            \coordinate[shift={(0.1cm,0)}] (topRight) at (outer box.north east);
            \draw[black] (outer box.south west) rectangle (topRight);
        }
    }
    \draw (current bounding box.north) pic{key};
}

\newcommandx{\strKeyXX}[1]{
    \readlist*\keyIn{#1}
    \draw (current bounding box.north) node[label, anchor=south]{
            \begin{tikzpicture}
                \coordinate (currentNode) at (0,0);
                \begin{scope}[local bounding box=outer box]
                    \foreach \i in {\keyInlen,...,1}{
                        \edef\nodeName{\keyIn[\i,1]}
                        \edef\nodeColour{\keyIn[\i,2]}
            
                        \coordinate (temp) at (currentNode);
                        \coordinate[shift={(0,0.5)}] (currentNode) at (temp);
    
                        \draw (currentNode) node[label, anchor=east](temp){\nodeName};
                        \draw[color=\nodeColour, line width=\hlWidth] (currentNode) --++(0.3,0);
                    }  
                \end{scope}
                \coordinate[shift={(0.1cm,0)}] (topRight) at (outer box.north east);
                \draw[black] (outer box.south west) rectangle (topRight);
            \end{tikzpicture}
    };
}

\newcommandx{\strLabel}[6][3=,4=,5=,6=1]{
    \tikzmath{
        \height=0.5*#6;
    }
    \coordinate[shift={(0,-\height)}] (labelTemp) at (#1);
    \strId{#1}[#3][#4][#5][#6]
    \draw (labelTemp) node[label]{#2};    
}

\newcommandx{\highlightSideX}[6][6=0]{
    \removeEmpty{#1}{\coloursI}        
    \listLength{\coloursI}{\colourslenI}
    \removeEmpty{#2}{\coloursII}
    \listLength{\coloursII}{\colourslenII}
    \removeEmpty{#3}{\coloursIII}
    \listLength{\coloursIII}{\colourslenIII}

    \halfPath[
        double,
        double=white,
        draw=white, 
        double distance=hlDblSep(#6), 
        line width=\hlWidth pt,
        rounded corners=0.5
    ]{#4}{#5}

    \foreach \c [count=\k] in \coloursI{
        \ifthenelse{\equal{\c}{}}{}{
        \halfPath[
            double,
            double=white,
            draw=\c, 
            double distance=hlDblSep(\colourslenI+\colourslenII+\colourslenIII+1-\k), 
            line width=\hlWidth pt,
            rounded corners=0.5
        ]{#4}{#5}
        }
    }
    \foreach \c [count=\k] in \coloursII{
        \ifthenelse{\equal{\c}{}}{}{
        \halfPath[
            double,
            double=white,
            draw=white, 
            double distance=(hlDblSep(\colourslenII+\colourslenIII+1-\k)), 
            line width=\hlWidth pt,
            rounded corners=0.5,
        ]{#4}{#5}
        \halfPath[
            double,
            double=white,
            draw=\c, 
            double distance=(hlDblSep(\colourslenII+\colourslenIII+1-\k)), 
            line width=\hlWidth pt,
            rounded corners=0.5,
            NTDash
        ]{#4}{#5}
        }
    }
    \foreach \c [count=\k] in \coloursIII{
        \ifthenelse{\equal{\c}{}}{}{
        \halfPath[
            double,
            double=white,
            draw=\c, 
            double distance=hlDblSep(\colourslenIII+1-\k), 
            line width=\hlWidth pt,
            rounded corners=0.5
            ]{#4}{#5}
        }
    }
}

\newcommandx{\highlightX}[9][8=,9=0;0]{
    \xdef\coloursTempI{#2}
    \xdef\coloursTempII{#3}
    \reverseList{#4}{\coloursTempIII}    
    \xdef\coloursTempIV{#5}

    \readlist*\whiteSpace{#9}
    \xdef\whiteSpaceI{\whiteSpace[1,1]}
    \xdef\whiteSpaceII{\whiteSpace[2,1]}
    \highlightSideX{\coloursTempI}{\coloursTempII}{\coloursTempIV}{#1}{#6}[\whiteSpaceI]
    \highlightSideX{\coloursTempI}{\coloursTempIII}{\coloursTempIV}{#1}{#7}[\whiteSpaceII]
    
    \listLengthX{#3}{\ntleftlen}
    \listLengthX{#4}{\ntrightlen}
    \tikzmath{
        integer \nt;
        \nt=\ntleftlen+\ntrightlen;
    }
    \ifnum \nt>0
        \draw[line width=\lineWidth,NTDash,#8] #1;
    \else
        \draw[line width=\lineWidth,#8] #1;
    \fi
}

\newcommandx{\highlightL}[6][5=,6=0;0]{
    \xdef\bfdInput{#2////}

    \readlist*\hlCol{\bfdInput}

    \xdef\colI{\hlCol[1,1]}
    \xdef\colII{\hlCol[1,2]}
    \xdef\colIII{\hlCol[1,3]}
    \xdef\colIV{\hlCol[1,4]}
    \highlightX{#1}{\colI}{\colII}{\colIII}{\colIV}{#3}{#4}[#5][#6]
}

\newcommandx{\drawSX}[5][4=,5=0;0]{
    \clipS{#1}{#2}
    \highlightL{(#1) to[out=-90, in=90] (#2)}{#3}{\leftClip}{\rightClip}[#4][#5]
}

\newcommandx{\strSXX}[3][3=]{
    \coordinate[shift={(#2)}] (temp) at (#1);
    \drawSX{#1}{temp}{#3}
    \coordinate (#1) at (temp);
}

\newcommandx{\drawIdToX}[4][3=,4=,usedefault]{
    \clipS{#1}{#2}
    \highlightL{(#1)--(#1|-#2)}{#3}{\leftClip}{\rightClip}[#4]
}

\newcommandx{\strIdToX}[4][3=,4=,usedefault]{
    \drawIdToX{#1}{#2}[#3][#4]
    \coordinate (#1) at (#1|-#2);
}

\newcommandx{\strIdX}[4][2=,3=1,4=,usedefault]{
    \coordinate[shift={(0,-#3)}] (temp) at (#1);
    \strIdToX{#1}{temp}[#2][#4]
}

\newcommandx{\strBraidXX}[6][3=,4=,5=1,6=;,usedefault]{
    \joinTarget{#1;#2}{0}{topLevel}
    \readlist*\options{#6}
    \xdef\optionI{\options[1,1]}
    \xdef\optionII{\options[2,1]}

    \coordinate[shift={(0,-#5)}] (overBot) at (#2|-topLevel);
    \coordinate[shift={(0,-#5)}] (undBot) at (#1|-topLevel);

    \gettikzxy{(#1)}{\overX}{\overY}
    \gettikzxy{(#2)}{\underX}{\underY}

    \ifdim \overX<\underX
        \xdef\overleft{1}
    \else
        \xdef\overleft{0}
    \fi

    \readlist*\undCol{#4////}
    \readlist*\overCol{#3////}

    \xdef\undNat{0}
    
    \ifnonempty{\undCol[1,2]}{
        \xdef\undNat{1}
    }
    \ifnonempty{\undCol[1,3]}{
        \xdef\undNat{1}
    }

    \ifnum \overleft=1
        \xdef\newUndFunCol{\overCol[1,1],\overCol[1,2]}        
        \reverseList{\undCol[1,3]}{\newOverFunCol}
        \ifnum \undNat=0
            \xdef\newUndFunCol{\newUndFunCol,\overCol[1,4]}        
        \fi
    \else
        \reverseList{\overCol[1,3]}{\newUndFunCol}
        \xdef\newUndFunCol{\overCol[1,1],\newUndFunCol}
        \xdef\newOverFunCol{\undCol[1,2]}
        \ifnum \undNat=0
            \xdef\newUndFunCol{\newUndFunCol,\overCol[1,4]}        
        \fi
    \fi

    \xdef\newUndCol{
        \newUndFunCol/\undCol[1,2]/\undCol[1,3]/\undCol[1,4]
    }
    \xdef\newOverCol{
        \overCol[1,1]/\overCol[1,2]/\overCol[1,3]/\newOverFunCol
    }

    \xdef\overPath{
        (#1) to[out=-90, in=90] (overBot)
    }
    \xdef\underPath{
        (#2) to[out=-90, in=90] (undBot)
    }

    \ifnum \overleft=1
        \clipS{#1}{overBot}[\botClip][\topClip]
    \else
        \clipS{#1}{overBot}[\topClip][\botClip]
    \fi

    \coloursLength{\overCol[1,1]/\overCol[1,2]/\overCol[1,4]}[\leftTop]
    \coloursLength{\overCol[1,1]/\overCol[1,3]/\overCol[1,4]}[\rightTop]
    \coloursLength{\overCol[1,1]/\overCol[1,2]/\newOverFunCol}[\leftBottom]
    \coloursLength{\overCol[1,1]/\overCol[1,3]/\newOverFunCol}[\rightBottom]

    \tikzmath{
        \leftWhiteSpace=1+max(\leftTop,\leftBottom);
        \rightWhiteSpace=1+max(\rightTop,\rightBottom);
    }

    \begin{scope}
        \begin{pgfinterruptboundingbox}
            \clip \overPath\topClip;
        \end{pgfinterruptboundingbox}
        \drawSX{#2}{undBot}{#4}[\optionII]
    \end{scope}
    \begin{scope}
        \begin{pgfinterruptboundingbox}
            \clip \overPath\botClip;
        \end{pgfinterruptboundingbox}
        \drawSX{#2}{undBot}{\newUndCol}[\optionII]
    \end{scope}
    
    \ifnum \undNat=0
        \drawSX{#1}{overBot}{white,#3}[\optionI]
    \else
        \ifnum \overleft=1
            \clipS{#2}{undBot}[\topClip][\botClip]
        \else
            \clipS{#2}{undBot}[\botClip][\topClip]
        \fi
        \begin{scope}
            \begin{pgfinterruptboundingbox}
                \clip \underPath\topClip;
            \end{pgfinterruptboundingbox}
            \drawSX{#1}{overBot}{#3}[\optionII][\leftWhiteSpace;\rightWhiteSpace]
        \end{scope}
        \begin{scope}
            \begin{pgfinterruptboundingbox}
                \clip \underPath\botClip;
            \end{pgfinterruptboundingbox}
            \drawSX{#1}{overBot}{\newOverCol}[\optionII][\leftWhiteSpace;\rightWhiteSpace]
        \end{scope}
    \fi
    \coordinate (#1) at (overBot);
    \coordinate (#2) at (undBot);
}

\newcommandx{\drawSCircX}[4][3=,4=0.5]{
    \joinTarget{#1;#2}{#4}{joinPoint}[][#4]
    \tikzmath{
        \a=sign(#4)*90;
        \shft=(3*(sign(#4)));
    }
    \highlightL{
        (#1) to[out=\a, in=180] (joinPoint) to[out=0, in=\a] (#2)
    }{#3}{
        |-++(3,-\shft)
        |-($(joinPoint)+(0,\shft)$)
        -|($(#1)+(-3,-\shft)$)
        -|(#1)
    }{
        --++(0,-\shft)
        -|(#1)
    }
}

\newcommandx{\strCupX}[4][3=,4=0.5,usedefault]{
    \drawSCircX{#1}{#2}[#3][-#4]
}

\newcommandx{\strCapX}[4][3=,4=0.5,usedefault]{
    \drawSCircX{#1}{#2}[#3][#4]
}

\newcommandx{\strCol}[4][4=1]{
    \coordinate[shift={(0,-#4*0.5)}] (temp) at (#1);
    \drawIdToX{#1}{temp}[#2]
    \coordinate[shift={(0,-#4*0.5)}] (#1) at (temp);
    \drawIdToX{temp}{#1}[#3]
    
    \readlist*\inLst{#2////}
    \readlist*\outLst{#3////}

    \listLengthX{\inLst[1,1],\inLst[1,2],\inLst[1,4]}{\lflen}
    \listLengthX{\inLst[1,1],\inLst[1,3],\inLst[1,4]}{\rtlen}
    \listLengthX{\outLst[1,1],\outLst[1,2],\outLst[1,4]}{\lflenII}
    \listLengthX{\outLst[1,1],\outLst[1,3],\outLst[1,4]}{\rtlenII}

    \tikzmath{
        \CCWidth=3+max(symWidth(0,\lflen,\rtlen),symWidth(0,\lflenII,\rtlenII));
    }
    \draw (temp) node[colChange={\CCWidth}]{};
}

\newcommandx{\strModX}[6][4=,5=0.5,6=0.5usedefault]{
    \coordinate[shift={(0,-#5)}] (joint) at (#1);

    \readlist*\modCol{#4////}
    \listLengthX{\modCol[1,1],\modCol[1,2],\modCol[1,4]}{\lflen}
    \listLengthX{\modCol[1,1],\modCol[1,3],\modCol[1,4]}{\rtlen}

    \tikzmath{
        \TCWidth=10+symWidth(0,\lflen,\rtlen);    
    }
    \draw (joint) node[mod={\TCWidth}](joint){#2};
    
    \ifdim #5 pt > 0 pt
        \coordinate (entry) at (#1);
        \xdef\modPathL{
            (entry)--(entry|-joint.north)-|(joint.south west)
        }
        \xdef\modPathR{
            (entry)--(entry|-joint.north)-|(joint.south east)
        }
    \else
        \coordinate (entry) at (joint.north);
        \xdef\modPathL{
            (joint.north)-|(joint.south west)
        }
        \xdef\modPathR{
            (joint.north)-|(joint.south east)
        }
    \fi

    \ifdim #6 pt > 0 pt
        \coordinate[shift={(0,-#6)}] (exit) at (joint);
        \coordinate (#3) at (exit);
        \xdef\modPathL{\modPathL-|(exit)}
        \xdef\modPathR{\modPathR-|(exit)}
    \else
        \coordinate (exit) at (joint.south);
        \xdef\modPathL{\modPathL-|(exit)}
        \xdef\modPathR{\modPathR-|(exit)}
    \fi

    \highlightSideX{\modCol[1,1]}{\modCol[1,2]}{\modCol[1,4]}{
        \modPathL
    }{
        --++(0,-3)
        -|($(joint.west)+(-3,0)$)
        |-($(entry)+(0,3)$)
        --(entry)
    }
    \draw[NTDash, rounded corners=0.5,line width=\lineWidth] \modPathL;

    \reverseList{\modCol[1,3]}{\rNTCol}
    \highlightSideX{\modCol[1,1]}{\rNTCol}{\modCol[1,4]}{
        \modPathR
    }{
        --++(0,-3)
        -|($(joint.east)+(3,0)$)
        |-($(entry)+(0,3)$)
        --(entry)
    }
    \draw[NTDash, rounded corners=0.5, line width=\lineWidth] \modPathR;
}

\newcommandx{\strTwoCellX}[7][4=,5=0,6=0.5,7=0.5,usedefault]{

    \bufferlist{\bufflst}{#1}
    \readlist*\lst{\bufflst}
    \bufferlist{\bufflst}{#3}
    \readlist*\outlst{\bufflst}
    \getNodes{#1}{\nodes}
    \joinTarget{\nodes}{-#6}{joint}[\nodesWidth]

    \listLengthX{\lst[1,2],\lst[1,3],\lst[1,5]}{\lflen}
    \listLengthX{\lst[-1,2],\lst[-1,4],\lst[-1,5]}{\rtlen}
    \listLengthX{\outlst[1,2],\outlst[1,3],\outlst[1,5]}{\lflenII}
    \listLengthX{\outlst[-1,2],\outlst[-1,4],\outlst[1,-5]}{\rtlenII}

    \tikzmath{
        \lshft=TCLShft(\lflen);
        \rshft=-TCLShft(\rtlen);
        \TCWidth=10+max(symWidth(\nodesWidth,\lflen,\rtlen),symWidth(((#5)*28.45276),\lflenII, \rtlenII));
    }
    \ifthenelse{\equal{#2}{-}}{
        \draw (joint) node[colChange={\TCWidth-7}](joint){};
    }{

        \draw (joint) node[twocellX={\TCWidth}](joint){#2};
    }

    \foreach \i in {1,...,\lstlen}{
        \ifnum \i=1
            \coordinate[xshift=\lshft] (hlcursor) at (\lst[1,1]);
            \clipS{hlcursor}{joint.west}
            \highlightSide{#4}{}{
                (hlcursor)
                --(hlcursor|-joint.north)
                -|(joint.west)
            }{\leftClip}
        \fi

        \ifnum \i=\lstlen
            \coordinate[xshift=\rshft] (hlcursor) at (\lst[-1,1]);
            \clipS{hlcursor}{joint.east}
            \highlightSide{#4}{}{
                (hlcursor)
                --(hlcursor|-joint.north)
                -|(joint.east)
            }{\rightClip}

            \xdef\fCol{\lst[-1,2]}
            \xdef\lntCol{\lst[-1,3]}
            \xdef\rntCol{\lst[-1,4]}
            \xdef\ifCol{\lst[-1,5]}

            \strIdToX{\lst[-1,1]}{joint}[\fCol/\lntCol/\rntCol/\ifCol]
        \else
            \listLengthX{\lst[\i,2],\lst[\i,4],\lst[\i,5]}{\rtlen}
            \listLengthX{\lst[\i+1,2],\lst[\i+1,3],\lst[\i+1,5]}{\lflen}

            \tikzmath{
                \rshft=-TCLShft(\rtlen);
                \lshft=TCLShft(\lflen);
            }

            \coordinate[xshift=\rshft] (hlcursor) at (\lst[\i,1]);
            \coordinate[xshift=\lshft] (hlcursor2) at (\lst[\i+1,1]);

            \highlightSide{#4}{}{
                (hlcursor)
                --(hlcursor|-joint.north)
                -|(hlcursor2)
            }{--++(0,3)
            -| (hlcursor)}

            \xdef\fCol{\lst[\i,2]}
            \xdef\lntCol{\lst[\i,3]}
            \xdef\rntCol{\lst[\i,4]}
            \xdef\ifCol{\lst[\i,5]}

            \strIdToX{\lst[\i,1]}{joint}[\fCol/\lntCol/\rntCol/\ifCol]
        \fi
    }    

    \getNodes{#3}{\outNodes}
    \splitTargets{joint}{\outNodes}{#5}{#7}

    \foreach \i in {1,...,\outlstlen}{
        \ifnum \i=1 then
            \listLengthX{\outlst[1,2],\outlst[1,3],\outlst[1,5]}{\lflen}

            \tikzmath{
                \lshft=TCLShft(\lflen);
            }
            \coordinate[xshift=\lshft] (hlcursor) at (\outlst[1,1]);
            \clipS{hlcursor}{joint.west}
            \highlightSide{#4}{}{
                (hlcursor)
                --(hlcursor|-joint.south)
                -|(joint.west)
            }{\leftClip}
        \fi

        \ifnum \i=\outlstlen then
            \listLengthX{\outlst[-1,2],\outlst[-1,4],\outlst[-1,5]}{\rtlen}

            \tikzmath{
                \rshft=-TCLShft(\rtlen);
            }
            \coordinate[xshift=\rshft] (hlcursor) at (\outlst[-1,1]);
            \clipS{hlcursor}{joint.east}
            \highlightSide{#4}{}{
                (hlcursor)
                --(hlcursor|-joint.south)
                -|(joint.east)
            }{\rightClip}

            \xdef\fCol{\outlst[-1,2]}
            \xdef\lntCol{\outlst[-1,3]}
            \xdef\rntCol{\outlst[-1,4]}
            \xdef\ifCol{\outlst[-1,5]}

            \drawSX{joint-|\outlst[\i,1]}{\outlst[\i,1]}{\fCol/\lntCol/\rntCol/\ifCol}
        \else
            \listLengthX{\outlst[\i,2],\outlst[\i,4],\outlst[\i,5]}{\rtlen}
            \listLengthX{\outlst[\i+1,2],\outlst[\i+1,3],\outlst[\i+1,5]}{\lflen}

            \tikzmath{
                \rshft=-TCLShft(\rtlen);
                \lshft=TCLShft(\lflen);
            }

            \coordinate[xshift=\rshft] (hlcursor) at (\outlst[\i,1]);
            \coordinate[xshift=\lshft] (hlcursor2) at (\outlst[\i+1,1]);

            \highlightSide{#4}{}{
                (hlcursor)
                --(hlcursor|-joint.south)
                -|(hlcursor2)
            }{--++(0,-3)
                -| (hlcursor)
            }

            \xdef\fCol{\outlst[\i,2]}
            \xdef\lntCol{\outlst[\i,3]}
            \xdef\rntCol{\outlst[\i,4]}
            \xdef\ifCol{\outlst[\i,5]}

            \drawSX{joint-|\outlst[\i,1]}{\outlst[\i,1]}{\fCol/\lntCol/\rntCol/\ifCol}
        \fi
    }
    \ifthenelse{\equal{#2}{-}}{
        \draw (joint) node[colChange={\TCWidth-7}](joint){};
    }{
        \draw (joint) node[twocellX={\TCWidth}](joint){#2};
    }
}

\newcommandx{\strLabelX}[5][3=,4=1,5=,usedefault]{
    \tikzmath{
        \height=0.5*#4;
    }
    \coordinate[shift={(0,-\height)}] (labelTemp) at (#1);
    \strIdX{#1}[#3][#4][#5]
    \draw (labelTemp) node[label]{#2};    
}

\newcommandx{\bufferlist}[2]{
    \readlist*\buffIn{#2}
    \xdef\newList{}
    \foreach \i in {1,...,\buffInlen}{
        \ifnum \i=1
            \xdef\newList{\buffIn[\i]/////}
        \else
            \xdef\newList{\newList;\buffIn[\i]/////}
        \fi
    }
    \xdef#1{\newList}
}

\newcommandx{\drawForkX}[4][2=-0.5,3=-0.5,4=]{
    \getNodes{#1}{\nodes}
    \joinTarget{\nodes}{#2}{join}[][#2]
    \coordinate[shift={(0,#3)}] (exit) at (join);
    
    \bufferlist{\prongstring}{#1}
    \readlist*\prong{\prongstring}

    \xdef\natTrans{0}

    \foreach \i in {1,...,\pronglen}{
        \ifnonempty{\prong[\i,3]}{
            \xdef\natTrans{1}
        }
    }

    \foreach \i in {1,...,\pronglen}{
        \ifnonempty{\prong[\i,4]}{
            \xdef\natTrans{1}
        }
    }
    \ifnum \natTrans>0
        \tikzset{
            temp/.style={
                NTDash,#4,line width=\lineWidth
            }
        }
    \else
        \tikzset{
            temp/.style={
                line width=\lineWidth,#4
            }
        }
    \fi

    \tikzmath{
        \s=sign(#2);
        \shft=(-3)*(sign(#2));
    }

    \pgfmathsetmacro\angle{
        -\s*joinAngle(1,\pronglen)
    }

    \clipS{\prong[1,1]}{exit}
    
    \highlightSideX{\prong[1,2]}{\prong[1,3]}{\prong[1,5]}   
    {
        (\prong[1,1])
        to[out=\s*90,in=\angle] (join) 
        to[out=\s*90, in=-\s*90] (exit)
    }{
        \leftClip
    }

    \draw[temp]      
    (\prong[1,1])
    to[out=\s*90,in=\angle] (join) 
    to[out=\s*90, in=-\s*90] (exit);

    \ifnum \pronglen>1
        \tikzmath{\finish=\pronglen-1;}
        \foreach \i in {1,...,\finish}{

            \xdef\currentNode{\prong[\i,1]}
            \xdef\nextNode{\prong[\i+1,1]}
            \reverseList{\prong[\i,4]}{\natCol}
            \pgfmathsetmacro\inAngle{
                -\s*joinAngle(\i,\pronglen)
            }
            \pgfmathsetmacro\outAngle{
                -\s*joinAngle(\i+1,\pronglen)
            }
            \highlightSideX{\prong[\i,2]}{\natCol}{\prong[\i,5]}{
                (\currentNode)
                to[out=\s*90,in=\inAngle] (join) 
                to[out=\outAngle, in=\s*90] (\nextNode)
            }{
                --++(0,\shft)
                -|(\currentNode)
            }

            \draw[temp]      
            (\currentNode)
            to[out=\s*90,in=\inAngle] (join);
            \draw[temp]
            (\nextNode)
            to[in=\outAngle, out=\s*90] (join);
       }
    \fi

    \reverseList{\prong[-1,4]}{\natCol}

    \pgfmathsetmacro\angle{
        -\s*joinAngle(\pronglen,\pronglen)
    }
    
    \clipS{\prong[-1,1]}{exit}

    \highlightSideX{\prong[-1,2]}{\natCol}{\prong[-1,5]}{
        (\prong[-1,1])
        to[out=\s*90,in=\angle] (join) 
        to[out=\s*90, in=-\s*90] (exit)
    }{\rightClip}

    \draw[temp]         
    (\prong[-1,1])
    to[out=\s*90,in=\angle] (join) 
    to[out=\s*90, in=-\s*90] (exit);
}

\newcommandx{\strJoinX}[5][3=0.5,4=0.5,5=]{
    \drawForkX{#1}[-#3][-#4][#5]
    \coordinate (#2) at (exit);
}

\newcommandx{\strSplitX}[6][3=1,4=0.5,5=0.5,6=]{
    \getNodes{#2}{\splitNodes}
    \tikzmath{
        \s=#4+#5;
    }
    \splitTargets{#1}{\splitNodes}{#3}{\s}
    \drawForkX{#2}[#5][#4][#6]
}

\newcommandx{\strModXX}[7][4=0,5=0.5,6=0.5,7=,usedefault]{
    \getNodes{#1}{\nodes}
    
    \joinTarget{\nodes}{-#5}{joint}[][-1]

    \bufferlist{\inputs}{#1/}
    \readlist*\lst{\inputs}

    \bufferlist{\outputs}{#3/}
    \readlist*\outlst{\outputs}

    \listLengthX{\lst[1,2],\lst[1,3],\lst[1,5]}{\lflen}
    \listLengthX{\lst[-1,2],\lst[-1,4],\lst[-1,5]}{\rtlen}
    \listLengthX{\outlst[1,2],\outlst[1,3],\outlst[1,5]}{\lflenII}
    \listLengthX{\outlst[-1,2],\outlst[-1,4],\outlst[-1,5]}{\rtlenII}

    \tikzmath{
        \lshft=TCLShft(\lflen);
        \rshft=-TCLShft(\rtlen);
        \TCWidth=5+max(symWidth(\nodesWidth,\lflen,\rtlen),symWidth(((#4)*28.45276),\lflenII, \rtlenII));
    }

    \draw (joint) node[mod2={\TCWidth}](joint){#2};
    
    \foreach \i in {1,...,\lstlen}{
        \ifnum \i<\lstlen
            \reverseList{\lst[\i,4]}{\natCol}
            \path (\lst[\i,1]) to node[pos=0.5,blank](halfway){} (\lst[\i+1,1]);
            
            \highlightSideX{\lst[\i,2]}{\natCol}{\lst[\i,5]}
            {
                (\lst[\i,1]) |- (halfway|-joint.north)
            }{
                |-(\lst[\i,1])
            }

            \highlightSideX{\lst[\i+1,2]}{\lst[\i+1,3]}{\lst[\i+1,5]}
            {
                (\lst[\i+1,1]) |- (halfway|-joint.north)
            }{
                |-(\lst[\i+1,1])
            }        
            \draw[NTDash] (\lst[\i,1]) |- (halfway|-joint.north);
            \draw[NTDash] (\lst[\i+1,1]) |- (halfway|-joint.north);
        \fi
    }
    
    \getNodes{#3}{\outNodes}
    \splitTargets{joint}{\outNodes}{#4}{#6}

    \foreach \i in {1,...,\outlstlen}{
        \ifnum \i<\outlstlen
            \reverseList{\outlst[\i,4]}{\natCol}
            \path (\outlst[\i,1]) to node[pos=0.5,blank](halfway){} (\outlst[\i+1,1]);

            \highlightSideX{\outlst[\i,2]}{\natCol}{\outlst[\i,5]}
            {
                (\outlst[\i,1]) |- (halfway|-joint.south)
            }{
                |-(\outlst[\i,1])
            }

            \highlightSideX{\outlst[\i+1,2]}{\outlst[\i+1,3]}{\outlst[\i+1,5]}
            {
                (\outlst[\i+1,1]) |- (halfway|-joint.south)
            }{
                |-(\outlst[\i+1,1])
            }        
            \draw[NTDash] (\outlst[\i,1]) |- (halfway|-joint.south);
            \draw[NTDash] (\outlst[\i+1,1]) |- (halfway|-joint.south);  
        \fi
    }

    \ifdim #5pt >0pt
        \xdef\leftEntry{\lst[1,1]}
        \xdef\rightEntry{\lst[-1,1]}
        \xdef\leftEntryColI{\lst[1,2]}
        \xdef\leftEntryColII{\lst[1,3]}
        \xdef\leftEntryColIII{\lst[1,5]}
        \xdef\rightEntryColI{\lst[-1,2]}
        \reverseList{\lst[-1,4]}{\rightEntryColII}
        \xdef\rightEntryColIII{\lst[-1,5]}
    \else
        \xdef\leftEntry{joint.north}
        \xdef\rightEntry{joint.north}
        \xdef\leftEntryColI{\outlst[1,2]}
        \xdef\leftEntryColII{\outlst[1,3]}
        \xdef\leftEntryColIII{\outlst[1,5]}
        \xdef\rightEntryColI{\outlst[-1,2]}
        \reverseList{\outlst[-1,4]}{\rightEntryColII}
        \xdef\rightEntryColIII{\outlst[-1,5]}
    \fi
    \ifnonempty{#3}{
        \xdef\leftExit{\outlst[1,1]}
        \xdef\rightExit{\outlst[-1,1]}
        \xdef\leftExitColI{\outlst[1,2]}
        \xdef\leftExitColII{\outlst[1,3]}
        \xdef\leftExitColIII{\outlst[1,5]}
        \xdef\rightExitColI{\outlst[-1,2]}
        \reverseList{\outlst[-1,4]}{\rightExitColII}
        \xdef\rightExitColIII{\outlst[-1,5]}
    }
    \ifempty{#3}{
        \xdef\leftExit{joint.south}
        \xdef\rightExit{joint.south}
        \xdef\leftExitColI{\lst[1,2]}
        \xdef\leftExitColII{\lst[1,3]}
        \xdef\leftExitColIII{\lst[1,5]}
        \xdef\rightExitColI{\lst[-1,2]}
        \xdef\rightExitColII{\lst[-1,3]}
        \xdef\rightExitColIII{\lst[-1,5]}
    }

    \clipS{\leftEntry}{joint.west}
    \highlightSideX{\leftEntryColI}{\leftEntryColII}{\leftEntryColIII}{
        (\leftEntry)|-(joint.north west)--(joint.west)
    }{
        \leftClip
    }
    \draw[NTDash] (\leftEntry)|-(joint.north west)--(joint.west);

    \clipS{\rightEntry}{joint.east}
    \highlightSideX{\rightEntryColI}{\rightEntryColII}{\rightEntryColIII}{
        (\rightEntry)|-(joint.north east)--(joint.east)
    }{
        \rightClip
    }
    \draw[NTDash] (\rightEntry)|-(joint.north east)--(joint.east);

    \clipS{\leftExit}{joint.west}
    \highlightSideX{\leftExitColI}{\leftExitColII}{\leftExitColIII}{
        (\leftExit)|-(joint.south west)--(joint.west)    
    }{
        \leftClip
    }
    \draw[NTDash] (\leftExit)|-(joint.south west)--(joint.west);    

    \clipS{\rightExit}{joint.east}
    \highlightSideX{\rightExitColI}{\rightExitColII}{\rightExitColIII}{
        (\rightExit)|-(joint.south east)--(joint.east)
    }{
        \rightClip
    }
    \draw[NTDash] (\rightExit)|-(joint.south east)--(joint.east);
}

\newcommandx{\strModXXX}[9][4=0,5=0.5,6=0.5,7=,8=,9=,usedefault]{
    \getNodes{#1}{\nodes}
    
    \joinTarget{\nodes}{-#5}{joint}[][-1]

    \bufferlist{\inputs}{#1/}
    \readlist*\lst{\inputs}

    \bufferlist{\outputs}{#3/}
    \readlist*\outlst{\outputs}

    \listLengthX{\lst[1,2],\lst[1,3],\lst[1,5]}{\lflen}
    \listLengthX{\lst[-1,2],\lst[-1,4],\lst[-1,5]}{\rtlen}
    \listLengthX{\outlst[1,2],\outlst[1,3],\outlst[1,5]}{\lflenII}
    \listLengthX{\outlst[-1,2],\outlst[-1,4],\outlst[-1,5]}{\rtlenII}

    \tikzmath{
        \lshft=TCLShft(\lflen);
        \rshft=-TCLShft(\rtlen);
        \TCWidth=10+max(symWidth(\nodesWidth,\lflen,\rtlen),symWidth(((#4)*28.45276),\lflenII, \rtlenII));
    }

    \draw (joint) node[mod2={\TCWidth}](joint){#2};
    
    \foreach \i in {1,...,\lstlen}{
        \ifnum \i<\lstlen
            \reverseList{\lst[\i,4]}{\natCol}
            \path (\lst[\i,1]) to node[pos=0.5,blank](halfway){} (\lst[\i+1,1]);
            
            \highlightSideX{#8,\lst[\i,2]}{\natCol}{\lst[\i,5]}
            {
                (\lst[\i,1]) |- (halfway|-joint.north)
            }{
                |-(\lst[\i,1])
            }

            \highlightSideX{#8,\lst[\i+1,2]}{\lst[\i+1,3]}{\lst[\i+1,5]}
            {
                (\lst[\i+1,1]) |- (halfway|-joint.north)
            }{
                |-(\lst[\i+1,1])
            }        
            \draw[NTDash] (\lst[\i,1]) |- (halfway|-joint.north);
            \draw[NTDash] (\lst[\i+1,1]) |- (halfway|-joint.north);
        \fi
    }
    
    \getNodes{#3}{\outNodes}
    \splitTargets{joint}{\outNodes}{#4}{#6}

    \foreach \i in {1,...,\outlstlen}{
        \ifnum \i<\outlstlen
            \reverseList{\outlst[\i,4]}{\natCol}
            \path (\outlst[\i,1]) to node[pos=0.5,blank](halfway){} (\outlst[\i+1,1]);

            \highlightSideX{#8,\outlst[\i,2]}{\natCol}{\outlst[\i,5]}
            {
                (\outlst[\i,1]) |- (halfway|-joint.south)
            }{
                |-(\outlst[\i,1])
            }

            \highlightSideX{#8,\outlst[\i+1,2]}{\outlst[\i+1,3]}{\outlst[\i+1,5]}
            {
                (\outlst[\i+1,1]) |- (halfway|-joint.south)
            }{
                |-(\outlst[\i+1,1])
            }        
            \draw[NTDash] (\outlst[\i,1]) |- (halfway|-joint.south);
            \draw[NTDash] (\outlst[\i+1,1]) |- (halfway|-joint.south);  
        \fi
    }

    \ifdim #5pt >0pt
        \xdef\leftEntry{\lst[1,1]}
        \xdef\rightEntry{\lst[-1,1]}
        \xdef\leftEntryColI{\lst[1,2]}
        \xdef\leftEntryColII{\lst[1,3]}
        \xdef\leftEntryColIII{\lst[1,5]}
        \xdef\rightEntryColI{\lst[-1,2]}
        \reverseList{\lst[-1,4]}{\rightEntryColII}
        \xdef\rightEntryColIII{\lst[-1,5]}
    \else
        \xdef\leftEntry{joint.north}
        \xdef\rightEntry{joint.north}
        \xdef\leftEntryColI{\outlst[1,2]}
        \xdef\leftEntryColII{\outlst[1,3]}
        \xdef\leftEntryColIII{\outlst[1,5]}
        \xdef\rightEntryColI{\outlst[-1,2]}
        \reverseList{\outlst[-1,4]}{\rightEntryColII}
        \xdef\rightEntryColIII{\outlst[-1,5]}
    \fi
    \ifnonempty{#3}{
        \xdef\leftExit{\outlst[1,1]}
        \xdef\rightExit{\outlst[-1,1]}
        \xdef\leftExitColI{\outlst[1,2]}
        \xdef\leftExitColII{\outlst[1,3]}
        \xdef\leftExitColIII{\outlst[1,5]}
        \xdef\rightExitColI{\outlst[-1,2]}
        \reverseList{\outlst[-1,4]}{\rightExitColII}
        \xdef\rightExitColIII{\outlst[-1,5]}
    }
    \ifempty{#3}{
        \xdef\leftExit{joint.south}
        \xdef\rightExit{joint.south}
        \xdef\leftExitColI{\lst[1,2]}
        \xdef\leftExitColII{\lst[1,3]}
        \xdef\leftExitColIII{\lst[1,5]}
        \xdef\rightExitColI{\lst[-1,2]}
        \xdef\rightExitColII{\lst[-1,3]}
        \xdef\rightExitColIII{\lst[-1,5]}
    }

    \clipS{\leftEntry}{joint.west}
    \highlightSideX{#8,\leftEntryColI}{\leftEntryColII}{\leftEntryColIII}{
        (\leftEntry)|-(joint.north west)--(joint.west)
    }{
        \leftClip
    }
    \draw[NTDash] (\leftEntry)|-(joint.north west)--(joint.west);

    \clipS{\rightEntry}{joint.east}
    \highlightSideX{#8,\rightEntryColI}{\rightEntryColII}{\rightEntryColIII}{
        (\rightEntry)|-(joint.north east)--(joint.east)
    }{
        \rightClip
    }
    \draw[NTDash] (\rightEntry)|-(joint.north east)--(joint.east);

    \clipS{\leftExit}{joint.west}
    \highlightSideX{#8,\leftExitColI}{\leftExitColII}{\leftExitColIII}{
        (\leftExit)|-(joint.south west)--(joint.west)    
    }{
        \leftClip
    }
    \draw[NTDash] (\leftExit)|-(joint.south west)--(joint.west);    

    \clipS{\rightExit}{joint.east}
    \highlightSideX{#8,\rightExitColI}{\rightExitColII}{\rightExitColIII}{
        (\rightExit)|-(joint.south east)--(joint.east)
    }{
        \rightClip
    }
    \draw[NTDash] (\rightExit)|-(joint.south east)--(joint.east);
}

\newcommandx{\coloursLength}[2][2=\colLen]{
    \readlist*\colLenColours{#1////}
    \listLengthX{\colLenColours[1,1]}{\clI}
    \listLengthX{\colLenColours[1,2]}{\clII}
    \listLengthX{\colLenColours[1,3]}{\clIII}
    \listLengthX{\colLenColours[1,4]}{\clIV}
    \tikzmath{
        integer #2;
        #2=\clI+\clII+\clIII+\clIV;
    }
}

\newcommandx{\strCrossCap}[6][4=,5=,6=0.5]{
    \path (#2) to node[pos=#6,blank](joint){} (#2|-#1);
    \readlist*\capCol{#4////}
    \readlist*\lineCol{#5////}

    \xdef\capColI{\capCol[1,1]}
    \xdef\capColII{\capCol[1,2]}
    \xdef\capColIII{\capCol[1,3]}
    \xdef\capColIV{\capCol[1,4]}

    \xdef\lineColII{\lineCol[1,2]}
    \xdef\lineColIII{\lineCol[1,3]}
    \xdef\lineColIV{\lineCol[1,4]}

    \strIdToX{#2}{joint}[#5]
    \strIdToX{#2}{#3}[\capColIII/\lineColII/\lineColIII/\lineColIV]

    \coloursLength{\capColI/\capColII/\capColIV}[\leftTop]
    \coloursLength{\capColI/\capColIII/\capColIV}[\leftBottom]
    \coloursLength{\capColI/\capColII/\lineColIII///}[\rightTop]
    \coloursLength{\capColI/\capColIII/\lineColIII///}[\rightBottom]

    \tikzmath{
        \topWhiteSpace=1+max(\leftTop,\rightTop);
        \botWhiteSpace=1+max(\leftBottom,\rightBottom);
    }

    \highlightL{
        (#1) to[out=90, in=180] (joint)
    }{
        #4
    }{
        -|++(3,3)
        -|($(#1)+(-3,-3)$)
        -|(#1)
    }{
        --++(3,-3)
        |-(#1)
    }[][\topWhiteSpace;\botWhiteSpace]
    \highlightL{
        (joint) to[out=0, in=90] (#3)
    }{
        \capColI/\capColII/\capColIII/\lineColIII
    }{
        -|++(3,3)
        -|($(#1)+(-3,-3)$)
        |-(#1)
    }{
        --++(-3,-3)
        |-(#1)
    }[][\topWhiteSpace;\botWhiteSpace]
}

\newcommandx{\strCrossCup}[7][4=,5=,6=2,7=0.5]{
    \path (#2) to node[pos=#7,blank](joint){} ++(0,-#6) node[blank](exit){};
    \readlist*\capCol{#4////}
    \readlist*\lineCol{#5////}

    \xdef\capColI{\capCol[1,1]}
    \xdef\capColII{\capCol[1,2]}
    \xdef\capColIII{\capCol[1,3]}
    \xdef\capColIV{\capCol[1,4]}

    \xdef\lineColII{\lineCol[1,2]}
    \xdef\lineColIII{\lineCol[1,3]}
    \xdef\lineColIV{\lineCol[1,4]}

    \strIdToX{#2}{joint}[#5]
    \strIdToX{#2}{exit}[\capColIII/\lineColII/\lineColIII/\lineColIV]

    \coloursLength{\capColI/\capColII/\capColIV}[\leftTop]
    \coloursLength{\capColI/\capColIII/\capColIV}[\leftBottom]
    \coloursLength{\capColI/\capColII/\lineColIII///}[\rightTop]
    \coloursLength{\capColI/\capColIII/\lineColIII///}[\rightBottom]

    \tikzmath{
        \topWhiteSpace=1+max(\leftTop,\rightTop);
        \botWhiteSpace=1+max(\leftBottom,\rightBottom);
    }

    \highlightL{
        (#1) to[out=-90, in=180] (joint)
    }{
        #4
    }{
        -|++(3,-3)
        -|($(#1)+(-3,3)$)
        -|(#1)
    }{
        -|($(joint|-#1)+(1,1)$)
        -|(#1)
    }[][\topWhiteSpace;\botWhiteSpace]
    \highlightL{
        (joint) to[out=0, in=-90] (#3)
    }{
        \capColI/\capColII/\capColIII/\lineColIII
    }{
        |-++(3,3)
        |-($(#1)+(-3,-3)$)
        |-(#1)
    }{
        |-($(#3|-joint)+(-1,1)$)
        |-(#1)
    }[][\topWhiteSpace;\botWhiteSpace]
}

\newcommandx{\strJoinAlt}[6][3=,4=,5=0.5,6=0.5]{    
    \readlist*\leftCol{#3///////}
    \readlist*\rightCol{#4////////}

    \place{joinLeft}{0}[-#5]{#1}
    \path (#1|-joinLeft) to node[pos=0.5,blank](join){} (#2|-joinLeft);
    
    \place{templeft}{0}[-#5-#6]{#1}
    \place{tempright}{0}[-#5-#6]{#2}
    
    \highlightSideX{\leftCol[1,1]}{\leftCol[1,2]}{\leftCol[1,3]}
    {(#1)--(templeft)}{
        --++(0,-1)
        -|($(#1)+(-1,1)$)
        -|(#1)
    }    
    
    \highlightSideX{\rightCol[1,4]}{\rightCol[1,5]}{\rightCol[1,6]}
    {(#2)--(tempright)}{
        --++(0,-1)
        -|($(#2)+(1,1)$)
        -|(#2)
    }   

    \highlightSideX{\rightCol[1,1]}{\rightCol[1,2]}{\rightCol[1,3]}
    {(#1)|-(join)-|(#2)}{
        --++(0,1)
        -|($(#1)+(0,1)$)
        -|(#1)
    }    
    \draw[line width=\lineWidth] (#1) to (templeft);
    \draw[line width=\lineWidth] (#2) to (tempright);

    \place{#1}{0}{templeft}
    \place{#2}{0}{tempright}
}

\newcommandx{\strLabelBraid}[6][3=,4=1,5=,6=,usedefault]{
    \strBraidXX{#1}{#2}[][][#4][#5;#6]
    \place{labelLoc}{0}[-0.5*#4]{topLevel}
    \node[blank, below] (labelLoc) at (labelLoc){\small{#3}};
}

\newcommandx{\strBraidCaseI}[6][3=,4=,5=1,6=;,usedefault]{
    \joinTarget{#1;#2}{0}{topLevel}
    \readlist*\options{#6}
    \xdef\optionI{\options[1,1]}
    \xdef\optionII{\options[2,1]}

    \coordinate[shift={(0,-#5)}] (overBot) at (#2|-topLevel);
    \coordinate[shift={(0,-#5)}] (undBot) at (#1|-topLevel);

    \gettikzxy{(#1)}{\overX}{\overY}
    \gettikzxy{(#2)}{\underX}{\underY}

    \ifdim \overX<\underX
        \xdef\overleft{1}
    \else
        \xdef\overleft{0}
    \fi

    \readlist*\undCol{#4////}
    \readlist*\overCol{#3////}

    \xdef\undNat{0}
    
    \ifnonempty{\undCol[1,2]}{
        \xdef\undNat{1}
    }
    \ifnonempty{\undCol[1,3]}{
        \xdef\undNat{1}
    }

    \ifnum \overleft=1
        \xdef\newUndFunCol{\overCol[1,1],\overCol[1,2]}        
        \reverseList{\undCol[1,3]}{\newOverFunCol}
        \ifnum \undNat=0
            \xdef\newUndFunCol{\newUndFunCol}
            \xdef\newOverFunCol{\newOverFunCol}
        \fi
    \else
        \reverseList{\overCol[1,3]}{\newUndFunCol}
        \xdef\newUndFunCol{\overCol[1,1],\newUndFunCol}
        \xdef\newOverFunCol{\undCol[1,2],\undercol[1,4]}
        \ifnum \undNat=0
            \xdef\newUndFunCol{\newUndFunCol,\overCol[1,4]}        
        \fi
    \fi

    \xdef\newUndCol{
        \newUndFunCol/\undCol[1,2]/\undCol[1,3]/\undCol[1,4]
    }
    \xdef\newOverCol{
        \overCol[1,1]/\overCol[1,2]/\overCol[1,3]/white,white,c3
    }

    \xdef\overPath{
        (#1) to[out=-90, in=90] (overBot)
    }
    \xdef\underPath{
        (#2) to[out=-90, in=90] (undBot)
    }

    \ifnum \overleft=1
        \clipS{#1}{overBot}[\botClip][\topClip]
    \else
        \clipS{#1}{overBot}[\topClip][\botClip]
    \fi

    \coloursLength{\overCol[1,1]/\overCol[1,2]/\overCol[1,4]}[\leftTop]
    \coloursLength{\overCol[1,1]/\overCol[1,3]/\overCol[1,4]}[\rightTop]
    \coloursLength{\overCol[1,1]/\overCol[1,2]/\newOverFunCol}[\leftBottom]
    \coloursLength{\overCol[1,1]/\overCol[1,3]/\newOverFunCol}[\rightBottom]

    \tikzmath{
        \leftWhiteSpace=1+max(\leftTop,\leftBottom);
        \rightWhiteSpace=1+max(\rightTop,\rightBottom);
    }

    \begin{scope}
        \begin{pgfinterruptboundingbox}
            \clip \overPath\topClip;
        \end{pgfinterruptboundingbox}
        \drawSX{#2}{undBot}{#4}[\optionII]
    \end{scope}
    \begin{scope}
        \begin{pgfinterruptboundingbox}
            \clip \overPath\botClip;
        \end{pgfinterruptboundingbox}
        \drawSX{#2}{undBot}{\newUndCol}[\optionII]
    \end{scope}
    
    \ifnum \undNat=0
        \drawSX{#1}{overBot}{white,#3}[\optionI]
    \else
        \ifnum \overleft=1
            \clipS{#2}{undBot}[\topClip][\botClip]
        \else
            \clipS{#2}{undBot}[\botClip][\topClip]
        \fi
        \begin{scope}
            \begin{pgfinterruptboundingbox}
                \clip \underPath\topClip;
            \end{pgfinterruptboundingbox}
            \drawSX{#1}{overBot}{#3}[\optionII][\leftWhiteSpace;\rightWhiteSpace]
        \end{scope}
        \begin{scope}
            \begin{pgfinterruptboundingbox}
                \clip \underPath\botClip;
            \end{pgfinterruptboundingbox}
            \drawSX{#1}{overBot}{\newOverCol}[\optionII][\leftWhiteSpace;\rightWhiteSpace]
        \end{scope}
    \fi
    \coordinate (#1) at (overBot);
    \coordinate (#2) at (undBot);
}

\newcommandx{\strBraidCaseII}[6][3=,4=,5=1,6=;,usedefault]{
    \joinTarget{#1;#2}{0}{topLevel}
    \readlist*\options{#6}
    \xdef\optionI{\options[1,1]}
    \xdef\optionII{\options[2,1]}

    \coordinate[shift={(0,-#5)}] (overBot) at (#2|-topLevel);
    \coordinate[shift={(0,-#5)}] (undBot) at (#1|-topLevel);

    \gettikzxy{(#1)}{\overX}{\overY}
    \gettikzxy{(#2)}{\underX}{\underY}

    \ifdim \overX<\underX
        \xdef\overleft{1}
    \else
        \xdef\overleft{0}
    \fi

    \readlist*\undCol{#4////}
    \readlist*\overCol{#3////}

    \xdef\undNat{0}
    
    \ifnonempty{\undCol[1,2]}{
        \xdef\undNat{1}
    }
    \ifnonempty{\undCol[1,3]}{
        \xdef\undNat{1}
    }

    \ifnum \overleft=1
        \xdef\newUndFunCol{\overCol[1,1],\overCol[1,2]}        
        \reverseList{\undCol[1,3]}{\newOverFunCol}
        \ifnum \undNat=0
            \xdef\newUndFunCol{\newUndFunCol}
            \xdef\newOverFunCol{\newOverFunCol}
        \fi
    \else
        \reverseList{\overCol[1,3]}{\newUndFunCol}
        \xdef\newUndFunCol{\overCol[1,1],\newUndFunCol}
        \xdef\newOverFunCol{\undCol[1,2],\undercol[1,4]}
        \ifnum \undNat=0
            \xdef\newUndFunCol{\newUndFunCol,\overCol[1,4]}        
        \fi
    \fi

    \xdef\newUndCol{
        \newUndFunCol,c2/\undCol[1,2]/\undCol[1,3]/\undCol[1,4]
    }
    \xdef\newOverCol{
        \overCol[1,1]/\overCol[1,2]/\overCol[1,3]/c2,c1,c3
    }

    \xdef\overPath{
        (#1) to[out=-90, in=90] (overBot)
    }
    \xdef\underPath{
        (#2) to[out=-90, in=90] (undBot)
    }

    \ifnum \overleft=1
        \clipS{#1}{overBot}[\botClip][\topClip]
    \else
        \clipS{#1}{overBot}[\topClip][\botClip]
    \fi

    \coloursLength{\overCol[1,1]/\overCol[1,2]/\overCol[1,4]}[\leftTop]
    \coloursLength{\overCol[1,1]/\overCol[1,3]/\overCol[1,4]}[\rightTop]
    \coloursLength{\overCol[1,1]/\overCol[1,2]/\newOverFunCol}[\leftBottom]
    \coloursLength{\overCol[1,1]/\overCol[1,3]/\newOverFunCol}[\rightBottom]

    \tikzmath{
        \leftWhiteSpace=1+max(\leftTop,\leftBottom);
        \rightWhiteSpace=1+max(\rightTop,\rightBottom);
    }

    \begin{scope}
        \begin{pgfinterruptboundingbox}
            \clip \overPath\topClip;
        \end{pgfinterruptboundingbox}
        \drawSX{#2}{undBot}{#4}[\optionII]
    \end{scope}
    \begin{scope}
        \begin{pgfinterruptboundingbox}
            \clip \overPath\botClip;
        \end{pgfinterruptboundingbox}
        \drawSX{#2}{undBot}{\newUndCol}[\optionII]
    \end{scope}
    
    \ifnum \undNat=0
        \drawSX{#1}{overBot}{white,#3}[\optionI]
    \else
        \ifnum \overleft=1
            \clipS{#2}{undBot}[\topClip][\botClip]
        \else
            \clipS{#2}{undBot}[\botClip][\topClip]
        \fi
        \begin{scope}
            \begin{pgfinterruptboundingbox}
                \clip \underPath\topClip;
            \end{pgfinterruptboundingbox}
            \drawSX{#1}{overBot}{#3}[\optionII][\leftWhiteSpace;\rightWhiteSpace]
        \end{scope}
        \begin{scope}
            \begin{pgfinterruptboundingbox}
                \clip \underPath\botClip;
            \end{pgfinterruptboundingbox}
            \drawSX{#1}{overBot}{\newOverCol}[\optionII][\leftWhiteSpace;\rightWhiteSpace]
        \end{scope}
    \fi
    \coordinate (#1) at (overBot);
    \coordinate (#2) at (undBot);
}

\newcommandx{\strModXXXX}[7][4=0,5=0.5,6=0.5,7=,usedefault]{
    \getNodes{#1}{\nodes}
    
    \joinTarget{\nodes}{-#5}{joint}[][-1]

    \bufferlist{\inputs}{#1/}
    \readlist*\lst{\inputs}

    \bufferlist{\outputs}{#3/}
    \readlist*\outlst{\outputs}

    \listLengthX{\lst[1,2],\lst[1,3],\lst[1,5]}{\lflen}
    \listLengthX{\lst[-1,2],\lst[-1,4],\lst[-1,5]}{\rtlen}
    \listLengthX{\outlst[1,2],\outlst[1,3],\outlst[1,5]}{\lflenII}
    \listLengthX{\outlst[-1,2],\outlst[-1,4],\outlst[-1,5]}{\rtlenII}

    \tikzmath{
        \lshft=TCLShft(\lflen);
        \rshft=-TCLShft(\rtlen);
        \TCWidth=5+max(symWidth(\nodesWidth,\lflen,\rtlen),symWidth(((#4)*28.45276),\lflenII, \rtlenII));
    }

    \draw (joint) node[mod2={\TCWidth}](joint){#2};
    
    \foreach \i in {1,...,\lstlen}{
        \ifnum \i<\lstlen
            \reverseList{\lst[\i,4]}{\natCol}
            \path (\lst[\i,1]) to node[pos=0.5,blank](halfway){} (\lst[\i+1,1]);
            
            \highlightSideX{\lst[\i,2]}{\natCol}{\lst[\i,5]}
            {
                (\lst[\i,1]) |- (halfway|-joint.north)
            }{
                |-(\lst[\i,1])
            }

            \highlightSideX{\lst[\i+1,2]}{\lst[\i+1,3]}{\lst[\i+1,5]}
            {
                (\lst[\i+1,1]) |- (halfway|-joint.north)
            }{
                |-(\lst[\i+1,1])
            }        
            \draw[NTDash] (\lst[\i,1]) |- (halfway|-joint.north);
            \draw[NTDash] (\lst[\i+1,1]) |- (halfway|-joint.north);
        \fi
    }
    
    \getNodes{#3}{\outNodes}
    \splitTargets{joint}{\outNodes}{#4}{#6}

    \foreach \i in {1,...,\outlstlen}{
        \ifnum \i<\outlstlen
            \reverseList{\outlst[\i,4]}{\natCol}
            \path (\outlst[\i,1]) to node[pos=0.5,blank](halfway){} (\outlst[\i+1,1]);

            \highlightSideX{\outlst[\i,2]}{\natCol}{\outlst[\i,5]}
            {
                (\outlst[\i,1]) |- (halfway|-joint.south)
            }{
                |-(\outlst[\i,1])
            }

            \highlightSideX{\outlst[\i+1,2]}{\outlst[\i+1,3]}{\outlst[\i+1,5]}
            {
                (\outlst[\i+1,1]) |- (halfway|-joint.south)
            }{
                |-(\outlst[\i+1,1])
            }        
            \draw[NTDash] (\outlst[\i,1]) |- (halfway|-joint.south);
            \draw[NTDash] (\outlst[\i+1,1]) |- (halfway|-joint.south);  
        \fi
    }

    \ifdim #5pt >0pt
        \xdef\leftEntry{\lst[1,1]}
        \xdef\rightEntry{\lst[-1,1]}
        \xdef\leftEntryColI{\lst[1,2]}
        \xdef\leftEntryColII{\lst[1,3]}
        \xdef\leftEntryColIII{\lst[1,5]}
        \xdef\rightEntryColI{\lst[-1,2]}
        \reverseList{\lst[-1,4]}{\rightEntryColII}
        \xdef\rightEntryColIII{\lst[-1,5]}
    \else
        \xdef\leftEntry{joint.north}
        \xdef\rightEntry{joint.north}
        \xdef\leftEntryColI{\outlst[1,2]}
        \xdef\leftEntryColII{\outlst[1,3]}
        \xdef\leftEntryColIII{\outlst[1,5]}
        \xdef\rightEntryColI{\outlst[-1,2]}
        \reverseList{\outlst[-1,4]}{\rightEntryColII}
        \xdef\rightEntryColIII{\outlst[-1,5]}
    \fi
    \ifnonempty{#3}{
        \xdef\leftExit{\outlst[1,1]}
        \xdef\rightExit{\outlst[-1,1]}
        \xdef\leftExitColI{\outlst[1,2]}
        \xdef\leftExitColII{\outlst[1,3]}
        \xdef\leftExitColIII{\outlst[1,5]}
        \xdef\rightExitColI{\outlst[-1,2]}
        \reverseList{\outlst[-1,4]}{\rightExitColII}
        \xdef\rightExitColIII{\outlst[-1,5]}
    }
    \ifempty{#3}{
        \xdef\leftExit{joint.south}
        \xdef\rightExit{joint.south}
        \xdef\leftExitColI{\lst[1,2]}
        \xdef\leftExitColII{\lst[1,3]}
        \xdef\leftExitColIII{\lst[1,5]}
        \xdef\rightExitColI{\lst[-1,2]}
        \xdef\rightExitColII{\lst[-1,3]}
        \xdef\rightExitColIII{\lst[-1,5]}
    }

    \clipS{\leftEntry}{joint.west}
    \highlightSideX{\leftEntryColI}{\leftEntryColII}{\leftEntryColIII}{
        (\leftEntry)|-(joint.north west)--(joint.west)
    }{
        \leftClip
    }
    \draw[NTDash] (\leftEntry)|-(joint.north west)--(joint.west);

    \clipS{\rightEntry}{joint.east}
    \highlightSideX{\rightEntryColI}{\rightEntryColII}{\rightEntryColIII}{
        (\rightEntry)|-(joint.north east)--(joint.east)
    }{
        \rightClip
    }
    \draw[NTDash] (\rightEntry)|-(joint.north east)--(joint.east);

    \clipS{\leftExit}{joint.west}
    \highlightSideX{\leftExitColI}{\leftExitColII}{\leftExitColIII}{
        (\leftExit)|-(joint.south west)--(joint.west)    
    }{
        \leftClip
    }
    \draw[NTDash] (joint.south)|-(joint.south west)--(joint.west);    

    \clipS{\rightExit}{joint.east}
    \highlightSideX{\rightExitColI}{\rightExitColII}{\rightExitColIII}{
        (\rightExit)|-(joint.south east)--(joint.east)
    }{
        \rightClip
    }
    \draw[NTDash] (joint.south)|-(joint.south east)--(joint.east);
}

\newcommandx{\strBraidCaseIII}[6][3=,4=,5=1,6=;,usedefault]{
    \joinTarget{#1;#2}{0}{topLevel}
    \readlist*\options{#6}
    \xdef\optionI{\options[1,1]}
    \xdef\optionII{\options[2,1]}

    \coordinate[shift={(0,-#5)}] (overBot) at (#2|-topLevel);
    \coordinate[shift={(0,-#5)}] (undBot) at (#1|-topLevel);

    \gettikzxy{(#1)}{\overX}{\overY}
    \gettikzxy{(#2)}{\underX}{\underY}

    \ifdim \overX<\underX
        \xdef\overleft{1}
    \else
        \xdef\overleft{0}
    \fi

    \readlist*\undCol{#4////}
    \readlist*\overCol{#3////}

    \xdef\undNat{0}
    
    \ifnonempty{\undCol[1,2]}{
        \xdef\undNat{1}
    }
    \ifnonempty{\undCol[1,3]}{
        \xdef\undNat{1}
    }

    \ifnum \overleft=1
        \xdef\newUndFunCol{\overCol[1,1],\overCol[1,2]}        
        \reverseList{\undCol[1,3]}{\newOverFunCol}
        \ifnum \undNat=0
            \xdef\newUndFunCol{\newUndFunCol}
            \xdef\newOverFunCol{\newOverFunCol}
        \fi
    \else
        \reverseList{\overCol[1,3]}{\newUndFunCol}
        \xdef\newUndFunCol{\overCol[1,1],\newUndFunCol}
        \xdef\newOverFunCol{\undCol[1,2],\undercol[1,4]}
        \ifnum \undNat=0
            \xdef\newUndFunCol{\newUndFunCol,\overCol[1,4]}        
        \fi
    \fi

    \xdef\newUndCol{
        c2,c1/c2,c1/
    }
    \xdef\newOverCol{
        /c2,c1,white,white
    }

    \xdef\overPath{
        (#1) to[out=-90, in=90] (overBot)
    }
    \xdef\underPath{
        (#2) to[out=-90, in=90] (undBot)
    }

    \ifnum \overleft=1
        \clipS{#1}{overBot}[\botClip][\topClip]
    \else
        \clipS{#1}{overBot}[\topClip][\botClip]
    \fi

    \coloursLength{\overCol[1,1]/\overCol[1,2]/\overCol[1,4]}[\leftTop]
    \coloursLength{\overCol[1,1]/\overCol[1,3]/\overCol[1,4]}[\rightTop]
    \coloursLength{\overCol[1,1]/\overCol[1,2]/\newOverFunCol}[\leftBottom]
    \coloursLength{\overCol[1,1]/\overCol[1,3]/\newOverFunCol}[\rightBottom]

    \tikzmath{
        \leftWhiteSpace=1+max(\leftTop,\leftBottom);
        \rightWhiteSpace=1+max(\rightTop,\rightBottom);
    }

    \begin{scope}
        \begin{pgfinterruptboundingbox}
            \clip \overPath\topClip;
        \end{pgfinterruptboundingbox}
        \drawSX{#2}{undBot}{#4}[\optionII]
    \end{scope}
    \begin{scope}
        \begin{pgfinterruptboundingbox}
            \clip \overPath\botClip;
        \end{pgfinterruptboundingbox}
        \drawSX{#2}{undBot}{\newUndCol}[\optionII]
    \end{scope}
    
    \ifnum \undNat=0
        \drawSX{#1}{overBot}{white,#3}[\optionI]
    \else
        \ifnum \overleft=1
            \clipS{#2}{undBot}[\topClip][\botClip]
        \else
            \clipS{#2}{undBot}[\botClip][\topClip]
        \fi
        \begin{scope}
            \begin{pgfinterruptboundingbox}
                \clip \underPath\topClip;
            \end{pgfinterruptboundingbox}
            \drawSX{#1}{overBot}{#3}[\optionII][\leftWhiteSpace;\rightWhiteSpace]
        \end{scope}
        \begin{scope}
            \begin{pgfinterruptboundingbox}
                \clip \underPath\botClip;
            \end{pgfinterruptboundingbox}
            \drawSX{#1}{overBot}{\newOverCol}[\optionII][\leftWhiteSpace;\rightWhiteSpace]
        \end{scope}
    \fi
    \coordinate (#1) at (overBot);
    \coordinate (#2) at (undBot);
}

\newcommandx{\strBraidCaseIV}[6][3=,4=,5=1,6=;,usedefault]{
    \joinTarget{#1;#2}{0}{topLevel}
    \readlist*\options{#6}
    \xdef\optionI{\options[1,1]}
    \xdef\optionII{\options[2,1]}

    \coordinate[shift={(0,-#5)}] (overBot) at (#2|-topLevel);
    \coordinate[shift={(0,-#5)}] (undBot) at (#1|-topLevel);

    \gettikzxy{(#1)}{\overX}{\overY}
    \gettikzxy{(#2)}{\underX}{\underY}

    \ifdim \overX<\underX
        \xdef\overleft{1}
    \else
        \xdef\overleft{0}
    \fi

    \readlist*\undCol{#4////}
    \readlist*\overCol{#3////}

    \xdef\newUndCol{
        c1,c2//c2,c1
    }
    \xdef\newOverCol{
        //white,white,c1,c2
    }

    \xdef\overPath{
        (#1) to[out=-90, in=90] (overBot)
    }
    \xdef\underPath{
        (#2) to[out=-90, in=90] (undBot)
    }

    \ifnum \overleft=1
        \clipS{#1}{overBot}[\botClip][\topClip]
    \else
        \clipS{#1}{overBot}[\topClip][\botClip]
    \fi

    \coloursLength{\overCol[1,1]/\overCol[1,2]/\overCol[1,4]}[\leftTop]
    \coloursLength{\overCol[1,1]/\overCol[1,3]/\overCol[1,4]}[\rightTop]
    \coloursLength{\overCol[1,1]/\overCol[1,2]/\newOverFunCol}[\leftBottom]
    \coloursLength{\overCol[1,1]/\overCol[1,3]/\newOverFunCol}[\rightBottom]

    \tikzmath{
        \leftWhiteSpace=1+max(\leftTop,\leftBottom);
        \rightWhiteSpace=1+max(\rightTop,\rightBottom);
    }

    \begin{scope}
        \begin{pgfinterruptboundingbox}
            \clip \overPath\topClip;
        \end{pgfinterruptboundingbox}
        \drawSX{#2}{undBot}{#4}[\optionII]
    \end{scope}
    \begin{scope}
        \begin{pgfinterruptboundingbox}
            \clip \overPath\botClip;
        \end{pgfinterruptboundingbox}
        \drawSX{#2}{undBot}{\newUndCol}[\optionII]
    \end{scope}
    
    \ifnum \undNat=0
        \drawSX{#1}{overBot}{white,#3}[\optionI]
    \else
        \ifnum \overleft=1
            \clipS{#2}{undBot}[\topClip][\botClip]
        \else
            \clipS{#2}{undBot}[\botClip][\topClip]
        \fi
        \begin{scope}
            \begin{pgfinterruptboundingbox}
                \clip \underPath\topClip;
            \end{pgfinterruptboundingbox}
            \drawSX{#1}{overBot}{#3}[\optionII][\leftWhiteSpace;\rightWhiteSpace]
        \end{scope}
        \begin{scope}
            \begin{pgfinterruptboundingbox}
                \clip \underPath\botClip;
            \end{pgfinterruptboundingbox}
            \drawSX{#1}{overBot}{\newOverCol}[\optionII][\leftWhiteSpace;\rightWhiteSpace]
        \end{scope}
    \fi
    \coordinate (#1) at (overBot);
    \coordinate (#2) at (undBot);
}

    \tikzstyle{finish}=[to path={(\tikztostart) -- (\tikztotarget.#1)}]
    \tikzstyle{start}=[to path={(\tikztostart.#1) -- (\tikztotarget)}]
    \tikzstyle{vertical}=[
      to path={
          (\tikztostart) 
          -- (\tikztostart|-\tikztotarget.north)\tikztonodes}
    ]

    \tikzstyle{verticalup}=[
      to path={
          (\tikztostart) 
          -- (\tikztostart|-\tikztotarget.south)\tikztonodes}
    ]

    \newenvironment{cdaligned}{\catcode`\&=4 \aligned}{\endaligned}

    \makeatletter
    \newcommand{\asti}{%
      \check@mathfonts
      \leavevmode
      {\ooalign{%
        \hidewidth
        \raisebox{.8ex}{\fontsize{\ssf@size}{0}\selectfont*}%
        \hidewidth\cr
        \i\cr
      }}%
    }
    \makeatother

    \newcommand{\eid}{\ecomp{\id}}
    \newcommand{\cll}{\operatorname{cl}}

\tikzexternaldisable

\renewbibmacro*{urldate}{
(visited on \printfield{urlday}/\printfield{urlmonth}/\printfield{urlyear})
}

\begin{document}
\pagenumbering{gobble}
\begin{titlepage}
  \newcommand*{\customtitle}{\begingroup%
  \vspace*{\baselineskip}
  \vfill
  \begin{figure}
    \begin{minipage}[t][0pt]{0.2\textwidth}
      \vspace{0pt}
      \rule{1pt}{\textheight}
    \end{minipage}
    \hspace{0.05\textwidth}
    \begin{minipage}[t][0pt]{0.7\linewidth}
      \vspace{0pt}
      \includegraphics[height=3cm]{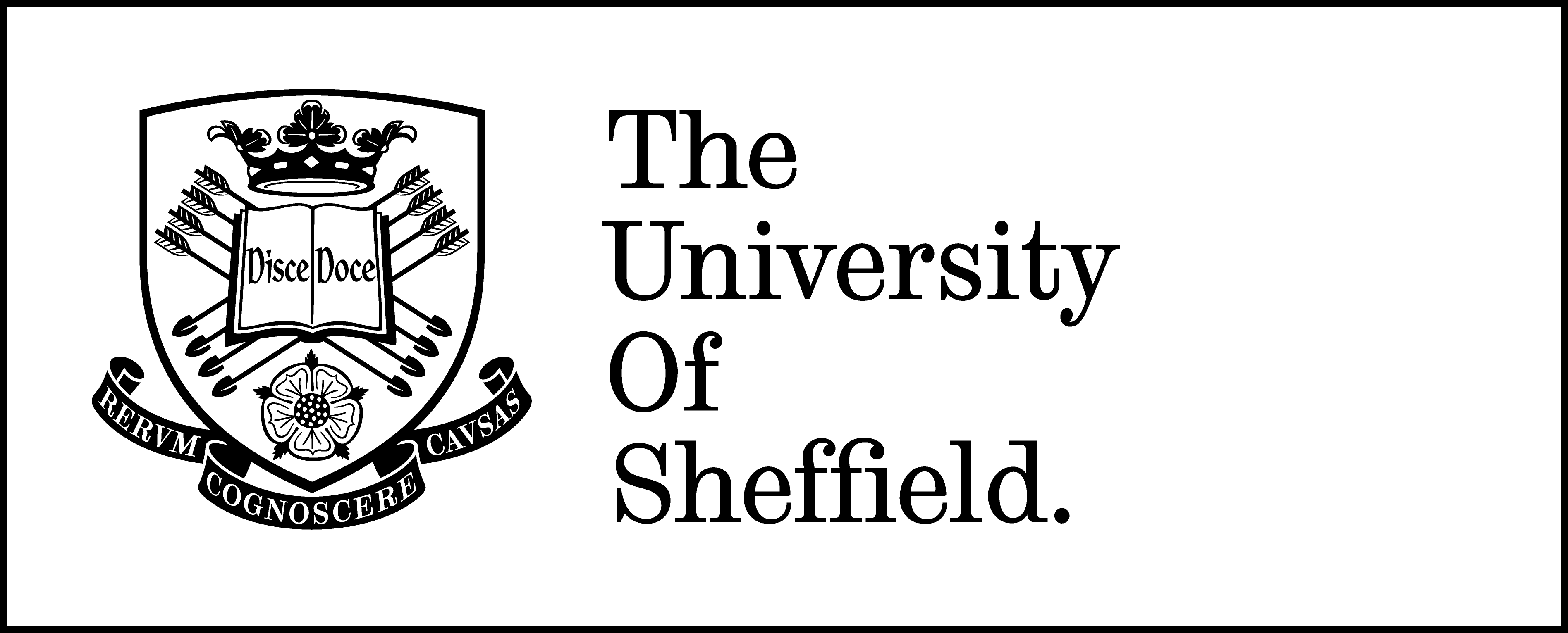}$ $\\[1.0\baselineskip]
      {\noindent \Large School of Mathematics and Statistics}\\[1\baselineskip]
      \parbox[b]{\textwidth}{
        \vspace{0.10\textheight}
        {\noindent \Huge\bfseries 
        Scalar Enrichment\\[0.5\baselineskip]
        \noindent and Cotraces\\[0.5\baselineskip]
        \noindent in Bicategories
        }\\[5\baselineskip]
        {\LARGE Callum William Reader}\\[\baselineskip]
        {\Large\itshape Supervised by Dr.\ Simon Willerton
        \vspace{0.30\textheight}

        \noindent Submitted for the
          degree of Doctor of Philosophy of
          Mathematics on 22\textsuperscript{nd} March 2023
          }
      }
    \end{minipage}
  \end{figure}
  \vfill
  \null
  \endgroup}

  \customtitle
\end{titlepage}
\begin{abstract}
  It is known that every monoidal bicategory has an associated braided monoidal category of scalars. In this thesis we show that every monoidal bicategory, which is closed both monoidally and compositionally, can be enriched over the monoidal 2-category of scalar-enriched categories. This enrichment provides a number of key insights into the relationship between linear algebra and category theory. 

  The enrichment replaces every set of 2-cells with a scalar, and we show that this replacement can be given in terms of the cotrace, first defined by Day and Street in the context of profunctors. This is analogous to the construction of the Frobenius inner product between linear maps, which is constructed in terms of the trace of linear maps. In linear algebra it is also possible to define the trace in terms of the Frobenius inner product. We show that the cotrace can be defined in terms of the enrichment, and in doing so we prove that the cotrace is an enriched version of the `categorical trace' studied by Ganter and Kapranov, and Bartlett. Thus, we unify the concept of a categorical trace with the concept of a cotrace.
  
  Finally, we study the relationship between the trace and the cotrace for compact closed bicategories. We show that the trace and cotrace have a structured relationship and share many of the properties of the linear trace including -- but not limited to -- dual invariance and linearity. Motivating examples are given throughout. We also introduce a decorated string diagram language to simplify some of the proofs.
\end{abstract}

{\let\cleardoublepage\clearpage%
\chapter*{Acknowledgements}
\setcounter{page}{1}
I would like to offer my thanks to all of those who have made the process of writing this thesis possible and even, at times, enjoyable.

Firstly, thank you to Simon Willerton, my supervisor, for the endless patience and encouragement, creative insights and tangential suggestions that led to this thesis being written. And thank you to Joseph Martin for letting me air my unfinished thoughts and giving me plenty of constructive feedback.

Thank you to the rest of the category theorists in the Sheffield department and beyond who have included me in an engaged and interested community. A particular thank you to those category theorists in the Edinburgh department who welcomed me to their seminars whenever I was in the city, and during lockdowns.

Thank you to Joe and Ceci for housing me when I first arrived in Sheffield. Thank you to Poppy and Ciaran for housing me when I had nowhere else to go. And thank you to Dan Graves for putting me up (and putting up with me) during the final months of my PhD. You have all made Sheffield feel like home.

Thank you to Ellie, Hope and Vinh for encouraging me to be silly when I've needed it. Thank you to James for sharing this journey with me and encouraging me to engage with algebraic topology. Thank you to Matt and Becky for sharing many trips when I've needed them, to Green Man and much further afield. Thank you to Huw and Erin for helping me see the world outside of mathematics, for bringing fresh eyes to old problems, and for keeping their door open to me at all times.

Thank you to my father who read my whole thesis despite having very little mathematical background, and thank you to my mother who will never read this thesis but will be proud of me for writing it. Without their constant love and support I would not have made it this far in my academic career. Thank you to John and Sally, and Henry and Lynda, who have always treated me like part of the family and taken an interest in my work.

Mostly I would like to thank my partner, Caz. They shared in my highs, kept me sane when things went wrong, taught me to enjoy reading again, supported me in all of my ill-advised creative endeavours, and did everything in their power to make the final days of writing bearable.
\setcounter{page}{0}

\tableofcontents
}
\afterpage{\blankpage}
\chapter*{Introduction}
\addcontentsline{toc}{chapter}{Introduction}
\setcounter{page}{1}
\pagenumbering{arabic}
In this thesis we will investigate how ideas from linear algebra have strong analogues in the theory of bicategories. In particular, we will focus on the role of scalars and the role of traces.

A bicategory consists of objects, 1-cells between objects, and, between any pair of 1-cells, a set of 2-cells. The category of inner product spaces has objects and 1-cells between them. But its closed structure gives us an extra piece of data. Between any pair of 1-cells there is a scalar. This scalar is given by taking the Frobenius inner product
\[
\langle f,g\rangle =\Tr(f^\dagger\circ g)   
\]
where $\Tr$ denotes the trace and $(-)^\dagger$ denotes the adjoint linear map. In this thesis we show that an analogous construction exists in any bicategory equipped with appropriate closed structures. That is to say, between any two 1-cells we can place a scalar and this scalar is defined in terms of a functor with trace-like properties.

This result has several consequences. Firstly, it highlights a canonical enrichment that gives many well-known bicategories their extra structure. 

Secondly, it unifies the cotrace as defined by Day and Street~\cite[def.~8]{street1997monoidal}, with the `categorical trace' as defined by Ganter and Kapranov~\cite[def.~3.1]{ganter} and the `2-trace' as defined by Bartlett~\cite[def.~7.8]{brucethesis} for 2-Hilbert spaces.

Finally, it provides a theoretical underpinning for Willerton's~\cite{willertonTalk,willertonUnpub} observation that the 2-trace seems to be somehow dual to the usual notion of trace. Willerton pointed out that, if we extend the definitions of trace and 2-trace to the context of a bicategory with duals, these two different traces always appear to give opposing results. For example, in the bicategory of profunctors the trace gives a coend, but the 2-trace gives an end; in the bicategory of bimodules, the trace gives coinvariants, but the 2-trace gives invariants. The problem with this observation was that the trace is a scalar -- that is, a map from the unit object to itself -- whereas the 2-trace is a \emph{set} of 2-cells. It was only after adding appropriate structure to the set that this ostensible duality made sense. Attempting to solve this problem gave life to this thesis, and the canonical scalar enrichment means that the trace and the enriched 2-trace -- what Day and Street call the cotrace -- live in the same category and can therefore be compared formally.

It should be noted here that there has been a lot of previous work in the literature on traces in bicategories, with some interesting applications in topology. Ponto~\cite[]{pontoOriginal} defined traces in bicategories with shadows, and later expanded the theory of shadows and traces alongside Shulman~\cite[]{pontoShulman}. This work is somewhat distinct from our own. Ponto's definition of a bicategorical trace relies on the bicategory coming equipped with a shadow structure, which is analogous to a symmetry structure for a monoidal category. The trace in this setting is an invariant of endo-2-cells. Our trace and cotrace do not require a shadow structure and are invariants of endo-1-cells. Interestingly, however, it would appear that the trace and cotrace do, in fact, provide shadow-like structures for the bicategory, but we do not explore that avenue in this thesis.

Let us explain how this scalar enrichment works. A scalar, in the context of linear algebra, is just a complex number. But we can also think of a scalar as a linear endomorphism on $\mathbb{C}$. This means that we can define set of scalars in a monoidal category to be morphisms from the unit object to itself. What's more, Kelly and Laplaza~\cite[prop.~6.1]{kelly1980coherence} pointed out that the set of scalars in a monoidal category automatically comes equipped with the structure of a commutative monoid, and that this commutative monoid then acts on every hom-set in the category. One way to view this scalar action is as follows. Let $\MM$ be a monoidal category, let $B$ be an object in $\MM$ and let $s\colon I \to I$ be a scalar. We can `spread' the scalar at $B$ by taking the composite
\[
B\xrightarrow{\sim} I\tensor B \xrightarrow{s\tensor B} I\tensor B \xrightarrow{\sim} B    
\] 
and then the action of scalars on the set $\MM(A,B)$ is given by postcomposition with the spread of $s$. In the particular case of finite dimensional Hilbert spaces, which form a closed category, this spread can be thought of as a linear morphism
\[
    \mathrm{Spr}\colon \mathbb{C} \to \FDHilb(B,B)   
\]
which has a linear adjoint. That linear adjoint is the trace.

In a monoidal bicategory $\BB$ a similar principle exists. The \emph{category} of scalars $\BB(I,I)$ is automatically equipped with a braided monoidal structure and acts on every hom-category. The action is given in terms of a functor that we call the spread functor
\[
    \Spr\colon \BB(I,I)\to \BB(B,B).
\]
Here we use the $\circ$ superscript to indicate that the spread functor is dual in some sense to another functor that we call the cospread functor, and to indicate that the spread functor is defined in terms of composition in the bicategory.

It can be easily shown that if the bicategory is right-monoidal-closed, and left-composition-closed, then the spread functor has a right adjoint. We call this right adjoint the cotrace functor,
\[
    \ctr\colon \BB(B,B)\to \BB(I,I)
\] 
since it coincides with the cotrace as defined by Day and Street~\cite[def.~8]{street1997monoidal}.

Here we use the overline and $\circ$ superscript to indicate that the cotrace functor is dual in some sense to the trace functor, and to indicate that the cotrace functor is defined in terms of lifts in the bicategory.

In our definition of the cotrace functor we have already assumed that the bicategory is left-composition-closed. That is to say, we have assumed that it has all lifts. Given $f,g\colon A\to B$ this means that we can always take the lift $f\lift g$ in the following diagram.
\[
\begin{tikzcd}[column sep=3.2em]
    B   
        & A
        \arrow[l, "f"{swap}]
            & A
            \arrow[l, dashed, "f\lift g"{swap}]
            \arrow[ll, "g", ""{swap,name=bottom}, bend left=65]
    \arrow[Leftarrow, to=\tikzcdmatrixname-1-2, from=bottom, verticalup]
\end{tikzcd}    
\]
In particular that means that given $f,g\colon A\to B$ we can always construct a scalar `between them' by taking the cotrace of the lift,
\[
    \ctr(f\lift g).
\]
In fact, this construction makes the action of scalars into a closed action. A result of Gordon and Power~\cite[thm.~3.7]{gordon1997enrichment} shows that closed actions give rise to enriched categories, and this leads us to our main theorem.
\begin{restatable*}{theorem}{scalarEnrichment}
    \label{thm:scalarEnrichment}
    Every left-composition-closed, right-monoidal-closed bicategory $\BB$ is the underlying bicategory of a $\cat{\BB(I,I)}$-enriched bicategory.
\end{restatable*}
Note that, by construction, the 2-hom-object $\underline{\BB}(A,B)(f,g)$ is given by taking the cotrace of the lift of $g$ through $f$. In the particular case that $f$ has a right adjoint $f^\dagger$, this means that this 2-hom-object is given by
\[
\ctr(f^\dagger \circ g)    
\]
analogous to the definition of Frobenius inner product.

Whilst the Frobenius inner product can be defined in terms of the trace, we can also define the trace in terms of the Frobenius inner product. If the Frobenius inner product had already been defined from first principles, the trace of $E\colon A\to A$ could be defined as
\[
\Tr(E)\coloneqq \langle \id, E\rangle.
\]
Similarly, we see that the cotrace can be defined in terms of the enrichment. For an endo-1-cell $f\colon A\to A$ we have that
\[
    \underline{\BB}(A,A)(\id, f) = \ctr(\id\lift f)\cong \ctr(f).
\]
This then unifies the cotrace with another trace-like construction. In studying categorical representation theory two sets of authors -- Ganter and Kapranov~\cite[def.~3.1]{ganter}, and Bartlett~\cite[def.~7.8]{brucethesis} -- have previously defined the trace of an endo-1-cell as being the hom-set \[\BB(A,A)(\id, f).\] Ganter and Kapranov refer to this as the `categorical trace', whereas Bartlett refers to this as the `2-trace'. In many contexts this hom-set comes with extra structure, but the scalar enrichment above proves that this extra structure is guaranteed, and that the 2-trace is the underlying set of the cotrace.

One categorical definition of trace that we have failed to address so far is the trace that can be defined in a compact-closed category. A compact-closed category, first defined by Kelly~\cite[p.~102]{kelly1972many}, is a symmetric monoidal category where every object has duals. That is, every $A$ has an associated $A^*$ and coevaluation and evaluation maps
\[
\coev\colon I\to A\tensor A^*\text{ and }  \ev\colon A^*\tensor A\rightarrow I,
\]
which in the string diagram language are written as caps and cups  \[\myinput{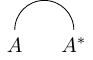}\qquad
\myinput{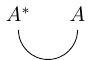}
\]
which satisfy the zigzag identities for adjunctions.
Given an endomorphism $E\colon A\to A$ in the category of finite dimensional vector spaces, we can define the trace to be the composite
\[
I\xrightarrow{coev} A\tensor A^* \xrightarrow{E\tensor A^*} A\tensor A^* \xrightarrow{\sim} A^*\tensor A\xrightarrow{ev} I
\]
which is given by the following string diagram.
\[
\myinput{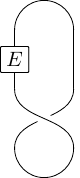}
\]
Of course, in compact-closed \emph{bicategories} -- first defined by Day and Street~\cite[def.~6]{street1997monoidal} -- it is also possible to define a trace, and in this context the trace gives a functor
\[
    \tr\colon \BB(A,A)\to \BB(I,I).  
\]
In the bicategory of profunctors, for example, the trace takes in a profunctor $P\colon \AA\profto \AA$ and returns the end of that profunctor. Rather interestingly, the cotrace gives the dual notion, the coend.

Not only does the cotrace provide a conceptual dual to the trace in many other examples, but we also see that in the case of bicategories, the two traces actually have remarkably similar definitions and interact nicely. They are both, for example, linear in an appropriate sense.

In the first chapter we give some useful definitions and propositions relating to bicategories. In particular, we cover composition-closedness, which is the first piece of structure necessary for scalar enrichment, and we give the definition of the 2-trace -- also known as the `categorical trace' -- which we will later see is the unenriched version of the cotrace. We also take a look at some motivating examples and what the 2-trace gives in each of these cases.

In the second chapter we cover pseudofunctors and pseudonatural transformations, and introduce an adaptation of the usual string diagram language which, in our opinion, makes proving some basic propositions much easier. We also cover some of the more technical details of bicategories including the coherence theorem and pseudoadjunctions.

In the third chapter we cover monoidal bicategories in detail. We give a proof, using our adapted string diagram language, that the monoidal category of scalars associated to a monoidal bicategory comes equipped with a braid. We also give a short account of monoidal-closed bicategories. This is the second piece of structure necessary for defining the cotrace.

In the fourth chapter we deal with categorical representations. This whole chapter hinges on what Garner calls the fundamental theorem of enriched category theory, first proved by Gordon and Street~\cite[thm.~3.7]{gordon1997enrichment} and later reproved by Janelidze and Kelly~\cite[sec.~6]{janelidze2000actions}. Given a monoidal category $\VV$, a $\VV$-representation is a category equipped with an oplaxly associative and unital action of $\VV$. By this we mean that there is a functor
\[
    \act\colon \VV\times \CC \to \CC
\]
which comes equipped with associator and unitor natural transformations satisfying certain coherence conditions. Analogously to the case for monoidal categories if, for every $C\in \CC$, the functor $(-)\odot C$ has a right adjoint, we call the representation closed. The fundamental theorem then says that a closed $\VV$-representation is the `same thing' as a copowered $\VV$-category. In order to prove our enrichment theorem in the following chapter, we give a slight variation on this result. We show that there is a 2-functor from the 2-category of \emph{all} closed representations to the 2-category of \emph{all} enriched categories.

In the fifth chapter we give our main theorem: the scalar enrichment of a monoidal bicategory. We begin by giving a short account of enriched bicategories and their underlying bicategories. We then go on to define the spread functor and the cotrace functor and show that they form an adjoint pair. These functors are then used to show that each $\BB(A,B)$ is a closed $\BB(I,I)$-representation and this, in essence, gives us our scalar enrichment. The following section shows that all of the additional data of bicategories -- composition, identities, associators, unitors -- is compatible with the representation structure. The consequence of this is that, not only are each of the hom-categories enriched, but $\BB$ is the underlying bicategory of a $\cat{\BB(I,I)}$-bicategory.

In the final chapter we look at more structured monoidal categories, starting with braided, moving on to symmetric and finishing with compact-closed. We give some basic propositions about compact-closed bicategories and take a look at the trace for compact-closed bicategories. We compare composition-closed compact-closed bicategories to dagger compact categories, to give some explanation as to why there is no cotrace for linear algebra -- the answer being that the cotrace and trace seem to coincide. We look at each of the main properties of the trace -- linearity, dual invariance, cyclicity, tensor preservation -- and show that analogous properties hold for both the cotrace and the trace. Finally, we study the codimension and dimension of objects in the bicategory and show that they always give a monoid-module pair.

\chapter{Bicategories}
In this introductory chapter we give an account of the basic definitions and properties of a bicategory, with a particular focus on certain motivating examples. We also give an account of the 2-trace and what it produces in each of our examples.

The first section includes definitions of key structures internal to a bicategory, such as adjoint equivalences. These will be important later on as we develop external structures for bicategories, such as a monoidal products. 

The following section focuses on a motivating definition for the thesis: the 2-trace. We give background and analogies with linear algebra that justify the definition, but show through examples that the 2-trace seems to be missing some structure.

In the final section we provide some of the basic theory of closed bicategories. Later on we will use the closed structure to define the cotrace functor, and it is this cotrace functor which ultimately gives us a canonical enrichment, providing the 2-trace with the structure that it otherwise lacks.
\section{Definitions and Examples}
A bicategory is a particular type of weak 2-category. By this we mean that a bicategory has objects, morphisms between objects and 2-cells between morphisms -- and that 2-cells can be composed unitally and associatively -- but composition of the morphisms is associative and unital only up to isomorphism.
\begin{definition}
A \textdef{bicategory} $\BB$ consists of:
\begin{itemize}
    \item a collection of \textdef{0-cells}, or objects, which, by abuse of notation, we denote $\BB$;
    \item for every pair of objects $A,B\in \BB$ a category $\BB(A,B)$, called the \textdef{hom-category}, the objects of which are called \textdef{1-cells}, and denoted $f\colon A\to B$, and the morphisms of which are called \textdef{2-cells} and denoted $\phi\colon f\Rightarrow g$;
    \item for every triple of objects $A,B,C\in \BB$ a functor
    \[\circ \colon \BB(B,C)\times \BB(A,B)\to \BB(A,C)\]
    called the \textdef{composition functor};
    \item for every object $A$ a 1-cell $\id\colon A\to A$ called the \textdef{identity} 1-cell;
    \item for every triple of 1-cells
      $f\colon A\to B, \;g\colon B\to C,\; h\colon C\to D $
     a natural isomorphism
     \[
        \comp{\alpha}_{f,g,h}\colon (h\circ g)\circ f\Rightarrow h\circ (g\circ f);
        \]
    called the \textdef{composition associator};
    \item for every 1-cell $f\colon A\to B$ natural isomorphisms
    \[
    \comp{\lambda}\colon \id_X\circ f\Rightarrow f\text{ and } \comp{\rho}\colon f\circ \id_X\Rightarrow f;
    \]
    called the \textdef{composition unitors}\footnote{
        It is not standard here to decorate our natural isomorphisms with $\circ$, and typically these natural isomorphisms are referred to simply as associators and unitors. Here we include the $\circ$ decoration and use the names `composition associator' and `composition unitor' in order to distinguish these natural isomorphisms from the associator and unitor for the tensor product in a monoidal bicategory.}.
\end{itemize}
The composition associator must satisfy the pentagon equation, given by the commutativity of the following diagram.
        \[
        \begin{tikzcd}[column sep=-1.5em, row sep=5em]
        \blank
            &&
                (k\circ h)\circ (g\circ f)
                \arrow[rrd, "\comp{\alpha}_{k,h,g\circ f}"]\\
        ((k\circ h)\circ g)\circ f
        \arrow[rru, "\comp{\alpha}_{k\circ h, g,f}"]
        \arrow[rd, "\comp{\alpha}_{k,g,h}\circ f"{swap}]
            &&
                && k\circ(h\circ(g\circ f))\\
        &   (k\circ (h \circ g)) \circ f
        \arrow[rr, "\comp{\alpha}_{k,h\circ g,f}"{swap}]
            &&
                k\circ ((h \circ g) \circ f)
                \arrow[ru, "k\circ \comp{\alpha}_{h,g,f}"{swap}]
        \end{tikzcd}
        \]
        The composition associators and unitors must satisfy the triangle equation given by the commutativity of the following diagram,
        \[
        \begin{tikzcd}
        (g\circ \id_B) \circ f
        \arrow[rd, "\comp{\rho}\circ \iota_f"{swap}]
        \arrow[rr, "\comp{\alpha}_{f,\id_B,g}"]
            && g\circ (\id_B\circ f)
            \arrow[ld, "\iota_g\circ \comp{\lambda}"]\\
        &g\circ f
        \end{tikzcd}
        \]
        for all $f\colon A\to B$, $g\colon B\to C$, where $\iota_f$ and $\iota_g$ are identity 2-cells.
\end{definition}
\begin{notation}
In a given bicategory $\BB$, unless stated otherwise, we will use the following notational conventions for the various cells and composites in $\BB$:
\begin{itemize}
    \item capital letters will represent the objects $A\in \BB$;
    \item lower case letters will represent the 1-cells $f\in \BB(A,B)$ and will typically be written $f\colon A\to B$;
    \item lower case Greek letters will represent the 2-cells $\alpha\in \BB(X,Y)(f,g)$, and will typically be written $\alpha\colon f\Rightarrow g$;
    \item the composite of two 2-cells, $\alpha\colon f\Rightarrow g$ and $\beta \colon g\Rightarrow h$ will be referred to as the \text{vertical composition} and written $\beta\cdot \alpha\colon f\Rightarrow h$;
    \item the composite of two 1-cells $f\colon X\to Y$ and $g\colon Y\to Z$ will be referred to as the horizontal composite and written $g\circ f\colon X\to Z$;
    \item the \text{horizontal composite} of two 2-cells $\alpha\colon f\Rightarrow g$ and $\beta \colon h\Rightarrow k$, will be written
    $\beta\circ \alpha \colon h\circ f \Rightarrow g\circ k$;
    \item we will write $\id_X\colon X\to X$ for the identity 1-cell;
    \item we will write $\iota_f\colon f\Rightarrow f$ for the identity 2-cell.
\end{itemize}
\end{notation}
\begin{definition}
    We call a bicategory a \textdef{2-category} if all of the composition unitors and associators are identity 2-cells.
\end{definition}
Throughout the thesis there will be several key examples of bicategory that we keep referring back to.
\begin{example}
 The bicategory $\Rel$ whose objects are sets, whose 1-cells are relations, and where there is a unique 2-cell $R\Rightarrow S$ if and only $R\subseteq S$. Horizontal composition is given by the usual composition of relations. That is, given $R\colon A\to B$ and $S\colon B\to C$ we have 
 \[
    S\circ R=\{(a,c)\mid \exists b\in B\text{ such that } (a,b)\in R\text{ and }(b,c)\in S\}.
 \]
\end{example}
\begin{example}
For some commutative ring $R$, the bicategory $\Bim_R$, whose objects are $R$-algebras, whose 1-cells $A\to B$ are $A$-$B$-bimodules, and whose 2-cells are bimodule homomorphisms. Horizontal composition is given by tensoring bimodules.
\end{example}
There is also a `locally derived' version of this bicategory.
\begin{example}
For some commutative ring $R$, the bicategory $\DGBim_R$ whose objects are $R$-algebras and whose hom-categories are given by the derived categories of bimodules. Horizontal composition is given by the derived tensor product -- sometimes called the total tensor product -- of chain complexes of modules. See for example Weibel's \cite[thm.~10.6.3]{Weibel} textbook.
\end{example}
\begin{example}
For some cosmos $\VV$ -- that is a complete, cocomplete, closed symmetrical monoidal category -- the bicategory $\VV\text{-}\Prof$ whose objects are $\VV$-categories, whose 1-cells are profunctors, and whose 2-cells are $\VV$-natural transformations. Horizontal composition of two profunctors, $P\colon \AA\profto \BB$ and $Q\colon \BB\profto \CC$ is given by a chosen coend
\[
\endint^{B\in \BB} Q(-,B)\tensor P(B,-).    
\]
\end{example}
Note that the above examples are all of a similar flavour. In fact, the first two examples are special cases of the fourth. Every set can be thought of as a discrete preorder. Every preorder can be thought of as a category enriched over the category of truth values $(F\rightarrow T)$ with monoidal product given by logical conjunction. Thus, the bicategory $\Rel$ can be viewed as the full sub-bicategory of $(F\rightarrow T)${-}$\Prof$,  consisting of just the discrete categories. 

Every algebra can be thought of as a one-object category enriched over the category of $R$-modules. Thus, the bicategory $\Bim_R$ can be identified with the full sub-bicategory of ($R$-Mod)-Prof consisting of just one-object categories.
\begin{example}
Given any category with pullbacks, $\CC$, we can construct $\Span(\CC)$ whose objects are those of $\CC$, whose 1-cells are spans, and whose 2-cells $\alpha\colon  f\Rightarrow g$ are maps in $\CC$ such that the following diagram commutes.
\[
\begin{tikzcd}
&S
\arrow[rd, "f_B"]
\arrow[ld, "f_A",swap]
\arrow[dd, "\alpha", description]\\
A
    &
        &B\\
&T
\arrow[lu, "g_A"]
\arrow[ru, "g_B"{swap}]
\end{tikzcd}
\]
The horizontal composition of two spans is given by taking a particular choice of pullback.
\end{example}
It typically helps to keep the bicategory $\Span(\Set)$ in mind. Note also that this can be thought of as the sub-bicategory of $\Set$-$\Prof$ consisting of only the discrete categories. This follows from the fact that every span $A\leftarrow S \rightarrow B$ can be thought of as a map $S\rightarrow A\times B$ and the fact that there is an equivalence between the slice category $\Set/(A\times B)$ and the functor category $\Cat(A\times B, \Set)$.
\begin{example}
    Given a topological space $T$ we can construct a bicategory $\Path(T)$ whose objects are the points of $T$, whose 1-cells are paths between points and whose 2-cells are homotopy classes of homotopies between paths. Horizontal composition is given by concatenation of paths.
\end{example}
As with monoidal categories, a major theme of bicategories is the idea of coherence. The term `coherence' was originally used by Mac Lane~\cite[p.~33]{maclane1963natural} in reference to the fact that any `structural' isomorphism in a monoidal category is unique. Suppose that we have a two composite 1-cells, $f$ and $g$, and suppose that $\gamma\colon f\Rightarrow g$ is an isomorphism composed entirely of associators and unitors. Coherence is the idea that any other isomorphism from $f$ to $g$ composed of associators and unitors should be equal to $\gamma$. The following useful results are sometimes referred to as `coherence results' in the spirit of the above.

\begin{proposition}
\label{prop:coherence}
In a bicategory $\AA$ the following diagrams commute for all objects $A,B,C$ and all morphisms $f\colon A\to B$, $g\colon B\to C$.
\[
\begin{tikzcd}
(\id_{C}\circ g)\circ f
\arrow[rr, "{\comp{\alpha}}"]
\arrow[rd, "{\comp{\lambda}\circ \iota}"{swap}]
    && \id_{C}\circ (g\circ f)
    \arrow[ld, "{\comp{\lambda}}"]\\
&   g\circ f
\end{tikzcd}
\]
\[
\begin{tikzcd}
    (g\circ f)\circ \id_A
    \arrow[rr,"\comp{\alpha}"]
    \arrow[rd, "\comp{\rho}"{swap}]
        && g\circ (f\circ\id_A)
        \arrow[ld, "\iota\circ \comp{\rho}"]\\
    &   g\circ f
\end{tikzcd}
\] 
\[
\begin{tikzcd}
    \id_A\circ \id_A
    \arrow[r,bend left, "\comp{\lambda}"]
    \arrow[r,bend right, "\comp{\rho}"{swap}]
        & \id_A
\end{tikzcd}    
\]
\end{proposition}
\begin{proof}
When Mac Lane~\cite[p.~42]{maclane1963natural} gave the first coherence results for monoidal categories, the analogues of these commuting diagrams were given as \emph{axioms} for monoidal categories. Theorems 6 and 7 of Kelly's~\cite[]{kelly1964maclane} response to this paper prove that they follow from the other axioms. The arguments in Kelly's paper are identical to the arguments for bicategories.
\end{proof}
In \autoref{section:strictification} we include the strictification theorem for bicategories that shows coherence does in fact hold. 

In the last chapter we will investigate how bicategories with extra structure relate to categories. It will be useful then, to have a way to turn a bicategory into a category.
\begin{proposition}
    For every bicategory $\BB$ there is a category, $\BB_0$, called the 2-skeleton of $\BB$, whose objects are those of $\BB$ and whose morphisms are isomorphism classes of those in $\BB$.
\end{proposition}
\begin{proof}
    We must show that $\BB_0$ is a category by defining a unital and associative composition. We define this composition in the only possible way:
    \[
    \hom{g}\circ \hom{f}\coloneqq \hom{g\circ f}. 
    \]
    To see that this is well-defined, note that if $f'\cong f$ and $g'\cong g$ then $g'\circ f'\cong g\circ f$ by horizontal composition of 2-cells. This composition is unital and associative by the fact that the composition unitors and associators are isomorphisms.
\end{proof}
\begin{remark}
This isn't the `underlying' category -- in fact an underlying category does not necessarily exist, since the composition of 1-cells isn't strictly associative. We call this the 2-skeleton since the hom-set of the 2-skeleton is the underlying set of the skeleton of the hom-category.
\end{remark}
We will often use string diagrams to describe 2-cells, which we read from right to left and top to bottom. To simplify these diagrams we will omit the composition associators and unitors. One way to view this is that the string diagrams provide proof schema -- associator and unitors can be inserted wherever necessary to provide the full proof. But the strictification theorem in \autoref{section:strictification} guarantees that any proof given using these string diagrams is a valid proof in a bicategory. The following definitions and propositions will be useful in later chapters.
\begin{definition}
    Let $f\colon A\to B$ and $g\colon B\to A$ be 1-cells in a bicategory.  We say that $(f,g,\eta,\epsilon)$ defines an \textdef{adjunction} between $A$ and $B$, if $\eta$ and $\epsilon$ are 2-cells
    \[
    \eta\colon \id_A \Rightarrow g\circ f \text{ and }
    \epsilon\colon f\circ g \Rightarrow \id_B,
    \]
    such that, when drawn as caps and cups
    \[
      \myinput{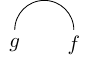}  \qquad  \myinput{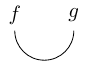}
    \] 
    the zigzag identities hold.
    \[
        \myinput{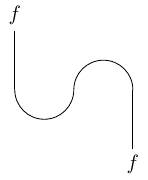}\quad=\quad\myinput{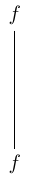} \qquad\myinput{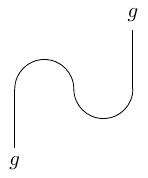}\quad=\quad\myinput{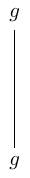}
    \]
    We call $\epsilon$ the \textdef{unit} and $\eta$ the \textdef{counit}. If there exist $\eta$ and $\epsilon$ such that $(f,g,\eta,\epsilon)$ defines an adjunction we write $f\dashv g$, and say that $f$ is \textdef{left adjoint} to $g$ and that $g$ is \textdef{right adjoint} to $f$.
\end{definition}
As with functors, a right or left adjoint is essentially unique.
\begin{proposition}
    If $g$ and $g'$ are both right adjoint to $f$ then there is a canonical isomorphism $g\cong g'$.
\end{proposition}
\begin{proof}[Sketch proof]
    This canonical isomorphism is given by the diagram below, where the cap on the left is the unit for the adjunction between $f$ and $g'$, and the cup on the right is the counit for the adjunction between $f$ and $g$. 
    \[\myinput{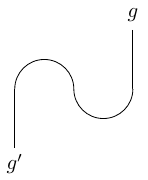}\]
    Details of the proof can be found in Johnson and Yau's~\cite[lem.~6.1.6]{yau2021bicat} book. 
\end{proof}
\begin{definition}
    Let $f\colon A\to B$ and $g\colon B\to A$ be 1-cells in a bicategory. We say that $(f,g,\eta,\epsilon)$ defines an \textdef{equivalence} between $A$ and $B$, if $\eta$ and $\epsilon$ are isomorphisms
    \[
        \eta\colon \id_A \Rightarrow g\circ f \text{ and }
        \epsilon\colon f\circ g \Rightarrow \id_B.
    \]
    We call $\eta$ the \textdef{unit} and $\epsilon$ the \textdef{counit}. If there is an equivalence between $A$ and $B$ we say that $A$ and $B$ are \textdef{equivalent}, and that $f$ and $g$ give an equivalence between $A$ and $B$. In the case that $\eta$ and $\epsilon$ are also the unit and counit of an adjunction $f\dashv g$ we say that $(f,g,\eta,\epsilon)$ defines an \textdef{adjoint equivalence}.
\end{definition}
If we draw $\eta$ and $\epsilon^{-1}$ as caps, and $\eta^{-1}$ and $\epsilon$ as cups
    \[
        \myinput{tikz_Chapter5_adjoints_counit}\qquad \myinput{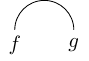}\qquad
        \myinput{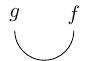}\qquad \myinput{tikz_Chapter5_adjoints_unit}
    \]
    then $(f,g,\eta,\epsilon)$ defines an equivalence if and only if the following identities hold.
    \[
        \myinput{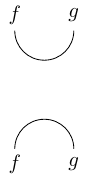}\quad=\quad\myinput{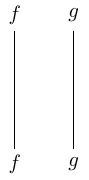}\qquad\myinput{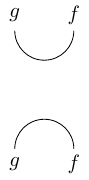}\quad=\quad\myinput{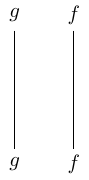}
    \]
    \[
        \myinput{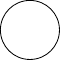}\quad=\quad
        \myinput{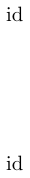}
    \]
\begin{proposition}
    If $(f,g,\eta,\epsilon)$ defines an adjoint equivalence then so does $(g,f,\epsilon^{-1},\eta^{-1})$. In particular if $f$ and $g$ give an adjoint equivalence between $A$ and $B$ then $f$ is left and right adjoint to $g$.
\end{proposition}
\begin{proof}[Sketch proof]
The fact that $(g,f,\epsilon^{-1}, \eta^{-1})$ gives an equivalence follows immediately from the fact that $\epsilon$ and $\eta$ are isomorphisms. To prove that this is an adjoint equivalence we need to show the zigzag identities. For $g$ this follows from the equalities between the following diagrams.
\begin{align*}
    \myinput{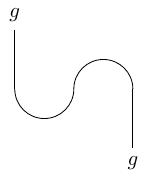}\;=\;\myinput{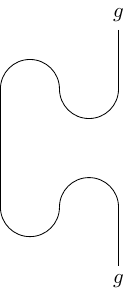}\;=\;\myinput{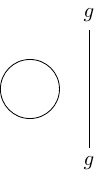}\;=\;\myinput{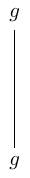}
\end{align*}
The first equality is by the zigzag identities for $f$ and $g$ and the following equalities follow from the fact that $\epsilon^{-1}\circ \epsilon=\iota_{f\circ g}$ and $\eta^{-1}\circ \eta=\iota_{\id}$. The proof of the zigzag identity for $f$ follows similarly.
\end{proof}
\begin{proposition}
    \label{prop:adjEquiv}
    If $(f,g,\eta,\epsilon)$ is an equivalence then there is an $\epsilon'$ such that $(f,g,\eta,\epsilon')$ is an adjoint equivalence.
\end{proposition}
\begin{proof}[Sketch Proof]
We simply define $\epsilon'$ to be the following 2-cell.
\[
    \myinput{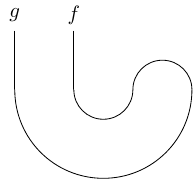}
\]
For more details, see for example, Johnson and Yau's~\cite[prop.~6.2.4]{yau2021bicat} book.
\end{proof}
Suppose that we have a sequence of 1-cells, $f_1,...,f_n$, and a sequence of 1-cells, $g_1,...,g_m$, and a 2-cell between them.
\[ 
    \fbox{\tikzsetnextfilename{Chapter5/adjoints/adjointEquivs/1}
\begin{tikzpicture}[baseline={(current bounding box.center)}, xscale=1
]

\drawStartNodes{
    f1/$f_1$/-1,
    dots/${\ldots}$/-1,
    fn/$f_n$/0
}
\strTwoCellX{fn;f1}{$\phi$}{gm;g1}[][2][1][1]

\place{dots}{1}{gm}

\drawEndNodes{
    gm/$g_m$,
    dots/${\ldots}$,
    g1/$g_1$
}

\end{tikzpicture}
}
\]
If each of $f_1,...,f_n$ and $g_1,...,g_m$ are adjoint equivalences we have another 2-cell given by the following diagram.
\[
    \fbox{\tikzsetnextfilename{Chapter5/adjoints/adjointEquivs/2}
\begin{tikzpicture}[baseline={(current bounding box.center)}, xscale=1
]

\drawStartNodes{
    f1/$g^\bullet_m$/-1,
    dots/$\ldots$/-1,
    fn/$g^\bullet_1$/0
}
\strTwoCellX{fn;f1}{$\phi^\bullet$}{gm;g1}[][2][1][1]

\place{dots}{1}{gm}

\drawEndNodes{
    gm/$f^\bullet_1$,
    dots/$\ldots$,
    g1/$f^\bullet_m$
}

\end{tikzpicture}
} \quad \coloneqq\quad \fbox{\tikzsetnextfilename{Chapter5/adjoints/adjointEquivs/3}
\begin{tikzpicture}[baseline={(current bounding box.center)}, xscale=1
]

\drawStartNodes{
    gm/$g^\bullet_m$/-1,
    dots/$\ldots$/-1,
    g1/$g^\bullet_1$/0
}
\strIdX{g1}[][3]
\strIdX{gm}[][3]
\place{f1}{-1}{g1}
\place{dots}{-1}{f1}
\node[blank] (dots) at (dots){$\ldots$};

\place{fn}{-5}{gm}
\place{f1A}{-5}{f1}
\place{fnA}{-1}{fn}

\place{dots}{-1}{fnA}
\node[blank] (dots) at (dots){$\ldots$};

\strCapX{f1A}{f1}[][2.5]
\strCapX{fnA}{fn}
\node[blank] (dots) at (dots){$\ldots$};

\strTwoCellX{fn;f1}{$\phi$}{gmA;g1A}[][2][1][1]
\place{dots}{1}{gmA}
\node[blank] (dots) at (dots){$\ldots$};

\strIdToX{g1}{g1A}
\strIdToX{gm}{gmA}

\strCupX{g1A}{g1}
\strCupX{gmA}{gm}[][2.5]
\place{dots}{1}{g1}
\node[blank] (dots) at (dots){$\ldots$};

\strIdX{f1A}[][5]
\strIdX{fnA}[][5]

\place{dots}{1}{f1A}

\drawEndNodes{
    f1A/$f^\bullet_1$,
    dots/$\ldots$,
    fnA/$f^\bullet_n$
}

\end{tikzpicture}}
\]
This might be thought of as the 2-cell $\phi$ rotated 180 degrees. Of course, if $\phi$ has an inverse then we can draw the inverse as in the following diagram.
\[
    \fbox{\tikzsetnextfilename{Chapter5/adjoints/adjointEquivs/4}
\begin{tikzpicture}[baseline={(current bounding box.center)}, xscale=1
]

\drawStartNodes{
    f1/$g_1$/-1,
    dots/${\ldots}$/-1,
    fn/$g_m$/0
}
\strTwoCellX{fn;f1}{$\phi^{-1}$}{gm;g1}[][2][1][1]

\place{dots}{1}{gm}

\drawEndNodes{
    gm/$f_n$,
    dots/${\ldots}$,
    g1/$f_1$
}

\end{tikzpicture}
}
\]
This might be thought of as the 2-cell $\phi$ reflected in the horizontal axis. If each of $f_1,...,f_n$ and $g_1,...,g_m$ are adjoint equivalences and $\phi$ is invertible then we can both reflect and rotate a 2-cell to get a new 2-cell $(\phi^\bullet)^{-1}=(\phi^{-1})^\bullet$ that we denote $\phi^\invbullet$, as in the following diagram.
\[
    \fbox{\tikzsetnextfilename{Chapter5/adjoints/adjointEquivs/5}
\begin{tikzpicture}[baseline={(current bounding box.center)}, xscale=1
]

\drawStartNodes{
    f1/$f^\bullet_n$/-1,
    dots/${\ldots}$/-1,
    fn/$f^\bullet_1$/0
}
\strTwoCellX{fn;f1}{$\phi^{\invbullet}$}{gm;g1}[][2][1][1]

\place{dots}{1}{gm}

\drawEndNodes{
    gm/$g_1^\bullet$,
    dots/${\ldots}$,
    g1/$g_m^\bullet$
}

\end{tikzpicture}}
\]
This might be thought of as the 2-cell $\phi$ reflected in the vertical axis.

Since we will be working with monoidal bicategories in later chapters, it is worth spelling out what we mean by a cartesian product of bicategories.
\begin{definition}
    Given bicategories $\AA$ and $\BB$ we define the \textdef{cartesian product} of $\AA$ and $\BB$, to be the bicategory, $\AA\times \BB$, whose
    \begin{itemize}
        \item objects are pairs $(A,B)$ for $A\in \AA$ and $B\in \BB$;
        \item hom-categories are given by
        \[
        (\AA\times \BB)((A,B), (A',B'))\coloneqq \AA(A,A')\times \BB(B,B');
        \]
        \item composition functor is given by
        \[
        \begin{tikzcd}[ampersand replacement=\&, cramped, every matrix/.append style={nodes={font=\scriptsize}} ]
            \begin{cdaligned}
                \AA(A',&A'')\times \BB(B',B'')\\
                &\times \AA(A,A')\times \AA(B,B')
            \end{cdaligned}
            \arrow[r, "\sim"]
            \& 
            \begin{cdaligned}
                \AA(A,&A')\times \AA(A, A')\\
                &\times \BB(B',B'')\times \BB(B,B')
            \end{cdaligned}
            \arrow[r, "{(\circ, \circ)}"]
            \&\AA(A,A'')\times \BB(B,B'');
        \end{tikzcd}
        \]
        \item identity 1-cells are given by $(\id,\id)\in \BB(A,A)\times \BB(B,B)$;
        \item composition associator, composition left unitor, and composition right unitor are given by
        \[
            (\comp{\alpha}, \comp{\alpha}), (\comp{\lambda},\comp{\lambda})\text{ and } (\comp{\rho}, \comp{\rho}).
        \]
    \end{itemize}
\end{definition}
In other words, objects are pairs of objects, 1-cells are pairs of 1-cells, 2-cells are pairs of 2-cells and composition is given pairwise.
\section{The 2-Trace}
\label{section:bicategoricaltrace}
Adjunctions are ubiquitous in category theory, and arguably form the foundation of the theory itself. But in true category theory style, the definition of an adjunction is lifted from elsewhere in mathematics -- namely linear algebra. Two linear maps, $F\colon V\rightleftarrows 
W\cocolon G$, between Hilbert spaces are adjoint whenever
\[
    \langle  Fv, w\rangle_{W}=\langle  v,Gw\rangle_V
\]
for all $v\in V$, $w\in W$, where $\langle-,- \rangle$ denotes the inner product on each space. Of course, two functors $F\colon \AA\rightleftarrows \BB\cocolon G$ are adjoint whenever
\[
    \BB(FA,B)\cong \AA(A,GB)    
\] 
naturally for all $A\in \AA$, $B\in \BB$.
Clearly these two definitions are formally similar. But this analogy runs much deeper. To understand the analogy we firstly need to understand some facts about Cauchy complete categories.
\begin{definition}
    Given a $\VV$-category $\DD$, a $\VV$-functor $F\colon \DD\to \CC$ and a copresheaf $W\colon \DD\to \VV$, the \textdef{weighted limit} $\lim^W F$, if it exists, is an object $\lim^W F$ such that there is a natural isomorphism
    \[
        \CC(c,\wlim^WF)\cong \hom{\DD,\VV}(W, \CC(c,F(-))).
    \]
    Given a $\VV$-category $\DD$, a functor $F\colon \DD\to \CC$ and a presheaf $W\colon \DD^{\op}\to \VV$, the \textdef{weighted colimit} $\colim^W F$, if it exists, is an object $\colim^W F$ such that there is a natural isomorphism
    \[
        \CC(\colim^WF, c)\cong \hom{\DD^{\op},\VV}(W, \CC(F(-),c)).
    \]
\end{definition}

    \begin{definition}
        A \textdef{conical limit} is a limit weighted by the copresheaf which is constant at the unit. In the case that $\VV$ is the category of sets, conical limits are the usual limits of category theory.

        A \textdef{conical colimit} is a colimit weighted by the presheaf which is constant at the unit of $\VV$. In the case that $\VV$ is the category of sets, conical limits are the usual colimits of category theory.
    \end{definition}
    \begin{definition}
            A weighted limit or colimit is called \textdef{absolute} if it is preserved by every functor.
    \end{definition}
    This definition of a Cauchy complete category is different, but equivalent, to the one initially suggested by Lawvere. Lawvere~\cite[p.~138]{lawvere1973} famously showed that every metric space can be thought of as an enriched category. In the same paper Lawvere~\cite[p.~163]{lawvere1973} also proved that a metric space $X$ is Cauchy complete (in the sense of metric spaces) if and only if the profunctors 
    \[
        X(x,-)\colon X\profto * \text{ and } X(-,x)\colon *\profto X
    \] 
    form an adjoint pair.

    A proof of Street~\cite[]{StreetAbsolute} shows that any limit, $\lim^W F$, weighted by a presheaf $W\colon \CC\to \Set$ is absolute if and only if the profunctor $W\colon \CC\profto *$ has a right adjoint, $\overline{W}$. Note that in such a case $\lim^W F\cong \colim^{\overline{W}} F$, and so this shows that absolute weighted limits and absolute weighted limits colimits coincide.
    \begin{definition}
    A category is called \textdef{Cauchy complete} if it has all absolute weighted limits and colimits. That is to say, it has any limit weighted by a presheaf which has a right adjoint as a profunctor. 
    \end{definition}
    In our analogy these absolute weighted (co)limits will play the role of weighted sums. Note that this is slightly at odds to the usual analogy which compares limits to products and colimits to sums. 

To highlight their similarities, given an indexing set $I$ and a weighting function $f\colon I\to \mathbb{C}$ we will write $\sum^f_{i\in I} x_i\coloneqq \sum_{i\in I} f(i)\cdot x_i$. We will write $\overline{f}$  for the pointwise conjugate of $f$. Consider the natural isomorphism between $\CC$ and the Hilbert space of linear endomorphisms from $\mathbb{C}$ to itself, $\operatorname{Hilb}(\mathbb{C},\mathbb{C})$. Every $\lambda\in \mathbb{C}$ has an adjoint given by its complex conjugate, $\overline{\lambda}$. So we will write $\overline{W}$ for the profunctor right adjoint to $W$. The idea behind this analogy explained by \autoref{table}.
\begin{figure}[h]
\renewcommand{\arraystretch}{1.5}
\begin{center}
    \begin{tabular}{c c c}
    \toprule
    &\textbf{Hilbert space $V$}
        &\textbf{Cauchy complete category $\CC$}\\
    \midrule
    Scalar
        & $c\in \mathbb{C}$
            & $S\in \Set$\\
    \midrule
    $0$-cell
        & Point $x$
            & Object $A$\\
    \midrule
    \multirow{2}{*}{Weighted sum}
        & Weighted sum
            &Absolute weighted (co)limit\\
        &$\sum^F_{i\in I} x_i$
            &$\lim^W_{i\in I} A_i$\\
    \midrule
    Sum size
        & Finite
            &Small\\
    \midrule
    1-cell
        & $\langle  x,y\rangle_V\in \mathbb{C}$
            & $\CC(A,B)\in \Set$\\
    \midrule
    \multirow{2}{*}{Sesquilinearity}
        &$\sum^g_{j\in J}\sum^f_{i\in I}\langle  v_j,w_i\rangle_V$
            &$\lim^V_{j\in J}\lim^W_{i\in I}\CC(A_j,B_i)$\\
        &$\qquad=\langle \sum^{\overline{g}}_{j\in J}v_j, \sum^f_{i\in I} w_i\rangle_V$
            &$\qquad\cong \CC(\lim^{\overline{V}}_{j\in J}A_j, \lim^{W}_{i\in I}B_i)$\\
    \midrule
    Map
        &Linear Map
            &Functor\\
    \midrule
    Binary operator
        &Sesquilinear form
            &Profunctor\\
    \bottomrule
    \end{tabular}
    \caption{The analogy between Hilbert spaces and Cauchy complete categories}
    \label{table}
\end{center}
\end{figure}

We tend to work simply with categories rather than Cauchy complete categories, since they occur more `naturally' but note that every category has a unique Cauchy completion, we take the subcategory of $\hom{\CC^{\op}, \Set}$ consisting only of presheaves that have adjoints in the bicategory of profunctors. See for example, Borceux and Dejean's~\cite[sec.~4]{borceux1986} survey on Cauchy complete categories. Similarly, we could work with an appropriate notion of `partial vector space' since any such space would be uniquely completable, but such objects are somehow much less natural than vector spaces.
\begin{remark}
Every Hilbert space has a natural metric structure induced by the inner product, and by definition the space is Cauchy complete. So every Hilbert space can be thought of as an enriched Cauchy complete category. Unfortunately this does not formalise the above analogy. Absolute weighted limits in metric spaces correspond directly to Cauchy limits, and have nothing to do with sums. In fact, it seems unlikely -- at least to this author -- that there is any way to recast inner products as enriched homs, or sums as absolute weighted limits.
\end{remark}

There is more to be said about this analogy than is immediately apparent. In both cases we also have a form of generalised weighted sum. In the case of Hilbert spaces we can take certain \emph{infinite} weighted sums. These are not preserved by arbitrary linear maps, but a linear map that preserves all infinite sums is exactly a continuous linear map, or equivalently a bounded linear map. And it is known that the bounded linear maps are exactly those with adjoints (see, for example, Debnat and Mikusi\'nski's~\cite[thm.~1.5.7]{debnath1999hilbert} book). 

In the case of categories we can take certain non-absolute weighted limits and weighted colimits. These are not preserved by arbitrary functors but rather by functors which have an adjoint.

It is worth noting that weighted sums in Hilbert spaces are generated by two separate processes: unweighted sums of vectors, and scalar products. This corresponds to the fact that, in an enriched category every weighted limit can be expressed using a combination of conical limits and particular weighted limits called powers. This can be found in, for example, Kelly's~\cite[thm.~3.73]{kelly1982basic} monograph. Thus, we can think of conical limits -- or conical colimits -- as being unweighted sums, and powers -- or copowers -- as being scalar products.

Delving deeper into linear algebra and category theory yields further comparisons, between the Riesz representation theorem and Freyd's representability theorem for example. But we are particularly interested in the case of `small' categories and Hilbert spaces. For categories `small' means small as in set-sized. For Hilbert spaces `small' means finite dimensional. 

The category of finite dimensional Hilbert spaces, $\FDHilb$ can be given a lot of extra structure. It can be given the structure of a compact-closed category, or that of a dagger category. For now, we're interested in the closed monoidal structure on $\FDHilb$. The finite dimensional Hilbert space of linear maps from $V$ to $W$ has inner product given by
\[
    \langle F,G\rangle\coloneqq \Tr(\overline{F}\circ G)
\]
where $\overline{G}$ is the adjoint of $G$. One consequence of this is that we can describe the trace of a linear endomorphism $E\colon A\to A$, or more generally a sesquilinear form, as the inner product
\[
    \langle \id,E\rangle.
\]
By analogy then, the 2-trace, as introduced by Ganter and Kapranov~\cite[def.~3.1]{ganter}, and Bartlett~\cite[def.~7.8]{brucethesis} is defined as follows.
\begin{definition}
Given an endomorphism $E\colon A\to A$ in some bicategory $\BB$, the \textdef{2-trace} is given by the set of 2-cells
\[
\BB(A,A)(\id,F).    
\] 
\end{definition}
 Note that this is distinct from the usual categorification of a trace, which relies on the compact-closed structure and duals. In \autoref{section:CompactClosed} will see how these two different perspectives on the trace interact. Most of the following examples are due to Willerton \cite{willertonTalk}.
\begin{example}
In $\Rel$, the identity relation on a set $A$ is given by the equality, or diagonal, relation
\[
    \Delta A=\{(a,a)\mid a\in A\}\subset A\times A,
\] so the 2-trace of an endo-relation $R\colon A\to A$ is $*$ if $R$ is reflexive and $\varnothing$ otherwise.
\end{example}
\begin{example}
The identity bimodule on an algebra $A$ is given by $A$ equipped with the bimodule structure given by the algebra multiplication. The 2-trace of an $A$-$A$-bimodule $M$ is the set of bimodule maps from $A$ to $M$. This is in natural bijection with the set of invariants, or what might be called the centre, of $M$:
\[
\{m\in M\mid \forall a\in A\; am=ma\}. 
\]
To see this, note that given a bimodule map $f\colon A\to M$ we know that $a\cdot f(1)= f(a\cdot 1)=f(a)=f(1\cdot a)=f(1)\cdot a$ and so $f(1)$ gives an invariant. Similarly, if $m$ is invariant then $f_m(a)\coloneqq m\cdot a=a\cdot m$ gives a bimodule map.
\end{example}
\begin{example}
The identity bimodule on an algebra $A$ is given by $A$ equipped with bimodule structure given by the algebra multiplication. The 2-trace of a DG-bimodule $M$ is then, by definition, given by the (underlying set of) $\mathrm{Ext}^\bullet(A,M)$. Cartan and Eilenberg~\cite[ch.~IX]{cartan1956homological} proved that this coincides with Hochschild \cite[]{hochschild} cohomology whenever $k$ is a field.
\end{example}
\begin{example}
The identity profunctor on a category $\AA$ is given by the Hom profunctor $\AA(-,-)$. The 2-trace of a profunctor $P\colon \AA\profto \AA$ is the set of natural transformations from $\AA(-,-)$ to $P$, given by the end
\[
    \endint_{A\in \AA}\VV(\AA(A,-), P(-,A)).
\]
By the Yoneda lemma this is equal to the underlying set of
\[
    \endint_{A\in \AA}P(A,A).
\]
\end{example}
\begin{example}
The identity span is given by the span $A\xleftarrow{\id} A \xrightarrow{\id} A$. The 2-trace of a span $A\xleftarrow{f} S \xrightarrow{g} A$ is given by the set of maps $\alpha\colon A\to S$ such that 
\[
    f\circ \alpha=g \circ \alpha=\id.
\]
In other words, it is the set of ``mutual sections'' of the span. 
\end{example}
\begin{example}
    The identity path at a point $x$ in a topological space $T$ is given by the constant loop. Then in $\Path(T)$, given a loop $p$ at $x$, the 2-trace is the set of homotopy classes of homotopies from the constant loop to the path $p$. In particular this is non-empty if and only if $p$ is nullhomotopic. 
\end{example}
It may be obvious from these examples that the 2-trace seems to lack some structure. In the first three cases it always happens to be that the 2-trace is the set underlying some other object: an $R$-module in the first case, a chain complex of $R$-modules in the second case, and a natural transformation object in $\VV$ in the third case. This is because our 2-trace lands in  the category of sets, but, like the linear trace, it ought to be a \emph{scalar} for our bicategory. We give an account of scalars in \autoref{section:Scalars}. In \autoref{section:cotrace} we will show that the 2-trace can always be replaced by a scalar called the cotrace. Furthermore, much as the linear trace gives rise to the Frobenius inner product, the cotrace gives rise to an enrichment of the bicategory in its category of scalars.
\section{Closed bicategories}
In classical category theory adjoint functors are often thought of as a sort of `weak inverse', and they underpin a lot of theory about algebraic structures. In this section we focus on a concept even weaker than adjoints -- Kan extensions. 

Typically, discussions of Kan extensions focus on functors, and such Kan extensions of functors are ubiquitous. As Mac Lane~\cite[p.~248]{maclane1971categories} remarks in \textbook{Categories for the Working Mathematician}: `the notion of Kan extensions subsumes all the other fundamental concepts of category theory'. But Kan extensions can be defined in any arbitrary bicategory. They will play an integral role in the scalar enrichment of our bicategories.
\begin{definition}
Let $\BB$ be a bicategory and let $f\colon A\to B$ and $g\colon A\to C$ be 1-cells. Then the \textdef{right extension} of $g$ along $f$, if it exists, consists of a 1-cell $g\extend f$ called the \textdef{extension} and a 2-cell $\eta$ called the \textdef{evaluation map}
\[
    \begin{tikzcd}
    C
    \arrow[rr, leftarrow, bend right=75, "g"{swap}, ""{name=target}]
        &B 
        \arrow[Rightarrow, to=target, "\eta"]
    \arrow[l, "g\extend f"{swap}]
            &A
            \arrow[l, "f"{swap}]
    \end{tikzcd}    
\]
such that for any other 1-cell $h\colon B\to C$ and 2-cell $\gamma\colon g\circ h\Rightarrow f$ there is a unique 2-cell $!\colon h\Rightarrow g\extend f$ such that the following two diagrams are equal.
\[
    \begin{tikzcd}
        C
        \arrow[rr, leftarrow, bend right=75, "g"{swap}, ""{name=target}]
            &B 
            \arrow[Rightarrow, to=target, "\gamma"]
        \arrow[l, "h"{swap}]
                &A
                \arrow[l, "f"{swap}]
    \end{tikzcd} 
    =   
    \begin{tikzcd}
    C
    \arrow[rr, leftarrow, bend right=75, "g"{swap}, ""{name=target}]
        &B 
        \arrow[l, bend right=75, "h"{name=source, swap}]
        \arrow[Rightarrow, to=target, "\eta"]
        \arrow[l, "g\extend f", ""{swap,name=target2}]
            &A
            \arrow[l, "f"]
            \arrow[Rightarrow, from=source, to=target2, "\;!"]
    \end{tikzcd}.
\]
\end{definition}
\begin{remark}
    We opt to use the lollipop, $\extend$, notation rather than any other notation for a number of reasons. Firstly it is one of the notations used for closed structures on monoidal categories, and as we are about to see extensions are a generalisation of that concept. Secondly the round end of the lollipop can be thought of like an arrowhead, indicating that $g\extend f$ extends from $f$ to $g$. Finally, the round end of the lollipop resembles the composition circle. This is useful for remembering that the evaluation map is given by $(g\extend f)\circ f\Rightarrow f$, since the $\circ$ always occurs on the opposite side to the $\circ$ at the end of the lollipop. Note also that the name of $g\extend f$ reads circle-side first, it is the extension of $g$ along $f$.
    \end{remark}
\begin{remark}
When dealing with extensions we will often opt to use pasting diagrams, rather than string diagrams. This is because generally speaking extensions, unlike adjoints, do not immediately lend themselves to a string diagram language. However, Baez and Stay~\cite[p.~30]{baez2010physics} have developed a string diagram language of clasps and bubbles for closed monoidal categories that would work equally well here.
\end{remark}
One example of extensions outside the 2-category of categories is extensions of profunctors. Given a functor $F\colon \mathscr{I}\to \AA$ and a profunctor $W\colon \mathscr{I}\profto \VV$, the weighted limit $\lim^{W}F\colon \VV\to \AA$ exists if and only if the following extension of profunctors
\[
    \begin{tikzcd}[column sep={5.4em}]
    \mathscr{V}
    \arrow[rr, leftarrow, "\shortmid"{marking}, bend right=75, "W"{swap, outer sep=2pt}, ""{name=target},]
        &\mathscr{A} 
        \arrow[Rightarrow, to=target, shorten <=10pt, shorten >=10pt, "\eta"]
    \arrow[l, "\shortmid"{marking}, "W\extend{\AA(F-,-)}"{swap,outer sep=2pt}]
            &\mathscr{I}
            \arrow[l, "{\AA(F-,-)}"{swap,outer sep=2pt}, "\shortmid"{marking}]
    \end{tikzcd}    
\]
exists and is representable, in which case $W\extend \AA(F-,-)=\AA(-,\lim^WF-)$. This simply follows from the definition of a weighted limit and the description of Kan extensions in the bicategory of profunctors given below.

As well as right extensions we also have left extensions. These are `vertically dual' to right extensions, in that a left extension in $\AA$ is a right extension in $\AA^{\operatorname{co}}$. We also have the notion of a right lift, that is `horizontally dual', given by taking a right extension in $\AA^{\op}$.
\begin{definition}
    Let $\BB$ be a bicategory and let $f\colon B\to C$ and $g\colon A\to C$ be 1-cells. Then the \textdef{right lift} of $g$ through $f$, if it exists, consists of a 1-cell $f\lift g$ called the \textdef{lift} and a 2-cell $\eta$ called the \textdef{evaluation map}
    \[
    \begin{tikzcd}
        C
        \arrow[rr, leftarrow, bend right=75, "g"{swap}, ""{name=target}]
            &B 
            \arrow[Rightarrow, to=target, "\eta"]
        \arrow[l, "f"{swap}]
                &A
                \arrow[l, "f\lift g"{swap}]
    \end{tikzcd} 
    \]
    such that for any other 1-cell $h\colon A\to B$ and 2-cell $\gamma\colon g\circ h\Rightarrow f$ there is a unique 2-cell $!\colon h\Rightarrow g\extend f$ such that the following two diagrams are equal.
\[
    \begin{tikzcd}
        C
        \arrow[rr, leftarrow, bend right=75, "g"{swap}, ""{name=target}]
            &B 
            \arrow[Rightarrow, to=target, "\gamma"]
        \arrow[l, "f"{swap}]
                &A
                \arrow[l, "h"{swap}]
    \end{tikzcd} 
    =   
    \begin{tikzcd}
    C
    \arrow[rr, leftarrow, bend right=75, "g"{swap}, ""{name=target}]
        &B 
        \arrow[Rightarrow, to=target, "\eta"]
        \arrow[l, "f"]
            &A
            \arrow[l, bend right=75, ""{name=source}]
            \arrow[l, "f\lift g", ""{swap,name=target2}]
            \arrow[Rightarrow, from=source, to=target2, "\;!"]
    \end{tikzcd}.
\]
\end{definition}
\begin{remark}
Once again the lollipop acts like an arrowhead indicating that the lift goes from $f$ to $g$, and once again the circles occur on opposite sides of the line in the evaluation map $f\circ f\lift g\Rightarrow f$. Similarly to extension we read this starting from the circle: $f\lift g$ denotes a lift of $g$ through $f$. In order to avoid confusion between extensions and lifts a useful mnemonic might be to note $\extend$ looks more like a lower-case `e' than $\lift$, and so $\extend$ denotes an \textbf{e}xtension.
\end{remark}
In this thesis we focus almost exclusively on \emph{right} lifts and extensions. As such, any reference to a ``\textdef{lift}'' may be assumed to be a reference to a right lift, and any reference to an ``\textdef{extension}'' may be assumed to be a reference to a right extension. In addition, whilst Kan extensions seem to occur more in `nature' than Kan lifts, it is right lifts that seem to make the most sense in our applications.
But it is worth noting that, conceptually, lifts are just extensions and by duality any results about lifts have an analogue for extensions.

It is well known that Kan extensions of functors and adjoint functors are intrinsically linked. If every right Kan extension along some functor $F$ exists, then precomposition with $F$ has a right adjoint. Similarly, if $F$ has a right adjoint $G$ then $G\circ H$ is the right Kan extension of $H$ along $F$ for any $H$. This is true inside any arbitrary bicategory.
\begin{proposition}
Let $\BB$ be a bicategory, let $X$, $A$ and $B$ be objects in $\BB$. If $g\colon A\to C$ is a 1-cell such that $f\lift g$ exists for all $f\colon B\to C$, then there is a functor
\[
-\lift g\colon \BB(B,C)^{\op}\to \BB(A,B).
\]
If $f\colon B\to C$ is a 1-cell such that $f\lift g$ exists for all $g\colon X\to B$, then there is a functor
\[
    f\lift -\colon \BB(A,C)\to\BB(A,B).
\]
Furthermore, the functor $f\lift -$ is right adjoint to the post-composition functor
\[
    f\circ- \colon \BB(A,B)\to \BB(X,B).
\]
\end{proposition}
\begin{proof}
    The argument is exactly the same as the argument given for Kan extensions of functors. See, for example, Riehl's~\cite[prop.~6.1.5]{riehl} book. 
\end{proof}
Clearly then, lifts are just a multi-object version of right closed structures for monoidal categories.
\begin{corollary}
If $f\colon B\to C$ is a 1-cell with a right adjoint $f^\dagger\colon C\to B$ then for any $g\colon A\to B$ the lift of $g$ through $f$ exists and $f\lift g=f^\dagger \circ g$. In particular there is a natural isomorphism $f\lift \id\cong f^\dagger$.
\end{corollary}
\begin{definition}
A bicategory is called \textdef{left-closed}, or \textdef{left-composition-closed} if it has all right lifts and \textdef{right-closed}, or \textdef{right-composition-closed}, if it has all right extensions. A bicategory that is both left- and right-closed is called \textdef{composition-closed}, or simply closed.
\end{definition}
\begin{remark}
We use the phrase ``composition-closed'' here to distinguish this closed structure from the monoidal-closed structure which will feature in later chapters.
\end{remark}
\label{section:CompClosed}
\begin{example}
The bicategory $\Rel$ is composition-closed. Given a relation $R\subseteq B\times C$ and a relation $S\subseteq A\times C$ the lift of $S$ through $R$ is given by
\[
    R\lift S=\{(a,b)\mid \forall c\in C,\; bRc\Rightarrow aSc\}=\{(a,b)\mid \forall c\in C,\; \neg bRc\lor aSc\}.
\]
\end{example}
\begin{remark}
Note that for all sets $A$ and $B$ there is a functor
\[
(-)^c\colon \Rel(A,B)\to \Rel(A,B)^{\op}    
\] 
that takes every relation $R$ to its complement $(A\times B)\setminus R$. The category $\Rel(A,B)$ has a monoidal structure given by taking unions and this complement functor gives $\Rel(A,B)$ the structure of a star-autonomous category. We also have a functor
\[
    (-)^T\colon \Rel(A,B)\to \Rel(B,A)   
\]
which takes a relation to its transpose. We can then give $R\lift S$ in terms of these two functors and the composition:
\begin{align*}
    R\lift S&=\{(a,b)\mid \forall c\in C\; \neg bRc\lor aSc\}\\
    &=\{(a,b)\mid \forall c\in C\; \neg (bRc\land \neg aSc)\}\\
    &=\{(a,b)\mid \neg \;\exists c\in C \;bRc\land \neg aSc\}\\
    &=(R\circ (S^{Tc}))^{Tc}.
\end{align*}
This is a closed structure similar to the closed structure for asymmetric star-autonomous categories, in the sense of Barr~\cite[]{Barr}. In other words, the closed structure is given in terms of a contravariant endofunctor, and the composition. Note that the relation $S^{Tc}$ is given by taking $S\lift (\Delta A)^c$, the lift of the anti-diagonal relation through $S$. So the composition-closed structure for $\Rel$ seems particularly structured, in that it is given by taking a kind of `horizontal dual', $(-)^T$ and a kind of `vertical dual' $(-)^c$.
\end{remark}
\begin{example}
    The bicategory $\Bim_R$ is composition-closed. Given a $B$-$C$-bimodule $M$ and an $A$-$C$-bimodule $N$ the set of right $C$-module homomorphisms $\Hom_C(M,N)$ can be given an $A$-$B$-bimodule structure by
    \[
    (a\cdot f\cdot b)(m)=a\cdot f(b\cdot m) 
    \]
    and this defines the right lift $M\lift N$. The right extension is defined similarly.
\end{example}
\begin{example}
    The bicategory $\DGBim_R$ is composition-closed. This follows from the derived tensor-hom adjunction
    \[
    \RHom(A\tensor^L B, C)\cong \RHom(A, \RHom(B,C)).    
    \]
    See, for example, Weibel's~\cite[thm.~10.8.7]{Weibel} textbook.
\end{example}
\begin{example}
The bicategory $\VV$-Prof is composition-closed. Given profunctors
$P\colon \BB\profto \CC$ and $Q\colon \AA\profto \CC$ the right lift is given by the coend
\[
  P\lift Q \cong \endint_{C\in \CC}  \VV(P(-,C), Q(-,C))
\]
and the right extension is given similarly.
\end{example}
\begin{remark}
    When $\VV$ is star-autonomous, the composition for $\VV$-Prof also resembles the closed structure for star-autonomous monoidal categories. Given $P\colon \BB\profto \CC$ and $Q\colon \AA\profto \CC$ the lift $P\lift Q$ is given by
    \begin{align*}
        \endint_{c\in \CC} \VV(P(-,C),Q(-,C))& \cong
        \endint_{c\in \CC} (P(-,C)\tensor Q(-,C)^*)^*\\ &\cong 
        \left(\endint^{c\in \CC} P(-,C)\tensor Q(-,C)^*\right)^*\\
        &\cong (P\circ Q^*)^*
    \end{align*}
    and in particular, $Q^*$ is given by taking the lift $Q^*\lift \Hom^*$, where $\Hom^*$ is the composite
    \[
      \AA\tensor \AA^{\op}\xrightarrow{\Hom^{\op}} \VV^{\op}\xrightarrow{(-)^*} \VV. 
    \]
    This is a generalisation of the case for $\Rel$, since the category $F\rightarrow T$ is star-autonomous with duals given by logical negation.
    \end{remark}
\begin{example}
    Day~\cite[prop.~4.1]{day1974spans} proved that $\Span(\CC)$ is composition-closed if and only if for every $x\in \CC$ the slice category, $\CC/x$, is cartesian closed. 
\end{example}
\begin{example}
    The bicategory $\Path(T)$ is composition-closed for any topological space $T$. Let $p$ be a path from $x$ to $y$, then the reverse path $\overline{p}$ is left and right adjoint to $p$. This follows from the fact that $p\circ \overline{p}$ and $\overline{p}\circ p$ are both contractible and so $p$ and $\overline{p}$ give an equivalence.
\end{example}

\chapter{Pseudofunctors and Pseudonatural Transformations}
There are several notions of a morphism between two bicategories. In this thesis we will be interested in the strong, but not strict, morphisms. These are pseudofunctors, where composition and identities are preserved up to invertible 2-cells. Between pseudofunctors live pseudonatural transformations and between pseudonatural transformations live modifications.

Trying to understand pseudofunctors, pseudonatural transformation and modifications via pasting diagrams can often be confusing and unintuitive. Since they form a tricategory they can be understood using surface diagrams, but 3-dimensional diagrams don't work particularly well on 2-dimensional paper. Instead, we develop a compromise approach. We encode the 3-dimensional data as colourings of strings. 

The first section is dedicated to defining pseudofunctors and pseudonatural transformations. We first give accounts of these via pasting diagrams before giving a decorated string diagram language to that allows us to interpret the axioms geometrically. The string diagram language represents pseudonaturality as a colour-changing braid.

The second section gives an account of whiskering for pseudonatural transformations. Whiskering of pseudonatural transformations has an impact on the colour changing of the braid. We also give a short proof that the usual braid equations, sometimes known as the Yang-Baxter equations, hold for our colour-changing braids.

The third section focuses on adjunctions internal to bicategories, and adjoint pairs of pseudonatural transformations. We show, using the decorated string diagram language, that a pair of pseudonatural transformations are adjoint if and only if their components are adjoint.

The fourth section gives the definition of pseudoadjunctions between bicategories, as well as some technical discussion about the necessity of coherence axioms.

The fifth and final section gives the definition of a biequivalence and the strictification theorem, which says that a bicategory is essentially the same thing as a 2-category. This gives a rigorous justification for the omission of associators and unitors from our diagrams.
\section{Decorated String Diagrams}
We begin with a definition of pseudofunctor, given in terms of pasting diagrams. This will ultimately be reinterpreted to give the foundation of our string diagram language.
\begin{definition}\label{def:pseudofunctor}
A \textdef{pseudofunctor} $F\colon \BB\to \CC$ between two bicategories consists of
\begin{itemize}
    \item a function \[
    F\colon \obj(\BB)\to \obj(\CC);
    \]
    \item for every pair of objects $A,B\in \BB$ a functor
    \[
    F_{A,B}\colon \BB(X,Y)\to \CC(F(X),F(Y));
    \]
    \item for every object $X$ an invertible 2-cell called the \textdef{identitor}
    \[
    F_X\colon \id_{F(X)}\Rightarrow F_{X,X}(\id_X); 
    \]
    \item for every triple of objects $X,Y,Z$ an invertible 2-cell, natural in $f$ and $g$, called the \textdef{compositor}
    \[
    F_{X,Y,Z}\colon F_{Y,Z}(g)\circ F_{X,Y}(f)\Rightarrow F_{X,Z}(g\circ f).
    \]
\end{itemize}
We will drop subscripts where they are clear from context. The above data adheres to the following axioms:
the identitor diagrams commute
\[
\begin{tikzcd}
F(f) \circ \id_{F(X)}
\arrow[r, "\comp{\lambda}"]
\arrow[d, "\iota_{F(f)}\circ F_X"{swap}]
    & F(f)\\
F(f)\circ F(\id_X)
\arrow[r, "F"{swap}]
    &
    F(f\circ \id_X)
    \arrow[u, "F(\comp{\lambda})"{swap}]
\end{tikzcd}
\]
\[
\begin{tikzcd}
\id_{F(Y)}\circ F(f)
\arrow[r, "\comp{\rho}"]
\arrow[d, "F_Y\circ \iota_{F(f)}"{swap}]
    & F\\
F(\id_Y)\circ F(f)
\arrow[r, "F"{swap}]
    &
    F(\id_Y\circ f)
    \arrow[u, "F(\comp{\rho})"{swap}]
\end{tikzcd}
\]
and the compositor diagram commutes
\[
\begin{tikzcd}
(F(h)\circ F(g))\circ F(f)
\arrow[r, "\comp{\alpha}"]
\arrow[d, "F\circ \iota_{F(f)}"{swap}]
    &F(h)\circ (F(g)\circ F(f))
    \arrow[d, "\iota_{F(h)}\circ F"]\\
F(h\circ g)\circ F(f)
\arrow[d, "F"{swap}]
    &F(h)\circ F(g\circ f)
    \arrow[d, "F"]\\
F((h\circ g)\circ f)
\arrow[r, "F(\comp{\alpha})"{swap}]
    &F(h\circ (g\circ f))
\end{tikzcd}
\]
\end{definition}
In order to keep track of pseudofunctors while working with string diagrams, we introduce new notational convention, inspired by the functorial boxes introduced by Melli\'es~\cite[]{mellies2006functorial}. This notational convention will rely on highlighting strings with various colours. The decision to use colours was made to avoid certain confusing ambiguities, but we realise that this may prove challenging for colour-blind readers. Due to size constraints, patterns could not be used in place of colours, however, in the interest of accessibility we have opted to use a palette based on the recommendation of Bang Wong~\cite[]{wong2011color}, and the interactive simulation tool provided by David Nichols~\cite[]{nichols}. The colours of this palette should be distinguishable for the majority of colour-blind readers.

Suppose that $F\colon \AA\to \BB$ is a pseudofunctor; $f\colon A\to B$ and $g\colon A\to B$ are 1-cells in $\AA$; and $\gamma\colon f\Rightarrow g$ is a 2-cell in $\AA$. We firstly assign a colour to our functor: let's say {\color{c1}$F$ is this colour}. The idea is that the functor ``wraps'' strings, which represent 1-cells, and beads, which represent 2-cells, in that colour.
\[
    \myinput{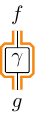}\quad\coloneq \quad
    \myinput{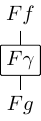}
\]
Note that this is compatible with the more traditional string diagram notation: we simply think of the identity functor as being the same as the background colour.

Suppose instead that we have a 2-cell $\xi$ in $\BB$, and it happens to be the case that the source of $\xi$ is $Ff$ and the target of $\xi$ is $Fg$, but $\xi$ isn't necessarily in the image of $F$. We wrap the string in the colour but leave the bead unwrapped.
\[
    \myinput{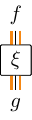}\quad \coloneq \quad      
    \myinput{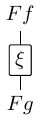}
\]
This indicates that the functor is applied to $f$ and $g$ but not to $\xi$. The identitor has an obvious expression in this language.
\[
    \myinput{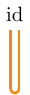}\quad \coloneq \quad      
    \myinput{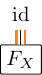}
\]
The compositor can be written in two difference ways.
\[
    \myinput{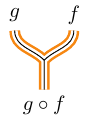}\quad \coloneq \quad \myinput{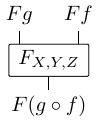}
\]
\[\myinput{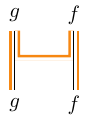}\quad \coloneq \quad \myinput{tikz_Chapter5_PFExample8}  
\]
Then for any 1-cells $f\colon A\to B$, $f'\colon A\to B$, $g\colon B\to C$, $g'\colon B\to C$; and any 2-cells $\alpha\colon f\Rightarrow f'$ and $\beta\colon g\to g'$, naturality of the compositor can be expressed by the following equality.
\[
    \myinput{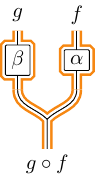}\quad=\quad
    \myinput{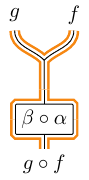}
\]
The unitor axioms can also be expressed by the following equalities.
\[
    \myinput{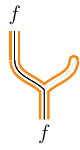}\quad = \quad \myinput{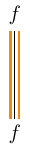} \quad = \quad \myinput{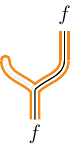}
\]
Using the proximity of strings to indicate bracketing, and using a change in proximity to denote the associator, we can write the compositor axiom as follows.
\[
    \myinput{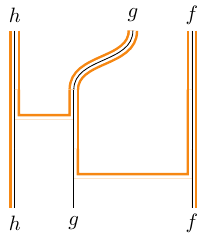} \quad=\quad  \myinput{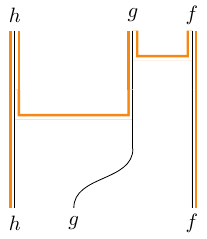}
\]
Or leaving the associator to be implicit we can instead write the following.
\[
   \myinput{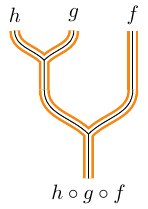}= \myinput{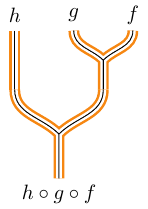}
\]
\begin{definition}
Suppose we have two pseudofunctors $F\colon \AA\to \BB$ and $G\colon \BB\to \CC$ then we can \textdef{compose} them to get $G\circ F\colon \AA\to \CC$ as follows. Firstly we compose their respective functions on objects to get
\[
    G\circ F\colon \obj(\AA)\to \obj(\CC);
\]
next we compose the respective 1-cell functors to get
\[
    (G\circ F)_X,Y\colon \AA(X,Y)\to \CC(G\circ F(X), G\circ F(Y));
\]
now we construct the identitor for $G\circ F$ as the composite 
\[
    \id_{GFX}\xrightarrow{G_{FX}} G(\id_{FX})\xrightarrow{GF_X}GF\id_X;
\]
and we construct the compositor for $G\circ F$ as the composite
\[
    G_{Y,Z}(F_{Y,Z}(g))\circ G(F_{X,Y}(f))\xrightarrow{G_{FX,FY,FZ}} G(F_{Y,Z}(g)\circ F_{X,Y}(f)) \xrightarrow{G(F_{X,Y,Z})} G(F(g\circ f)).
\]
\end{definition}
To convey a composite functor in our string diagram language we firstly wrap strings and beads in $F$ and then wrap the whole thing in $G$. For example, given a functor {\color{c1}$F\colon \AA\to \BB$ in this colour}; a functor {\color{c2}$G\colon \BB\to \CC$ in this colour}; a pair of 1-cells $f,g\colon A\to B$ in $\AA$; and a 2-cell $\gamma\colon f\Rightarrow g$ in $\AA$, the following two diagrams denote the same 2-cell.
\[
    \myinput{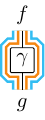} \coloneq \quad \myinput{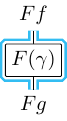}
\]
\begin{remark}
    Composing pseudofunctors is not strictly associative. Suppose that we had a functor ${\color{c3}H\colon \CC\to \DD}$. We choose the following convention. 
    \[
        \myinput{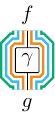}\quad\coloneq \quad   
        \myinput{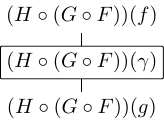}
    \]
    Whilst we \emph{could} use changes in the proximity of the wrappings to denote $(H\circ G)\circ F$, we have no need to here and doing so would greatly increase the size of our diagrams. In fact, one consequence of the coherence theorem that we see later is that we can treat pseudofunctors as if associativity holds on the nose.
\end{remark}
The identitor and compositor for $G\circ F$ can be expressed by the following diagrams.
\[
    \myinput{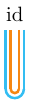}\quad \coloneqq\quad \myinput{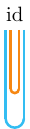}
\]
\[
    \myinput{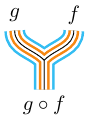}\quad \coloneqq\quad \myinput{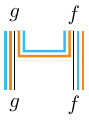}\quad\coloneqq\quad \myinput{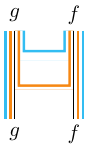} 
\]
\begin{definition}
Given two pseudofunctors $F,G\colon \BB\to \CC$, a \textdef{pseudonatural transformation} $n\colon F\Rightarrow G$ consists of:
\begin{itemize}
    \item a family of 1-cells $(n_A\colon FA\to GA)_A$ indexed by the objects of $\AA$
    \item a family of invertible 2-cells,
    \[
        \left(
    \begin{tikzcd}
    FA
    \arrow[r, "F(f)"]
    \arrow[d, "n_A"{swap}]
        &FB
        \arrow[dl, Rightarrow, "n_f"]
        \arrow[d, "n_B"]\\
    GA
    \arrow[r, "G(f)"{swap}]
        &GB
    \end{tikzcd}
    \right)_f
    \]
    called \textdef{naturalisors}, indexed by the 1-cells of $\AA$;
\end{itemize}
such that for every $f\colon A\to B$, $g\colon B\to C$, $f'\colon A\to B$ and $\gamma\colon f\Rightarrow f'$ the following composites are equal.
\[
\begin{tikzcd}[column sep=7.5ex, row sep=7.5ex]
FA
\arrow[d, "n_A"{swap}]
\arrow[r, "F(g\circ f)"]
    &FC
    \arrow[dl, Rightarrow, "n_{g\circ f}"{description}]
    \arrow[d, "n_C"]\\
GA
\arrow[r, "G(g\circ f)"{swap}]
    &GC
\end{tikzcd}=
\begin{tikzcd}[column sep=7.5ex, row sep=7.5ex]
FA
\arrow[d, "n_A"{swap}]
\arrow[r, "F(f)"]
    &FB
    \arrow[dl, Rightarrow, "n_{f}"{description}]
    \arrow[d, "n_B"{description}]
    \arrow[r, "F(g)"]
        &FC
        \arrow[dl, Rightarrow, "n_g"{description}]
        \arrow[d,"n_C"]
        \arrow[ll, bend right, "F(g\circ f)"{swap}, ""{name=top}]
        \\
GA
\arrow[r, "G(f)"{swap}]
    &GB
    \arrow[r, "G(g)"{swap}]
        &GC
        \arrow[ll, bend left, "G(g\circ f)", ""{swap,name=bottom}]
        \arrow[phantom, from=1-2, to=top, "\cong"]
        \arrow[phantom, from=2-2, to=bottom, "\cong"]
\end{tikzcd}
\]
\[
    \begin{tikzcd}[column sep=10ex, row sep=10ex]
        FA
        \arrow[r, "F(f)", ""{swap, name=f}]
        \arrow[d, equal]
            &FB
            \arrow[d,equal]\\
        FA
        \arrow[r, "F(f')"{description}, ""{name=ff}]
        \arrow[Rightarrow, "F(\gamma)"{description}, from=f, to=ff]
        \arrow[d, "n_A"{swap}]
            &FB
            \arrow[dl, Rightarrow, "n_{f'}"{description}]
            \arrow[d, "n_B"]\\
        GA
        \arrow[r, "G(f')"{swap}]
            &GB
        \end{tikzcd}
        =
        \begin{tikzcd}[column sep=10ex, row sep=10ex]
        FA
            \arrow[r, "F(f)"]
            \arrow[d, "n_A"{swap}]
                &FB
                \arrow[dl, Rightarrow, "n_f"{description}]
                \arrow[d, "n_B"]\\
            GA
            \arrow[r, "G(f)"{description}, ""{swap, name=f}]
            \arrow[d, equal]
                &GB
                \arrow[d,equal]\\
            GA
            \arrow[r, "G(f')"{swap}, ""{name=ff}]
            \arrow[Rightarrow, "G(\gamma)"{description}, from=f, to=ff]
                &GB
        \end{tikzcd}
\]
\[
    \begin{tikzcd}[column sep=10ex, row sep=10ex]
        FA
        \arrow[r, "F(\id_A)"]
        \arrow[d, "n_A"{swap}]
            &FA
            \arrow[dl, Rightarrow, "n_{id_A}"{description}]
            \arrow[d, "n_A"]\\
        GA
        \arrow[r, "G(\id_A)"{description}, ""{swap, name=id}]
        \arrow[d, equal]
            &GA
            \arrow[d, equal]\\
        GA
        \arrow[r, "\id_{GA}"{swap}, ""{name=id2}]
        \arrow[Rightarrow, from=id, to=id2, "G_A"{description}]
            &GA
    \end{tikzcd}
    =
    \begin{tikzcd}[column sep=10ex, row sep=10ex]
        FA
        \arrow[r, "F(\id_A)", ""{swap, name=id}]
        \arrow[d, equal]
            &FA
            \arrow[d, equal]\\
        FA
        \arrow[r, "\id_{FA}"{description}, ""{name=id2}]
        \arrow[Rightarrow, "F_A"{description}, from=id, to=id2]
        \arrow[d, "n_A"{swap}]
        \arrow[dr, "n_A"{description}, ""{name=mid}, ""{swap, name=mid2}]
            &FA
            \arrow[Rightarrow, to=mid, "\sim", sloped]
            \arrow[d, "n_A"]\\
        GA
        \arrow[Leftarrow, to=mid2, "\sim", sloped]
        \arrow[r, "\id_{GA}"{swap}, ""{name=id2}]
            &GA
    \end{tikzcd}
\]
\end{definition}
Pseudonatural transformations also fit into our string diagram language. Suppose we have two pseudofunctors $F,G\colon \AA\to \BB$, and we assign {\color{c1} $F$ this colour} and {\color{c2} $G$ this colour}. We draw a pseudonatural transformation $n\colon {\color{c1}F}\Rightarrow {\color{c2}G}$ as a dashed string highlighted on the right by the colours for {\color{c1}$F$} and highlighted on the left by the colours for {\color{c2}$G$}.
\[
    \myinput{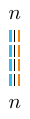}\quad\coloneqq\quad \myinput{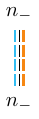}\quad\coloneq\quad \myinput{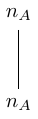}
\]
For any given 1-cell $f$, the naturalisor is drawn as a braid, with $F(f)$ passing underneath $n$. The wrapping around a 1-cell $f$ changes from the colour of ${\color{c1}F}$ to the colour of ${\color{c2}G}$ as it passes underneath.
\[
    \myinput{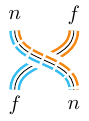}\quad\coloneqq\quad \myinput{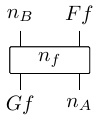}
\]
\[
    \myinput{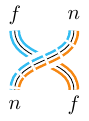}\quad\coloneqq \quad \myinput{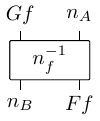}
\]
Note then that only a string wrapped in {\color{c1}this colour} can pass underneath a dashed line of {\color{c1}this colour}.
Given any 2-cell, $\gamma\colon f\Rightarrow g$, in $\AA$ the naturalisor axioms are then given by the following equalities.
\[
    \myinput{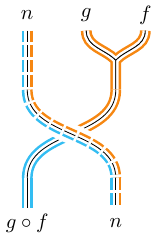}=
    \myinput{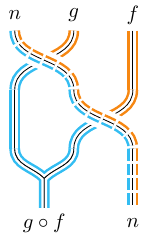}
\]
\[
    \myinput{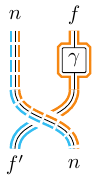}=
    \myinput{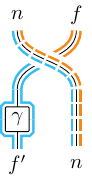}
\]
\[
    \myinput{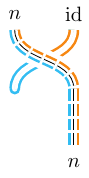}=
    \myinput{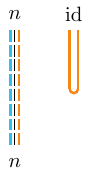}
\]
Drawing the naturalisor as a braid is not a new concept: previous authors have given similar accounts. However, including the colour wrappings helps to make clear when a particular braiding is valid.

Intuitively anything wrapped in the right-hand colour can pass underneath the braid and be dyed the left-hand colour. This is true from any direction, as proven by the following three propositions.

\begin{proposition}
    \label{prop:altNatAxioms}
    Each of the following equalities is equivalent to the first naturalisor axiom.
    \[
        \myinput{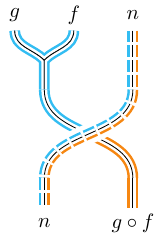}=\myinput{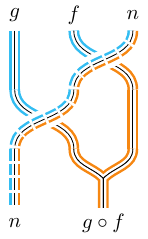}
    \]
    \[
        \myinput{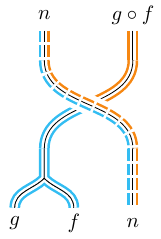}=\myinput{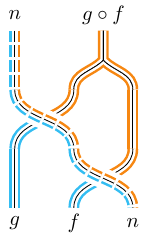}
    \]
    \[
        \myinput{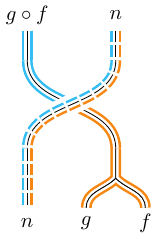}=\myinput{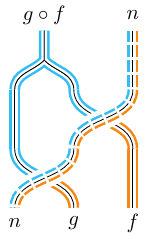}
    \]
    \end{proposition}
    \begin{proof}
    Each of these equivalences essentially follows from the fact that every 2-cell in the axiom has an inverse. For example, suppose that the first naturalisor axiom holds. Then the first equality above must hold by the following string of equalities:
    \[
        \myinput{tikz_PNat_Axioms_1A}=\fbox{\tikzsetnextfilename{PNat/Axioms/AxiomsXFlipProof/1}
\begin{tikzpicture}[baseline={(current bounding box.center)}, xscale=1
]

\drawStartNodes{
    n/$n$/-1,
    f/$f$/-1,
    g/$g$/0
}
\strBraidXX{n}{f}[/c2/c1][c2]
\strIdX{g}[c2]
\strBraidXX{n}{g}[/c2/c1][c2]
\strIdX{f}[c1][2]
\strBraidXX{n}{g}[/c2/c1][c1]
\strBraidXX{n}{f}[/c2/c1][c1]
\strIdX{g}[c2]

\strIdX{n}[/c2/c1]
\strJoinX{g/c2;f/c2}{gof}

\strBraidXX{n}{gof}[/c2/c1][c2][1.5]

\strIdX{n}[/c2/c1][0.5]
\strIdX{gof}[c1][0.5]

\drawEndNodes{
    gof/$g\circ f$,
    n/$n$
}

\end{tikzpicture}}=\fbox{\tikzsetnextfilename{PNat/Axioms/AxiomsXFlipProof/2}
\begin{tikzpicture}[baseline={(current bounding box.center)}, xscale=1
]

\drawStartNodes{
    n/$n$/-1,
    f/$f$/-1,
    g/$g$/0
}
\strBraidXX{n}{f}[/c2/c1][c2]
\strIdX{g}[c2]
\strBraidXX{n}{g}[/c2/c1][c2]
\strIdX{f}[c1]
\strJoinX{g/c1;f/c1}{gof}
\strIdX{n}[/c2/c1]
\strBraidXX{n}{gof}[/c2/c1][c1][1.5]
\strBraidXX{n}{gof}[/c2/c1][c2][1.5]

\strIdX{n}[/c2/c1][0.5]
\strIdX{gof}[c1][0.5]

\drawEndNodes{
    gof/$g\circ f$,
    n/$n$
}

\end{tikzpicture}}=\myinput{tikz_PNat_Axioms_2A}
    \]
    where the first equality holds by the fact that we are adding the identity 2-cell, the second equality holds by the first naturalisor axiom, and the third equality holds since we are removing the identity 2-cell. The other proofs follow in a similar fashion.
    \end{proof}
    \begin{proposition}
    \label{prop:altNatAxioms1}
        The following equality is equivalent to the second naturalisor axiom.
        \[
            \myinput{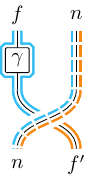}=\myinput{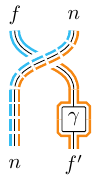}
        \]
    \end{proposition}
    \begin{proof}
    The proof follows similarly to the above proposition.
    \end{proof}
    \begin{proposition}
    \label{prop:altNatAxioms2}
        Each of the following equalities is equivalent to the third naturalisor axiom.
        \[
            \myinput{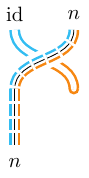}=\myinput{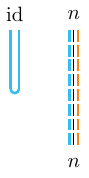}
        \]
        \[
            \myinput{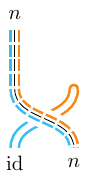}=\myinput{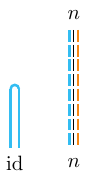}
        \]
        \[
            \myinput{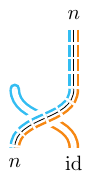}=\myinput{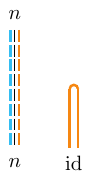}
        \]
        \end{proposition}
        \begin{proof}
            Again, this follows similarly to the above.
        \end{proof}
As with natural transformations, the pointwise composite of pseudonatural transformations yields another pseudonatural transformation. Suppose we have pseudofunctors \[{\color{c1}F},{\color{c2}G},{\color{c3}H}\colon \AA\to\BB\] and pseudonatural transformations, $n\colon {\color{c1}F}\Rightarrow {\color{c2}G}$, $m\colon {\color{c2}G}\Rightarrow {\color{c3}H}$. We can define the naturalisor for $m_A\circ n_A$ as follows,
\[
    \myinput{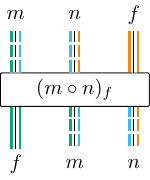}\coloneqq
    \myinput{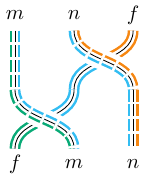}    
\]
and clearly all of the necessary axioms hold since they hold for the two individual braidings. This also means that the decorations for pseudonatural transformations are compatible with composition, in that it makes perfect sense to draw the following diagram for the identity 2-cell from $m\circ n$ to $m\circ n$.
\[
        \myinput{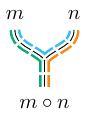}
\]
\begin{definition}
    Given a pair of pseudofunctors $F\colon \AA\to \BB$ and $G\colon \AA\to \BB$ and a pair of pseudonatural transformations $n\colon F\Rightarrow G$ and $m \colon F \Rightarrow G$, a \textdef{modification} $\mu\colon F \Rrightarrow G$ consists of a collection of 2-cells $(\mu_A\colon n_A\Rightarrow m_A)_A$ indexed by the objects of $\AA$ such that for any $f\colon A\to B$ in $\AA$ the following composite 2-cells are equal.
\[\begin{tikzcd}[column sep=7.5ex, row sep=7.5ex]
	FA & FA & FB \\
	GA & GA & GB
	\arrow[from=1-2, to=1-3, "F(f)"]
	\arrow[""{name=0, anchor=center, inner sep=0}, from=1-2, to=2-2, "n_A"{description}]
	\arrow[""{name=1, anchor=center, inner sep=0}, from=1-1, to=2-1, "m_A"{swap}]
	\arrow[equal, from=1-1, to=1-2]
	\arrow[equal, from=2-1, to=2-2]
	\arrow[from=2-2, to=2-3, "G(f)"{swap}]
	\arrow[from=1-3, to=2-3, "n_B"]
	\arrow["{n_F}"{swap}, shorten <=8pt, shorten >=8pt, from=1-3, to=2-2, Rightarrow]
	\arrow["{\mu_A}"', shorten <=12pt, shorten >=12pt, Rightarrow, from=0, to=1]
\end{tikzcd}
        =
        \begin{tikzcd}[column sep=7.5ex, row sep=7.5ex]
            FA & FB & FB \\
            GA & GB & GB
            \arrow[equal,from=1-2, to=1-3]
            \arrow[""{name=0, anchor=center, inner sep=0}, from=1-3, to=2-3, "n_B"]
            \arrow[""{name=1, anchor=center, inner sep=0}, from=1-2, to=2-2, "m_B"{description}]
            \arrow[from=1-1, to=1-2, "F(f)"]
            \arrow[from=2-1, to=2-2, "G(f)"{swap}]
            \arrow[from=2-2, to=2-3, equal]
            \arrow[from=1-1, to=2-1, "m_A"{swap}]
            \arrow["{m_F}"{swap}, shorten <=8pt, shorten >=8pt, from=1-2, to=2-1, Rightarrow]
            \arrow["{\mu_B}"', shorten <=12pt, shorten >=12pt, Rightarrow, from=0, to=1]
        \end{tikzcd}
    \]
\end{definition}
Suppose that {\color{c1} $F$ is this colour} and {\color{c2} $G$ is this colour}. We will draw modifications as beads with dashed lines, highlighted by the colours of their source and target pseudonatural transformations.
\[
\myinput{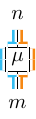}\quad \coloneqq\quad
\myinput{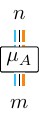}\quad\coloneqq \quad\myinput{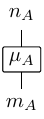}
\]
This means that we can write the modification axiom in the following way.
\[
    \myinput{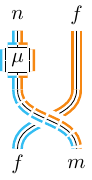}\quad =
\quad    \myinput{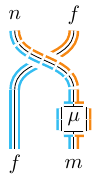}
\]
Intuitively, we can ``slide'' a modification over a braid because of the protective bubble around it, but we can't necessarily do this with an arbitrary 2-cell. As with sliding under the braid, this is true in any direction.
\begin{proposition}
The following equality is equivalent to the modification axiom.
\[
    \myinput{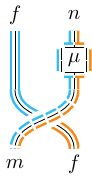}=
    \myinput{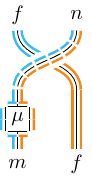}
\]
\end{proposition}
\begin{proof}
This proof follows analogously to the above propositions.
\end{proof}

Like pseudonatural transformations modifications behave well with respect to composition. Suppose that $n,m,l\colon {\color{c1}F}\Rightarrow {\color{c2}G}$ are pseudonatural transformations, $\mu \colon n\Rrightarrow m$ and $\nu\colon m\Rrightarrow l$. If we take the pointwise, vertical composite of two modifications we have the following equalities, and so $\nu\cdot \mu$ also defines a modification. 
\[
    \myinput{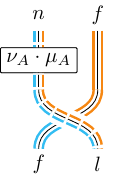}=
    \myinput{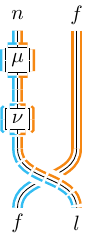}=
    \myinput{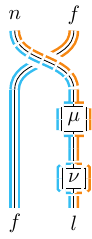}=
    \myinput{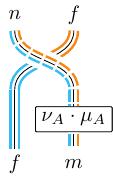}
\]
Similarly, suppose that $n,n'\colon {\color{c1} F}\Rightarrow {\color{c2}G}$ and that $m,m'\colon {\color{c2} G}\Rightarrow {\color{c3}H}$ are pseudonatural transformations, and that $\mu\colon n\Rrightarrow n'$, $\nu\colon m\Rrightarrow m'$ are modifications. The pointwise horizontal composite of $\mu$ and $\nu$ also gives a modification since the following equalities hold.
\[
    \myinput{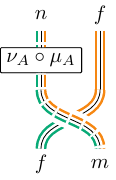}=
    \myinput{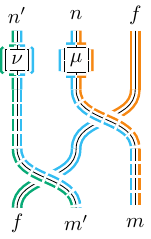}=
    \myinput{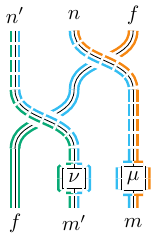}=
    \myinput{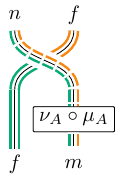}
\]
\begin{remark}
    Given any monoidal category we can construct a braided monoidal category called the monoidal centre -- sometimes referred to as the Drinfeld centre -- independently discovered by Joyal and Street \cite[def.~3]{JOYAL1991tortile}, and Drinfeld. Drinfeld realised the centre as a special case of Majid's~\cite[ex.~3.4]{shahn1991} duals for functored monoidal categories. Another way to define the monoidal centre is to first think of a monoidal category $\MM$ as a one-object bicategory, and then take the category of pseudonatural transformations from the identity on $\MM$ to itself:
    \[\operatorname{Bicat}(\MM,\MM)(\id,\id).\]
    
    One advantage of the decorated string diagram notation is that it makes immediately clear how the monoidal centre is braided.
    
    For the sake of illustration we first consider the centre of a monoid. Let $\MM$ be a monoid viewed as a category. Now consider the set $\Cat(\MM,\MM)(\Id,\Id)$, of natural transformations from the identity functor to itself, which has a canonical monoid structure given by composition. Every natural transformation $(n_A)_{A\in \MM}$ is a family of morphisms in $\MM$, indexed by the objects of $\MM$ such that for every $f$ the following diagram commutes.
    \[
    \begin{tikzcd}
        A
        \arrow[r, "f"]
        \arrow[d, "n_A"{swap}]
            &B
            \arrow[d,"n_B"]\\
        A
        \arrow[r,"f"{swap}]
            &B
    \end{tikzcd}
    \]
    However, since we know that $\MM$ has exactly one object, each natural transformation is simply a map $n_*$ in $\MM$ such that the following diagram commutes for every $f$,
    \[
    \begin{tikzcd}
        *
        \arrow[r, "f"]
        \arrow[d, "n_*"{swap}]
            &*
            \arrow[d,"n_*"]\\
        *
        \arrow[r,"f"{swap}]
            &*
    \end{tikzcd}  
    \]
    and so, in other words, $n_*$, is a natural transformation if and only if it is a member of the centre of $\MM$. Now consider the bicategorical equivalent of this statement. Let $\MM$ be a monoidal category viewed as a one-object bicategory, and consider the category \[\operatorname{Bicat}(\MM,\MM)(\id,\id)\]
    whose objects are pseudonatural transformations from the identity pseudofunctor to the identity pseudofunctor, and whose morphisms are modifications. This category has a canonical monoidal structure given by composing the natural transformations. We can immediately see that this category is braided monoidal, since all objects and morphisms are of the forms given below.
    \[
    \myinput{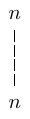}    
    \qquad
    \myinput{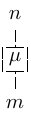}
    \]
    Since the only wrapping is the transparent wrapping of the identity functor, every string can braid over every other string. And since every morphism is a modification, any morphism can pass over or under any braid. In other words, the braiding for this monoidal structure is just given by the naturalisor for each pseudonatural transformation.
\end{remark}
It is worth noting that modifications interact well with inverses.
\begin{proposition}
Suppose that $F,G\colon \AA\to \BB$ are pseudofunctors, $n,m\colon F\Rightarrow G$ are pseudonatural transformations, and $\mu\colon n\Rrightarrow m$ is a modification from $n$ to $m$. If, for every $A\in \AA$, the 2-cell $\mu_A$ has an inverse, then the $\AA$-indexed family $(\mu^{-1}_A)_{A\in \AA}$ defines a modification $\mu^{-1}\colon m\Rrightarrow n$.
\end{proposition}
\begin{proof}
The proof of this follows similarly to the proof of \autoref{prop:altNatAxioms}.
\end{proof}
Finally, we should also note that our modifications often have composites of natural transformations as their sources and targets. In the interest of saving space we introduce the following definition for modifications $\nu\colon m\circ n\Rrightarrow m'\circ n'$,
\[
    \myinput{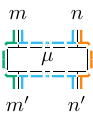}\quad \coloneqq \quad \myinput{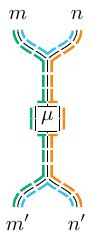}
\]
and similarly for composites of arbitrary numbers of pseudonatural transformations.
\section{Whiskering}  
As with natural transformations we can compose pseudonatural transformations not just with other pseudonatural transformations, but also with pseudofunctors. This is sometimes known as whiskering.
\begin{definition}
Let $F,G\colon\AA\to \BB$ and ${H},{K}\colon \BB\to \CC$ be pseudofunctors, and let $n\colon {F}\Rightarrow{G}$ and $m\colon H\Rightarrow K$  be pseudonatural transformations. 

The \textdef{left}, or \textdef{outer}, \textdef{whiskering} of $H$ with $n$ is the pseudonatural transformation whose 1-cells and naturalisors are given by \[\left(H(n_A)\right)_{A\in \AA}\text{ and }(H(n_f))_{{\morphindx{f}{A}{B}{\AA}}}.\]

The \textdef{right}, or \textdef{inner}, \textdef{whiskering} of $F$ with $m$ is the pseudonatural transformation whose 1-cells and naturalisors are given by \[(m_{FA})_{A\in \AA}\text{ and }(m_{Ff})_{\morphindx{f}{A}{B}{\AA}}.\]
\end{definition}
Let's assign these functors the following colours:  ${\color{c1}F}$, ${\color{c2}G}$, ${\color{c3}H}$, and ${\color{c4}K}$. Then the 1-cells of an outer whiskering have two obvious diagrammatic representations.
\[
    \fbox{\tikzsetnextfilename{PNat/Axioms/Whiskering/1}
\begin{tikzpicture}[baseline={(current bounding box.center)}, xscale=1
]

\drawStartNodes{
    n/$n$/0
}
\strIdX{n}[c3/c2/c1][2]

\drawEndNodes{
    n/$n$
}

\end{tikzpicture}
}\quad = \quad \fbox{\tikzsetnextfilename{PNat/Axioms/Whiskering/2}
\begin{tikzpicture}[baseline={(current bounding box.center)}, xscale=1
]

\drawStartNodes{
    n/${Hn}$/0
}
\strIdX{n}[/c3,c2/c1,c3][2]

\drawEndNodes{
    n/${Hn}$
}

\end{tikzpicture}
}
\]
In the left-hand diagram we view $Hn$ as a pseudonatural transformation $n\colon {\color{c1}F}\Rightarrow {\color{c2}G}$ with ${\color{c3}H}$ applied. In the right-hand diagram we view $Hn$ as a pseudonatural transformation $Hn\colon {\color{c3}H}\circ {\color{c1}F}\Rightarrow {\color{c3}H}\circ {\color{c2}G}$ in its own right. This means that the naturalisors also have two diagrammatic representations.
\[
    \fbox{\tikzsetnextfilename{PNat/Axioms/Whiskering/3}
\begin{tikzpicture}[baseline={(current bounding box.center)}, xscale=1
]

\drawStartNodes{
    n/$n$/2,
    f/$f$
}
\strBraidXX{n}{f}[c3/c2/c1][c3,c1][2]
\drawEndNodes{
    n/$n$,
    f/$f$
}

\end{tikzpicture}
}\quad\coloneqq\quad\fbox{\tikzsetnextfilename{PNat/Axioms/Whiskering/4}
\begin{tikzpicture}[baseline={(current bounding box.center)}, xscale=1
]

\drawStartNodes{
    n/$Hn$/2,
    f/$f$
}
\strBraidXX{n}{f}[/c3,c2/c1,c3][c3,c1][2]
\drawEndNodes{
    n/$Hn$,
    f/$f$
}

\end{tikzpicture}
}
\]
In the left diagram, we simply depict the naturalisor for $Hn$ as a braid where everything is wrapped in the colour of ${\color{c3}H}$. On the right-hand side we depict the naturalisor for $Hn$ as a braid that changes the ${\color{c3}H}\circ {\color{c1}F}$ wrapping to the ${\color{c3}H}\circ {\color{c2}G}$ wrapping.

Analogously the 1-cells of an inner whiskering have two diagrammatic representations.
\[
    \fbox{\tikzsetnextfilename{PNat/Axioms/Whiskering/5}
\begin{tikzpicture}[baseline={(current bounding box.center)}, xscale=1
]

\drawStartNodes{
    m/$m$/0
}
\strIdX{m}[/c4/c3/c1][2]

\drawEndNodes{
    m/$m$
}

\end{tikzpicture}
}\quad\coloneqq\quad\fbox{\tikzsetnextfilename{PNat/Axioms/Whiskering/6}
\begin{tikzpicture}[baseline={(current bounding box.center)}, xscale=1
]

\drawStartNodes{
    m/$m_{F-}$/0
}
\strIdX{m}[/c4,c1/c1,c3][2]

\drawEndNodes{
    m/$m_{F-}$
}

\end{tikzpicture}
}
\]
The naturalisors also have two diagrammatic representations.
\[
    \fbox{\tikzsetnextfilename{PNat/Axioms/Whiskering/7}
\begin{tikzpicture}[baseline={(current bounding box.center)}, xscale=1
]

\drawStartNodes{
    m/$m$/2,
    f/$f$/0
}
\strBraidXX{m}{f}[/c4/c3/c1][c3,c1][2]
\drawEndNodes{
    m/$m$,
    f/$f$
}

\end{tikzpicture}
}\quad\coloneqq\quad\fbox{\tikzsetnextfilename{PNat/Axioms/Whiskering/8}
\begin{tikzpicture}[baseline={(current bounding box.center)}, xscale=1
]

\drawStartNodes{
    m/$m_{F-}$/2,
    f/$f$
}
\strBraidXX{m}{f}[/c4,c1/c1,c3][c3,c1][2]
\drawEndNodes{
    m/$m_{F-}$,
    f/$f$
}

\end{tikzpicture}
}
\]
Given a whiskering of a pseudonatural transformation with a pseudofunctor, it will often be useful to consider both diagrammatic representations in the same diagram. To indicate a change of perspective we will change the colouring on the string by way of a black box. For example, we might draw the following diagram.
\[
    \myinput{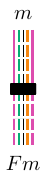}
\] 
Note how the outer colour changes from a solid line to a dashed line. Of course the diagram above simply represents the identity 2-cell from $Fm$ to $Fm$, but explicitly changing perspective in this way can be helpful in certain proofs and constructions. For other changes of perspective we will also use this same black box. For example in the following diagram we change our perspective, forgetting that $Fm_A$ is a component in a pseudonatural transformation and changing it to be a plain old 1-cell.
\[
    \myinput{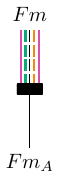}
\] 
There is also special case for whiskering that is slightly more complicated. That is when we have a pseudonatural transformation passing over another pseudonatural transformation. For example, suppose we want to pass $m_{FA}$ over the top of $H(n_A)$. Since $m$ is a pseudonatural transformation, and $n_A$ is a 1-cell, we have a 2-cell in the following diagram,
\[
    \begin{tikzcd}[column sep=10ex, row sep=10ex]
    HFA
    \arrow[r, "H(n_A)"]
    \arrow[d, "m_{FA}"{swap}]
        &HGA
        \arrow[d, "m_{GA}"]
        \arrow[dl, Rightarrow, "m_{n_A}"{description}]\\
    KFA
    \arrow[r, "K(n_A)"{swap}]
        &KGA
    \end{tikzcd}
\]
and since this 2-cell is a naturalisor, it adheres to all of the braiding axioms. But notice that on the right we have $m_{GA}$ and on the left we have $m_{FA}$, so our inner whiskering has changed. This is due to a sort of interchange law. We have swapped the outer whiskering of $n$ by the properties of $m$ and swapped the inner whiskering of $m$ by the properties of $n$. If we write this out in our graphical language we get the following braid.
\[
  \fbox{\tikzsetnextfilename{PNat/Axioms/Whiskering/9}
\begin{tikzpicture}[baseline={(current bounding box.center)}, xscale=1
]

\drawStartNodes{
    m/$m$/2,
    n/$n$
}
\strBraidXX{m}{n}[/c4/c3/c2][c3/c2/c1][2]
\drawEndNodes{
    n/$n$,
    m/$m$
}

\end{tikzpicture}
}
\]
These special braids have some interesting properties. For a start they adhere to a multi-object (or multi-coloured) version of the braid relations, the relations that govern the Yang-Baxter equations.
\begin{proposition}
Given any pseudofunctors ${\color{c1}F},{\color{c2}G}\colon\AA\to \BB$ and ${\color{c3}H},{\color{c4}K}\colon \BB\to \CC$; any pseudonatural transformations $n\colon {\color{c1}F}\Rightarrow{\color{c2}H}$ and $m\colon {\color{c3}H}\Rightarrow {\color{c4}K}$; and any morphism $f\colon A\to B$ in $\AA$, we have the following equality.
\[
\myinput{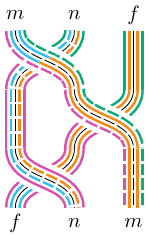}=\myinput{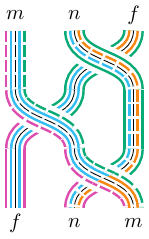}
\]
\end{proposition}
\begin{proof}
This is a simple application of the braiding laws that follows from viewing the braiding of $n$ over $f$ as the 2-cell $H(n_f)$.
\[
\myinput{tikz_PNat_Yang-Baxter_2}=\myinput{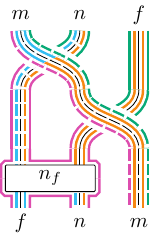}=\myinput{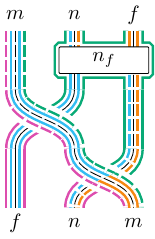}=\myinput{tikz_PNat_Yang-Baxter_1}
\]
\end{proof}
Another way of interpreting this proposition is that braiding a pseudonatural transformation over another pseudonatural transformation yields a modification.
\section{Adjunctions}
If every morphism in a natural transformation has an inverse, then the natural transformation itself has an inverse. In this section we prove that the same holds for adjoint pseudonatural transformations.
\begin{proposition}
Let $F,G\colon \AA\to \BB$ be pseudofunctors, and let $n\colon F\Rightarrow G$ be a pseudonatural transformation. 
\begin{itemize}
\item   If $n_A$ has a left adjoint $n_A^L$ for every object $A\in \AA$, then there is a pseudonatural transformation $n^L\colon G\Rightarrow F$ whose 1-cells are given by are $(n^L_A)_{A\in \AA}$. 
\item If $n_A$ has a right adjoint $n_A^R$ for every object $A\in \AA$ then there is a pseudonatural transformation $n^R\colon G\Rightarrow F$ whose 1-cells are given by $(n^R_A)_{A\in \AA}$. 
\end{itemize}
\end{proposition}
\begin{proof}
To prove this we need, for every 1-cell $f$ in $\AA$, to be able to define a 2-cell $n^L_f$ and $n^R_f$ such that the pseudonatural transformation axioms hold. The definitions of these 2-cells are given by \autoref{adjointprop}.
\begin{figure}[h]
\[
    (n^L_f)^{-1}\coloneqq \myinput{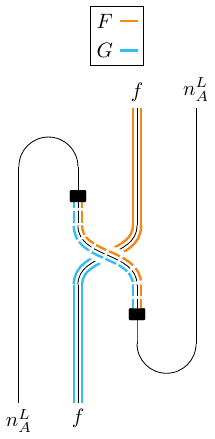}\qquad
    n^R_f\coloneqq \myinput{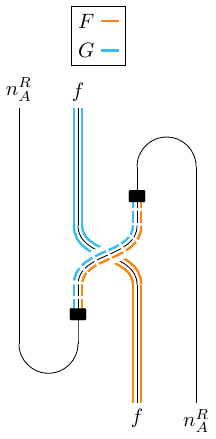}
\]
\caption{Definition of $(n^L_f)^{-1}$ and $(n^R_f)$.}
\label{adjointprop}
\end{figure}
To prove that these 2-cells adhere to the naturalisor axioms we appeal to \autoref{prop:altNatAxioms}, \autoref{prop:altNatAxioms1} and \autoref{prop:altNatAxioms2}.
\end{proof}

The above result allows us to include cups and caps in our decorated string diagram language. If $l\colon F\Rightarrow G$ is left adjoint to $r\colon G\Rightarrow F$ it makes perfect sense to colour cups and caps in the following way.
\[
\myinput{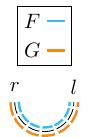} \qquad \myinput{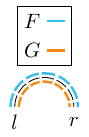}
\]
Our aim now is to prove that these pointwise adjoints give adjoints in the bicategory of pseudonatural transformations and modifications. To do this we need to show that the pointwise cups and caps give modifications. This is proven by the following sequence of propositions.
\begin{proposition}
\label{prop:cupsandbraids}
Let ${\color{c1}F},{\color{c2}G}\colon \AA\to \BB$ be pseudofunctors, and let $l\colon F\Rightarrow G$ and $r\colon G\Rightarrow F$ be pseudonatural transformations such that $l$ is the left adjoint of $r$. Then for any 1-cell $f\colon A\to B$ in $\AA$ we have the following equalities.
\[
    \myinput{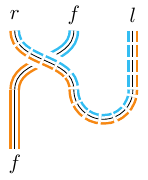}=
    \myinput{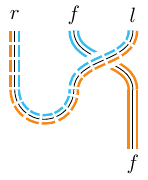}
\] 
\[ 
    \myinput{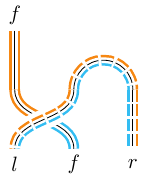}=
    \myinput{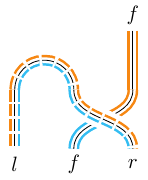}
\]
\end{proposition}
\begin{proof}
We have the following two equalities, firstly by definition and then by the yanking condition.
\[
    \myinput{tikz_PNat_Adj_6}=
    \myinput{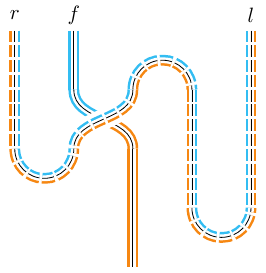}=
    \myinput{tikz_PNat_Adj_7}
\]
The case for caps follows similarly.
\end{proof}
One result of this proposition is that we can unambiguously define the following crossing diagrams.
\[
    \myinput{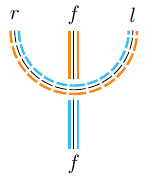}\coloneqq \myinput{tikz_PNat_Adj_4}=\myinput{tikz_PNat_Adj_5}
\]
\[
    \myinput{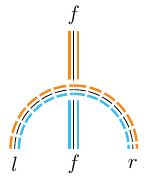}\coloneqq \myinput{tikz_PNat_Adj_4}=\myinput{tikz_PNat_Adj_5}
\]
This result also lets us prove that a cup or a cap can be pulled over a string.
\begin{proposition}
    Let ${\color{c1}F},{\color{c2}G}\colon \AA\to \BB$ be pseudofunctors, and let $l\colon F\Rightarrow G$ and $r\colon G\Rightarrow F$ be pseudonatural transformations such that $l$ is the left adjoint of $r$. Then for any 1-cell $f\colon A\to B$ in $\AA$ we have the following equalities.
\[
    \myinput{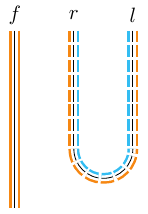}=
    \myinput{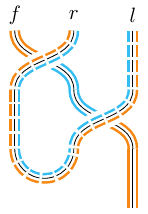}
\]
\[
    \myinput{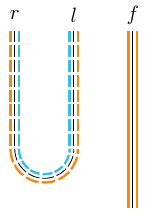}=
    \myinput{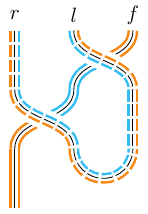}
\]
\[
    \myinput{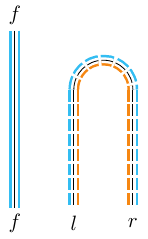}=
    \myinput{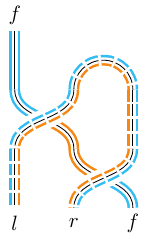}
\]
\[
    \myinput{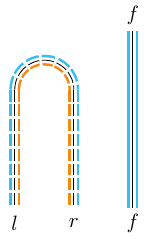}=
    \myinput{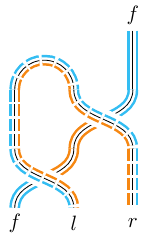}
\]
\end{proposition}
\begin{proof}
To prove the first equality note that the following holds by definition and then cancellation.
    \[
        \myinput{tikz_PNat_Adj_10}=
        \myinput{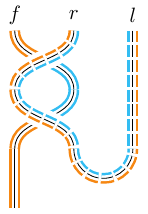}=
        \myinput{tikz_PNat_Adj_9}
    \]
The proofs of the other equalities are analogous.
\end{proof}
\begin{corollary}
    Let $\AA$ and $\BB$ be bicategories, and $\Bicat(\AA,\BB)$ be the bicategory whose objects are pseudofunctors, whose 1-cells are pseudonatural transformations and whose 2-cells are modifications. A pseudonatural transformation $n\colon F\Rightarrow G$ has a left adjoint in $\Bicat(\AA,\BB)$ if and only if $n_A$ has a left adjoint for every $A\in \AA$.
\end{corollary}
\begin{proof}
    This follows from the fact that the cup and cap above form modifications.
\end{proof}
\section{Biequivalences and Strictification}
In this section we explore the strictification theorem, which essentially allows us to ignore composition associators and unitors in bicategories. This justifies their omission in string diagrams and will simplify certain proofs later on. First we need some notion of what it means for two bicategories to be the same.
\label{section:strictification}
\begin{definition}
    A \textdef{biequivalence} between bicategories $\AA$ and $\BB$ is a pseudofunctor $F\colon \AA\to \BB$ that is \begin{itemize}
        \item essentially surjective: for every $B\in \BB$ there is an $A\in \AA$ with $F(A)$ equivalent to $B$;
        \item a local equivalence: the hom-functor $F\colon \AA(A,B)\to \BB(FA,FB)$ is an equivalence for all $A,B\in \AA$.
    \end{itemize}
    We say that two bicategories are \textdef{biequivalent} if there is a biequivalence between them.
\end{definition}
This definition is often easier to work with but, as with equivalences of categories, there is an alternative definition in terms of pseudo-inverses, see for example Johnson and Yau's \cite[def.~6.2.10, thm~7.4.1]{yau2021bicat} monograph on bicategories. It follows from this alternative definition that biequivalence of bicategories defines an equivalence relation. Now that we have a notion of what it means for two bicategories to be the same in a weak sense, we can express the strictification theorem for bicategories.
    \begin{theorem}
    Every bicategory is biequivalent to a 2-category.
    \end{theorem}
    The proof of this essentially follows from an argument given by Mac Lane and Par\'e~\cite[p.~61]{maclane1985coherence}, which itself is a refinement and generalisation of Mac Lane's~\cite[thm.~5.1]{maclane1963natural} own coherence theorem for monoidal categories.
    However, it seems to have been Power~\cite[p.~172]{power1989general} who first wrote this result in terms of a biequivalence. Gordon, Power and Street~\cite[thm.~1.4]{street1995coherence} also pointed out that this follows as an immediate consequence of the Yoneda lemma for bicategories.

    It is worth noting that, although there is a standard method to construct such a 2-category, Lack~\cite[ex.~1.5]{lack2010} has pointed out that ``naturally occurring bicategories tend to be biequivalent
    to naturally occurring 2-categories''.
    \label{section:NaturalStrict}
    \begin{example}
    Carboni and Johnstone~\cite[lem.~4.8]{carboni1995connected} showed that the bicategory of profunctors $\Prof$ is biequivalent to the sub-2-category of $\Cat$ consisting of presheaf categories and cocontinuous functors between them. This essentially follows from the fact that a profunctor
    \[
        P\colon \BB^\op\times \AA\to \Set
    \]
    can be curried to give a functor
    \[
        P\colon \AA\to \hom{\BB^\op, \Set}
    \]
    and the left Kan extension of $P$ along the Yoneda embedding is the unique cocontinuous extension of $\AA$ to $\hom{\AA^{\op}, \Set}$. The same construction works for $\VV$-$\Prof$, so long as $\VV$ is a cosmos.
    \end{example}
    \begin{example}
        Lack~\cite[ex.~1.5]{lack2010} pointed out that the bicategory of spans over some category with pullbacks is equivalent to the sub-2-category of $\Cat$ consisting of slice categories of $\CC$, $\CC/C$, and functors between them
        \[
            \CC/C\xrightarrow{f^*} \CC/S\xrightarrow{g_!} \CC/D
        \] 
        where $f^*$ is the pullback functor for some $f\colon S\to C$, and $g_!$ is the postcomposition functor for $g\colon S\to D$. The correspondence is given by identifying the functor $g_!\circ f^*$ with the span $C\xleftarrow{f} S \xrightarrow{g} D$.
    \end{example}
    Strictification gives us a formal grounding for omitting identities and associators from our string diagrams: it simply doesn't matter where we put them.

\section{Pseudoadjunctions}
In this section we explore a particular bicategorical analogue for the notion of adjoint. As pointed out by Gray~\cite[ch.~6]{gray1974formal} there are many, many ways to construct something resembling an adjunction between bicategories. Since we are working with pseudofunctors and not lax functors, for our purposes we settle on a relatively strong definition. 
\begin{definition}
Given two pseudofunctors $L\colon \AA\to \BB$ and $R\colon \BB\to \AA$, we say that $L$ is \textdef{left pseudoadjoint} to $R$ if there are functors
\[
E_{A,B}\colon \BB(LA,B)\rightleftarrows \AA(A,RB)\cocolon I_{A,B}     
\]
natural in both $A$ and $B$, that form an equivalence. When $L$ is left pseudoadjoint to $R$ we say that $R$ is \textdef{right pseudoadjoint} to $L$ and that $L$ and $R$ form a \textdef{pseudoadjunction}.
\end{definition}
Note that, pseudoadjunctions are sometimes referred to as \emph{biadjunctions}. It is worth being explicit about what this means, and we will express this definition in two ways. The first way is more succinct. Firstly note that we have functors
\[
{\color{c1} \BB(L-,-)} \colon \AA^{\op}\times\BB\to \Cat \text{ and }{\color{c2} \AA(-,R-)}\colon \AA^{\op}\times \BB\to \Cat.
\]
Then $L$ being left pseudoadjoint to $R$ means that there are a pair of pseudonatural equivalences in the bicategory $\Cat$, expressed by the following diagrams.
\[
\myinput{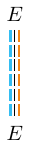}\qquad \myinput{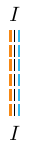}
\]
The second way to understand this definition is less succinct but perhaps more familiar, since it deals with 1-cells and 2-cells in $\AA$ and $\BB$. Assigning {\color{c1} $L$ this colour} and {\color{c2} $R$ this colour}, $L$ being left pseudoadjoint to $R$ means that there are pseudonatural transformations, called the unit and counit expressed by the following diagrams,
\[
    \myinput{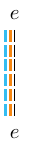}\qquad \myinput{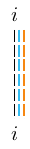}
\]
such that the following pairs of whiskerings form pseudonatural adjoint equivalences.
\[
    \myinput{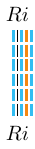}\qquad \myinput{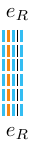}\qquad\qquad\qquad    \myinput{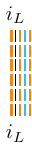}\qquad\myinput{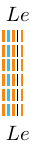}
\]

To seasoned category theorists this definition may seem na\"ive. It seems as though we've simply replaced the word ``isomorphism" with ``adjoint equivalence", and typically in higher category theory that isn't enough. We usually also need coherence diagrams. 

\begin{figure}[h]
\[
    \myinput{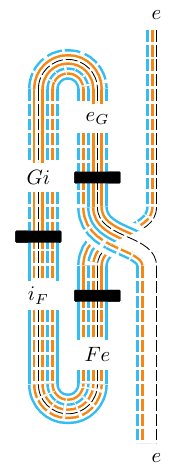}\quad=\quad\myinput{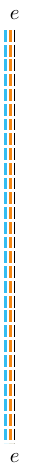}
\]
\caption{A string diagram depiction of one of the swallowtail conditions.}
\label{swallowtail}
\end{figure}

As pointed out by Verity~\cite[lem.~1.3.9]{verity} the case of pseudoadjoints is unusual in that the coherence diagrams come for free. By \autoref{prop:adjEquiv} we can always upgrade our equivalence of categories
\[
    E_{A,B}\colon \BB(LA,B)\rightleftarrows \AA(A,RB)\cocolon I_{A,B}    
\] 
to an adjoint equivalence of categories. When $E$ and $I$ form an adjoint equivalence, the associated unit $e$ and counit $i$ have two coherence conditions, sometimes called the swallowtail conditions. For strict 2-functors the first of these conditions is given by \autoref{swallowtail}, and the second is given by a similar equality for $i$. Since the swallowtail conditions aren't strictly necessary for pseudofunctors we make no use of them in this thesis.
\begin{lemma}
    If $G$ and $G'$ are both right pseudoadjoint to $F$ then there is a pseudonatural adjoint equivalence from $G'$ to $G$.
\end{lemma}
\begin{proof}
This follows as a consequence of the bicategorical Yoneda lemma that can be found, for example, in Johnson and Yau's~\cite[lem.~8.3.12]{yau2021bicat} book on bicategories, and is given as a result in Verity's~\cite[lem.~1.3.9]{verity} thesis. However, we give the sketch of a more constructive approach here. Let $e$ and $i$ be the unit and counit for the $F,G$ pseudoadjunction and let $e', i'$ be the unit and counit for the $F,G'$ adjunction. Then the 1-cells of the equivalence are given by
\[
G'\xrightarrow{e_{G'}} GFG' \xrightarrow{G(i')} G    
\]
and 
\[
G \xrightarrow{e'_{G}}G'FG \xrightarrow{G'(i)} G'.
\]
The proof is essentially the same as the proof that right adjoint 1-cells are unique up to isomorphism, but equalities are replaced with naturalisors. To be explicit, we have a 2-cell given by the pasting diagram in \autoref{pastingAdjUn}, or equivalently the string diagram in \autoref{stringAdjUn}.
\begin{figure}
\[
\begin{tikzcd}[column sep=3.2em, row sep=3.2em]
G'
\arrow[rr, "e_{G'}", ""{swap, name=ttl}]
\arrow[dr, "e'_{G}"]
\arrow[dddrrr, equal, bend right=45, ""{name=outerequal}]
    &\blank
        &GFG'
        \arrow[r, "G(i')", ""{swap,name=ttr}]
        \arrow[d, "e'_{GFG'}"]
            &G'
            \arrow[d, "e'_{G}"]\\
\blank
    &G'FG'
    \arrow[rd, equal, ""{swap,name=innerequal}, bend right=30]
    \arrow[r, "G'Fe_{G'}"{name=btl}]
        &G'FGFG'
        \arrow[phantom, "\cong", to=innerequal]
        \arrow[r,"G'FG(i')"{name=btr}]
        \arrow[d, "G'(i_{FG'})"{name=lbr}, ""{swap, name=lbr2}]
            &GFG'
            \arrow[dd, "G'(i)"{name=rbr}, ]\\
\blank
    &\blank
        &G'FG'
        \arrow[dr,"G(i')"]
            &\blank\\
\blank
    &\blank
        &\blank
            &G'
\arrow[phantom, "\cong", from=ttl, to=btl]
\arrow[phantom, "\cong", from=ttr, to=btr]
\arrow[phantom, "\cong", from=lbr2, to=rbr]
\arrow[phantom, "\cong", from=innerequal, to=outerequal]
\end{tikzcd}    
\]
\caption{}
\label{pastingAdjUn}
\end{figure}
\begin{figure}
\[
    \myinput{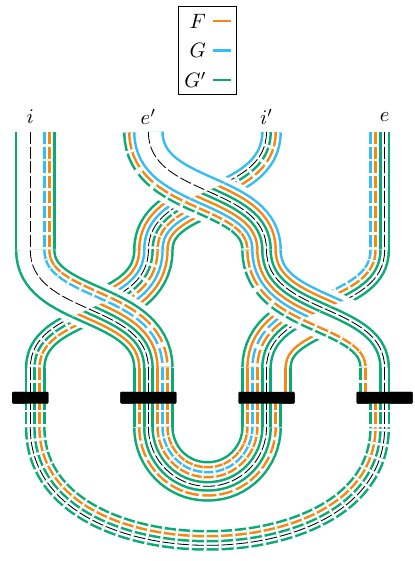}
\]
\caption{}
\label{stringAdjUn}
\end{figure}
This 2-cell is just a composite of naturalisors and the equivalence conditions for the units and counits. The 2-cell giving the isomorphism between the identity and the reverse composite is given similarly. Of course here we are treating pseudofunctors as if they compose strictly associatively, and so this particular `proof' technically only holds for 2-functors.
\end{proof}

\chapter{Monoidal Bicategories}
The aim of this thesis is to show that bicategories with certain closed structures have a canonical enrichment over their monoidal category of scalars. To understand scalars we must first understand monoidal bicategories.

In the first section we give the definition of a monoidal bicategory and give a monoidal structure for each of our motivating examples.

In the second section we focus on the role of the interchangerator. This is a 2-cell that follows from the existence of the compositor for the tensor pseudofunctor. Much of this section is dedicated to the braid-like properties of the interchangerator.

In the third section we reprove a known result: the fact that the monoidal category of scalars is braided monoidal. We prove this using our decorated string diagram language, and without appealing to any results that strictify the monoidal structure. This follows from the braid-like nature of the interchangerator.

In the fourth section we give an account of semi-strictification for bicategories, a result which shows that every monoidal bicategory is equivalent to one where all of the 1-cell and 2-cell data that define the monoidal bicategory, except the interchangerator, are degenerate.

In the final section we give the definition of a right-monoidal-closed bicategory and give a closed structure for each of our motivating examples. This monoidal-closed structure is the second type of closed structure that our bicategories require in order to be equipped with a cotrace functor.
\section{Definitions and Examples}
\begin{definition}
A \textdef{monoidal bicategory} is a bicategory $\BB$ equipped with a pseudofunctor $\tensor \colon \BB\times \BB\to \BB$, a unit object $I\in \BB$ along with pseudonatural adjoint equivalences:
\begin{itemize}
    \item the \textdef{associator} $a_{A,B,C}\colon (A\tensor B)\tensor C\to A\tensor (B\tensor C)$;
    \item the \textdef{left unitor} $l_A\colon I\tensor A\to A$;
    \item the \textdef{right unitor} $r_A\colon A\tensor I\to A$;
\end{itemize}
as well as \textdef{coherence modifications}:
\begin{itemize}
    \item the \textdef{2-associator} whose components are given by the following diagram;
    \[
        \myinput{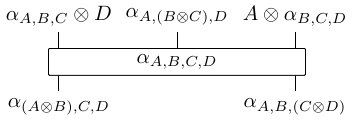}   
    \]
    \item the \textdef{left}, \textdef{middle} and \textdef{right 2-unitors} whose components are given by the following diagrams.
    \[
        \myinput{tikz_MonBicat_leftunitor}\qquad 
        \myinput{tikz_MonBicat_middleunitor} \qquad 
        \myinput{tikz_MonBicat_rightunitor}
    \]
\end{itemize}
These coherence 2-cells adhere to their own coherence axioms, which we choose to omit here but can be found in, for example, Stay's~\cite[def.~4.4]{stay2016} paper on compact-closed bicategories.
\end{definition}
\begin{remark}
Note that we did not use any coloured strings for our modifications above but instead defined them by their components. This is because, as we point out in \autoref{prop:narypseudofunctor}, a binary pseudofunctor such as $-\tensor -$ can be thought of as a collection of unary pseudofunctors $A\tensor -$ and $-\tensor A$ for each object $A\in \BB$. For example, depending on the context, the 1-cell $a_{A,B,C}\colon (A\tensor B)\tensor C\to A\tensor (B\tensor C)$ might be thought of as a pseudonatural transformation
\[a_{-,B,C}\colon(-\tensor C)\circ (-\tensor B)\Rightarrow-\tensor (B\tensor C),\] or a pseudonatural transformation
\[
    a_{A,B,-}\colon (A\tensor B)\tensor -\Rightarrow (A\tensor -)\circ (B\tensor -)
\]
Different perspectives are useful in different arguments, and so we choose the decorations for our diagrams depending on the context.
\end{remark}
We do not include the coherence axioms here since they are lengthy and complicated, and we never appeal to them directly. To understand what the axioms achieve consider the following. Suppose that we have five objects $A,B,C,D,E$ in $\MM$ and we tensor them together to get $(((A\tensor B)\tensor C)\tensor D)\tensor E$. There are two distinct ways in which we can use a string of associators to re-bracket and end up with $A\tensor (B\tensor (C\tensor (D\tensor E)))$. These two distinct strings, or paths, are naturally isomorphic, or homotopic, via a composite of 2-associators. However, a priori, there are also two distinct composites of 2-associators. One axiom asserts that these two composites of 2-associators are identical, thus forcing the two paths of associators to be \emph{uniquely} isomorphic. The other axioms give similar assertions but with regard to composites of 2-unitors.
\begin{example}
    The bicategory $\Rel$ has a monoidal bicategory structure where the tensor product is the cartesian product of sets. 
\end{example}
\begin{example}
    The bicategory $\Bim_R$ has a monoidal bicategory structure where the tensor product is given by taking the tensor product of $R$-algebras. The unit is given by $R$.
\end{example}
\begin{example}
    The bicategory $\DGBim_R$ has a monoidal bicategory structure where the tensor product is given by taking the tensor product of $R$-algebras. The unit is given by $R$.
\end{example}
\begin{example}
    The bicategory $\VV$-Prof has a monoidal bicategory structure where the tensor product is given by taking the tensor product of categories. The unit is given by the one object category $*$ with hom-object $I\in \VV$.
\end{example}
\begin{example}
    If $\CC$ has finite limits then $\Span(\CC)$ has a monoidal bicategory structure where the tensor product is given by taking the product of objects in $\CC$.
\end{example}
\begin{example}
    If $M$ is a topological monoid then $\Path(M)$ has a monoidal bicategory structure given by taking the monoidal product of points in $M$.
\end{example}

\begin{remark}
    Recall that in \autoref{section:NaturalStrict} we gave two examples of `naturally strictifiable' bicategories. These two examples were motivating examples for an observation, made by Lack~\cite[ex.~1.5]{lack2010}, that many bicategories seem to be naturally strictifiable. Interestingly, the natural strictification for both of these examples arises from what might be called the `3-trace'. For every monoidal bicategory there is a pseudofunctor 
    \[
        \BB\xrightarrow{\BB(I,-)} \Cat.
    \]
    We can think of $\BB$ as a one-object tricategory where the unit is the identity 1-cell. Thus, the pseudofunctor above is defined analogously to the 2-trace defined in \autoref{section:bicategoricaltrace}.
    For profunctors and spans the natural strictification is given by this pseudofunctor.

    Firstly consider the bicategory $\VV$-$\Prof$, where the unit object is given by the point $\VV$-category $*$. Given a category $\AA$, then we have that \[\VV\mhyphen\Prof(*,\AA)\cong \cat{\VV}(\AA^{\op}, \VV)\] i.e. $\AA$ is taken to the presheaf category on $\AA$. For a profunctor $P\colon A\profto B$, the functor $\VV\mhyphen\Prof(*,P)$ is given by postcomposition by $P$. Consider a presheaf $F\colon \AA^{\op}\to \VV$. It is known that $F$ is some weighted colimit of representables $F\cong \colim^W \AA(-,A)$. The postcomposition applied to $F$, $\VV\mhyphen \Prof(*,P)(F)$, is then given by
    \[
        \endint^{A\in \AA} \colim^W \AA(-,A)\tensor P(A,-)\cong \colim^W\endint^{A\in \AA} \AA(-,A)\tensor P(A,-)\cong \colim^W P(A,-)
    \]
    and so $\VV$-$\Prof(*, P)$ is the cocontinuous extension of $P\colon A\to \hom{B^{\op}, \VV}$.

    Now let us turn to the case for spans. Here the unit object is given by the terminal object $*$, and the category $\Span(\CC)(*,C)$ is isomorphic to the slice category $\CC/C$. Given a span\[ 
    \begin{tikzcd}
        & S
        \arrow[ld, "f"{swap}]
        \arrow[rd, "g"]\\
    C   && D
    \end{tikzcd}
    \]
    the functor $\Span(\CC)(*,S)$ is postcomposition with $S$, which corresponds to first pulling back along $f$ and then composing with $g$. In other words, $\Span(\CC)(*,S)$ is given by the composite
    \[
       \CC/C\xrightarrow{f^*}\CC/S\xrightarrow{g_!} \CC/D.
    \]
\end{remark}
As with monoidal categories, it will often be useful to think of the tensor product of a monoidal category $\AA$ as a collection of unary pseudofunctors, indexed by the objects of $\AA$.
\begin{proposition}
\label{prop:narypseudofunctor}
Let $F\colon \AA\times \BB\to \CC$ be a pseudofunctor. For every $A\in \AA$ there is a pseudofunctor
\[
    F(A,-)\colon \BB\to \CC
\]
whose definition on 1-cells and 2-cells is given by the following diagram,
\[
    \fbox{\tikzsetnextfilename{MonBiCat/BinaryPFunctor/2cell/1}
\begin{tikzpicture}[baseline={(current bounding box.center)}, xscale=1
]
\drawStartNodes{
    f/${f}$/0
}

\strTwoCellX{f}{$\gamma$}{g}[c1][][0.75][0.75]

\strKeyXX{
    $F(A,-)$/c1
}
\drawEndNodes{
    g/$g$
}
\end{tikzpicture}
}\coloneqq\fbox{\tikzsetnextfilename{MonBiCat/BinaryPFunctor/2cell/2}
\begin{tikzpicture}[baseline={(current bounding box.center)}, xscale=1
]
\drawStartNodes{
    f/${(\id_A,f)}$/0
}

\strTwoCellX{f}{$(\iota,\gamma)$}{g}[c2][][0.75][0.75]

\strKeyXX{
    $F(-,-)$/c2
}

\drawEndNodes{
    g/${(\id_A,g)}$
}

\end{tikzpicture}
}
\]
and whose compositor and identitor are defined by the following two diagrams.
\[
    \fbox{\tikzsetnextfilename{MonBiCat/BinaryPFunctor/Compositor/1}
\begin{tikzpicture}[baseline={(current bounding box.center)}, xscale=1
]
\drawStartNodes{
    f/${f}$/-1.5,
    g/$g$/0
}

\strJoinX{g/c1;f/c1}{gf}[1.25][1.25]

\strKeyXX{
    $F(A,-)$/c1
}
\drawEndNodes{
    gf/$g\circ f$
}

\end{tikzpicture}
}\coloneqq\fbox{\tikzsetnextfilename{MonBiCat/BinaryPFunctor/Compositor/2}
\begin{tikzpicture}[baseline={(current bounding box.center)}, xscale=1
]
\drawStartNodes{
    f/${(\id_A,f)}$/-1.5,
    g/${(\id_A,g)}$/0
}

\strJoinX{g/c2;f/c2}{gf}[0.75][0.25]
\strTwoCellX{gf}{$(\comp{\lambda},\iota)$}{gf}[c2][][][1]

\strKeyXX{
    $F(-,-)$/c2
}
\drawEndNodes{
    gf/${(\id_A,g\circ f)}$
}

\end{tikzpicture}
}
\]
\[
    \fbox{\tikzsetnextfilename{MonBiCat/BinaryPFunctor/Identitor/1}\begin{tikzpicture}[baseline={(current bounding box.center)}, xscale=1
]

\drawStartNodes{
    id/$\id$/0
}

\strId{id}[c1][][][1][white]
\idCup{id}{c1}

\strKeyXX{
    $F(A,-)$/c1
}
\end{tikzpicture}
}\coloneqq \fbox{\tikzsetnextfilename{MonBiCat/BinaryPFunctor/Identitor/2}\begin{tikzpicture}[baseline={(current bounding box.center)}, xscale=1
]

\drawStartNodes{
    id/$\id$/0
}

\strId{id}[c2][][][1][white]
\idCup{id}{c2}

\strKeyXX{
    $F(-,-)$/c2
}
\end{tikzpicture}
}
\]
For every $B\in \BB$ there is a pseudofunctor
\[
    F(-,B)\colon \AA\to \CC    
\]
whose definition on 1-cells and 2-cells is given by the following diagram,
\[
    \fbox{\tikzsetnextfilename{MonBiCat/BinaryPFunctor/2cell/3}
\begin{tikzpicture}[baseline={(current bounding box.center)}, xscale=1
]
\drawStartNodes{
    f/${f}$/0
}

\strTwoCellX{f}{$\gamma$}{g}[c3][][0.75][0.75]

\strKeyXX{
    $F(-,B)$/c3
}
\drawEndNodes{
    g/$g$
}
\end{tikzpicture}
}\coloneqq\fbox{\tikzsetnextfilename{MonBiCat/BinaryPFunctor/2cell/4}
\begin{tikzpicture}[baseline={(current bounding box.center)}, xscale=1
]
\drawStartNodes{
    f/${f}$/0
}

\strTwoCellX{f}{$(\gamma,\iota)$}{g}[c2][][0.75][0.75]

\strKeyXX{
    $F(-,-)$/c2
}
\drawEndNodes{
    g/$g$
}
\end{tikzpicture}
}
\]
and whose compositor and identitor are defined by the following two diagrams.
\[
    \fbox{\tikzsetnextfilename{MonBiCat/BinaryPFunctor/Compositor/3}
\begin{tikzpicture}[baseline={(current bounding box.center)}, xscale=1
]
\drawStartNodes{
    f/${f}$/-1.5,
    g/$g$/0
}

\strJoinX{g/c3;f/c3}{gf}[1.25][1.25]

\strKeyXX{
    $F(-,B)$/c3
}
\drawEndNodes{
    gf/$g\circ f$
}

\end{tikzpicture}
}\coloneqq\fbox{\tikzsetnextfilename{MonBiCat/BinaryPFunctor/Compositor/4}
\begin{tikzpicture}[baseline={(current bounding box.center)}, xscale=1
]
\drawStartNodes{
    f/${(f,\id_B)}$/-1.5,
    g/${(g,\id_B)}$/0
}

\strJoinX{g/c2;f/c2}{gf}[0.75][0.25]
\strTwoCellX{gf}{$(\iota,\comp{\lambda})$}{gf}[c2][][][1]

\strKeyXX{
    $F(-,-)$/c2
}
\drawEndNodes{
    gf/${(g\circ f,\id_B)}$
}

\end{tikzpicture}
}
\]
\[
    \fbox{\tikzsetnextfilename{MonBiCat/BinaryPFunctor/Identitor/3}\begin{tikzpicture}[baseline={(current bounding box.center)}, xscale=1
]

\drawStartNodes{
    id/$\id$/0
}

\strId{id}[c3][][][1][white]
\idCup{id}{c3}

\strKeyXX{
    ${F(-,B)}$/c3
}
\end{tikzpicture}
}\coloneqq \fbox{\tikzsetnextfilename{MonBiCat/BinaryPFunctor/Identitor/4}\begin{tikzpicture}[baseline={(current bounding box.center)}, xscale=1
]

\drawStartNodes{
    id/$\id$/0
}

\strId{id}[c2][][][1][white]
\idCup{id}{c2}

\strKeyXX{
    $F(-,-)$/c2
}
\end{tikzpicture}
}
\]
\end{proposition}
\begin{proof}
The proof follows from the axioms of pseudofunctors, and the naturality of the compositor.
\end{proof}
Then, by induction, we can fix any variable in a pseudonatural transformation to get a new pseudonatural transformation.
\begin{proposition}
    Suppose that $F,G\colon \AA_1\times...\times\AA_n\to \BB$ are pseudofunctors and $n\colon F\Rightarrow G$ a pseudonatural transformation. Then for any $i$ and any $A\in \AA_i$, 
    \[
        n_{-,...,A,...,-}\colon F(-,...,A,...,-)\Rightarrow G(-,...,A,...,-)
    \] 
    is a pseudonatural transformation.
\end{proposition}
\begin{proof}
    This follows previous proposition by induction.
\end{proof}

\section{The Interchangerator}
The definition of a monoidal bicategory is a clear generalisation of the definition of a monoidal category. This view does, however, obfuscate a rather fundamental property. In the same way that bicategories can be thought of as `multi-object monoidal categories', monoidal bicategories might be thought of as `multi-object \emph{braided} monoidal categories'. In particular a one-object monoidal bicategory is the same thing as a braided monoidal category. This is due to the presence of an invertible 2-cell that we call the \textdef{interchangerator}.

\begin{definition}
    Let $\BB$ be a monoidal bicategory and let $f\colon A_0\to A_1$ and $g\colon B_0\to B_1$ be 1-cells and, by an abuse of notation, let $A$ denote $\id_A$ for any object $A$. The \textdef{interchangerator} is an invertible 2-cell in the diagram below,
    \[
        \begin{tikzcd}
        A_0\tensor B_0
        \arrow[r, "f\tensor B_0"]
        \arrow[d, "A_0\tensor g"{swap}]
            &A_1\tensor B_0
            \arrow[d, "A_1\tensor g"]
            \arrow[ld, Rightarrow, shorten <=10pt, shorten >=10pt, "\chi_{f,g}"{decoration}]\\
        A_0\tensor B_1 
        \arrow[r, "f\tensor B_1"{swap}]
            &A_1\tensor B_1
        \end{tikzcd}
    \]
    that is given by the following composite.
\[
  \myinput{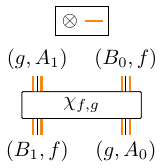}\quad\coloneqq\quad \myinput{tikz_MonBicat_BraidDef_1}
\]
\end{definition}
Note that here before the braid, on the left we right-tensor $f$ by the source of the 1-cell $g$, and on the right we left-tensor $g$ by the target of the 1-cell $f$. After the braid we right-tensor $f$ by the target of the 1-cell $g$, and we left-tensor $g$ by the source of the 1-cell $f$. 

\begin{figure}[h]
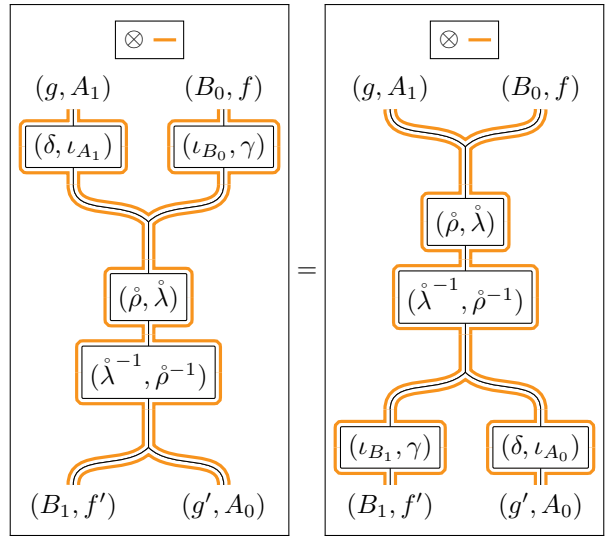

    \[
        \myinput{tikz_MonBicat_BraidNatProof_1}=
        \myinput{tikz_MonBicat_BraidNatProof_2}
    \]
    \caption{Naturality for the interchangerator.}
    \label{naturalityFigure}
    \end{figure}

There is something slightly subtle going on here. In light of strictification from the previous section it may seem that this interchangerator might as well be thought of as an identity 2-cell. The unitors are essentially identity 2-cells. But the join and split are \emph{not} inverses since they are indexed differently from one another. Using the labelling convention from \autoref{def:pseudofunctor}, the join is given by
\[
    \tensor_{(B_1,A_1),(B_0,A_1),(B_0,A_0)}\colon \tensor(g,A_1)\circ\tensor(B_0,f)\Rightarrow \tensor(g\circ B_0, A_1\circ f),
\]
but the split is given by
\[
    \tensor^{-1}_{(B_1,A_1),(B_1,A_0),(B_0,A_0)}\colon \tensor(g\circ B_0, A_1\circ f) \Rightarrow \tensor(B_1,f)\circ \tensor(g,A_0).
\]
Consider for a second a monoidal category. In a monoidal category we have the interchange \emph{law}. There is a sense in which this law is a sort of multi-object commutativity condition. It says that, so long as $f$ and $g$ are in different `lanes', `do $f$ then $g$' is the same thing as `do $g$ then $f$'. The interchangerator plays a similar role, and we will show that it adheres to braid-like conditions.
\begin{proposition}
    The interchangerator is invertible.
\end{proposition}
\begin{proof}
    This follows from the fact that $\tensor$ is a pseudofunctor, and so the compositor is invertible, as are the composition unitors.
\end{proof}
\begin{proposition}
    The interchangerator is natural in $f$ and $g$.
\end{proposition}
\begin{proof}
    Naturality says that, given any $\gamma\colon f\Rightarrow f'$ and $\delta\colon g\Rightarrow g'$, the equality in \autoref{naturalityFigure} holds.
    But this simply follows from the naturality of the compositor, the naturality of the unitors, and then the naturality of the compositor again.
\end{proof}

\begin{proposition}
    Given a pseudofunctor $F\colon \AA\to\BB$, and maps $f\colon A\to B$, $g\colon B\to C$ and $h\colon C\to D$ in $\AA$, the following 2-cells are equal.
    \[
        \myinput{tikz_MonBicat_BraidProof_1}=
        \myinput{tikz_MonBicat_BraidProof_2}
    \]
\end{proposition}
Note that here that on the left we omit the associator for $\BB$, but on the right we include the associator for $\AA$ with the pseudofunctor applied.
\begin{proof}
    Consider the associator diagram for pseudofunctors. By inverting the appropriate 2-cells we get the commuting diagram below.
    \[
    \begin{tikzcd}
    F(h)\circ (F(g)\circ F(f))
    \arrow[r, "\comp{\alpha}"]
        & (F(h)\circ F(g))\circ F(f)
        \arrow[d, "F \circ \iota_{F(f)}"]\\
    F(h)\circ F(g\circ f)
    \arrow[u, "\iota_{F(h)}\circ F^{-1}"]
    \arrow[d, "F"{swap}]
        & F(h\circ g)\circ F(f)\\
    F(h\circ (g\circ f))
    \arrow[r, "F(\comp{\alpha})"{swap}]
        & F((h\circ g)\circ f)
        \arrow[u, "F^{-1}"{swap}]
    \end{tikzcd}
    \]
    The string diagram on the left is the composite around the top of the diagram and the string diagram on the right is the composite around the bottom of the diagram.
    \end{proof}
    \begin{lemma}
        \label{lemma:interchangebraid}
        Let $\AA$ be a monoidal bicategory. For every quadruple of 1-cells $f\colon A_0\to A_1$, $g\colon A_1\to A_2$, $h\colon B_0\to B_1$ and $k\colon B_1\to B_2$ in $\AA$ the interchangerator adheres to the following generalised braid equations.
        \begin{figure}[H]
        \[
            \myinput{tikz_MonBicat_BraidProof_1}=
            \myinput{tikz_MonBicat_BraidProof_2}
        \]
        \end{figure}
        \begin{figure}[H]
        \[
            \myinput{tikz_MonBicat_BraidProof_1A}=
            \myinput{tikz_MonBicat_BraidProof_2A}
        \]
        \end{figure}
        \end{lemma}
    \begin{proof}
    We only show that the first generalised braid equality holds, the proof of the second is analogous. By definition, we know that the left-hand side of the equation is equal to \autoref{intBraid1}. By the above proposition, and the fact that $(\comp{\alpha},\comp{\alpha})$ is the associator for $\AA\times \AA$, we know that \autoref{intBraid1} and \autoref{intBraid2} are equal.

    Using naturality of the compositor we can move every 2-cell to the middle, so \autoref{intBraid2} equals \autoref{intBraid3}. By \autoref{prop:coherence}, as well as the unitor axiom, we can compose the three middle 2-cells, so we know that \autoref{intBraid3} is equal to \autoref{intBraid4}.

    Next, we can use the compositor axiom for pseudofunctors to swap which side the first join and the final split are on, and so \autoref{intBraid4} is equal to \autoref{intBraid5}. Then we can use the unitor axiom and \autoref{prop:coherence} to compose associators with unitors, so we know that \autoref{intBraid5} is equal to \autoref{intBraid6}.

    Finally, splitting the first and final centre 2-cells, and using the compositor axiom we have that \autoref{intBraid6} is equal to \autoref{intBraid7}. Now rearranging our composites we have that \autoref{intBraid7} is equal to \autoref{intBraid8}. And \autoref{intBraid8} is equal to the right-hand side of the equation by definition.

    \begin{figure}
        \centering
        \begin{minipage}{.5\textwidth}
            {\[\myinput{tikz_MonBicat_BraidProof_3}\]}
            \caption{}
            \label{intBraid1}
        \end{minipage}%
        \begin{minipage}{.5\textwidth}
        {\[\myinput{tikz_MonBicat_BraidProof_4}\]}
        \caption{}
        \label{intBraid2}
        \end{minipage}
    \end{figure}

    \begin{figure}
        \begin{minipage}{.5\textwidth}
            \[\myinput{tikz_MonBicat_BraidProof_5}\]
            \caption{}
            \label{intBraid3}
        \end{minipage}%
        \begin{minipage}{.5\textwidth}
        \[\myinput{tikz_MonBicat_BraidProof_6A}\]
        \caption{}
        \label{intBraid4}
        \end{minipage}
    \end{figure}

    \begin{figure}
        \begin{minipage}{.5\textwidth}
            \[\myinput{tikz_MonBicat_BraidProof_7A}\]
            \caption{}
            \label{intBraid5}
        \end{minipage}%
        \begin{minipage}{.5\textwidth}
        \[\myinput{tikz_MonBicat_BraidProof_8A}\]
        \caption{}
        \label{intBraid6}
        \end{minipage}
    \end{figure}

    \begin{figure}
        \begin{minipage}{.5\textwidth}
            \[\myinput{tikz_MonBicat_BraidProof_9}\]
            \caption{}
            \label{intBraid7}
        \end{minipage}%
        \begin{minipage}{.5\textwidth}
        \[\myinput{tikz_MonBicat_BraidProof_10}\]
        \caption{}
        \label{intBraid8}
        \end{minipage}
    \end{figure}
\end{proof}
\newpage
\section{Scalars}
\label{section:Scalars}
There is a classical result in the theory of monoidal categories, seemingly first proved by Kelly and Laplaza~\cite[prop.~6.1]{kelly1980coherence}, which tells us that every monoidal category gives rise to a commutative monoid. Given a monoidal category, $\AA$, we define the monoid of scalars in $\AA$ to be given by the endomorphisms at the unit, $\AA(I,I)$, equipped with composition. To understand this nomenclature note that in the category of Hilbert spaces this monoid is isomorphic to $\mathbb{C}$ equipped with multiplication.

The result of Kelly and Laplaza says that the monoid of scalars is always commutative, and the proof of this is a very straightforward application of the Eckmann-Hilton argument. Every pair of scalars $f,g\colon I\to I$ can be composed in two ways. We can multiply horizontally
\[
    I\cong I\tensor I\xrightarrow{f\tensor g} I\tensor I\cong I,
\] 
or vertically
\[
    I\xrightarrow{f}I\xrightarrow{g}I.
\]
But the functoriality of the tensor product guarantees that \[(f'\tensor g')\circ (f\tensor g)=(f'\circ f)\tensor (g'\circ g),\] and so, by the Eckmann-Hilton argument, the two multiplications agree and are commutative.

Monoidal bicategories also have scalars, but the collection of scalars now form a monoidal category. Let $\BB$ be a monoidal bicategory, then we define the monoidal category of scalars to be the category $\BB(I,I)$ with tensor product given by composition. In the previous section we saw how monoidal bicategories carry a sort of multi-object braid in the form of the interchangerator. In this section we see that the interchangerator gives rise to an actual braid structure on the monoidal category of scalars. One way to prove this is to repeat the strategy for the lower dimensional case. This would involve using a 2-dimensional version of the Eckmann-Hilton argument. Such a strategy does work, and this 2-dimensional version of the Eckmann-Hilton argument can be found in \textbook{Braided Monoidal Categories} by Joyal and Street~\cite[prop.~3]{joyal1986braided}. We will also see later that this fact follows almost immediately from a semi-strictification result as proved by Gordon, Power and Street~\cite[cor.~8.7]{street1995coherence}. Here we prove that the braiding exists directly, using the string diagram language developed in the previous chapter.
\begin{definition}
    Given a monoidal bicategory $\BB$, the \textdef{monoidal category of scalars} for $\BB$ is the category $\BB(I,I)$ equipped with the tensor product given by composition.
\end{definition}
\begin{example}
    The category of scalars for $\Rel$ consists of the two sets $*$ and $\varnothing$. This monoidal category is monoidally equivalent to the category of boolean truth values $
        (F\rightarrow T)
    $ equipped with logical conjunction.
\end{example}
\begin{example}
    The category of scalars for $\Bim_R$ is the category of $R$-$R$-bimodules. But since $R$ is commutative, this is just $R$-Mod with the usual tensor product.
\end{example}
\begin{example}
    The category of scalars for $\DGBim_R$ is the category of chain complexes of $R$-Modules. When $R$ is a field, this is equivalent to the category of graded vector spaces.
\end{example}
\begin{example}
    The bicategory $\VV$-Prof has the category
    \[
    \Prof(*,*)\cong \mathrm{Func}(*,\VV)\cong \VV    
    \]
    as its category of scalars. The composition of profunctors from $*$ to $*$ is just given by the tensor for $\VV$.
\end{example}
\begin{example}
    If $\CC$ has finite limits then the scalar category associated to the bicategory $\Span(\CC)$ consists of spans $*\leftarrow A\rightarrow *$, and since maps to $*$ are unique, this means the category of scalars is just the category $\CC$. Pulling back $A\to *\leftarrow B$ just gives the product $A\times B$ and so the monoidal structure is given by taking products.
\end{example}
\begin{example}
    If $M$ is a topological monoid then the scalars for $\Path(M)$ are loops at the unit $e$ and the monoidal category structure comes from concatenation of loops. Then the category of scalars is the monoidal category whose skeleton is the fundamental group of $M$.
\end{example}

In a monoidal category it is possible to prove that the left and right unitor agree at the unit. That is to say $I\tensor I\xrightarrow{l} I$ is equal to $I\tensor I\xrightarrow{r} I$. In a monoidal bicategory we can construct an invertible 2-cell that encodes this data. The following construction is given as a string diagram by Garner and Shulman~\cite[lem.~2.1]{garner2016enriched}, but here we give the construction using our decorated string diagram language.
\begin{proposition}
    In every monoidal bicategory there is an invertible 2-cell $\theta\colon l_I\Rightarrow r_I$.
\end{proposition}
\begin{proof}
Firstly note that there is an invertible modification in the diagram below,
\[
\begin{tikzcd}
I\tensor A
\arrow[rr, "l_A"]
\arrow[d, "r^\bullet_{I\tensor A}"{swap}]
    &\blank
    \arrow[d, "\theta_0", Rightarrow]
        &A\\
I\tensor (A\tensor I)
\arrow[r, "a_{I,A,I}"{swap}]
    &(I\tensor A)\tensor I
    \arrow[r, "l_{A\tensor I}"{swap}]
        &A\tensor I
        \arrow[u, "r_A"{swap}]
\end{tikzcd}    
\]
given by the following composite.
\[
  \myinput{tikz_MonBicat_theta1}
\] 
There are also invertible modifications in the following two diagrams,
\[
\begin{tikzcd}
I\tensor (I\tensor B)
\arrow[rr,"{l_{I\tensor B}}", bend left, ""{name=top,swap}]
\arrow[rr, "I\tensor l_B"{swap}, bend right, ""{name=bot}]
    && I\tensor l_B
    \arrow[from=top, to=bot, Rightarrow, "\theta_l"]
\end{tikzcd} 
\qquad
\begin{tikzcd}
    C\tensor I
    \arrow[rr,"{r^\bullet_{C\tensor I}}", bend left, ""{name=top,swap}]
    \arrow[rr, "r^\bullet_C\tensor I"{swap}, bend right, ""{name=bot}]
        && (C\tensor I)\tensor I
        \arrow[from=top, to=bot, Rightarrow, "\theta_r"]
    \end{tikzcd} 
\]
given by the following two composites
\[
\myinput{tikz_MonBicat_theta3} \qquad\myinput{tikz_MonBicat_theta2}
\]
Taking each of these modifications at their $I$'th component, and composing with the $I$'th component of the middle unitor $\mu$, we have an invertible 2-cell, that we call $\theta$, given below.
\[
\myinput{tikz_MonBicat_theta}
\]
\end{proof}
We now define the scalar braid morphism as follows. Firstly let {\color{c3}$-\tensor -$ be this colour}, let {\color{c1}$I\tensor -$ be this colour} and let {\color{c2}$-\tensor I$ be this colour}. Note that when we are composing scalars we can draw the interchangerator as in the diagram below.
\[\myinput{tikz_MonBicat_ScalarBraidDef_2}\qquad=\qquad\myinput{tikz_MonBicat_ScalarBraidDef_1}\]
Then the braid for scalars is defined by the following diagram.
\[\altchi_{f,g}\quad\coloneqq\quad\myinput{tikz_MonBicat_ScalarBraidDef_3}\] 
\begin{remark}
    Note that this particular example highlights why colouring pseudofunctors can be helpful. Without the colourings we see that $\altchi$ is given by \autoref{badBraid}.
    \begin{figure}
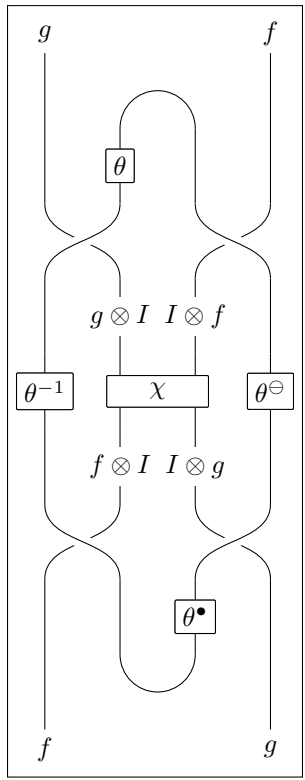

        \[\myinput{tikz_MonBicat_ScalarBraidDef_4}\]
        \caption{A confusing way to define the braid}
        \label{badBraid}
    \end{figure} 
    This diagram breaks geometric intuition. It appears as though there is a contractible loop on top of two separate strings. Once we contract the loop we're just left with $\chi$. But this makes no sense, since there is not an interchangerator morphism from $g\circ f$ to itself.
\end{remark}
\begin{proposition}
    \label{pseudoinvertible}
    Let ${\color{c1}F},{\color{c2}G}\colon \AA\to \BB$ be pseudofunctors, and let $n\colon {\color{c1}F}\Rightarrow {\color{c2}G}$ be a pseudo-invertible natural transformation. Then the following equality holds.
    \[
        \fbox{\tikzsetnextfilename{MonBiCat/ScalarBraidProof/Prop/1}
\begin{tikzpicture}[baseline={(current bounding box.center)}, xscale=1
]
\drawStartNodes{
    nB/${n^\bullet}$/1,
    f/{$f$}/1,
    n/{$n$}/0
}

\strBraidXX{n}{f}[/c2/c1][c2]

\strIdToX{nB}{f}[/c1/c2]
\strCupX{nB}{n}[/c1/c2]

\place{nB}{1}[-2]{f}
\place{n}{1}{nB}
\strCapX{nB}{n}[/c1/c2]

\strIdToX{f}{nB}[c1]
\strBraidXX{nB}{f}[/c1/c2][c1]

\strIdToX{n}{nB}[/c2/c1]

\drawEndNodes{
    nB/$n^\bullet$,
    n/$n$,
    f/$f$
}

\end{tikzpicture}}\quad=\quad\fbox{\tikzsetnextfilename{MonBiCat/ScalarBraidProof/Prop/3}
\begin{tikzpicture}[baseline={(current bounding box.center)}, xscale=1
]
\drawStartNodes{
    nB/${n^\bullet}$/1,
    f/{$f$}/1,
    n/{$n$}/0
}

\strIdX{nB}[/c1/c2][4]
\strIdX{f}[c2][4]
\strIdX{n}[/c2/c1][4]

\drawEndNodes{
    nB/$n^\bullet$,
    n/$n$,
    f/$f$
}

\end{tikzpicture}
}
    \]
\end{proposition}
\begin{proof}
    The left-hand side of the equation is equal to the following diagram by \autoref{prop:cupsandbraids}.
    \[\fbox{\tikzsetnextfilename{MonBiCat/ScalarBraidProof/Prop/2}
\begin{tikzpicture}[baseline={(current bounding box.center)}, xscale=1
]
\drawStartNodes{
    nB/${n^\bullet}$/1,
    f/{$f$}/1,
    n/{$n$}/0
}

\strBraidXX{n}{f}[/c2/c1][c2]

\strIdToX{nB}{f}[/c1/c2]
\strCupX{nB}{n}[/c1/c2]

\place{n}{-1}[-2]{f}
\place{nB}{-1}{n}
\strCapX{nB}{n}[/c1/c2]

\strIdToX{f}{n}[c1]
\strBraidXX{n}{f}[/c2/c1][c1]

\strIdToX{nB}{n}[/c1/c2]

\drawEndNodes{
    nB/$n^\bullet$,
    n/$n$,
    f/$f$
}

\end{tikzpicture}}
    \]    
    But this diagram is equal to the right-hand side of the equation by the definition of a pseudo-invertible 1-cell, and the invertibility of the naturalisor.
\end{proof}
\begin{lemma}
The family of 2-cells $(\altchi_{f,g})$ is a braiding for the monoidal category of scalars.
\end{lemma}
\begin{proof}
We must show that $\altchi$ is invertible, natural and satisfies the braid equations. Invertibility simply follows from the fact that $\altchi$ is a composite of invertible 2-cells. Naturality follows from the axioms of natural transformations and the fact that $\chi_{f,g}$ is natural in $f$ and $g$. 

We only show that one of the braid equations hold. The other follows analogously. Firstly note that by the above proposition, and by the following two equalities 
\[
    \myinput{tikz_MonBicat_ScalarBraidProof_A}=
    \myinput{tikz_MonBicat_ScalarBraidProof_B}\qquad
    \myinput{tikz_MonBicat_ScalarBraidProof_C}=
    \myinput{tikz_MonBicat_ScalarBraidProof_D}
\]
we know that \autoref{braidDiag1} and \autoref{braidDiag2} are equal.
\begin{figure}
    \begin{minipage}{.5\textwidth}
        \[\myinput{tikz_MonBicat_ScalarBraidProof_1}\]
        \caption{}
        \label{braidDiag1}
    \end{minipage}%
    \begin{minipage}{.5\textwidth}
    \[\myinput{tikz_MonBicat_ScalarBraidProof_2}\]
    \caption{}
    \label{braidDiag2}
    \end{minipage}
\end{figure}
Now by \autoref{lemma:interchangebraid} we know that \autoref{braidDiag2} is equal to \autoref{braidDiag3}. By \autoref{prop:narypseudofunctor} and the axioms of natural transformations we know that \autoref{braidDiag3} is equal to \autoref{braidDiag4}.
\begin{figure}
    \begin{minipage}{.5\textwidth}
        \[\myinput{tikz_MonBicat_ScalarBraidProof_3}\]
        \caption{}
        \label{braidDiag3}
    \end{minipage}%
    \begin{minipage}{.5\textwidth}
    \[\myinput{tikz_MonBicat_ScalarBraidProof_4}\]
    \caption{}
    \label{braidDiag4}
    \end{minipage}
\end{figure}
Thus, the braid equation holds.
\end{proof}
This is fortunate for us. Our ultimate aim is to show that every left-composition-closed, right-monoidal-closed bicategory can be enriched over $\cat{\BB(I,I)}$. Note that in $\cat{\VV}$ the tensor product is given by $\AA\tensor \BB$ whose objects are pairs of objects $(A,B)$ and whose hom-objects are given by \[
    (\AA\tensor \BB)((A,B), (A',B'))\coloneqq \AA(A,A')\tensor \BB(B,B').
\]
But note then that this definition only makes sense when $\VV$ is braided. The source of the composition morphism must be
\[
\AA(A',A'')\tensor \BB(B',B'')\tensor \AA(A,A')\tensor \BB(B,B')    
\]
and so without braiding the centre two hom-objects there is no canonical choice of composition.
\section{Semi-Strictification}
In the previous chapter we saw how every bicategory can be thought of as a 2-category, and so composition unitors and composition associators can really be thought of as equalities. Something similar holds for monoidal bicategories -- all the tensor unitors and tensor associators might as well be equalities -- but crucially the interchangerator is not degenerate. This gives us a sort of semi-strictification result. It turns out that every monoidal bicategory can be thought of as a Gray monoid.
\begin{definition}
    \label{def:GrayMonoid}
A \textdef{Gray monoid}, $\MM$, is a 2-category equipped with a unit object $I$, and for every object $A\in \MM$, a left-tensor 2-functor $L_A\colon \MM\to \MM$ and a right-tensor 2-functor $R_A\colon \MM\to \MM$ such that
        \begin{itemize}
            \item tensors of \emph{objects} are unambiguous: $L_A(B)=R_B(A)$;
            \item tensoring is unital: $L_I=R_{I}$ and they are both the identity pseudofunctor;
            \item tensoring is associative: $L_{L_A(B)}=L_A\circ L_B,\quad  R_{R_A(B)}=R_B\circ R_A,\quad R_B\circ L_A=L_A\circ R_B.$
        \end{itemize}
    We often write $A\tensor -$ for $L_A$ and $-\tensor A$ for $R_A$. By the above axioms, for any pair of sequences, $(A_1,...,A_n)$ and $(B_1,...,B_m)$, of \emph{objects} in $\MM$ we have an unambiguous 2-functor
    \[
    A_1\tensor ... \tensor A_m\tensor (-)\tensor B_1\tensor ... \tensor B_m\colon \MM\to \MM.
    \]
Additionally, for all objects, $A$ and $B$, and all arrows, $f\colon A_0\to A_1$ and $g\colon B_0\to B_1$, $\MM$ comes equipped with an invertible 2-cell, $\chi_{g,f}$, called the \textdef{interchangerator} given by the following diagram.
\[
\myinput{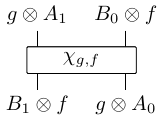}
\]
The interchangerator is unital in the sense of the following equalities.
\begin{equation}
    \label{def:GrayId}
    \myinput{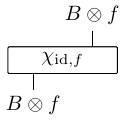}\quad = \quad \myinput{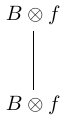}\qquad
    \myinput{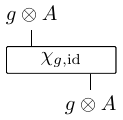}\quad=\quad     \myinput{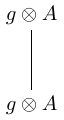}
\end{equation}
The interchangerator is natural in both variables, in the sense that, given $\phi\colon f\to f'$ and $\psi\colon g\to g'$ the following composites are equal.
\begin{equation}
    \label{def:GrayNat}
    \myinput{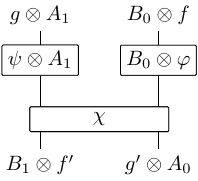}\quad = \quad \myinput{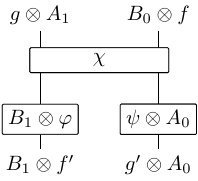}
\end{equation}
That is compatible with associativity in the sense that, given any object $A$, the following equalities hold.
\begin{equation}
    \label{def:GrayAss1}
    \myinput{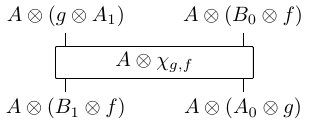}\quad = \quad\myinput{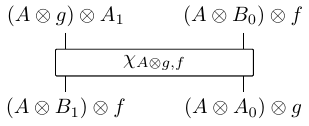}
\end{equation}
\begin{equation}
    \label{def:GrayAss2}
    \myinput{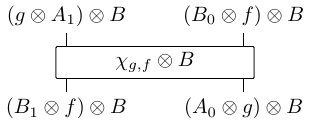}\quad = \quad    \myinput{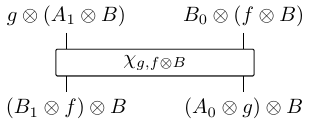}
\end{equation}
And finally, that adheres to the generalised braid laws, in the sense that, given 1-cells $f\colon A_0\to A_1$, $g\colon A_1\to A_2$, $h\colon B_0\to B_1$ and $k\colon B_1\to B_2$ the following composites are equal.
\begin{equation}
    \label{def:GrayBraid1}
    \myinput{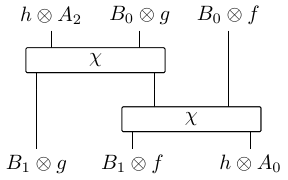}\quad = \quad \myinput{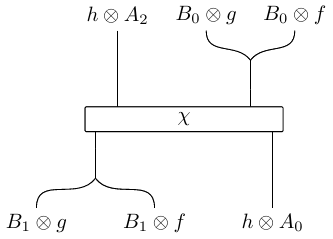}
\end{equation}
\begin{equation}
    \label{def:GrayBraid2}
    \myinput{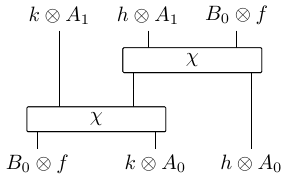} \quad = \quad  \myinput{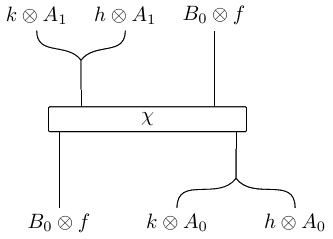}
\end{equation}
\end{definition}
\begin{remark}
Note that the definition of a Gray monoid as given by Day and Street~\cite[def.~1]{street1997monoidal} has a small error. It claims that $\chi_{f,g}$ is the identity 2-cell if $f$ and $g$ are both identity 1-cells. Of course, this should hold if $f$ \emph{or} $g$ are the identity.
\end{remark}
From the data of a Gray-monoid we can define a single pseudofunctor
\[
    \overline{\tensor} \colon \MM\times \MM\to \MM
\]
where 
\begin{itemize}
    \item given objects $A,B\in \MM$, we define $A\overline{\tensor} B\coloneqq A\tensor B$;
    \item given 1-cells $f\colon A_0\to A_1$ and $g\colon B_0\to B_1$ we define $f\overline{\tensor} g\coloneqq (g\tensor A_1)\circ (B_0\tensor f)$;
    \item given 2-cells $\phi\colon f\Rightarrow f'$ and $\psi\colon g\Rightarrow g'$, where $f,f'\colon A_0\to A_1$ and $g, g'\colon B_0\to B_1$ we define $\psi\tensor \phi\coloneqq (\psi\tensor A_1)\circ (B_0\tensor \phi)$;
    \item given 1-cells, $f\colon A_0\to A_1$, $f'\colon A_1\to A_2$, $g\colon B_0\to B_1$, $g'\colon B_1\to B_2$ we define the compositor to be the following 2-cell.
    \[
    \myinput{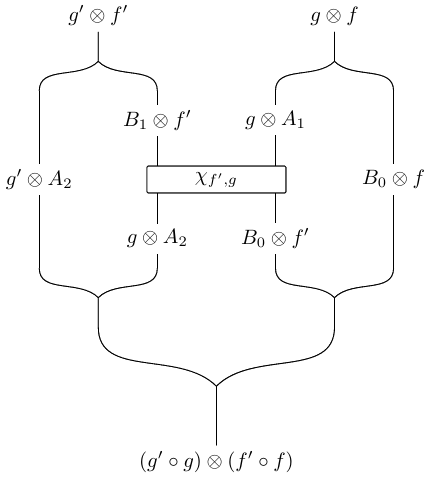}    
    \]
\end{itemize}
This gives every Gray monoid the structure of a monoidal bicategory where the unitors and associator are identity 1-cells. It is worth noting however, that this is not unique or even canonical. We chose, for example, to define $f\overline{\tensor} g$ as $(A_1\tensor g)\circ (f\tensor B_0)$, but we could just have easily defined it to be $(f\tensor B_0)\tensor (A_1\tensor g)$. This too defines a monoidal bicategory. 

Gordon, Power and Street~\cite[cor.~8.7]{street1995coherence} proved that every monoidal bicategory is essentially a Gray monoid.
\begin{theorem}
    Every monoidal bicategory is monoidally biequivalent to a Gray monoid.
\end{theorem}
As a result of this theorem we will prove certain results for Gray monoids rather than monoidal bicategories, but the proof for monoidal bicategories can always be recovered. This will greatly simplify certain diagrams that would otherwise be unreasonably large.
\section{Closed Monoidal Bicategories}
\label{section:ClosedMonoidal}
We've already seen that bicategories can be closed with respect to their composition. But just like monoidal categories, monoidal bicategories can also be closed with respect to their monoidal product. The definition of a monoidal-closed bicategory is deceptively simple.
\begin{definition}
A \textdef{right-monoidal-closed bicategory}, $\AA$ is a monoidal bicategory $\AA$ equipped with a pseudofunctor $\hom{-,-}\colon \AA^{\op}\times \AA\to \AA$ such that for any $A\in \AA$ there is a pseudoadjunction
\[
    -\tensor A\colon \AA\rightleftarrows \AA\cocolon \hom{A,-}.
\]
\end{definition}
\begin{remark}
    Note that bicategories with this structure are often referred to as right closed monoidal. However, since we are interested in bicategories that are closed compositionally \emph{and} monoidally, we reorder the adjectives and introduce hyphens to remove any potential ambiguities. For the sake of consistency we also apply this convention to categories.
\end{remark}
\begin{example}
    The bicategory $\Rel$ is right-monoidal-closed with $\hom{A,B}$ given by $B\times A$. The unit and counit for the pseudoadjunction are relations
    \[
    R\colon  B\to B\times (A\times A) \text{ and } S\colon  (B\times A)\times A\to B    
    \]
    which are given by the permuted diagonal relations.
\end{example}
\begin{example}
    The bicategory $\Bim_R$ is right-monoidal-closed with $\hom{A,B}$ given by $B\tensor A^{\op}$. The existence of the unit and counit for the pseudoadjunction follows from the fact that an $A$-$B$-bimodule is a left $A\tensor B^{\op}$ module and a right $A^{\op}\tensor B$ module. We need bimodules
    \[
    U\colon B\to (B\tensor A)\tensor A^{\op} \text{ and } E\colon (B\tensor A^{\op})\tensor A\to B
    \]
    which are given by the identity bimodule for $B\tensor A$, thought of as a $B$-$((B\tensor A)\tensor A^{\op})$ bimodule, and the identity bimodule for $B\tensor A^{\op}$, thought of as a $((B\tensor A^{\op})\tensor A)$-$B$ bimodule.
\end{example}
\begin{example}
    The bicategory $\DGBim_R$ is a right-monoidal-closed with $\hom{A,B}$ given by $A\tensor B^{\op}$. The unit and counit are constructed analogously to the above.
\end{example}
\begin{example}
    The bicategory $\VV$-Prof is right-monoidal-closed with $\hom{\AA,\BB}$ given by $\BB\tensor \AA^{\op}$. The unit and counit are given by the identity profunctor
    \[
    \Hom\colon \BB\tensor \AA \profto \BB\tensor \AA
    \]
    thought of as profunctors
    \[
    \Hom\colon \BB\profto (\BB\tensor \AA)\tensor \AA^{\op} \text{ and } \Hom\colon (\BB\tensor \AA^{\op})\tensor \AA\profto \BB.
    \]
\end{example}
\begin{example}
    If $\CC$ has all limits then the bicategory $\Span(\CC)$ is right-monoidal-closed with $\hom{A,B}=B\times A$. The unit and counit are given by spans
    \[
        \begin{tikzcd}
        &B\times A
        \arrow[ld, "p_B"{swap}]
        \arrow[rd, "\id\times \Delta"]\\
    B&& B\times A \times A
    \end{tikzcd}
    \text{ and }
    \begin{tikzcd}
        &B\times A
        \arrow[rd, "p_B"]
        \arrow[ld, "\Delta\times \id"{swap}]\\
    B\times A\times A && B
    \end{tikzcd}
    \]
    where $\Delta$ is the diagonal map, and where we've ignored associators to save space.
\end{example}
\begin{example}
    If $G$ is a topological group then $\Path(G)$ is a right-monoidal-closed bicategory with $\hom{x,y}=yx^{-1}$. The unit and counit are the identity paths since $\hom{x,yx}=yxx^{-1}=y$ and $\hom{x,y}x=yx^{-1}x=y$.
\end{example}
Notice that in all of the examples given above, the closed structure is given in relation to the tensor product. This is because the examples we are interested in constitute compact-closed bicategories, which we study in \autoref{section:CompactClosed}. Right-monoidal-closed bicategories have many similar properties to right-monoidal-closed categories. One particular feature we will make heavy use of is the ability to replace a given 1-cell with its `name'.
\begin{proposition}
    \label{prop:nameRealisation}
    In a right-monoidal-closed bicategory $\BB$, for every $A,B\in \BB$ there is a pseudonatural adjoint equivalence
    \[
        \name{(-)}\colon \BB(A,B) \rightleftarrows \BB(I,\hom{A,B}) \cocolon \unname{(-)}.
    \]
    We call $\name{f}$ the \textdef{name} of $f$ and $\unname{g}$ the \textdef{realisation} of $g$.
\end{proposition}
\begin{proof}
Let $E$ and $U$ be the adjoint equivalences associated to the adjunction given by the closed structure. Then the name functor, $\name{(-)}$, is defined by the composite
\[
    \BB(A,B)\xrightarrow{\BB(l, B)}\BB(I\tensor A, B) \xrightarrow{E} \BB(I, \hom{A,B})
\]
and the realisation functor $\unname{(-)}$ is defined by the composite
\[
    \BB(I,\hom{A,B})\xrightarrow{U}\BB(I\tensor A,B)\xrightarrow{\BB(l^\bullet,B)} \BB(A,B). 
\]
These composites give an adjoint equivalence since $l$ and $l^\bullet$ give an adjoint equivalence and so do $E$ and $U$.
\end{proof}
\begin{definition}
    If a bicategory is both right-monoidal-closed and composition-closed then we will refer to it as a \textdef{monoidal-closed composition-closed bicategory}.
\end{definition}

\chapter{Monoidal Actions}
Let's take a break from monoidal bicategories for a while to consider good, old-fashioned monoidal categories. In the next chapter we will prove one of our main results: every left-composition-closed, right-monoidal-closed bicategory $\BB$ can be enriched in $\cat{\BB(I,I)}$. We haven't yet given an account of what an enriched bicategory actually is, but in this case it involves replacing every category $\BB(A,B)$ with a $\BB(I,I)$-category, replacing the composition and identity functors with $\BB(I,I)$-functors and replacing the unitors and associators with $\BB(I,I)$-natural transformations. 

 Given some braided monoidal $\VV$, working with $\VV$-categories can often be unwieldy, especially when definitions are given indirectly. The aim of this chapter is to simplify this process by working with categories \emph{acted on} by $\VV$, called $\VV$-representations, instead of $\VV$-categories. If a monoidal category $\VV$ is like a higher-dimensional monoid, then a $\VV$-representation is like a higher dimensional $\VV$-set. Representations of this kind seem to have been first suggested by B\'enabou~\cite[def.~2.3]{benabou1967introduction} but, as with all good definitions, have been rediscovered independently many times. They also go by several different names, one that seems to have gained traction in the applied category theory community is ``actegory'' coined by McCrudden~\cite[sec.~3]{mccrudden2000categories}.

As proved by Gordon and Power~\cite[thm.~3.7]{gordon1997enrichment}, closed, strong $\VV$-representations are the `same thing' as copowered $\VV$-categories via a 2-equivalence of 2-categories
\[
    \Repp{\VV}^{\str}_{\cll} \xrightarrow{\sim} \cat{\VV}_{\operatorname{cop}}.
\]
The \emph{idea}, then, is that, given a bicategory $\BB$, we firstly show that every hom-category $\BB(A,B)$ is a closed $\BB(I,I)$-representation. We then show that all of the data of $\BB$ -- the composition, the identities, the associators and unitors --  are all compatible with the representation structure. Then we use the above 2-functor to change the base from representations to enriched categories. 

In reality this does not work, the reason being that enrichment requires a monoidal structure for $\Repp{\VV}$. Whilst $\Repp{\VV}$ has an obvious tensor product, the tensor of two closed $\VV$-representations may not be closed. 

If we were enriching over monoidal \emph{categories} this would be easily rectifiable. We could instead define a multicategory structure on $\Repp{\VV}_{\cll}$ -- induced by the tensor product on $\Repp{\VV}$ -- and then enrich over the multicategory, in the sense of Lambek~\cite[p.~106]{lambek1969deductive}. If we showed that the functor above could be given a multifunctor structure then this would give a base-change functor from closed representations to enriched categories. We are, however, trying to enrich \emph{bicategories} over monoidal bicategories. A similar approach to the above is likely possible using Leinster's highly abstract machinery of enriched generalised multicategories \cite{leinsterGeneral}. We instead use a more direct approach.

Firstly we show that there is an enriching 2-functor from \emph{all} closed representations to \emph{all} enriched categories:
\[
    \underline{(-)}\colon \Rep_{\cll} \to \EnCat.
\]
We then introduce two monoidal 2-categories: the 2-category of closed $\VV$-iterated representations includes all categories with a $\VV^n$ action for some $n\in \N$; the 2-category of $\VV$-iterated categories includes all categories enriched over $\VV^n$ for some $n\in \N$. These 2-categories are intended to behave similarly to the free monoidal 2-categories on $\Repp{\VV}_{\cll}$ and $\cat{\VV}$, and the restriction of the enriching 2-functor to these 2-categories is monoidal. Finally, we introduce the `collapsing 2-functor' which sends every $\VV^n$-category to its associated $\VV$-category via the tensor product for $\VV$, and show that this 2-functor is monoidal. As a result, when we enrich in the next chapter, it suffices to show that our enriching data is all given in terms of $\VV$-iterated representations, before collapsing this data to $\VV$-categories.
\section{Categorical Representations}
Given a monoid $M$, we can define an $M$-set to be a set on which $M$ acts associatively and unitally. Given a monoidal category $\VV$, we can define a $\VV$-representation to be a category on which $\VV$ acts oplaxly associatively and oplaxly unitally. In this sense the definition of a $\VV$-representation is a categorification of the definition of an $M$-set.
\begin{definition}
    Let $\VV$ be a monoidal category. A category $\CC$ \textdef{acted on by} $\VV$, or a $\VV$\textdef{-representation}, is a category $\CC$ equipped with a functor
    \[
        \act\colon \VV\times \CC \to \CC
    \]
    and natural transformations -- the \textdef{associator} $a$, and the \textdef{unitor} $u$ -- whose components are
    \[
        a\colon (V\tensor W)\act C \to V\act (W\act C) \text{ and } u\colon I\act C\to C. 
    \]
    These natural transformations adhere to coherence conditions given by the following commuting diagrams.
    \[
    \begin{tikzcd}
        ((U\tensor V)\tensor W)\act A
        \arrow[r, "\alpha\act A"]
        \arrow[d, "a"]
            &
            (U\tensor (V\tensor W))\act A
            \arrow[r, "a"]
                    &U\act ((V\tensor W)\act A)
                    \arrow[d, "U\act a"]\\
        (U\tensor V)\act (W\act A)     
        \arrow[rr, "a"]
                &\blank
                    & (U\act (V\act (W\act A)))
    \end{tikzcd}
    \]
    \[
        \begin{tikzcd}
        (I\tensor U)\act A
        \arrow[rd,"\lambda\act A"{swap}]
        \arrow[rr, "a"]
                && I\act (U\act A)
                \arrow[dl, "u"]\\
        \blank
            &U\act A
        \end{tikzcd}  
    \]
    \[
        \begin{tikzcd}
        (U\tensor I)\act A
        \arrow[rd, "\rho\act A"{swap}]
        \arrow[rr, "a"]
            &\blank
                &U\act(I\act A)
                \arrow[ld, "U\act u"]\\
        \blank
            &U\act A.
        \end{tikzcd}
    \]
    We call a representation \textdef{strong} if the associator and unitor are isomorphisms.
\end{definition}
\begin{remark}
We describe the action as oplaxly associative and unital because, as pointed out to the author by Leinster~\cite[]{leinsterPrivate}, if we curry the action to get a functor
\[
    \hat{\act}\colon \VV\to \hom{\CC,\CC}
\]
then it is easy to see that $\act$ is an action if and only if $\hat{\act}$ is an oplax functor.
\end{remark}
\begin{remark}
For the sake of space and readability we will quite often drop bracketing for actions. For example, we might write
\[
    V\tensor W \act C\text{ or } V\act W\act C    
\]
since the meaning is unambiguous.
\end{remark}
\begin{example}
    Any monoidal category $\VV$ is a $\VV$-representation where the action is given by the tensor product.
\end{example}
\begin{example}
     If $\BB$ is a bicategory, then for any pair of objects $A$ and $B$, $(\BB(B,B),\circ)$ forms a monoidal category and $\BB(A,B)$ is a $\BB(B,B)$-representation via post-composition. The associator and unitor for the representation are just given by the composition unitor and associator.
\end{example}
\begin{example}
    If $\BB$ is a bicategory, then for any pair of objects $A$ and $B$, $(\BB(A,A),;)$ forms a monoidal category, where $(f;g)\coloneqq (g\circ f)$, and $\BB(A,B)$ is a $\BB(B,B)$-representation via pre-composition. The associator and unitor for the representation are just given by the composition unitor and associator.
\end{example}
If we have an $N$-Set, $A$, and a monoid homomorphism, $f\colon M\to N$, then we can equip $A$ with an $M$-Set structure, called the pullback action. Something similar holds for oplax monoidal functors and $\VV$-representations.
\begin{proposition}
    \label{prop:pullbackRep}
    Let $\CC$ be a $\WW$-representation. If $F\colon \VV\to \WW$ is an oplax monoidal functor then there is a $\VV$-representation $\CC_F$ whose underlying category is $\CC$, and whose action is given by $F(V)\act_{\CC}C$ for all $C\in \CC$ and $V\in \VV$.
\end{proposition}
\begin{proof}
    This is an immediate consequence of the above remark. If we think of our action as being an oplax functor
    \[
    \hat{\act}\colon \WW\to \hom{\CC,\CC}    
    \]
    then the composite $\hat{\act}\circ F$ will also be an oplax functor.
\end{proof}
If a $\VV$-representation is a higher-dimensional $M$-set, then the next definition might be thought of as a higher-dimensional torsor, as observed by Willerton~\cite{willerton2013torsors}. Not only is there an action of the monoidal category on each $C\in \CC$, but between any $C,D\in \CC$ there lies an object in $\VV$.
\begin{definition}
    We call a $\VV$-representation $\CC$ \textdef{closed} if, for every $A\in \CC$ there is a functor $\close{A,-}\colon \CC\to \VV$ such that $\close{A,-}$ is right-adjoint to $-\act A$.
\end{definition}
A representation being closed does not guarantee that its pullback will be, unless the pullback is along a functor with a right adjoint.
\begin{proposition}
    \label{prop:pullbackAdj}
    Let $\CC$ be a closed $\WW$-representation. If $F\colon \VV\to \WW$ is an oplax monoidal functor with a right adjoint $G$, then $\CC_F$ is closed with the closed structure given by $
    G\circ \close{C,-}    
    $
    for all $C\in \CC$.
\end{proposition}
\begin{proof}
    This follows from the sequence of natural isomorphisms
    \[
    \CC(F(V)\odot C,D)\cong \CC(F(V),\close{C,D})\cong \CC(V, G(\close{C,D})),
    \]
    induced by the adjunction given by the closed structure, and the fact that $G$ is right adjoint to $F$.
\end{proof}
It is known that a monoidal-closed category can be enriched over itself. The following theorem generalises this idea to closed representations.
\begin{theorem}
    \label{lemma:fundamentalTheorem}
    If $\CC$ is a closed strong $\VV$-representation then there is a $\VV$-enriched category $\underline{\CC}$ whose underlying category is $\CC$, such that $\underline{\CC}(A,B)=\close{A,B}$ for all $A,B\in \CC$. Furthermore, $\underline{\CC}$ is copowered over $\VV$ and the underlying category of any copowered $\VV$-category is a $\VV$-representation.
\end{theorem}
The above theorem is a corollary of a result due to Gordon and Power~\cite[thm.~3.7]{gordon1997enrichment} but was also later proved more explicitly in the appendix of a paper by Janelidze and Kelly~\cite[sec.~6]{janelidze2000actions}. We will make heavy use of this result, so it is worth giving the constructions from the proof.
\begin{proof}[Sketch Proof]
Given a $\VV$-representation the enriched composition for $\underline{\CC}$ is defined as the adjunct to
\[
    \close{B,C}\tensor \close{A,B}\act A\xrightarrow{a}\close{B,C}\act \close{A,B}\act A\xrightarrow{\close{B,C}\tensor \epsilon} \close{B,C}\act B\xrightarrow{\epsilon} C    
\]
and the enriched identity for $\underline{\CC}$ is defined as the adjunct 
\[
I\act C\xrightarrow{u} C.
\]
Proving the axioms hold is a case of taking the adjuncts of the diagrams which give the axioms for representations, and showing that these adjuncts correspond to the axioms of an enriched category. 
\end{proof}
\begin{remark}
This theorem tells us that closed strong representations correspond to enriched categories with copowers. However, the proof also shows that all closed representations give rise to an enriched category, but the resulting enriched category is copowered if and only if the representation is strong.
\end{remark}
\begin{definition}
    \label{def:linearFunctor}
    Let $\CC$ be a $\VV$-representation, and let $\DD$ be a $\WW$-representation. A \textdef{linear functor} consists of a functor $F\colon \CC\to \DD$, a lax monoidal functor $(G,n,i)\colon \VV\to \WW$ and a natural transformation
    \[
        m\colon G(-)\act F(-)\to F(-\act -)    
    \]
    such that the following associativity and unitality diagrams both commute.
    \[
        \begin{tikzcd}[column sep=4.8em]
            G(V)\tensor G(W)\act F(C) 
            \arrow[r, "n\act F(C)"]
            \arrow[d, "a"{swap}]
                &G(V\tensor W)\act F(C)
                \arrow[d, "m"]\\
            G(V)\act G(W)\act F(C)
            \arrow[d, "G(V)\act m"{swap}]
                &F(V\tensor W\act C)
                \arrow[d, "F(a)"]\\
            G(V)\act F(W\act C)
            \arrow[r, "m"{swap}]
                &F(V\act W\act C)
        \end{tikzcd}  
    \]
    \[
        \begin{tikzcd}[column sep=4.8em]
            I\act F(C)
            \arrow[d,"u"{swap}]
            \arrow[r, "i\act F(C)"]
                &G(I)\act F(C)
                \arrow[d, "m"]\\
            F(C)
                &F(I\act C)
                \arrow[l, "F(u)"]
        \end{tikzcd}
    \]
    The composite of two linear functors $(F,G,m)\colon \AA\to \BB$ and $(F',G',m')\colon \BB\to \CC$ is given by the functor $F\circ F'$, the lax monoidal functor $G\circ G'$ and the composite natural transformation
    \[
        G'G(V)\act F'F(C)\xrightarrow{m'}F'(G(V)\act F(C))\xrightarrow{F'(m)} F'F(V\act C).  
    \]
    We call $F$ a $\VV$\textdef{-linear functor} if $G$ is the identity.
\end{definition}
Elsewhere in the literature, for example Capucci and Gavroni\'c's~\cite[p.~28]{capucci2022actegories} survey, linear functors are made up of \emph{strong} monoidal functors. Given that the representations in this thesis can be defined as oplax monoidal functors, and that oplax functors induce pullback representations, it may seem odd that the definition above includes a \emph{lax} monoidal functor instead. This is because oplax functors are the `wrong' choice for closed representations. We want to think of closed representations as enriched categories, and base changes of enriched categories only exist for \emph{lax} monoidal functors. Another justification for using lax, rather than oplax monoidal functors, is internal modules.
\begin{definition}
    Let $\AA$ be a $\VV$ representation, and let $(v,\mu,\iota)$ define a monoid in $\VV$. A $v$\textdef{-module} internal to $\AA$ is an object $k\in \AA$ equipped with a morphism $v\odot k\xrightarrow{\nu} k$, called the \textdef{action morphism}, such that the following associativity and unitality diagrams both commute.
    \[
    \begin{tikzcd}
    v\tensor v\odot k
    \arrow[d, "\mu"{swap}]
    \arrow[r, "a"]
        & v\odot v\odot k
        \arrow[r, "v\odot \nu"]
            &v\odot k
            \arrow[d, "\nu"]\\
    v\odot k
    \arrow[rr, "\nu"{swap}]
        &
            &k
    \end{tikzcd}    
    \]
    \[
    \begin{tikzcd}
        I\odot k
        \arrow[rr, "\iota\odot k"]
        \arrow[rd, "u"{swap}]
            && v\odot k
            \arrow[ld, "\nu"]\\
        &k
    \end{tikzcd}    
    \]
\end{definition}
\begin{example}
    Any module internal to a monoidal category is a module in the sense above, where the monoidal category acts on itself via the tensor product.
\end{example}
The above definition and the following lemma are both instances of what Baez and Dolan~\cite[p.~156]{BaezMicrocosm} call the `microcosm principle'. The idea being that if we take some sort of algebraic structure that can be defined on sets, and categorify it, then we can define that same algebraic structure internal to the categorification. What's more, any structure preserving functor will preserve the internal algebraic structure.
\begin{lemma}
    \label{lemma:linearPreservesModules}
    Let $\AA$ be a $\VV$-representation and $\BB$ be a $\WW$-representation. If $(V,\mu,\iota)$ is a monoid in $\VV$, $(K,\nu)$ is a $V$-module in $\AA$, and
    \[
    (F,G,m)\colon \AA\to \BB    
    \]
    is a linear functor, then $F(K)$ is a $G(V)$ module with action map given by the composite
    \[
    G(V)\odot F(K)\xrightarrow{m} F(V\odot K)\xrightarrow{\mu} F(K).
    \]
\end{lemma}
\begin{proof}
    Recall that if $G\colon \VV\to \WW$ is a lax monoidal functor, and $(V,\mu,\iota)$ defines a monoid in $\VV$, then $G(V)$ defines a monoid in $\WW$ with multiplication map
    \[
        G(V)\tensor G(V)\xrightarrow{n} G(V\tensor V)\xrightarrow{G(\mu)} G(V) 
    \]
    and unit map
    \[
        I\xrightarrow{i} G(I) \xrightarrow{G(\iota)} G(V).
    \]
    Then we must prove that $F(K)$ and the action map satisfy the unitality and associativity axioms. To see that unitality holds, notice that the internal polygons of the following diagram commute by the unitality axiom for $F$, the naturality of $m$ and the unitality for $K$.
    \[
    \begin{tikzcd}[column sep=4.2em, row sep=4.2em]
        I\odot F(K)
        \arrow[rrdd, bend right=70, "u"]
        \arrow[r, "i\odot F(K)"]
            &G(I)\odot F(K)
            \arrow[r, "G(\iota)\odot F(K)"]
            \arrow[d, "m"]
                &G(V)\odot F(K)
                \arrow[d, "m"]\\
            &F(I\odot K)
            \arrow[r, "F(\iota\odot K)"]
            \arrow[rd, "F(u)"{swap}]
                &F(V\odot K)
                \arrow[d, "F(\nu)"]\\
            & 
                &F(K)
    \end{tikzcd}    
    \]
    But the outer morphisms of this diagram give the unitality diagram for $F(\nu)\circ m$. To see that associativity holds, notice that the interior rectangles of \autoref{ch4d1} commute,
    \begin{figure}
    \[
    \rotatebox{90}{
    \begin{tikzcd}[every matrix/.append style={nodes={font=\scriptsize}}, column sep=4.2em, row sep=4.2em, ampersand replacement=\&]
        G(V)\tensor G(V)\odot F(K)
        \arrow[r, "a"]
        \arrow[d, "n\odot F(K)"{swap}]
            \&G(V)\odot G(V)\odot F(K)
            \arrow[r, "G(V)\odot m"]
                \&G(V)\odot F(V\odot K)
                \arrow[d, "m"]
                \arrow[r, "G(V)\odot \nu"]
                    \&G(V)\odot F(K)
                    \arrow[d, "m"]\\
        G(V\tensor V)\odot F(K)
        \arrow[d, "G(\mu)\odot F(K)"{swap}]
        \arrow[r, "m"]
            \&F(V\tensor V\odot K)
            \arrow[d, "F(\mu\odot K)"]
            \arrow[r, "F(a)"]
                \&F(V\odot V\odot K)
                \arrow[r, "F(V\odot \nu)"]
                    \&F(V\odot K)
                    \arrow[d, "F(\nu)"]\\
        G(V)\odot F(K)
        \arrow[r, "m"{swap}]
            \&F(V\odot K)
            \arrow[rr, "F(\nu)"{swap}]
                \&\blank
                    \&F(K)
    \end{tikzcd}    
    }
    \]
    \caption{}
    \label{ch4d1}
\end{figure}
    by associativity for $F$, naturality of $m$ and the associativity of $\nu$. But the outer morphisms of this diagram give the associativity diagram for $F(\nu)\circ m$. Thus, the associativity axiom holds, and so $F(K)$ is a $G(V)$ module with action morphism $F(\nu)\circ m$.
\end{proof}
As well as linear functors, we also have a notion of linear natural transformation.
\begin{definition}
\label{def:linearNT}
Given a $\VV$-representation $\AA$, a $\WW$-representation $\BB$ and linear functors
\[
(F,G,m), (F',G',m')\colon \AA\to \BB    
\]
a \textdef{linear natural transformation} $(\nu,\mu)\colon (F,G,m)\Rightarrow (F',G',m')$ consists of a natural transformation $\nu\colon F\Rightarrow F'$ and a monoidal natural transformation $\mu\colon G\Rightarrow G'$ such that the following diagram commutes for all $V\in \VV$ and $A\in \AA$.
\[
\begin{tikzcd}
    G(V)\act F(A)
    \arrow[r, "\mu\act \nu"]
    \arrow[d, "m"{swap}]
        & G'(V)\act F'(A)
        \arrow[d, "m'"]\\
    F(V\act A)
    \arrow[r, "\nu"{swap}]
        &F'(V\act A)
\end{tikzcd}
\]
Both vertical and horizontal composition of linear natural transformations are given component-wise. We call $\nu$ a $\VV$\textdef{-linear natural transformation} if $\mu$ is $\id\colon \Id\to \Id$. 
\end{definition}
\begin{definition}
These definitions give rise to a number of useful 2-categories:
\begin{itemize}
    \item For every monoidal $\VV$ we denote by $\VV\operatorname{-Rep}$ the 2-category of all $\VV$-representations, $\VV$-linear functors and $\VV$-linear natural transformations. 
    \item For every right-monoidal-closed $\VV$ we denote by $\VV\operatorname{-Rep}_{\operatorname{cl}}$ the 2-category of all closed $\VV$-representations, $\VV$-linear functors and $\VV$-linear natural transformations.
    \item We denote by $\operatorname{Rep}_{\operatorname{cl}}$ the 2-category of all closed representations over right-closed categories, linear functors and linear natural transformations.
\end{itemize}
\end{definition}
\begin{remark}
    As mentioned above our 2-category $\operatorname{Rep}$ differs from other categories of representations (or actegories) since our 1-cells use lax monoidal functors. In, for example, Capucci and Gavranovi\'c's review, their 2-category of all actegories uses strong monoidal functors. This is because they first define a 2-functor
    \[
    \rep{(-)}\colon\operatorname{MonCat}^{\operatorname{coop}}_{\operatorname{str}}\to \cat{2}    
    \]
    using \autoref{prop:pullbackRep} to send every strong monoidal functor $G\colon \VV\to \WW$ to a pullback functor
    \[
    G^*\colon \rep{\WW}\to \rep{\VV}.    
    \] 
    Then they use the bicategorical Grothendieck construction of Carrasco, Cegarra and Garz\'on~\cite[def.~7.2]{carrasco2009nerves} to define a 2-category of all representations.

    Using this construction it would also be possible to define a category of all representations where the 1-cells contain oplax monoidal functors. However, for a lax monoidal functor $G\colon \VV\to \WW$ there is no pullback functor
    \[
    G^*\colon \WW\to \VV.
    \]
    Thus, there is no 2-functor and no 2-categorical Grothendieck construction.
    \end{remark}
It is worth mentioning here the stronger version of the fundamental theorem of enriched category theory, \autoref{lemma:fundamentalTheorem}, proved by Gordon and Power~\cite[thm.~3.7]{gordon1997enrichment}. 
\begin{theorem}
    Let $\VV$ be a right-monoidal-closed category and let $\cat{\VV}_{\operatorname{cop}}$ be the 2-category of copowered $\VV$-categories, $\VV$-functors and $\VV$-natural transformations. Then there is an equivalence of 2-categories between $\VV\operatorname{-Rep}_{\operatorname{cl}}$ and $\cat{\VV}_{\operatorname{cop}}$.
\end{theorem}
\begin{remark}
In its original form this result was given for actions of bicategories and categories enriched \emph{over} bicategories. These are not to be confused with the enriched bicategories in the next chapter, which are enriched over monoidal bicategories.
\end{remark}
\section{Enrichment via Actions}
The theorem above gives slightly more categorical weight to the claim that copowered $\VV$-categories are `the same thing as' strong, closed $\VV$-representations. Unfortunately for our purposes this is not sufficient. Let $\BB$ be a left-composition-closed, right-monoidal-closed bicategory. In the next chapter our first aim will be to show that, for all $A,B\in \BB$, $\BB(A,B)$ can be replaced by a scalar-enriched category $\underline{\BB}(A,B)$. We will show that the category of scalars, $\BB(I,I)$, acts on $\BB(A,B)$ and that this action is closed. However, our second aim will be to show that we can replace the horizontal composition functor with an enriched functor. This would mean equipping the composition functor
\[
\circ \colon \BB(B,C)\times \BB(A,B)  \to \BB(A,C)
\]
with a $\BB(I,I)$-linear functor structure. The problem here is that there is no obvious choice of \emph{closed} action by $\BB(I,I)$ on $\BB(B,C)\times \BB(A,B)$. There \emph{is}, however, an obvious choice of closed action by $\BB(I,I)\times \BB(I,I)$. And there is also a monoidal functor \[\circ\colon \BB(I,I)\times \BB(I,I)\to \BB(I,I).\] This allows us to exploit the fact that any monoidal functor gives rise to a base change for enriched categories. In this section we give a 2-functor 
\[
    \underline{(-)}\colon \operatorname{Rep}_{\operatorname{cl}}\to \operatorname{EnCat}
\]from the 2-category of \emph{all} closed representations over right-monoidal-closed categories, to the 2-category of \emph{all} enriched categories.
\begin{remark}
As mentioned in the previous section, if we define our linear functors to be those with an oplax monoidal component, we can show that $\operatorname{Rep_{\cll}}$ is the result of a bicategorical Grothendieck construction. We can similarly construct the category of all enriched categories via a bicategorical Grothendieck construction: we take the 2-functor
\[
    {\cat{(-)}} \colon \operatorname{MonCat}_{\operatorname{lax}}\to \cat{2}  
\]
which sends every $\VV$ to the 2-category $\cat{\VV}$, every lax monoidal functor $F\colon \VV\to \WW$ to its associated base-change functor $(-)_F\colon \cat{\VV}\to\cat{\WW}$, and every monoidal natural transformation $\mu\colon F\to G$ to a pseudonatural transformation, made up of functors $\underline{\mu}\colon \CC_F\to \CC_G$ which are the identity on objects and whose associated morphisms are given by
\[\mu\colon F(\CC(A,B))\to G(\CC(A,B)).\]

This hints at a 2-functor arising from these Grothendieck constructions. There is, however, a variance issue: note that the sources of each of the 2-functors do not agree
\[
(-)\mhyphen\Rep_{\cll}\colon \MonCat_{\operatorname{oplax}}^{\operatorname{coop}}\to \cat{2} \text{ and } \cat{(-)}\colon \MonCat_{\operatorname{lax}}\to \cat{2}.
\]
However, the right adjoint of every oplax monoidal functor can be given the structure of a lax monoidal functor. So consider, instead, the 2-category $\operatorname{ClMonCat}_{\operatorname{LAdj}}$ whose objects are right-monoidal-closed categories, whose 1-cells $A\to B$ are given by pairs \[(L\colon A\to B,R\colon B\to A)\] where $L$ is oplax monoidal and $R$ is the lax monoidal functor right adjoint to $L$. Consider also the 2-category $\operatorname{ClMonCat}_{\operatorname{RAdj}}$ whose objects are right-monoidal-closed categories, whose 1-cells $A\to B$ are given by pairs \[(R\colon A\to B,L\colon B\to A)\] where $R$ is lax monoidal and $L$ is the oplax monoidal functor left adjoint to $R$. We have a 2-equivalence
\[
\operatorname{ClMonCat}_{\operatorname{LAdj}}^{\operatorname{coop}} \sim \operatorname{ClMonCat}_{\operatorname{RAdj}}.
\]
It seems likely then, using the theorem above and \autoref{prop:pullbackAdj}, that we could define a pseudonatural adjoint equivalence in the following diagram.
\[
\begin{tikzcd}
    \operatorname{ClMonCat}_{\operatorname{RAdj}}
    \arrow[r, Rightarrow, "\sim"]
    \arrow[d, "\cat{(-)_{\operatorname{cop}}}"{swap}, ""{name=left}]
        &\operatorname{ClMonCat}_{\operatorname{LAdj}}^{\operatorname{coop}}
        \arrow[d, "\clrep{(-)}", ""{swap, name=right}]\\
    \cat{2}
    \arrow[r, equal]
        &\cat{2}.
    \arrow[phantom, "\cong", from=right, to=left]
\end{tikzcd}
\]
Presumably this would induce a 2-functor between the two associated bicategorical Grothendieck constructions. However, we instead define our 2-functor directly, since we would like to work with linear functors where the lax monoidal functor component does not necessarily have an adjoint.
\end{remark}
\begin{definition}
    We define the 2-category of all enriched categories, $\EnCat$, as follows.
    \begin{itemize}
        \item The objects of $\EnCat$ are pairs $(\CC,\VV)$ where $\VV$ is a braided, right-monoidal-closed category, and $\CC$ is a $\VV$-enriched category.
        \item The 1-cells of $\EnCat$ are pairs $(F,G)\colon (\CC,\VV)\to (\DD,\WW)$ where $G\colon \VV\to \WW$ is a lax monoidal functor, $F\colon \CC_G\to \DD$ is an enriched functor, and where $(-)_G$ denotes the base-change 2-functor associated to $G$.
        \item The 2-cells of $\EnCat$ are given by pairs of the following form:
        \[
            \begin{tikzcd}[row sep=0.6em, column sep=3.2em]
            (\CC,\VV)
            \arrow[rr, bend left, "{(F,G)}", ""{swap, name=top}]
            \arrow[rr, bend right, "{(F',G')}"{swap}, ""{name=bottom}]
                &\phantom{\CC_{G'}}
                    &(\DD,\WW)
                    \arrow[Rightarrow, from=top, to=bottom, "{(\nu,\mu)}", to path={
                        (\tikztostart) 
                        -- (\tikztostart|-\tikztotarget.north) \tikztonodes}]         
            \end{tikzcd}
        \]
        where $\mu\colon G\Rightarrow G'$ is a monoidal natural transformation, and where $\nu$ is a $\WW$-natural transformation in the following diagram.
        \[
            \begin{tikzcd}[row sep=0.6em, column sep=3.2em]
            {\CC}_G
            \arrow[rr, bend left=75, "F", ""{swap, name=top}]
            \arrow[r, "\underline{\mu}"{swap}]
                &
                {\CC}_{G'}
                \arrow[Leftarrow, to=top, "\nu"{swap}, to path={
                    (\tikztostart) 
                    -- (\tikztostart|-\tikztotarget.south) \tikztonodes}]
                \arrow[r, "F'"{swap}]
                    &{\DD}
            \end{tikzcd}
        \]
        Here $\underline{\mu}$ is the functor which is the identity on objects and has associated functor morphism $
            \mu\colon G(\CC(A,B))\to G'(\CC(A,B))
        $
        for all objects $A$ and $B$ in $\CC$.
    \end{itemize}
    Composition is defined by the following rules.
    \begin{itemize}
        \item The composite of 1-cells,
        \[
        (\BB,\UU)\xrightarrow{(F,G)}(\CC,\VV)\xrightarrow{(F',G')}(\DD,\WW),
        \]
        is given by
        \[
        (\BB_{G})_{G'}\xrightarrow{F_{G'}} \CC_{G'}\xrightarrow{F} \DD
        \]
        in $\cat{\WW}$.
        \item The vertical composite of 2-cells in the following diagram
        \[
            \begin{tikzcd}[column sep=3.2em]
                (\CC,\VV)
                \arrow[rr, bend left=75, "{(F,G)}", ""{swap, name=top}]
                \arrow[rr, "{(F',G')}"{description}, ""{name=midtop}, ""{swap,name=midbot}]
                \arrow[rr, bend right=75, "{(F'',G'')}"{swap}, ""{name=bot}]
                    &&(\DD,\WW)
                \arrow[Rightarrow, from=top, to=midtop, "{(\nu,\mu)}", vertical]
                \arrow[Rightarrow, from=midbot, to=bot, "{(\nu',\mu')}", vertical]
            \end{tikzcd}
            \]
            is given by the following composite in $\cat{\WW}$.
            \[
                \begin{tikzcd}[row sep=0.6em, column sep=0.6em]
                {\CC}_G
                \arrow[rr, "\underline{\mu}"]
                \arrow[rrrr, bend left=75, "\underline{F}", ""{swap, name=top}]
                \arrow[rrrr, phantom, ""]
                    &&\underline{\CC}_{G'}
                    \arrow[Rightarrow, from=top, vertical, "{\nu}"]
                    \arrow[rd, "\underline{\mu'}"{swap}, bend right]\arrow[rr, "{F'}", ""{swap,name=midtop}]
                        &&{\DD}\\
                \blank
                    &\hphantom{{\CC}_{G''}}
                        &\blank
                            &{\CC}_{G''}
                            \arrow[Rightarrow, "{\nu'}"{outer sep=2pt}, from=midtop]
                            \arrow[ru, "{F''}"{swap}, bend right]
                                &\blank
                \end{tikzcd}
            \]
       \item The left whiskering in the following diagram
        \[
            \begin{tikzcd}[row sep=0.6em, column sep=3.2em]
            (\BB,\UU)
            \arrow[r, "{(F,G)}"]
                &(\CC,\VV)
            \arrow[rr, bend left, "{(F',G')}", ""{swap, name=top}]
            \arrow[rr, bend right, "{(F'',G'')}"{swap}, ""{name=bottom}]
                &\phantom{\CC_{G'}}
                    &(\DD,\WW)
                    \arrow[Rightarrow, from=top, to=bottom, "{(\nu,\mu)}", to path={
                        (\tikztostart) 
                        -- (\tikztostart|-\tikztotarget.north) \tikztonodes}]         
            \end{tikzcd}
        \]
        is given by the following left whiskering in $\cat{\WW}$.
        \[
            \begin{tikzcd}[row sep=0.6em, column sep=3.2em]
            (\BB_{G})_{G'}
                \arrow[r, "{F_{G'}}"]
                &\CC_{G'}
                \arrow[rr, bend left=75, "F'", ""{swap, name=top}]
                \arrow[r, "\underline{\mu}"{swap}]
                &
                {\CC}_{G''}
                \arrow[Leftarrow, to=top, "\nu"{swap, outer sep=2pt}, to path={
                    (\tikztostart) 
                    -- (\tikztostart|-\tikztotarget.south) \tikztonodes}]
                \arrow[r, "F''"{swap}]
                    &{\DD}
            \end{tikzcd}
        \]
        \item The right whiskering in the following diagram
        \[
            \begin{tikzcd}[row sep=0.6em,column sep=3.2em]
            (\BB,\UU)
            \arrow[rr, bend left, "{(F,G)}", ""{swap, name=top}]
            \arrow[rr, bend right, "{(F',G')}"{swap}, ""{name=bottom}]
                &\phantom{\CC_{G'}}
                    &(\CC,\VV)
                    \arrow[Rightarrow, from=top, to=bottom, "{(\nu,\mu)}", to path={
                        (\tikztostart) 
                        -- (\tikztostart|-\tikztotarget.north) \tikztonodes}]
                        \arrow[r, "{(F'',G'')}"] 
                         & (\DD,\WW)
            \end{tikzcd}
        \]
        is given by the following right whiskering in $\cat{\WW}$.
        \[
            \begin{tikzcd}[row sep=0.6em,column sep=3.2em]
            ({\BB}_{G})_{G''}
            \arrow[rr, bend left=75, "F_{G''}", ""{swap, name=top}]
            \arrow[rr, phantom, ""{name=bottom}]
            \arrow[r, "\underline{\mu}_{G''}"{swap}]
                &
                ({\BB}_{G'})_{G''}
                \arrow[Leftarrow, to=top, "\nu_{G''}"{swap}, to path={
                    (\tikztostart) 
                    -- (\tikztostart|-\tikztotarget.south) \tikztonodes}]
                \arrow[r, "F'_{G''}"{swap}]
                    &{\CC}_{G''}
                    \arrow[r, "F''"]
                        &\DD
            \end{tikzcd}
        \]
    \end{itemize}
\end{definition} 
The following result is a generalisation of the fact that if $F\colon \AA\to \BB$ is a lax monoidal functor between monoidal-closed categories, then $F$ is lax closed. That is to say that whenever $F$ is lax monoidal, there are natural transformations
\[
F(\hom{A,B})\to \hom{(FA,FA)} \text{ and } I\to F(I)
\]
which adhere to certain coherence axioms. The proof follows similarly.
\begin{lemma}
    \label{lemma:actionFuncIsEnrichedFunc}
    If $\CC$ is a closed $\VV$-representation, $\DD$ is a closed $\WW$-representation and $(F,G,m)\colon \CC\to \DD$ is a linear functor between them, then $F$ is the underlying functor of a $\WW$-functor $\underline{F}\colon\underline{\CC}_G\to \underline{\DD}$, where $\underline{\CC}_G$ is the $\WW$-enriched category given by base-change under $G$.
\end{lemma}
\begin{proof}
In order to save space, let $\close{A,B}$ be denoted by $B^A$ for both the closed structure in $\CC$ and the closed structure in $\DD$. Additionally, let $\eta$ and $\epsilon$ be the unit and counit of the adjunction that gives the right closed structure. We define $\underline{F}$ to be equal to $F$ on objects. Since $F$ preserves actions we have a map
\[
G(B^A)\act FA\xrightarrow{m} F(B^A\act A) \xrightarrow{F(\epsilon)} FB.
\]
The adjunct of this map gives a map $\underline{F}\colon G(B^A)\to  FB^{FA}$ which will be functor morphism associated to our enriched functor. We firstly show that the diagram for identities, \autoref{diag:enrichedId}, commutes. Since $F$ is action preserving we already know that \autoref{unitorAction} commutes, and we will show that the adjunct of \autoref{unitorAction} is \autoref{diag:enrichedId}.
\begin{figure}[h]
\begin{minipage}{0.5\textwidth}
    \[
        \begin{tikzcd}
        I
        \arrow[dd, "\underline{\id}"{swap}]
        \arrow[dr, "\underline{\id}"]\\
        \blank  &G(C^C)
        \arrow[dl, "\underline{F}"]\\
        FC^{FC}
        \end{tikzcd}    
        \]
    \caption{}
    \label{diag:enrichedId}
\end{minipage}    
\begin{minipage}{0.5\textwidth}
\[\begin{tikzcd}[column sep=4.8em]
    I\act F(C)
    \arrow[d,"u"]
    \arrow[r, "i\act F(C)"]
        &G(I)\act F(C)
        \arrow[d, "m"]\\
    F(C)
        &F(I\act C)
        \arrow[l, "F(u)"]
\end{tikzcd}
\]
\caption{}
\label{unitorAction}
\end{minipage}
\end{figure}

Firstly note that the adjunct of $u$ is 
$
I\xrightarrow{\epsilon}(I\act FC)^{FC}\xrightarrow{(u\act FC)^{FC}}FC^{FC} 
$
which is, by definition, the enriched identity map for $FC$. So the adjunct of the left-hand side of \autoref{unitorAction} is the left-hand side of \autoref{diag:enrichedId}. 

Now consider the right-hand side of \autoref{unitorAction}. The adjunct of the right-hand side is equal to the composite
\[
I
\xrightarrow{\eta}
(I\act FC)^{FC}
\xrightarrow{(i\act FC)^{FC}}
(GI\act FC)^{FC}
\xrightarrow{m^{FC}}
(F(I\act C))^{FC}
\xrightarrow{Fu^{FC}}
FC^{FC},
\]
which is the right-hand path in \autoref{ch4diag2}, where we have omitted certain functors on morphisms for the sake of readability.

\tikzstyle{start}=[to path={(\tikztostart.#1) -- (\tikztotarget)}]
\tikzstyle{finish}=[to path={(\tikztostart) -- (\tikztotarget.#1)}]
\begin{figure}[ht]
\[
\begin{tikzcd}[ampersand replacement=\&, column sep=1.2em]
I
\arrow[r, "\eta"]
\arrow[d, "i"{swap}]
    \&(I\act FC)^{FC}
    \arrow[d, "i"]\\
GI
\arrow[d, "\eta"{swap}]
    \&(GI\act FC)^{FC}
    \arrow[r,"m"]
    \arrow[d, "\eta"]
        \&F(I\act C)^{FC}
        \arrow[d, "\eta"]
        \arrow[rd, start=0,equal]\\
G((I\act C)^C)
\arrow[d, "u"{swap}]
    \&(G((I\act C)^C)\act FC)^{FC}
    \arrow[d, "u"]
        \& (F((I\act C)^{C}\act C))^{FC}
        \arrow[r,"\epsilon"]
        \arrow[d, "u"]
            \&(F(I\act C))^{FC}
            \arrow[d, "u"]\\
G(C^C)
    \arrow[r, "\eta"]
        \& (G(C^C)\act FC)^{FC}
        \arrow[r,"m"]
            \&(F(C^C\act C))^{FC}
            \arrow[r, "\epsilon"]
                \& FC^{FC}
\end{tikzcd}
\]
\caption{}
\label{ch4diag2}
\end{figure}
The leftmost rectangle of \autoref{ch4diag2} commutes by naturality of $\eta$, the middle rectangle commutes by naturality of $m$, the rightmost rectangle commutes by naturality of $\epsilon$, and the triangle commutes by the zigzag identity for adjunctions. But the leftmost path is $\underline{\id}$ and the bottom path is $\underline{F}$. Thus, the adjunct of \autoref{unitorAction} is \autoref{diag:enrichedId}, and since \autoref{unitorAction} commutes so does \autoref{diag:enrichedId}.

Now we need to show that the composition square, \autoref{diag:enrichedComp}, commutes, where $\underline{\circ}$ denotes the enriched composition morphism. Since $F$ is linear we already know that \autoref{diag:functorandcomposite} commutes. We will show that the adjunct of \autoref{diag:functorandcomposite} is \autoref{diag:enrichedComp}.
\begin{figure}
    \[
    \begin{tikzcd}
    G(B^A)\tensor G(C^B)
    \arrow[d, "\underline{F}\tensor \underline {F}"{swap}]
    \arrow[r, "\underline{\circ}"]
        &G(C^A)
        \arrow[d, "\underline{F}"]\\
    FB^{FA}\tensor FC^{FB}
    \arrow[r, "\underline{\circ}"{swap}]
        &FC^{FA}
    \end{tikzcd}
    \]
    \caption{}
    \label{diag:enrichedComp}
\end{figure}
\begin{figure}
\[\begin{tikzcd}[column sep=3.6em]
G(C^B)\tensor G(B^A) \act FA
\arrow[r, "n \act FA"]
\arrow[d, "a"{swap}]
    &G(C^B\tensor B^A)\act FA
    \arrow[d,"m"]\\
G(C^B)\act G(B^A) \act FA
\arrow[d, "G(C^B)\act m"{swap}]
    &F(C^B\tensor B^A\act A)
    \arrow[d, "F(a)"]\\
G(C^B)\act F((B^A)\act A)
\arrow[d, "m"{swap}]
    &F(C^B\act B^A\act A)
    \arrow[d,"F(C^B\act \epsilon)"]\\
F(C^B\act B^A\act A)
\arrow[d, "F(C^B\act \epsilon)"{swap}]
    &F(C^B\act B)
    \arrow[d, "F(\epsilon)"]\\
F(C^B\act B)
\arrow[r, "F(\epsilon)"{swap}]
    &F(C)
\end{tikzcd}   
\] 
\caption{}
\label{diag:functorandcomposite}
\end{figure}
Recall that, by definition, the composition morphism in the enriched category $\underline{\CC}_G$ is given by
\begin{align*}
    G(C^B)\tensor G(B^A)\xrightarrow{n} G(C^B\tensor B^A)&\xrightarrow{G(\eta)} G((C^B\tensor B^A\act A)^A) \xrightarrow{G((a)^A)} \\G((C^B\act B^A\act A)) &\xrightarrow{G((C^B\act \epsilon)^A)} G((C^B\act B)^A)\xrightarrow{G(\epsilon^A)} G(C^A).
\end{align*}
Then the right-hand side of \autoref{diag:enrichedComp} is adjunct to the left-hand path of \autoref{ch4diag3}.
But \autoref{ch4diag3} commutes by naturality of $m$, $\epsilon$ and the zigzag laws for identities. Thus, the adjunct of the right-hand path of \autoref{diag:enrichedComp} is the right-hand path of \autoref{diag:functorandcomposite}. 

\begin{figure}[p]
    \[
    \begin{tikzcd}
        G(C^B)\tensor G(B^A)\act FA
        \arrow[d, "n"]\\
        G(C^B\tensor B^A)\act FA
        \arrow[r, "m"]
        \arrow[d, "\eta"]
            &F(C^B\tensor B^A\act A)
            \arrow[d, "\eta"]
            \arrow[dr, equal, start=0]\\
        G((C^B\tensor B^A\act A)^A)\act FA 
        \arrow[d, "a"]
            &F((C^B\tensor B^A\act A)^A\act A) 
            \arrow[r, "\epsilon"]
            \arrow[d,"a"]
                &F(C^B\tensor B^A\act A) 
                \arrow[d, "a"]\\
        G((C^B\act B^A\act A)^A)\act FA 
        \arrow[d, "\epsilon"]
            &F((C^B\act B^A\act A)^A\act A) 
            \arrow[d, "\epsilon"]
                &F(C^B\act B^A\act A) 
                \arrow[d, "\epsilon"]\\
        G((C^B\act B)^A)\act FA 
        \arrow[d, "\epsilon"]
            &F((C^B\act B)^A\act A) 
            \arrow[d,"\epsilon"]
                &F(C^B\act B) 
                \arrow[d,"\epsilon"]\\
        G(C^A)\act FA
        \arrow[r, "m"]
            &F(C^A\act A)
            \arrow[r, "\epsilon"]
                &FC
    \end{tikzcd}
    \]
    \caption{}
    \label{ch4diag3}
\end{figure}

Now we show that the left-hand path of \autoref{diag:enrichedComp} is adjunct to the left-hand path of \autoref{diag:functorandcomposite}. Recall that the composition morphism in the enriched category $\underline{\DD}$ is given by the adjunct of
\[
C^B\tensor B^A\act A \xrightarrow{a} C^B\act B^A \act A \xrightarrow{C^B\act \epsilon} C^B\act B\xrightarrow{\epsilon} C.\] Then the left-hand side of \autoref{diag:enrichedComp} is the adjunct of the left-hand side of \autoref{ch4diag4}.

But \autoref{ch4diag4} commutes by naturality of $a$ and $\epsilon$, and the zigzag identities for adjunctions. Thus, the adjunct of the left-hand path of \autoref{diag:functorandcomposite} is the left-hand path of \autoref{diag:enrichedComp}. Hence, the adjunct of \autoref{diag:functorandcomposite} is \autoref{diag:enrichedComp} and the latter commutes since the former does.
 
Finally, to prove that the underlying functor, $\underline{F}_0$, of $\underline{F}$ is $F$ we need to show that for every $f\colon A \to B$, $\underline{F}_0(f)=F(f)$. Note that by definition $\underline{F}_0(f)$ is given by the adjunct of
\[
    I\xrightarrow{i} G(I)\xrightarrow{G(\eta)}G((I\act A)^A)\xrightarrow{G(u^A)} G(A^A)\xrightarrow{G(f^A)} G(B^A) \xrightarrow{\underline{F}} FB^{FA},
\]
which is the bottom path in \autoref{ch4diag5}. On the other hand, $F(f)$ is given by the adjunct of
\[
    I\xrightarrow{i} GI\xrightarrow{\eta} (GI\act FA)^{FA} \xrightarrow{(F(m))^{FA}} (F(I\act A))^{FA}\xrightarrow{(F(u))^{FA}}FA^{FA}\xrightarrow{(Ff)^{FA}}FB^{FA},
\]
which is the top path in \autoref{ch4diag5}. But the rectangles in \autoref{ch4diag5} commute by the naturality of $\eta$, $m$ and $\epsilon$. The triangle commutes by the zigzag laws for adjunctions. Thus, $\underline{F}_0=F$ and we are done. 
\end{proof}

\begin{figure}[h]
    \[
    \begin{tikzcd}[column sep=1.2em]
    I
    \arrow[d, "i"]\\
    GI
    \arrow[r, "\eta"]
    \arrow[d, "\eta"]
        & (GI\act FA)^{FA}
        \arrow[r, "m"]
        \arrow[d, "\eta"]
            & (F(I\act A))^{FA}
            \arrow[d, "\eta"]
            \arrow[rd, equal, start=0]\\
    G((I\act A)^A)
    \arrow[d, "u"]
        &(G((I\act A)^A))^{FA}\act FA
        \arrow[d, "u"]
            & (F((I\act A)^A\act A))^{FA}
            \arrow[d, "u"]
            \arrow[r, "\epsilon"]
                &(F(I\act A))^{FA}
                \arrow[d, "u"]\\
    G(A^A)
    \arrow[d, "f"]
        &(G(A^A)\act FA)^{FA}
        \arrow[d,"f"]
            &(F(A^A\act A))^{FA}
            \arrow[d, "f"]
                &FA^{FA}
                \arrow[d, "f"]\\
    G(B^A)
    \arrow[r, "\eta"]
        & (G(B^A)\act FA)^{FA}
        \arrow[r, "m"]
            & (F(B^A\act A))^{FA}
            \arrow[r, "\epsilon"]
                & FB^{FA}
    \end{tikzcd}    
    \]
    \caption{}
    \label{ch4diag5}
    \end{figure}
\begin{figure}[p]
\[
\rotatebox{90}{
\begin{tikzcd}[ampersand replacement=\&, every matrix/.append style={nodes={font=\scriptsize}}, cramped]
G(C^B)\tensor G(B^A)\act FA
\arrow[r, "a"]
\arrow[dd, "\eta\tensor \eta"{swap}]
    \& G(C^B)\act G(B^A)\act FA
    \arrow[d, "\eta"]
    \arrow[equal, ddrr, start=0]\\
\blank
    \& G(C^B)\act (G(B^A\act FA)^{FA})\act FA
    \arrow[d, "\eta"]\\
(G(C^B) \act FB)^{FB}\tensor (G(B^A)\act FA)^{FA}\act FA
\arrow[dd,"m\tensor m"{swap}]
    \&(G(C^B) \act FB^{FB})\act (G(B^A)\act FA)^{FA}\act FA
    \arrow[d, "m"]
    \arrow[r, "\epsilon"]
        \&\bullet
        \arrow[r, "\epsilon"]
            \& G(C^B)\act G(B^A)\act FA
            \arrow[d, "m"]\\
\blank
    \&
    (G(C^B)\act FB)^FB\act (F(B^A\act A))^{FA} \act FA
    \arrow[d, "m"]
        \&\blank
            \& G(C^B)\act F(B^A\act A)
            \arrow[d, "m"]\\
(F(C^B\act B))^{FB}\tensor (F(B^A\act A))^{FA} \act FA
\arrow[dd, "\epsilon\tensor \epsilon"{swap}]
    \& (F(C^B\act B))^{FB}\act (F(B^A\act A))^{FA} \act FA
    \arrow[d, "\epsilon"]
        \&\blank
            \& F(C^B\act B^A\act A)
            \arrow[d, "\epsilon"]\\
\blank
    \& (F(C^B\act B))^{FB}\act FB^{FA} \act FA
    \arrow[d, "\epsilon"]
        \&\blank
            \& F(C^B\act B)
            \arrow[d, "\epsilon"]\\
FC^{FB}\tensor FB^{FA} \act FA
\arrow[r, "a"{swap}]
    \&FC^{FB}\act FB^{FA} \act FA
    \arrow[r, "\epsilon"{swap}]
        \&FC^FB\act FB
        \arrow[r, "\epsilon"{swap}]
            \&FC
\end{tikzcd}
}
\]
\caption{}
\label{ch4diag4}
\end{figure}
\newpage
\begin{lemma}
    \label{lemma:actionNTIsEnrichedNT}
    If $(\CC,\VV)$ and $(\DD,\WW)$ are closed representations over right-monoidal-closed categories, and we have a linear natural transformation
    \[
        \begin{tikzcd}[row sep=0.6em, column sep=3.2em]
        (\CC,\VV)
        \arrow[rr, bend left, "{(F,G,m)}", ""{swap, name=top}]
        \arrow[rr, bend right, "{(F',G',m')}"{swap}, ""{name=bottom}]
            &\phantom{\CC_{G'}}
                &(\DD,\WW)
                \arrow[Rightarrow, from=top, to=bottom, "{(\mu,\nu)}", to path={
                    (\tikztostart) 
                    -- (\tikztostart|-\tikztotarget.north) \tikztonodes}]         
        \end{tikzcd}
    \]
    then $\nu$ is the underlying natural transformation of an enriched natural transformation 
    \[
        \begin{tikzcd}[row sep=0.6em, column sep=3.2em]
        {\underline{\CC}}_G
        \arrow[rr, bend left=75, "\underline{F}", ""{swap, name=top}]
        \arrow[r, "\underline{\mu}"{swap}]
            &
            {\underline{\CC}}_{G'}
            \arrow[Leftarrow, to=top, "\nu"{swap}, to path={
                (\tikztostart) 
                -- (\tikztostart|-\tikztotarget.south) \tikztonodes}]
            \arrow[r, "\underline{F'}"{swap}]
                &{\underline{\DD}}
        \end{tikzcd}
    \]
    where $\underline{\mu}$ is the enriched functor that is the identity on objects and has functor morphism \[
        G(\CC(A,B))\xrightarrow{\mu} G'(\CC(A,B)).
    \]
\end{lemma}
\begin{proof}
We will prove that the adjunct of $\nu$, that is to say the map $I\to F'A^{FA}$,
gives $\underline{\nu}$. In order to show that $\underline{\nu}$ defines an enriched natural transformation we must show that \autoref{diag:enrichedNT} commutes. Since $(\nu,\mu)$ defines a linear natural transformation we already know that \autoref{diag:linNT} commutes. We will show that \autoref{diag:enrichedNT} is the adjunct of \autoref{diag:linNT}. 
\begin{figure}[H]
\[\begin{tikzcd}
G(B^A)
\arrow[r, "r^{-1}"]
\arrow[d, "l^{-1}"{swap}]
    &G(B^A)\tensor I
    \arrow[r, "(\underline{F'}\circ \underline{\mu})\tensor \underline{\nu}"]
        &F'B^{F'A}\tensor F'A^{FA}
        \arrow[d, "\underline{\circ}"]\\
I\tensor G(B^A)
\arrow[r, "\underline{\nu}\tensor \underline{F}"{swap}]
    &F'B^{FB}\tensor FB^{FA}
    \arrow[r, "\underline{\circ}"{swap}]
        & F'B^{FA}
\end{tikzcd}\]
\caption{}
\label{diag:enrichedNT}
\end{figure}
\begin{figure}[H]
\[\begin{tikzcd}
    G(B^A)\act FA
    \arrow[r,"\nu\act \mu"]
    \arrow[d, "m"{swap}]
        &G'(B^A)\act F'A
        \arrow[d, "m'"]\\
    F(B^A\act A)
    \arrow[d, "\nu"{swap}]
        &F'(B^A)
        \arrow[d, "F'(\epsilon)"]\\
    F'(B^A\act A)
    \arrow[r, "F'(\epsilon)"{swap}]
        &F'(A)
\end{tikzcd}\]
\caption{}
\label{diag:linNT}
\end{figure}
Firstly note that the right-hand side of \autoref{diag:linNT} is adjunct to the topmost path in \autoref{ch4.1} which commutes for the following reasons.

\tikzstyle{start}=[to path={(\tikztostart.#1) -- (\tikztotarget)}]
\begin{figure}[p]
\[
\rotatebox{90}{
\begin{tikzcd}[ampersand replacement=\&, every matrix/.append style={nodes={font=\scriptsize}, column sep=1.2em}, cramped]
    G(B^A)
    \arrow[d, "r^{-1}"]
    \arrow[r, "\eta"]
        \&(G(B^A)\act FA)^{FA}
        \arrow[d, "l^{-1}"]
        \arrow[dddrr, rounded corners, to path={
            (\tikztostart.0) -- 
            (\tikzcdmatrixname-1-5.center)-- 
            (\tikzcdmatrixname-4-5.center)-- 
            (\tikztotarget.0)}, equal]
            \&\blank
                \&\blank
                    \&\blank\\
    G(B^A)\tensor I
    \arrow[d, "\eta"]
    \arrow[ddd, "\id\tensor \underline{\nu}"{description}, 
    bend right=80]
        \&(G(B^A)\tensor I\act FA)^{FA}
        \arrow[r, "a"]
        \arrow[d, "\eta"]
            \&(G(B^A)\act I\act FA)^{FA}
            \arrow[d, "\eta"]
            \arrow[dr, equal, start=0]\\
    G(B^A)\tensor (I\act FA)^{FA} 
    \arrow[d, "u"]
        \&\bullet
        \arrow[d, "u"]
            \&\bullet
            \arrow[d, "u"]
            \arrow[r, "\epsilon"]
                \& (G(B^A)\act I\act FA)^{FA}
                \arrow[d, "u"]
                    \&\blank\\
    G(B^A)\tensor FA^{FA} 
    \arrow[d, "\nu"]
        \&\bullet
        \arrow[d, "\nu"]
            \&\bullet
            \arrow[d, "\nu"]
                \& (G(B^A)\odot FA)^{FA}
                \arrow[d, "\nu"]
                    \&\blank\\
    G(B^A)\tensor F'A^{FA}
    \arrow[dddd, "(\underline{F'}\circ \underline{\mu})\tensor \id"{description, pos=0.4}, bend right=80]
    \arrow[d, "\mu"]
        \&\bullet
        \arrow[d, "\mu"]
            \&\bullet
            \arrow[d, "\mu"]
               \& (G(B^A)\odot F'A)^{FA}
                \arrow[d, "\mu"]
                \\
    G'(B^A)\tensor F'A^{FA} 
    \arrow[d, "\eta"]
        \&\bullet
        \arrow[d, "\eta"]
            \&\bullet
            \arrow[d, "\eta"]
                \&(G'(B^A)\odot F'A)^{FA}
                \arrow[rd, equal, start=0]
                \arrow[d, "\eta"]
                \\
    (G'(B^A)\odot F'A)^{F'A}\tensor F'A^{FA} 
    \arrow[d,"m"]
        \&\bullet
        \arrow[d,"m"]
            \&\bullet
            \arrow[d,"m"]
                \&\bullet
                \arrow[d,"m"]
                \arrow[r, "\epsilon"]
                    \& (G'(B^A)\odot F'A)^{FA}
                    \arrow[d, "m"]
                    \\
    (F'(B^A\act A))^{F'A}\tensor F'A^{FA} 
    \arrow[d, "\epsilon"]
        \&\bullet
        \arrow[d, "\epsilon"]
            \&\bullet
            \arrow[d, "\epsilon"]
                \&\bullet
                \arrow[d, "\epsilon"]
                    \&(G(B^A)\odot FA)^{FA}
                    \arrow[d, "\epsilon"]
                        \\
    F'B^{F'A}\tensor F'A^{FA} 
    \arrow[rrrr, bend right=20, "\underline{\circ}"{swap}]
    \arrow[r, "\eta"]
        \&\left(
            \begin{cdaligned}
            FA\tensor &F'A^{FA}\\
            &\act F'B^{F'A}
        \end{cdaligned}
        \right)^{FA}
        \arrow[r, "a"]
            \&\left(
                \begin{cdaligned}
                FA\act &F'A^{FA}\\
                &\act F'B^{F'A}
            \end{cdaligned}
            \right)^{FA}
            \arrow[r, "\epsilon"]
                \&(F'A\act F'B^{F'A})^{FA}
                \arrow[r, "\epsilon"]
                    \&F'B^{FA}
\end{tikzcd}
}
\]
\caption{}
\label{ch4.1}
\end{figure}
The top right pentagon commutes because of the unit axiom for representations; the triangles commute because of the zigzag identity for adjunctions; the squares commute because of the naturality of $\eta$, $a$ and $\epsilon$; and the outer globules commute by definition.

Thus, the top right path of \autoref{diag:linNT} is adjunct to the top right path of \autoref{diag:enrichedNT}. 
An analogous commutative diagram proves that the bottom left path of \autoref{diag:linNT} is adjunct to the bottom left path of \autoref{diag:enrichedNT}. Thus, \autoref{diag:linNT} is adjunct to \autoref{diag:enrichedNT}, and since \autoref{diag:linNT} commutes, so does \autoref{diag:enrichedNT}. Thus, $\underline{\nu}$ defines an enriched natural transformation. The fact that the underlying natural transformation of $\underline{\nu}$ is $\nu$ follows immediately from the definition of $\underline{\nu}$.
\end{proof}
We now have an assignment $\underline{(-)}\colon \Rep_{\cll}\to \EnCat$. The fact that this assignment preserves identities follows immediately from the definition. All that remains is to prove that this assignment preserves composition. For the sake of readability, all referenced commutative diagrams can be found at the end of the section.
\begin{proposition}
Let $\BB$ be a $\UU$-representation, $\CC$ be a $\VV$-representation and $\DD$ be a $\WW$-representation. If $(F,G,m)\colon \BB\to \CC$ and $(F',G',m')\colon \CC\to \DD$ are both linear functors, then $\underline{F'\circ F}=\underline{F'}\circ \underline{F}_{G'}$ where $\underline{F}_{G'}$ is $\underline{F}$ after the base change induced by $G'$.
\end{proposition}
\begin{proof}
By definition $\underline{F'\circ F}$ and $\underline{F'}\circ \underline{F}_{G'}$ agree on objects. Thus, it suffices to prove that their underlying associated functor morphisms are equal. The functor morphism for $\underline{F}_{G'}$ is given by the top most path \autoref{ch4.2}, and the functor morphism for $\underline{F'}$ is given by the right-hand path. The functor morphism for $\underline{F'\circ F}$ is given by the left-hand and bottom path. But the rectangles in the diagram commute by naturality of $\eta$, $m'$ and $\epsilon$, and the triangle commutes by the zigzag identity for adjunctions. Thus, the functor morphisms for $\underline{F'}\circ \underline{F}_{G'}$ and $\underline{F'\circ F}$ are equal.
\tikzstyle{finish}=[to path={(\tikztostart) -- (\tikztotarget.#1)}]
\end{proof}

We now move on to the composition of linear natural transformations, beginning with vertical composition.
\begin{lemma}
    Let $\CC$ be a closed $\VV$-representation and $\DD$ be a closed $\WW$-representation. For every pair of linear natural transformations
    \[
    \begin{tikzcd}[column sep=3.2em]
        (\CC,\VV)
        \arrow[rr, bend left=75, "{(F,G,m)}", ""{swap, name=top}]
        \arrow[rr, "{(F',G',m')}"{description}, ""{name=midtop}, ""{swap,name=midbot}]
        \arrow[rr, bend right=75, "{(F'',G'',m'')}"{swap}, ""{name=bot}]
            &&(\DD,\WW)
        \arrow[Rightarrow, from=top, to=midtop, "{(\nu,\mu)}"]
        \arrow[Rightarrow, from=midbot, to=bot, "{(\nu',\mu')}"]
    \end{tikzcd}
    \]
    we have that the enriched natural transformation
    \[
        \begin{tikzcd}[row sep=0.6em,column sep=3.2em]
        {\underline{\CC}}_G
        \arrow[rr, bend left=75, "\underline{F}", ""{swap, name=top}]
        \arrow[r, "\underline{\mu'\cdot \mu}"{swap}]
            &
            {\underline{\CC}}_{G'}
            \arrow[Leftarrow, to=top, "\underline{\nu'\cdot \nu}"{swap}, to path={
                (\tikztostart) 
                -- (\tikztostart|-\tikztotarget.south) \tikztonodes}]
            \arrow[r, "\underline{F''}"{swap}]
                &{\underline{\DD}}
        \end{tikzcd}
    \]
    is equal to the following vertical composite in $\cat{\WW}$.
    \[
        \begin{tikzcd}[row sep=0.6em, column sep=0.6em]
        \underline{\CC}_G
        \arrow[rr, "\underline{\mu}"]
        \arrow[rrrr, bend left=75, "\underline{F}", ""{swap, name=top}]
        \arrow[rrrr, phantom, ""]
            &&\underline{\CC}_{G'}
            \arrow[Rightarrow, from=top, vertical, "\underline{\nu}"]
            \arrow[rd, "\underline{\mu'}"{swap}, bend right]\arrow[rr, "\underline{F'}", ""{swap,name=midtop}]
                &&\underline{\DD}\\
        \blank
            &\hphantom{\underline{\CC}_{G''}}
                &\blank
                    &\underline{\CC}_{G''}
                    \arrow[Rightarrow, "\underline{\nu'}", from=midtop, vertical]
                    \arrow[ru, "\underline{F''}"{swap}, bend right]
                        &\blank
        \end{tikzcd}
    \]
    \end{lemma}
    \begin{proof}
        By definition of the whiskering of enriched natural transformations, we know that $(\underline{\mu}\circ\underline{\nu'})$ is given by components
        \[
            I \xrightarrow{\underline{\nu'}_{\underline{\mu}(A)}} F''(\underline{\mu}A)^{F'(\underline{\mu}A)}
        \]
        but since $\underline{\mu}$ is the identity on objects this means that $(\underline{\mu}\circ\underline{\nu'})$ is given by
        \[
            I\xrightarrow{\underline{\nu'}_A} F''A^{F'A}.
        \]
        Then note that, by definition, the composite
        $(\mu\circ \underline{\nu'})\cdot(\underline{\nu})$ is given by
        \[
            I\xrightarrow{l^{-1}} I\tensor I\xrightarrow{\nu\tensor \nu'} F''^{F'A}\tensor F'A^{FA}\xrightarrow{\underline{\circ}} F''^{FA}  
        \]
        in $\WW$, which is the top right-hand path in \autoref{4.14}. This diagram commutes by naturality, the zigzag identities of adjunctions and the unit axioms of $\WW$-representations. But the bottom leftmost path is, by the naturality of $u$, equal to
    \[
        I\xrightarrow{\eta} (I\odot FA)^{FA}\xrightarrow{u} FA^{FA} \xrightarrow{\nu'\cdot\nu} F''A^{FA} 
    \]
    which is by definition $\underline{\nu'\cdot \nu}$, completing the proof.
    \end{proof}
\begin{lemma}
Let $\BB$ be a closed $\UU$-representation, let $\CC$ be a closed $\VV$-representation and let $\DD$ be a closed $\WW$-representation where $\UU$, $\VV$ and $\WW$ are all right-closed. For every linear natural transformation given by the whiskering
        \[
            \begin{tikzcd}[row sep=0.6em,column sep=3.2em]
            (\BB,\UU)
            \arrow[r, "{(F,G,m)}"]
                &(\CC,\VV)
            \arrow[rr, bend left, "{(F',G',m')}", ""{swap, name=top}]
            \arrow[rr, bend right, "{(F'',G'',m'')}"{swap}, ""{name=bottom}]
                &\phantom{\CC_{G'}}
                    &(\DD,\WW)
                    \arrow[Rightarrow, from=top, to=bottom, "{(\nu,\mu)}", to path={
                        (\tikztostart) 
                        -- (\tikztostart|-\tikztotarget.north) \tikztonodes}]         
            \end{tikzcd}
        \]
        we have that the enriched natural transformation
        \[
            \begin{tikzcd}[row sep=0.6em, column sep=3.2em]
            {\underline{\BB}}_{G'\circ G}
            \arrow[rr, bend left=75, "\underline{F'\circ F}", ""{swap, name=top}]
            \arrow[r, "\underline{\mu\circ G}"{swap}]
                &
                {\underline{\BB}}_{G''\circ G}
                \arrow[Leftarrow, to=top, "\underline{\nu\circ F}"{swap}, to path={
                    (\tikztostart) 
                    -- (\tikztostart|-\tikztotarget.south) \tikztonodes}]
                \arrow[r, "\underline{F''\circ F}"{swap}]
                    &{\underline{\DD}}
            \end{tikzcd}
        \]
        is equal to the enriched natural transformation given by the following whiskering in $\cat{\WW}$.
        \[
            \begin{tikzcd}[row sep=0.6em,column sep=3.2em]
            (\underline{\BB}_{G})_{G'}
                \arrow[r, "{\underline{F}_{G'}}"]
                &\underline{\CC}_{G'}
                \arrow[rr, bend left=75, "\underline{F}'", ""{swap, name=top}]
                \arrow[r, "\underline{\mu}"{swap}]
                &
                {\underline{\CC}}_{G''}
                \arrow[Leftarrow, to=top, "\underline{\nu}"{swap}, to path={
                    (\tikztostart) 
                    -- (\tikztostart|-\tikztotarget.south) \tikztonodes}]
                \arrow[r, "\underline{F''}"{swap}]
                    &\underline{\DD}
            \end{tikzcd}
        \]
\end{lemma}
\begin{proof}
Firstly note that $\underline{\BB}_{G'\circ G}=(\BB_{G})_{G'}$. Now note that by definition $\underline{\nu\circ F}$ is given by the composite
\[
I\xrightarrow{\eta} I\odot F'FA^{F'FA} \xrightarrow{u} F'FA^{F'FA} \xrightarrow{\nu_{F}} F''FA^{F'FA}   
\]
but by the definition of left whiskering in $\cat{\WW}$ we know that $\underline{\nu}\circ \underline{F}_G$ is given by 
\[
I \xrightarrow{\underline{\nu}_F} F''FA^{F'FA}
\]
which is, by definition, the composite above. Thus, the two natural transformations are equal.
\end{proof}
\begin{lemma}
    Let $\BB$ be a closed $\UU$-representation, let $\CC$ be a closed $\VV$-representation and let $\DD$ be a closed $\WW$-representation where $\UU$, $\VV$ and $\WW$ are all right-closed. For every linear natural transformation given by the whiskering
        \[
            \begin{tikzcd}[row sep=0.6em,column sep=3.2em]
            (\BB,\UU)
            \arrow[rr, bend left, "{(F,G)}", ""{swap, name=top}]
            \arrow[rr, bend right, "{(F',G')}"{swap}, ""{name=bottom}]
                &\phantom{\CC_{G'}}
                    &(\CC,\VV)
                    \arrow[Rightarrow, from=top, to=bottom, "{(\nu,\mu)}", to path={
                        (\tikztostart) 
                        -- (\tikztostart|-\tikztotarget.north) \tikztonodes}]
                        \arrow[r, "{(F'',G'')}"] 
                         & (\DD,\WW)
            \end{tikzcd}
        \]
            we have that the enriched natural transformation
            \[
                \begin{tikzcd}[row sep=0.6em, column sep=3.2em]
                {\underline{\BB}}_{G''\circ G}
                \arrow[rr, bend left=75, "\underline{F''\circ F}", ""{swap, name=top}]
                \arrow[r, "\underline{G''\circ \mu}"{swap}]
                    &
                    {\underline{\BB}}_{G''\circ G'}
                    \arrow[Leftarrow, to=top, "\underline{F''\circ \nu}"{swap}, to path={
                        (\tikztostart) 
                        -- (\tikztostart|-\tikztotarget.south) \tikztonodes}]
                    \arrow[r, "\underline{F''\circ F'}"{swap}]
                        &{\underline{\DD}}
                \end{tikzcd}
            \]
            is equal to the enriched natural transformation given by the following whiskering in $\cat{\WW}$.
            \[
                \begin{tikzcd}[row sep=0.6em, column sep=3.2em]
                (\underline{\BB}_{G})_{G''}
                \arrow[rr, bend left=75, "{\underline{F}_{G''}}", ""{swap, name=top}]
                \arrow[rr, phantom, ""{name=bottom}]
                \arrow[r, "{\underline{\mu}_{G''}}"{swap}]
                    &
                    (\underline{\BB}_{G'})_{G''}
                    \arrow[Leftarrow, to=top, "{\underline{\nu}_{G''}}"{swap}, to path={
                        (\tikztostart) 
                        -- (\tikztostart|-\tikztotarget.south) \tikztonodes}]
                    \arrow[r, "{\underline{F'}_{G''}}"{swap}]
                        &\underline{\CC}_{G''}
                        \arrow[r, "{\underline{F''}}"]
                            &\underline{\DD}
                \end{tikzcd}
            \]
    \end{lemma}
\begin{proof}
Note that, by definition, the composite $\underline{F''}\circ \underline{\nu}_{G''}$ is given by
\[
I\xrightarrow{i} G'(I)\xrightarrow{G'(\underline{\nu})} G'(F'A^{FA})\xrightarrow{G'(\underline{F''})} G'(F''F'A^{F''FA}),
\]
which is the right-hand path in \autoref{4.15}.  
This diagram commutes by naturality, the zigzag identities for adjunctions and the unitor axioms for representations. But the bottom leftmost path of \autoref{4.15} is, by definition, $\underline{F''\circ \nu}$. Hence, the two natural transformations are equal.
\end{proof}

Thus, we have shown that the assignment given by $\underline{(-)}$ preserves all composition, and so it gives a 2-functor. Furthermore, by construction the assignment preserves underlying categories, functors, and natural transformations. Thus, we have the following theorem. 
\begin{theorem}
    \label{theorem:enrichingfunctor}
    There is a 2-functor, called the enriching 2-functor, $\underline{(-)}\colon \Rep_{\cll}\to \EnCat$ such that the following diagram commutes
    \[
    \begin{tikzcd}
        \Rep_{\cll}
        \arrow[dd, "\underline{(-)}"{swap}]
        \arrow[rd, "(-)_0"]\\
        \blank
            & \Cat\\
        \EnCat
        \arrow[ur, "(-)_0"{swap}]
    \end{tikzcd}    
    \]
\end{theorem}
\begin{proof}
The first two lemmata of the section give a way to define the 2-functor. The remaining lemmata prove that this definition does, in fact, preserve vertical and horizontal composition. The diagram commutes by definition of the enrichment.
\end{proof}

\section{The Base Change 2-Functor}
The aim of this chapter is to provide the necessary machinery to enrich a bicategory in the monoidal 2-category of $\VV$-categories. In order to do this we want to first enrich over closed $\VV$-representations and then use a monoidal 2-functor to change the base to $\VV$-categories. Unfortunately this is not possible, since the product of two closed $\VV$-representations is not necessarily closed, and so $\Repp{\VV}_{\cll}$ does not have a monoidal structure. If we forget the 2-categorical data, then $\Repp{\VV}_{\cll}$ has a multicategory structure induced by the monoidal structure of $\Repp{\VV}$. 

A result of Hermida~\cite[thm.~7.2]{Hermida} tells us that there is a 2-adjunction
\[
    F\colon \mathrm{MultiCat} \rightleftarrows \mathrm{MonCat}\cocolon U.
\]
Note that enriching in the multicategory $U(\cat{\VV})$ is the same thing as enriching in $\cat{\VV}$. Thus, we could first enrich in the multicategory $\Repp_{\VV}{\cll}$ and then use the enriching 2-functor to construct a multifunctor,
\[
    \Repp{\VV}_{\cll}\xrightarrow{\underline{(-)}} U(\cat{\VV}),
\]
to change our base of enrichment. Alternatively, if we wanted to avoid enriching over multicategories then we could instead construct the adjunct of this multifunctor. Explicitly we could construct the monoidal functor
\[
    F(\Repp{\VV}_{\cll})\xrightarrow{F\underline{(-)}} FU(\cat{\VV})\xrightarrow{\epsilon} \cat{\VV}
\]
where $\epsilon$ is the counit of the 2-adjunction above. In this section we construct monoidal 2-categories that behave similarly to $F(\Repp{\VV}_{\cll})$ and $FU(\cat{\VV})$, we show that the enriching 2-functor induces a monoidal 2-functor between them, and we give a `collapsing' monoidal 2-functor analogous to the counit above. These 2-categories and 2-functors give the machinery necessary for the next chapter.

The following monoidal 2-category is slightly larger than the free monoidal category, $F(\Repp{\VV}_{\cll})$, in that it contains additional objects and 1-cells, but is easier to describe for our purposes.

In what follows we will invoke Mac Lane's coherence theorem and treat $\VV$ as a strict, unbiased monoidal category, so that $\VV^n$ is identified with any \[\prod_{i=1}^p\VV^{n_i}\] such that $\sum_{i=1}^p n_i=n$. Note also that we include $0$ in the natural numbers, $\VV^0$ is given by the terminal category, $*$, and $\tensor^0\colon *\to \VV$ is given by the functor that picks out the unit.
    \begin{definition}
        The 2-category, $\Repp{\VV}_{\cll}^\times$, of \textdef{closed $\VV$-iterated representations}, is the 2-category described as follows.
        \begin{itemize}
            \item The objects are given by pairs $(\AA,n)$, where $n$ is a natural number and $\AA$ is a closed $\VV^{n}$-representation.
            \item  The 1-cells $(\AA,n)\to (\BB,p)$ are given by linear functors
            \[
                (F,G,m)\colon (\AA,\VV^n)\to (\BB,\VV^p)
            \]
            where $G$ is of the form
            \[
            \VV^n \xrightarrow{\prod_{i=1}^p \tensor^{n_i}} \VV^p 
            \]
            and where $\sum_{i=1}^p n_i=n$.
            \item The 2-cells are given by linear natural transformations where the monoidal natural transformation is a coherence isomorphism.
        \end{itemize}
        This 2-category has a monoidal structure described as follows. 
        \begin{itemize}
            \item The tensor product of objects $(\AA,n)$ and $(\BB,p)$ is given by
            \[
            (\AA,n)\tensor (\BB,p)\coloneqq (\AA\times\BB,n+p),    
            \]
            where the action on $\AA\times \BB$ is given by
            \[
                (\VV^{n+p})\times (\AA\times \BB)\xrightarrow{\sim} (\VV^n\times \AA)\times (\VV^p\times \BB) \xrightarrow{\act_\AA\times \act_\BB} \AA\times \BB,
            \] and the closed structure of $\AA\times \BB$ is given by
            \[
                (\AA\times \BB)^{\op}\times (\AA\times \BB)\xrightarrow{\sim} (\AA^{\op}\times \AA)\times (\BB^{\op}\times \BB)\xrightarrow{\close{-,-}\times \close{-,-}} \VV^n\times \VV^p \xrightarrow{\sim}\VV^{n+p}.
            \]
            \item The tensor product of 1-cells \[(F,G,m)\colon (\AA,n)\to (\BB,p)\text{ and }(F',G',m')\colon (\AA',n')\to (\BB',p')\] is given by
                \[
                    (F,G,m)\times (F',G',m')\coloneq (F\times F', G\times G',m\times m'),
                \]
            where we have identified $\VV^n\times \VV^{n'}$ with $\VV^{n+n'}$ for the sake of brevity.
            \item The tensor product of 2-cells is given by the cartesian product. 
            \item The unitors and associators are given by the unitors and associators for categories.
            \item The 2-unitors and 2-associators are identities.
        \end{itemize}
    \end{definition}
    We can then define the monoidal 2-category of $\VV$-iterated categories analogously.
    \begin{definition}
        The 2-category, $\cat{\VV}^\times$, of \textdef{$\VV$-iterated categories} is the 2-category described as follows.
        \begin{itemize}
            \item The objects are given by pairs $(\AA,n)$, where $n$ is a natural number and $\AA$ is a $\VV^{n}$-category.
            \item  The 1-cells $(\AA,n)\to (\BB,p)$ are given by enriched functors
            \[
                (F,G)\colon (\AA,\VV^n)\to (\BB,\VV^p)
            \]
            where $G$ is of the form
            \[
            \VV^n \xrightarrow{\prod_{i=1}^p \tensor^{n_i}} \VV^p 
            \]
            and where $\sum_{i=1}^p n_i=n$.
            \item The 2-cells are given by enriched natural transformations where the monoidal natural transformation is a coherence isomorphism.
        \end{itemize}
        This 2-category has a monoidal structure described as follows. 
        \begin{itemize}
            \item The tensor product of objects $(\AA,n)$ and $(\BB,p)$ is given by
            \[
            (\AA,n)\tensor (\BB,p)\coloneqq (\AA\tensor \BB,n+p)
            \]
            \item The tensor product of 1-cells \[(F,G)\colon (\AA,n)\to (\BB,p)\text{ and }(F',G')\colon (\AA',n')\to (\BB',p')\] is given by
                \[
                    (F,G)\times (F',G')\coloneq (F\times F', G\times G'),
                \]
            where we have identified $\VV^n\times \VV^{n'}$ with $\VV^{n+n'}$ for the sake of brevity.
            \item The tensor product of 2-cells is given by the cartesian product. 
            \item The unitors and associators are given by the unitors and associators for categories.
            \item The 2-unitors and 2-associators are identities.
        \end{itemize}
    \end{definition}
    In order to translate between closed $\VV$-iterated representations and $\VV$-iterated categories we need to use monoidal 2-functors. A definition of monoidal 2-functors can be found in Day and Street's \cite[def.~2]{street1997monoidal} \textbook{Monoidal bicategories and Hopf algebroids}. Given Gray monoids $\AA$ and $\BB$, a monoidal 2-functor $F\colon \AA\to \BB$ comes equipped with a pseudonatural transformation 
    \[
        m\colon F(-)\tensor F(-)\Rightarrow F(-\tensor -)
    \]
    as well as several coherence modifications. We don't give an account of these since the monoidal 2-functors we are interested in adhere to the same coherence laws as monoidal functors. In other words, the coherence modifications are identities.
    \begin{lemma}
        There is a monoidal 2-functor, $\underline{(-)}\colon \Repp{\VV}_{\cll}^\times\to \cat{\VV}^\times,$
        that we also call the \textdef{enriching 2-functor}, which takes every closed representation to its associated enriched category.
        \label{enrichingFuncMon}
    \end{lemma}
    \begin{proof}
        The construction of this 2-functor is completely analogous to the construction of the 2-functor in \autoref{theorem:enrichingfunctor}. Note that the monoidal structure is entirely given by identity functors: the unit in $\Repp{\VV}_{\cll}^\times$ is $*$ with trivial action and closed structure given by the unit. The enriching functor takes this to the one-object category with hom-object given by $I$, which is the unit in $\cat{\VV}$. Note also that, by construction, the following diagram commutes.
        \[
            \begin{tikzcd}[column sep=4em]
                \Repp{\VV}_{\cll}^\times \times \Repp{\VV}_{\cll}^\times
                    &\cat{\VV}^\times\times \cat{\VV}^\times\\
                \Repp{\VV}_{\cll}^\times
                    & \cat{\VV}^\times
                \arrow[from=1-1, to=2-1, "\tensor"{swap}]
                \arrow[from=1-1, to=1-2, "\underline{(-)}\times \underline{(-)}"]
                \arrow[from=2-1, to=2-2, "\underline{(-)}"{swap}]
                \arrow[from=1-2, to=2-2, "\tensor"]
            \end{tikzcd}
        \]
    Thus we can take both coherence maps to be the identity functors. Coherence modifications are given by identity natural transformations.
    \end{proof}
    \begin{definition}
        We define the \textdef{collapsing 2-functor}, $
        C\colon \cat{\VV}^\times \to \cat{\VV}$, 
        to be the 2-functor where:
        \begin{itemize}
            \item  every object $(\AA,n)$ is sent to the $\VV$-category $\AA_{\tensor^n}$;
            \item every 1-cell $(F,G)\colon (\AA,n)\to (\BB,p)$ is sent to the map
                \[
                    \AA_{\tensor^n} \xrightarrow{\underline{c}} (\AA_{G})_{\tensor^p}\xrightarrow{F_{\tensor^p}} \BB_{\tensor^p},
                \]
            where here $\underline{c}$ denotes the functor which is the identity on objects and has hom-morphism given by a coherence map; 
            \item every enriched natural transformation, $(\nu,\mu)$,
            \[
                \begin{tikzcd}[row sep=0.6em, column sep=3.2em]
                (\AA,n)
                \arrow[rr, bend left, "{(F,G,m)}", ""{swap, name=top}]
                \arrow[rr, bend right, "{(F',G',m')}"{swap}, ""{name=bottom}]
                    &\phantom{\CC_{G'}}
                        &(\BB,p)
                        \arrow[Rightarrow, from=top, to=bottom, "{(\mu,\nu)}", to path={
                            (\tikztostart) 
                            -- (\tikztostart|-\tikztotarget.north) \tikztonodes}]         
                \end{tikzcd}
            \]
            is sent to the enriched natural transformation in the following diagram.
            \[
                \begin{tikzcd}[row sep=0.6em, column sep=3.2em]
                \AA_{\tensor^n}
                \arrow[r, "\underline{c}"{swap}]
                \arrow[rr, bend right=75, "\underline{c'}"{swap}, ""{ name=bottom1}]
                    &
                ({\AA}_G)_{\tensor^p}
                \arrow[from=bottom1, phantom, "="]
                \arrow[rr, bend left=75, "\underline{F}_{\tensor^p}", ""{swap, name=top}]
                \arrow[r, "\underline{\mu}_{\tensor^p}"{swap}]
                    &
                    ({\AA}_{G'})_{\tensor^p}
                    \arrow[Leftarrow, to=top, "\nu"{swap}, to path={
                        (\tikztostart) 
                        -- (\tikztostart|-\tikztotarget.south) \tikztonodes}]
                    \arrow[r, "\underline{F'}_{\tensor^p}"{swap}]
                        &{\BB_{\tensor^p}}
                \end{tikzcd}
            \]
        \end{itemize}
    \end{definition}
    \begin{lemma}
        The collapsing 2-functor has a strong monoidal structure.
        \label{collapsingFuncMon}
    \end{lemma}
    \begin{proof}
        To define the strong monoidal structure we need coherence 1-cells. Firstly, this means that we must define an invertible pseudonatural $\VV$-functor 
        \[
            M\colon C(-)\tensor C(-)\to C(-\tensor -).
        \]
        Note that on objects, $C(\AA,n)\tensor C(\BB,p)$ is equal to $C(\AA\tensor \BB,n+p)$. If $\AA$ is a $\VV^{ n}$ category then $\AA(A,A')$ is of the form $(\AA(A,A')_1,....,\AA(A,A')_n)$ where each $\AA(A,A')_i$ is in $\VV$. Now suppose that $\AA$ is a $\VV^{n}$-category and $\BB$ is a $\VV^{p}$-category. Letting 
        \[
            X_i\coloneqq\begin{cases}
                \AA(A,A')_i &\text{if }i\leq n,\\
                \BB(B,B')_i  &\text{otherwise,}
            \end{cases}  
        \]
        we know that
        \[
        (C(\AA,n)\tensor C(\BB,p))((A,A'),(B,B'))= (\tensor_{i=1}^n X_i) \tensor (\tensor_{j=n+1}^p X_i).    
        \]
        On the other hand, we have
        \[
        C((\AA,n)\tensor (\BB,p))((A,A'), (B,B'))=\tensor_{i=1}^{n+m}X_i.
        \]
        But there is a coherence natural isomorphism in $\VV$,
        \[
            \mu\colon (\tensor_{i=1}^n X_i)\tensor (\tensor_{j=n+1}^m X_j)\xrightarrow{\sim} \tensor_{i=1}^{m} X_i,
        \]
        so we define $M$ to be the identity-on-objects functor with underlying functor morphism $\mu$. We know that $M$ is invertible since $\mu$ is. The naturalisors for $M$ are identity 2-cells. This follows from the fact that $\mu$ is natural.

        Note that the collapse of the unit $\VV$-iterated category, $C((*,0))$, is just the one-object $\VV$-category with hom-object give by $I$. This is the unit of $\cat{\VV}$, and we can take the second coherence 1-cell to be the identity. 
    
        The rest of the structure for a strong monoidal functor consists of coherence modifications. Since $M$ is the identity on objects, and $\mu$ satisfies coherence laws by construction, these modifications are all identities.
        \end{proof}
        \tikzstyle{finish}=[to path={(\tikztostart) -- (\tikztotarget.#1)}]
\begin{figure}[p]
\[
    \rotatebox{90}{
    \begin{tikzcd}[ampersand replacement=\&, cramped]
    G'G(B^A)
    \arrow[r, "\eta"]
    \arrow[d, "\eta"]
        \&G'((G(B^A)\act FA)^{FA})
        \arrow[r,"m"]
            \&G'(F(B^A\act A)^{FA})
            \arrow[r, "\epsilon"]
                \&G'(FB^{FA})
                \arrow[d, "\eta"]\\
    (G'G(B^A)\act F'FA)^{F'FA}
    \arrow[d,"m'"]
    \arrow[r,"\eta"]
        \&\bullet
        \arrow[r, "m"]
            \&\bullet
            \arrow[r, "\epsilon"]
                \&(G'(FB^{FA})\act F'FA)^{F'FA}
                \arrow[d, "m'"]\\
    (F'(G(B^A\act FA)))^{F'FA}
    \arrow[dr, equal, finish=180]
    \arrow[r, "\eta"]
        \&\bullet
        \arrow[r, "m"]
        \arrow[d, "\epsilon"]
            \&\bullet
            \arrow[r, "\epsilon"]
                \&(F'(FB^{FA}\act FA))^{F'FA}
                \arrow[d, "\epsilon"]\\
    \blank
        \& (F'(G(B^A\act FA)))^{F'FA}
        \arrow[r, "m"]
            \&(F'F(B^A\act A))^{F'FA}
            \arrow[r, "\epsilon"]
                \&(F'F(B))^{F'FA}
    \end{tikzcd}
    }
\]
\caption{}
\label{ch4.2}
\end{figure}        
        \begin{figure}[p]
    \[
    \rotatebox[]{90}{
    \begin{tikzcd}[ampersand replacement=\&, column sep=1.2em, every matrix/.append style={nodes={font=\scriptsize}}, cramped]
        I
        \arrow[r, "l^{-1}"]\
        \arrow[d,"\eta"]
            \&I\tensor I
            \arrow[r, "\eta"]
                \&I\tensor (I\odot FA)^{FA}
                \arrow[r, "u"]
                    \& I\tensor FA^{FA}
                    \arrow[r, "\nu"]
                        \& I\tensor F'A^{FA}
                        \arrow[r, "\eta"]
                            \&
                            \begin{cdaligned}(I\act &F'A)^{F'A}\\
                                &\tensor F'A^{FA}
                            \end{cdaligned}
                            \arrow[r, "\nu'\cdot u"]
                                    \&F''A^{F'A}\tensor F'A^{FA}
                                    \arrow[d, "\eta"]\\
        (I\odot FA)^{FA}
        \arrow[r, "l^{-1}"]
        \arrow[rrrdd, equal, rounded corners, 
        to path={
            (\tikztostart.270) -- 
            (\tikzcdmatrixname-5-1.center)-- 
            (\tikzcdmatrixname-5-4.center)-- 
            (\tikztotarget.270)},]
            \& 
            \left(\begin{cdaligned}
            I\tensor &I\\
            &\odot FA
            \end{cdaligned}\right)^{FA}
            \arrow[d, "a"]
            \arrow[r,"\eta"]
                \&\bullet
                \arrow[r, "u"]
                    \&\bullet
                    \arrow[r, "\nu"]
                        \&\bullet
                        \arrow[r, "\eta"]
                            \&\bullet
                            \arrow[r, "\nu'\cdot u"]                                \&\left(\begin{cdaligned}
                                        F''A^{F'A}\tensor F'&A^{FA}\\
                                        &\odot FA
                                    \end{cdaligned}\right)^{FA}
                                    \arrow[d, "a"]\\
        \blank
            \&
                I\odot I \odot FA
            \arrow[r,"\eta"]
            \arrow[rd, equal, start=270]
                \&
                \bullet
                \arrow[r, "u"]
                \arrow[d, "\epsilon"]
                    \&\bullet
                    \arrow[r, "\nu"]
                        \&\bullet
                        \arrow[r, "\eta"]
                            \&\bullet
                            \arrow[r, "\nu'\cdot u"]                                \&\left(\begin{cdaligned}
                                        F''A^{F'A}\odot F'&A^{FA}\\
                                        &\odot FA
                                    \end{cdaligned}
                                    \right)^{FA}
                                    \arrow[d, "\epsilon"]\\
        \blank
            \&\blank
                \&(I\odot I\odot FA)^{FA}
                \arrow[r, "u"]
                    \&
                    (I\odot FA)^{FA}
                    \arrow[r, "\nu"]
                        \&(I\odot F'A)^{FA}
                        \arrow[r, "\eta"]
                        \arrow[rd, equal, start=270]
                            \&\bullet
                            \arrow[r, "\nu'\cdot u"]                        \arrow[d, "\epsilon"]
                                    \&(F''A^{F'A}\odot F'A)^{FA}
                                    \arrow[d, "\epsilon"]\\
        \blank
            \&\blank
                \&\blank
                    \&\blank
                        \&\blank
                            \&(I\odot F'A)^{FA}
                            \arrow[r, "\nu'\cdot u"]
                                    \&F''A^{FA}
    \end{tikzcd}
    }
    \]
    \caption{}
    \label{4.14}
    \end{figure}        
        \begin{figure}[p]
\[
\rotatebox{90}{
\begin{tikzcd}[column sep=1.2em, ampersand replacement=\&, cramped]
    I
    \arrow[d, "\eta"]
    \arrow[r, "i"]
        \&G''(I)
        \arrow[r, "\eta"]
            \&G''((I\odot FA)^{FA})
            \arrow[r, "u"]
                \&G''(FA^{FA})
                \arrow[r, "\nu"]
                    \&G''(F'A^{FA})
                    \arrow[d, "\eta"]\\
    (I\odot F''FA)^{F''FA}
    \arrow[r, "i"]
    \arrow[rrrdd, bend right=45, "u"{swap}]
        \& (G''I\odot F''FA)^{F''FA}
        \arrow[d, "m"]
        \arrow[r, "\eta"]
            \&\bullet
            \arrow[r, "u"]
                \&\bullet
                \arrow[r, "\nu"]
                    \&(G''(F'A^{FA})\odot F''FA)^{F''FA}
                    \arrow[d, "m"]\\
\blank
    \&(F''(I\odot FA))^{F''FA}
    \arrow[r, "\eta"]
    \arrow[rd, equal, finish=180]
        \&(F''(I\odot FA)^{FA}\odot FA)^{F''FA}
        \arrow[d, "\epsilon"]
        \arrow[r, "u"]
            \&\bullet
            \arrow[r, "\nu"]
                \&(F''(F'A^{FA}\odot FA))^{F''FA}
                \arrow[d, "\epsilon"]\\
\blank
    \&\blank
        \&(F''(I\odot FA))^{F''FA}
        \arrow[r, "u"]
            \&(F''FA)^{F''FA}
            \arrow[r, "\nu"]
                \&F''F'A^{F''FA}\\
\end{tikzcd}
}
\]
\caption{}
\label{4.15}
\end{figure}

\chapter{The Cotrace and Scalar Enrichment}
In this chapter we show one of the main results of the thesis. If $\BB$ is a left-composition-closed right-monoidal-closed bicategory, then $\BB$ is naturally enriched over the monoidal bicategory of scalar enriched categories, $\cat{\BB(I,I)}$. What this means is that every hom-category in $\BB$ can be enriched over $\BB(I,I)$, that horizontal composition and identities can be given the structure of enriched functors, and that the associator and unitors for composition can be given the structure of enriched natural transformations.

In the first section we give the definition of enriched bicategories and their underlying bicategories. The concept of a bicategory enriched over a monoidal bicategory seems to have been first introduced by Hoffnung~\cite[def.~6]{hoffnungthesis} in his thesis, but has subsequently been reintroduced by Garner and Shulman~\cite[sec.~2]{garner2016enriched}. 

In the second section we define the scalar spread functor and show that this induces an action by $\BB(I,I)$ on all hom-categories. We show that this functor has a right adjoint which we call the cotrace, and that this right adjoint turns our action into a closed action.

In the final section we show that the composition and identity functors can be given the structure of linear functors, and that the composition unitors and associators can be given the structure of linear natural transformations. This, combined with the closed action of $\BB(I,I)$ on the hom-categories is enough to prove our main theorem, after which we see how the scalar enrichment works for each of our motivating examples. We show that this enrichment is given by the cotrace in the same way that the Frobenius inner product is given by the trace. We also show that in the context of this enrichment, the cotrace is the enriched version of the  2-trace given in \autoref{section:bicategoricaltrace}, and we prove that a lot of the structure of bicategories, such as the composition-closed structure, remains in the associated scalar-enriched bicategory.
\section{Enriched Bicategories}
\begin{figure}
    \[\myinput{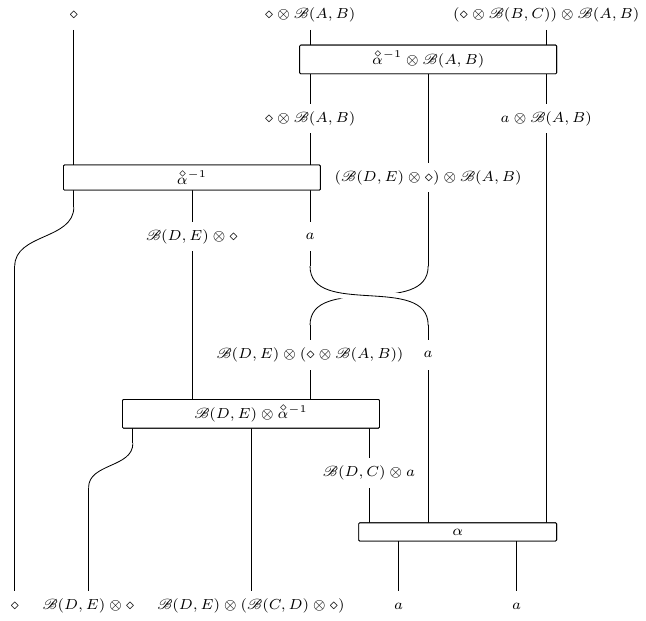} \]
    \caption{}
    \label{eBicatAx1}
\end{figure}

\begin{figure}
    \[\myinput{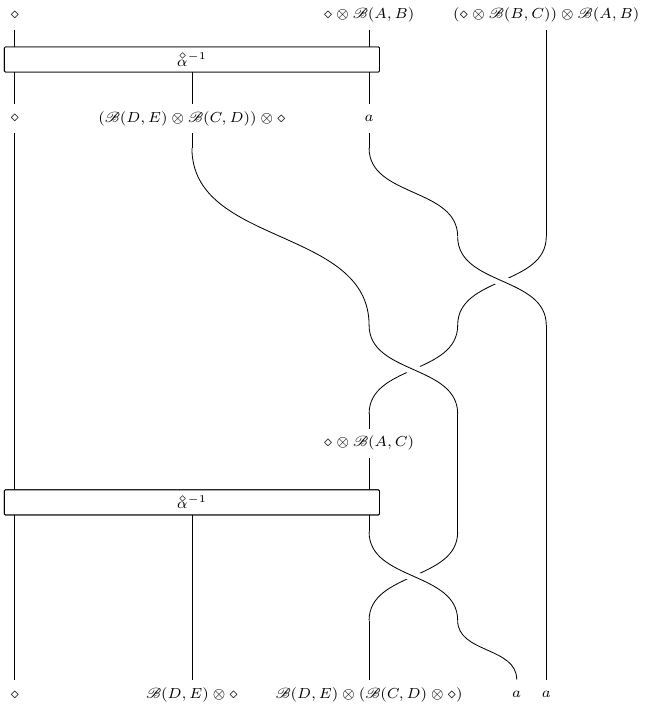}\]
    \caption{}
    \label{eBicatAx2}
\end{figure}

\begin{figure}
    \[\myinput{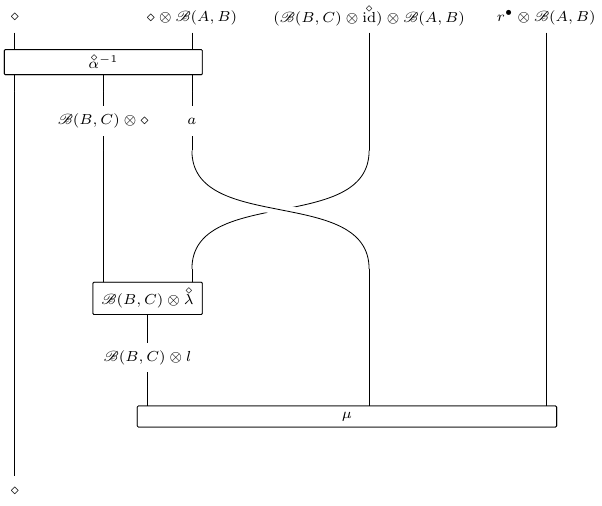}\]
    \caption{}
    \label{eBicatAx3}
\end{figure}

\begin{figure}
    \[\myinput{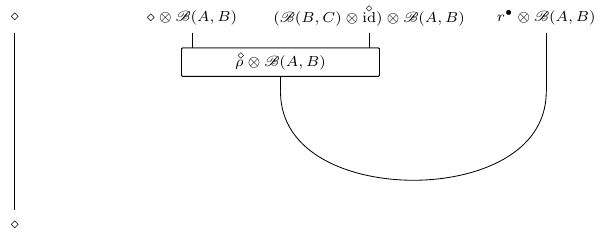}\]
    \caption{}
    \label{eBicatAx4}
\end{figure}

In this section we give a brief account of enriched bicategories, a concept which is core to our main theorem.
\begin{definition}
\label{def:enrichedBicat}
Let $\VV$ be a monoidal bicategory. A bicategory $\BB$ \textdef{enriched in $\VV$} consists of
\begin{itemize}
    \item a collection of objects $\obj(\BB)$ which, by abuse of notation, we denote $\BB$;
    \item for every pair of objects, $A,B\in \BB$, an object, $\BB(A,B)\in \VV$, which we call the hom-object;
\end{itemize}
as well as a horizontal category structure:
\begin{itemize}
    \item for every object, $A\in \BB$, a 1-cell, $\eid\colon I\to \BB(A,A)$, called the enriched identity 1-cell;
    \item for every triple of objects, $A,B,C\in \BB$, a 1-cell, 
    \[
        \ecirc_{A,B,C}\colon \BB(B,C)\tensor \BB(A,B)\to \BB(A,C),
    \]
    called enriched horizontal composition;
\end{itemize}
and finally natural isomorphisms that encode the associativity and unitality of composition:
\begin{itemize}
    \item for every quadruple of objects, $A$, $B$, $C$, $D$, an invertible 2-cell called the enriched composition associator, given below;
    \[
        \myinput{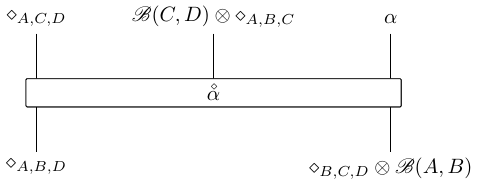}
    \]
    \item for every pair of objects $A$ and $B$, a pair of natural isomorphisms called the left enriched unitor and the right enriched unitor, given below;
    \[
        \myinput{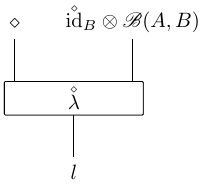}\qquad \myinput{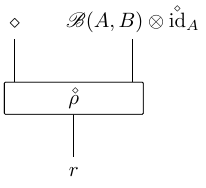} 
    \]
\end{itemize}
such that \autoref{eBicatAx1} is equal to \autoref{eBicatAx2}, and \autoref{eBicatAx3} is equal to \autoref{eBicatAx4}.
\end{definition}

\begin{remark}
    Note that, as mentioned previously, this is distinct from the notion of a category enriched over a bicategory. When enriching over a bicategory, we treat the bicategory as a generalised (or multi-object) monoidal category. By this we mean that if $\AA$ is enriched over a bicategory $\BB$, then $\AA(A,B)$ is a 1-cell in $\BB$ and composition in $\AA$ is given by a 2-cell in $\BB$. This is a one-dimensional construction in the sense that if $\AA$ is enriched over a bicategory then it has an underlying \emph{category}. 
    
    This chapter deals with enriched bicategories, which are enriched over \emph{monoidal} bicategories. By this we mean that if $\AA$ is enriched over a bicategory $\BB$, then $\AA(A,B)$ is an object in $\BB$ and composition is given by a 1-cell in $\BB$. This is a two-dimensional construction in the sense that if $\AA$ is enriched over a monoidal bicategory then it has an underlying \emph{bicategory}.
    \end{remark}
\begin{example}
Just as every category is a category enriched over the category of sets, every bicategory is a bicategory enriched over the 2-category of categories. In this case, the first axiom corresponds to the commutativity of the following pentagon;
\[
\begin{tikzcd}[column sep=-1.5em, row sep=5em]
\blank
    &&
        (k\circ h)\circ (g\circ f)
        \arrow[rrd, "\comp{\alpha}_{k,h,g\circ f}"]\\
((k\circ h)\circ g)\circ f
\arrow[rru, "\comp{\alpha}_{k\circ h, g,f}"]
\arrow[rd, "\comp{\alpha}_{k,g,h}\circ f"{swap}]
    &&
        && k\circ(h\circ(g\circ f))\\
&   (k\circ (h \circ g)) \circ f
\arrow[rr, "\comp{\alpha}_{k,h\circ g,f}"{swap}]
    &&
        k\circ ((h \circ g) \circ f)
        \arrow[ru, "k\circ \comp{\alpha}_{h,g,f}"]
\end{tikzcd}
\]
and the second axiom corresponds to the commutativity of the following triangle.
\[
\begin{tikzcd}
(g\circ \id_B) \circ f
\arrow[rd, "\comp{\rho}\circ \iota_f"{swap}]
\arrow[rr, "\comp{\alpha}_{f,\id_B,g}"]
    && g\circ (\id_B\circ f)
    \arrow[ld, "\iota_g\circ \comp{\lambda}"]\\
&g\circ f
\end{tikzcd}
\]
\end{example}
\begin{definition}
The \textdef{underlying} bicategory $\BB_{0}$ of a $\VV$-bicategory $\BB$ is a bicategory with
\begin{itemize}
    \item the same objects as $\BB$;
    \item hom-categories given by $\VV(I,\BB(A,B))$ for $A,B\in \BB$;
    \item identities given by $\eid\in \VV(I, \BB(A,A))$;
    \item composition given by the following functor;
    \[\begin{tikzcd}
        \VV(I,\BB(B,C))\times \VV(I,\BB(A,B))
        \arrow[d, "\tensor"{swap}]\\
        \VV(I\tensor I, \BB(B,C)\tensor \BB(A,B))
        \arrow[d, "l^\bullet"{swap}]\\
        \VV(I, \BB(B,C)\tensor \BB(A,B))
        \arrow[d, "\ecirc"{swap}]\\
        \VV(I, \BB(A,C))
    \end{tikzcd}\]
\end{itemize}  
and associator and unitor modifications constructed similarly.
\end{definition}
\section{Scalar Actions and the Cotrace}
In 1980 Kelly and Laplaza~\cite[prop.~6.1]{kelly1980coherence} pointed out that for every monoidal category, $\MM$, the commutative monoid of scalars $\MM(I,I)$ acts on every hom-set via the composite
\[
    A\xrightarrow{\rho^{-1}} A\tensor I \xrightarrow{f\tensor s} B\tensor I \xrightarrow{\rho} B.
\]
An analogous construction holds for scalars in monoidal bicategories. We show that, for any $B\in \AA$, there is a monoidal functor called the spread functor \[\Act\colon (\BB(I,I),\circ)\to (\BB(B,B),\circ).\] 
Since $\BB(B,B)$ acts on $\BB(A,B)$ by post-composition, $\BB(I,I)$ must also act on $\BB(A,B)$. Next we show that this functor has a right adjoint, called the cotrace functor. Since we have left-composition-closedness, the action of $\BB(I,I)$ on $\BB(A,B)$ is closed via the composition-closed structure and the cotrace functor. Thus, we can enrich each hom-category $\BB(A,B)$ over $\BB(I,I)$. 
\begin{definition}
    In any monoidal-closed, composition-closed bicategory $\BB$, for any $A\in \BB$ there is a functor, called the \textdef{spread},
    \[
    \Act_A\colon \BB(I,I)\to \BB(A,A)    
    \]
    where $\Act_A(s)$ is given by the following composite
    \[
        A\xrightarrow{l^\bullet} I\tensor A\xrightarrow{s\tensor A} I\tensor A \xrightarrow{l} A.
    \]
\end{definition}
\begin{remark}
    We call this the spread for two reasons. Firstly, in the case of finite dimensional Hilbert spaces, the spread takes in a map $s\colon \mathbb{C}\to \mathbb{C}$ and returns the map $s\cdot(-)\colon A\to A$. If $s$ is thought of as linearly distorting $\mathbb{C}$, then $s\odot (-)$ `spreads' this distortion evenly across every dimension.

    Secondly, in the case of finite dimensional Hilbert spaces, the spread is the map that is linearly adjoint to the trace map, and we are about to see that our spread is adjoint to what we will define as the cotrace functor. It is often convenient to think of an adjoint as a conceptual inverse. In English, a trace is a small, concentrated remnant left behind, indicating the prior presence of something larger. Thus, taking the trace is indicative of taking a small amount, representative of something much larger. To spread is to take a concentrated mass and distribute it over a large area. In this sense, we argue that taking a spread is somehow conceptually inverse to the idea of taking a trace.
\end{remark}
\begin{remark}
    Note, also, the inclusion of the composition circle in the symbol $\dSpr$. This is used to distinguish this functor from the cospread functor that we will define later.
\end{remark}
\begin{proposition}
\label{spreadStrongMonoidal}
The spread functor can be given the structure of a strong monoidal functor.
\end{proposition}
\begin{proof}
The associated multiplication natural isomorphism is given by \autoref{ch5.1} whilst the associated unit natural isomorphism is given by \autoref{ch5.2}. The fact that these adhere to the monoidal axioms follows immediately from the axioms of pseudofunctors and string diagrams.
\begin{figure}
\[
\myinput{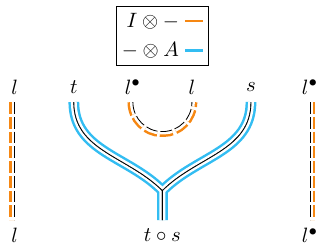}    
\]
\caption{}
\label{ch5.1}
\end{figure}
\begin{figure}
\[
\myinput{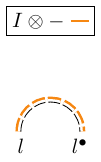}
\]
\caption{}
\label{ch5.2}
\end{figure}
\end{proof}
\begin{corollary}
    \label{homCatAction}
    For every pair of objects $A,B$ in a left-composition-closed right-monoidal-closed bicategory the category $\BB(A,B)$ can be given the structure of a $\BB(I,I)$-representation, where, for every scalar, $s$, and every $f\colon A\to B$, we define the action to be
    \[
        s\odot f\coloneqq \Spr_B(s)\circ f.
    \]
\end{corollary}
\begin{proof}
This follows from the fact that postcomposition gives an action, that oplax monoidal functors give restricted representations, and that the spread functor is strong (and therefore oplax) monoidal.
\end{proof}
\begin{example}
    In the bicategory $\Rel$ the monoidal category of scalars can be thought of as truth values. The spread of $*$, thought of as `true', is the identity and the spread of $\varnothing$, thought of as `false', is the empty relation. Thus, 
    \[
    s\odot R=\begin{cases} R &\text{ if } s=*\\
        \varnothing &\text{ if } s=\varnothing.
    \end{cases}  
    \]
\end{example}
\begin{example}
    In the bicategory $\Bim_R$ the monoidal category of scalars is given by $R$-Mod. The spread of an $R$-module $M$ gives the $A$-$A$-bimodule $M\tensor A$. The action is defined as
    \[
    a'\cdot (m\tensor a)\cdot a'' = m\tensor (a'\cdot a\cdot a'').
    \]
    Then for an $A$-$B$-bimodule $N$, $M\odot N$ is given by $M\tensor N$, where the action is defined as
    \[
        a\cdot (m\tensor n) \cdot b= m\tensor (a\cdot n \cdot b).        
    \]
\end{example}
\begin{example}
    In the bicategory $\DGBim$ the monoidal category of scalars is given by chain complexes in $R$-Mod. The action of chain complexes of $R$-modules is defined similarly to the above but with the derived tensor.
\end{example}
\begin{example}
    In the bicategory $\VV$-Prof the monoidal category of scalars is given by $\VV$. The spread of $v\in \VV$ at some category $\CC$ is the constant profunctor $v \colon \CC\profto \CC$ given by the functor
    \[
    v\colon \CC^\op\tensor \CC\to \VV,    
    \]
    and so for some profunctor $Q\colon \CC\profto \DD$ the action is defined by
    \[
    v\odot Q(C,D)= Q(C,D)\tensor v.
    \]
\end{example}
\begin{example}
    If $\CC$ has finite limits then the spread of $C\in \CC$ is given by the span
    \[
        \begin{tikzcd}
        &C\times A
        \arrow[ld, "p_A"{swap}]
        \arrow[rd, "p_A"]\\
        A
            &&A
        \end{tikzcd}
    \]
    where $p_A$ is the projection map. Then by properties of limits, given a span $S= A\xleftarrow{f} S\xrightarrow{g} B$ we have that the span $C\odot S$ is given by
    \[
    \begin{tikzcd}
        &C\times S\arrow[d, "p_S"{swap}]\\
        &S
        \arrow[ld, "f"{swap}]
        \arrow[rd, "g"]\\
    A
        &&B
    \end{tikzcd}
    \]
\end{example}
\begin{example}
    Given a topological monoid $M$, scalars in $\Path(M)$ are loops at $e$. Given a loop $p$ at $e$, the scalar spread of $p$ at $x$ is the loop $x\cdot p$. This then acts on paths $x$ to $y$ by concatenation.
\end{example}
It will be useful to have a number of different perspectives on this action, and the next two propositions show that two of the other obvious ways to define an action of $\BB(I,I)$ on $\BB(A,A)$ are essentially the same as this one.
\begin{lemma}
The $\BB(I,I)$-representation $\BB(A,B)$, where, for a scalar, $s$, and a 1-cell $f\colon A \to B$, the action is given by
\[s\odot f\coloneqq \Spr_{B}(s)\circ f,\] is isomorphic as a $\BB(I,I)$-representation to the $\BB(I,I)$-representation $\BB(A,B)$ with action given by \[s\odot f\coloneqq f\circ \Spr_{A}{(s)}.\]
\end{lemma}
\begin{proof}
In order to prove this we will show that there are linear functors given by 
\[
(\id,\id, m)\colon \BB(A,B)\rightleftarrows \BB(A,B)\cocolon (\id,\id,m^{-1}), 
\]
and so we just need to prove that there is a natural isomorphism
\[
    m\colon \Spr_B(s)\circ f \Rightarrow f\circ \Spr_A(s)
\]
which adheres to the axioms of linear functors. Using semi-strictification and treating $\BB$ as a Gray monoid, $m$ is defined as the interchangerator,
\[
    \myinput{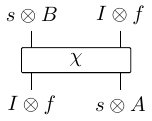}
\]
and the axioms of linear functors hold as a direct consequence of the axioms for the interchangerator. Specifically the associativity diagram commutes because of \autoref{def:GrayAss2} and the unitality diagram commutes because of \autoref{def:GrayId}.
\end{proof}
Recall the name and realisation functors defined in \autoref{prop:nameRealisation}.
\begin{lemma}
    The $\BB(I,I)$-representation $\BB(A,B)$, where, for a scalar, $s$, and a 1-cell $f\colon A \to B$, the action is given by \[s\odot f\coloneqq \Spr_{B}(s)\circ f,\] is equivalent, as a $\BB(I,I)$-representation, to the $\BB(I,I)$-representation $\BB(I, \hom{A,B})$ with action given by \[s\odot g\coloneqq \name{g}\circ s,\] 
    for $g\colon I\to \hom{A,B}$.
    \label{prop:altAction}
\end{lemma}
\begin{proof}
To show the second equivalence holds we need to show that the name and realisation functors can each be equipped with a $\BB(I,I)$-linear functor structure
\[
    (\name{(-)}, m)\colon \BB(A,B) \rightleftarrows \BB(I, \hom{A,B})\cocolon (\unname{(-)}, m')
\]
that give an equivalence between $\BB(A,B)$ with action given by $s\odot f\coloneqq f\circ \Spr_A(s)$ and $\BB(I,\hom{A,B})$. In order to define $m$ note that the 2-cell in \autoref{ch5.3} has $\name{f\circ \Act_A(s)}$ as its source, and $\name{f}\circ s$ as its target.
\begin{figure}
\[
\myinput{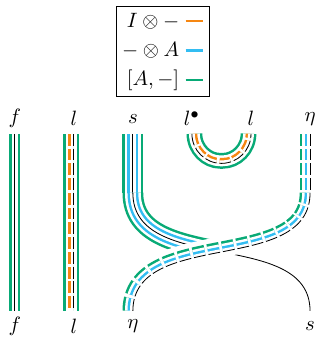}    
\]
\caption{}
\label{ch5.3}    
\end{figure}
Let us once again treat $\BB$ as a Gray monoid. Then, \autoref{ch5.3} can be simplified to give the diagram below.
\[
\myinput{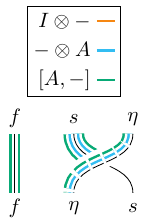}    
\]
The fact that this defines a linear functor follows from the axioms of pseudonatural transformations. To define $m'$ note that, given $g\colon I \to \hom{A,B}$ and $s\colon I \to I$, $\unname{g\circ s}$ is given by the composite
\[
    A=I\tensor A\xrightarrow{s\tensor A} I\tensor A \xrightarrow{g\tensor A} \hom{A,B}\tensor A \xrightarrow{\epsilon} B,
\]
which is equal to $\unname{g}\circ \Spr_A(s)$. Thus, we can let $m'$ be the identity, and $m'$ immediately adheres to the linear functor axioms.

All that remains is to show the natural transformations which give the equivalence are also $\BB(I,I)$-linear. That means we need to show that, for $f\colon A\to B$ and $g\colon I\to \hom{A,B}$, the following two diagrams commute.
\[
    \begin{tikzcd}
        \unname{\name{f}}\circ \Spr_A(s)
        \arrow[r, "\sim"]
        \arrow[d, equal]
            & f\circ \Spr_A(s)\\
        \unname{\name{f}\circ s}
        \arrow[r, "m"{swap}]
            &\unname{\name{f\circ \Spr_A(s)}}
            \arrow[u, "\sim"{sloped,swap}]
    \end{tikzcd}
    \qquad
    \begin{tikzcd}
        \name{\unname{g}}\circ s
        \arrow[r, "\sim"]
        \arrow[d, "m"{swap}]
            & g\circ s\\
        \name{\unname{g}\circ \Spr_A(s)}
        \arrow[r, equal]
        &\name{\unname{g\circ s}}
        \arrow[u, "\sim"{sloped,swap}]
    \end{tikzcd} 
\]
The first diagram commutes because the top morphism is given by \autoref{ch5.4} and the bottom morphism is given by \autoref{ch5.5}, which are equal by the properties of pseudonatural transformations. The second diagram commutes by a similar application of pseudonatural transformation axioms.
\begin{figure}
\[
\myinput{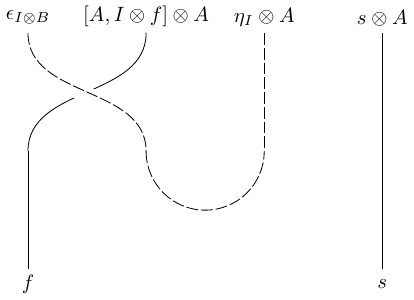}
\]
\caption{}
\label{ch5.4}
\end{figure}
\begin{figure}
    \[
    \myinput{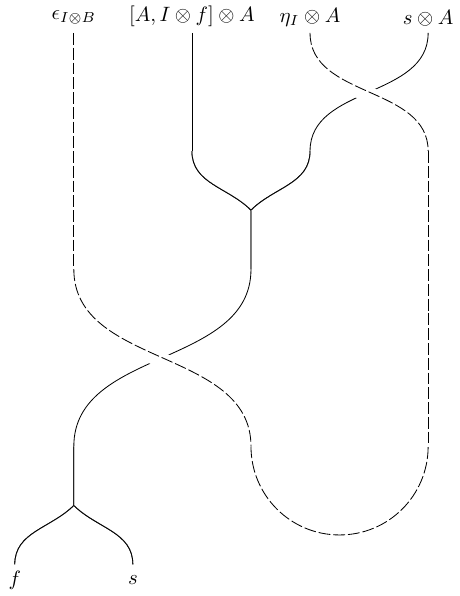}
    \]
    \caption{}
    \label{ch5.5}
    \end{figure}
\end{proof}
One advantage of these other perspectives is that it makes clear the fact that $\BB(A,B)$ has a \emph{closed} action by scalars.
\begin{proposition}
    \label{spreadAdj}
    For every $A\in \BB$ the spread functor has a right adjoint
    \begin{align*}
        \dTr_A\colon \BB(A,A)&\to \BB(I,I),
    \end{align*}
    that we call the \textdef{cotrace functor}, given by the composite   \[
        \BB(A,A)\xrightarrow{\name{(-)}}\BB(I,\hom{A,A})\xrightarrow{\name{\id_A}\lift (-)} \BB(I,I).  
    \]
    \end{proposition}
    \begin{proof}
    By \autoref{prop:altAction} we have natural isomorphisms
    \begin{align*}
        \BB(A,A)(\Act(s), f)& \cong \BB(A,A)(\id\circ \Spr(s), f)\\
        &\cong \BB(A,A)\left(\unname{\name{\id_A}\circ s}, f\right)\\
        &\cong \BB(I,\hom{A,A})\left(\name{\id_A}\circ s, \name{f}\right)\\
        &\cong \BB(I, I)\left(s, \name{\id_A}\lift \name{f}\right)
    \end{align*}
    and so the functor $\name{\id_A}\lift \name{(-)}$ is right adjoint to the spread functor.
    \end{proof}
We call this functor the cotrace functor since it gives the cotrace as defined by Day and Street~\cite[def.~8]{street1997monoidal}. The choice of notation $\dTr$ is to highlight that the cotrace functor is constructed using a lift, and to distinguish this functor from the trace functor that we define later.
\begin{example}
    In the bicategory $\Rel$, given a relation $R\colon A\to A$, its name $\name{R}$ is given by $*\times R\subseteq *\times (A\times A)$. Lifting the name of $R$ through $\name{\id}$ then gives the relation
    \[
        \name{\id}\lift \name{R}=\{(*,*)\mid \forall (a,a')\in A\times A\; a=a'\Rightarrow R(a,a')\}=\begin{cases}
            * &\text{ if } R \text{ is reflexive};\\
            \varnothing &\text{ otherwise} 
        \end{cases}
    \]
\end{example}
\begin{example}
    In the bicategory $\Bim_R$, given an $A$-$A$ bimodule, $M$, $\name{M}$ is the bimodule $M$ thought of as a right $A\tensor A^{\op}$ module and $\name{\id}$ is $A$ thought of a right $A\tensor A^{\op}$ bimodule. The lift of these gives
    \[
        \name{\id}\lift \name{M}=\Hom_{A\tensor A^{\op}}(A, M)
    \]
    which is the $R$-$R$-bimodule, or rather the $R$-module, of homomorphisms from $A$ to $M$.  
    Then, as with the 2-trace, there is an isomorphism between the cotrace of $M$ and the invariants, or centre, of $M$.
\end{example}
\begin{example}
    In the bicategory $\DGBim_R$, when $R$ is a field, the cotrace of a differential graded $A$-$A$-bimodule $M$ is given by the Hochschild cohomology of $A$ with coefficients in $M$. This follows from the fact that the cotrace of $M$ is given by $\RHom(A,M)$ whose $n$'th component is $\Ext_{A\tensor A^{\op}}^n(A,M)$. As shown by Cartan and Eilenberg\cite[ch.~IX]{cartan1956homological} this is equivalent to the definition of the Hochschild cohomology of $A$ with coefficients in $M$.
\end{example}
\begin{example}
    In the bicategory $\VV$-Prof, given a profunctor $P\colon \AA\profto \AA$ the names of $P$ and $\id$ are the profunctors
    \[
        \name{P}\colon *\profto \AA\tensor \AA^{\op}\text{ and } \name{\Hom}\colon *\profto \AA\tensor \AA^{\op}
    \]
    which are isomorphic to $P$ and $\Hom$ as functors. This means that the cotrace is given by
    \[
    \name{\Hom}\lift \name{P}= \endint_{(A,A')\in \AA^{\op}\tensor \AA} \VV(\Hom(A,A'), P(A,A'))= \endint_{A\in \AA}P(A,A)
    \]
    where the second equality follows from the Yoneda lemma for enriched categories.
\end{example}
\begin{example}
    In $\Span(\Set)$ the cotrace is identical to the  2-trace: it is the set of `mutual sections' of the span.
\end{example}
\begin{example}
    In the bicategory $\Path(G)$ for a topological group $G$, the name of a path $p\colon x\to y$ is the path $p\cdot y^{-1}$. Thus, if $p$ is a loop at $x$ then the cotrace of $p$ is just the loop $p\cdot x^{-1}$.
\end{example}
It should be clear from the examples that in each case the underlying set of the cotrace is the 2-trace. The following results explain why this is the case. 
\begin{lemma}
    \label{lemma:closedStructureDef}
    In a left-composition-closed, right-monoidal-closed bicategory $\BB$, every hom-category $\BB(A,B)$ has the structure of a closed $\BB(I,I)$-representation. The action is given by the functor
    \[
        \BB(I,I)\times \BB(A,B) \xrightarrow{\spr(-)\circ -} \BB(A,B)  
    \]
    and the closed structure is given by the functor
\[
\BB(A,B) \times \BB(A,B) \xrightarrow{\dTr((-) \lift (-))} \BB(I,I).
\]
\end{lemma}
\begin{proof}
This follows from the fact that $\BB(B,B)$ acts on $\BB(A,B)$ by composition, the fact that this action is closed via taking right lifts and the following three results: \autoref{spreadStrongMonoidal}, which says that the spread functor is strong monoidal; \autoref{spreadAdj}, which says that the cotrace functor is right adjoint to the spread functor; and \autoref{prop:pullbackAdj}, which says that oplax monoidal functors with right adjoints preserve closed representations.
\end{proof}
\begin{corollary}
\label{lemma:enrichmentDef}
In a left-composition-closed, right-monoidal-closed bicategory $\BB$, every hom-category $\BB(A,B)$ is the underlying category of a $\BB(I,I)$-category, $\underline{\BB}(A,B)$ whose hom-objects are given by
\[
\underline{\BB}(A,B)(f,g):= \dTr(f\lift g).
\]
\end{corollary}
\begin{proof}
This is a direct consequence of the above lemma and \autoref{lemma:fundamentalTheorem}.
\end{proof}
In other words, between any two 1-cells we have a 2-cell-object that lives in the braided monoidal category of scalars. What's more, this 2-cell-object is analogous to the construction of the Frobenius inner product.

Recall that the Frobenius inner product is defined in terms of the trace via
\[
\langle f,g\rangle\coloneqq \Tr(f^\dagger\circ g)
\]
where $f^\dagger$ denotes the linear adjoint to $f$. 
Suppose that $f\colon A\to B$ and $g\colon A\to B$ are 1-cells in $\BB$, and suppose that $f$ has a right adjoint $f^\dagger$. In such a case $f\lift g$ is exactly $f^\dagger \circ  g$ and so the 2-cell-object is defined as
\[\underline{\BB}(A,B)(f,g)=\dTr(f^\dagger\circ g).\]
In particular, we know that $\BBB(A,A)(\id, f)\cong \ctr(\id\lift f)\cong \ctr(f)$. This proves in terms of enriched categories that the 2-trace is always the underlying set of the cotrace.

Whilst the cotrace gives a good conceptual basis for our enrichment, the above lemma doesn't make it particularly clear what these 2-cell-objects look like for a category theorist. The next proposition gives us a different perspective.
\begin{proposition}
    \label{prop:altEnrichment}
Given 1-cells $f\colon A\to B$ and $g\colon A\to B$ in a left-composition-closed, right-monoidal-closed bicategory there is an isomorphism, natural in both $f$ and $g$,
\[
\name{\id_A}\lift \name{f\lift g}\cong \name{f}\lift \name{g}.
\]
\end{proposition}
\begin{proof}
Let $s$ be a scalar. By \autoref{prop:altAction} we know that there is a sequence of natural isomorphisms
\begin{align*}
    \BB(I,I)\left(s,\name{\id}\lift \name{f\lift g}\right)&\cong \BB(I,\hom{A,A})\left(\name{\id}\circ s, \name{f\lift g}\right)\\
    &\cong \BB(A,A)\left(\unname{\name{\id}\circ s}, f\lift g\right)\\
    &\cong \BB(A,A)\left(\id \circ \Spr(s), f\lift g\right)\\
    &\cong \BB(A,A)\left(\Spr(s), f\lift g\right)\\
    &\cong \BB(A,B)\left(f\circ \Spr(s), g\right)\\
    &\cong \BB(I,\hom{A,B})\left(\name{f\circ \Spr(s)},\name{g}\right)\\
    &\cong \BB(I,\hom{A,B})\left(\name{f}\circ s, \name{g}\right)\\
    &\cong \BB(A,B)\left(s, \name{f}\lift \name{g}\right)
\end{align*}
and so by the Yoneda lemma we know that $\name{f}\lift \name{g}\cong \name{\id_A}\lift \name{f\lift g}$. 
\end{proof}
In other words, given a pair of 1-cells, $f,g\colon A\to B$, the scalar that lives between them is the lift in the following diagram.
\[
\begin{tikzcd}[column sep=3.2em]
    \hom{A,B}
        & I
        \arrow[l, "\name{f}"{swap}]
            & I
            \arrow[l, dashed]
            \arrow[ll, "{\name{g}}", ""{swap,name=bottom}, bend left=65]
    \arrow[Leftarrow, to=\tikzcdmatrixname-1-2, from=bottom, verticalup]
\end{tikzcd}    
\] 
\section{Scalar Enrichment}
In the previous section we showed that in a left-composition-closed, right-monoidal-closed bicategory every hom-category $\BB(A,B)$ can be replaced by a $\BB(I,I)$-category $\underline{\BB}(A,B)$. In this section we will show that all of $\BB$ can be replaced by a $\cat{\BB(I,I)}$-bicategory. In other words, all of the other bicategorical data can be replaced by $\BB(I,I)$-functors and $\BB(I,I)$-natural transformations. In order to achieve this we will show that each of these components can be equipped with a linear structure, and by the results of the previous chapter this immediately corresponds to an enriched structure.

We end the section by showing that most of the additional structure on $\BB$, such as the composition-closed structure, is carried over to our enriched bicategory $\underline{\BB}$.
\begin{proposition}
There is a linear functor from the $(\BB(I,I)\times \BB(I,I))$-representation $\BB(B,C)\times \BB(A,C)$ to the $\BB(I,I)$-representation $\BB(A,C)$, for which the associated functor is composition,
\[
    \circ\colon \BB(B,C)\times \BB(A,B) \to \BB(A,C),
\]
for which the associated monoidal functor is composition,
\[
    \circ \colon \BB(I,I)\times \BB(I,I)\to \BB(I,I),  
\]
and for which the associated natural transformation is given by the diagram below for every $s,t\in \BB(I,I)$, every $f\in \BB(A,B)$ and every $g\in \BB(B,C)$.
\[
    \myinput{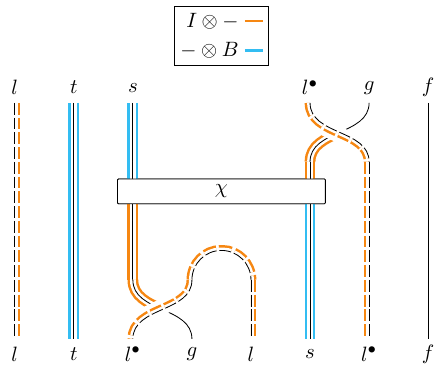}
\]
\end{proposition}
\begin{proof}
We begin by proving that the first diagram of \autoref{def:linearFunctor} commutes. In order to simplify the proof and save space we will prove this result using the semi-strictification theorem for monoidal bicategories. That is, we will treat $\BB$ as a Gray monoid. In this case the natural transformation $m$ reduces to the following 2-cell.
\[
\myinput{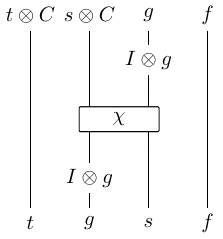}
\]
In the first diagram of \autoref{def:linearFunctor} let $V=(t',s')$, $W=(t,s)$ and $C=(g,f)$. Since $n$ is given by the braid in $\BB(I,I)$, and since in a Gray monoid this is simply the interchangerator, the left-hand path around the diagram is given by \autoref{compLinear:LHS}, and the right-hand path around the diagram is given by \autoref{compLinear:RHS}.
\begin{figure}[H]
    \[\myinput{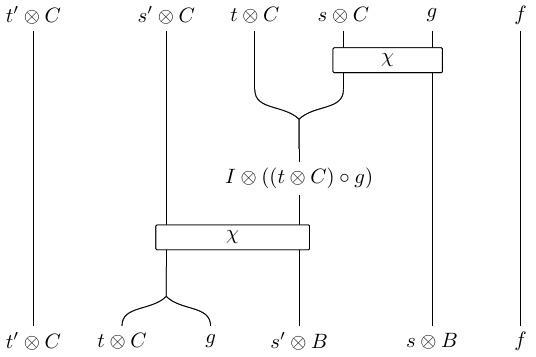}\]
    \caption{}
    \label{compLinear:LHS}
\end{figure}
\begin{figure}[H]
\[\myinput{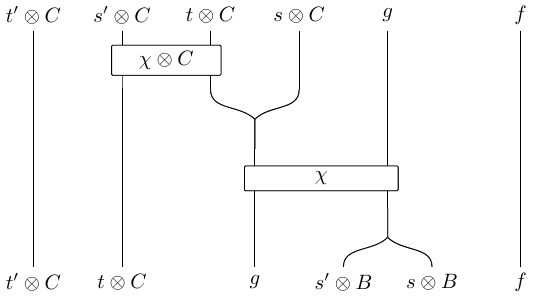}\]  
\caption{}
\label{compLinear:RHS}
\end{figure}
Note however, that by the axioms of a Gray monoid -- \autoref{def:GrayMonoid}, namely (\ref{def:GrayBraid1}), (\ref{def:GrayAss2}) and (\ref{def:GrayBraid2}) -- \autoref{ch5.6}, \autoref{ch5.7}, \autoref{ch5.8}, and \autoref{ch5.9} are equal. Thus, \autoref{compLinear:LHS} is equal to \autoref{compLinear:RHS}, and so the first axiom of linear functors holds.
\begin{figure}
\[\myinput{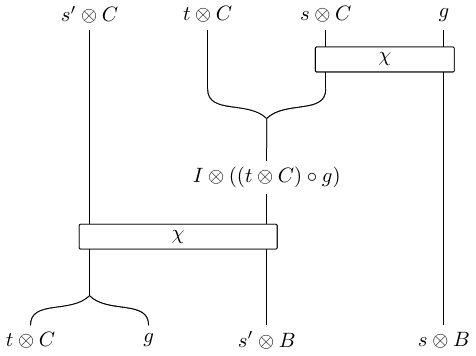}\]
\caption{}
\label{ch5.6}
\end{figure}
\begin{figure}
\[\myinput{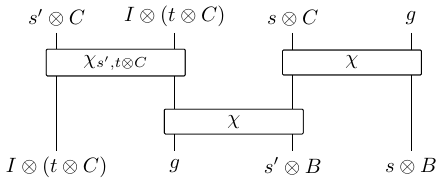}\]
\caption{}
\label{ch5.7}
\end{figure}
\begin{figure}
\[\myinput{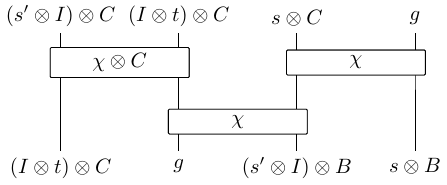}\]
\caption{}
\label{ch5.8}
\end{figure}
\begin{figure}
\[\myinput{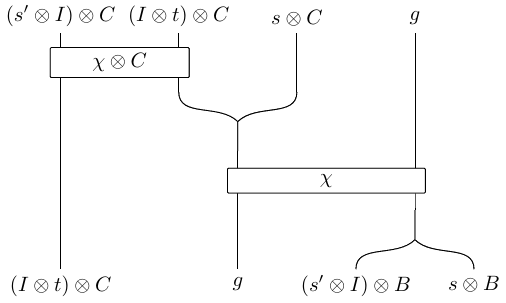}\]
\caption{}
\label{ch5.9}
\end{figure}
To prove that the second diagram commutes note that the right-hand side of the diagram is given by \autoref{ch5.1.1} and the left-hand side of the diagram is given by \autoref{ch5.2.1}.  These two diagrams are equal by properties of adjoint equivalences and pseudonatural transformations. Thus, both axioms of \autoref{def:linearFunctor} hold and so we are done.
\begin{figure}
\[
    \myinput{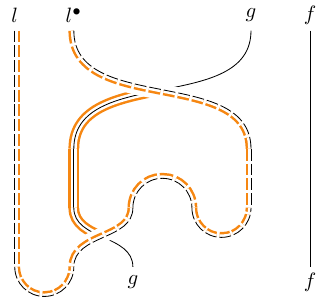}
\]
\caption{}
\label{ch5.1.1}
\end{figure}
\begin{figure}
\[
    \myinput{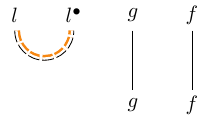}
\]
\caption{}
\label{ch5.2.1}
\end{figure}
\end{proof}
\newpage
\begin{proposition}
    Let $*$ be the one-object, one-morphism category. There is a linear functor from the trivial $*$-representation, $*$, to the $\BB(I,I)$-representation $\BB(A,A)$ for which the associated functor is the functor that picks out the identity,
    \[
        \id_A\colon * \to \BB(A,A)
    \]
    for which the associated monoidal functor is the functor that picks out the identity
    \[
        \id_I \colon *\to \BB(I,I)    
    \]
    and for which the associated natural transformation is given by the diagram below for every $s,t\in \BB(I,I)$, every $f\in \BB(A,B)$ and every $g\in \BB(B,C)$.
    \[\myinput{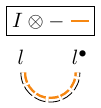}\]
    \end{proposition}
    \begin{proof}
    We firstly show that the first diagram of \autoref{def:linearFunctor} commutes. The path around the right-hand side is given by the following diagram.
        \[\myinput{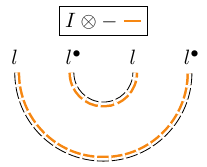}\]
    The path around the left-hand side is given by the following diagram.
        \[\myinput{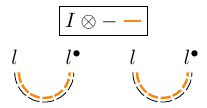}\]
        But note that both of the above 2-cells are equal to the 2-cell below.
        \[\myinput{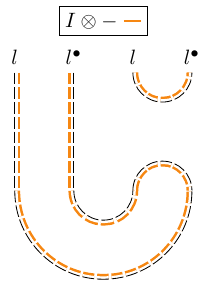}\]
    Now consider the second diagram of \autoref{def:linearFunctor}. Since we are treating composition unitors as identity 2-cells, both the right-hand path and the left-hand path are given by the diagram below.
    \[\myinput{tikz_ScalarActions_IdentityLinear_m}\]
    Thus, the diagram commutes and so $(\id,\id,m)$ is a linear functor.
    \end{proof}

    These linear functors will, of course, be transformed into the $\BB(I,I)$-functors necessary for our bicategorical enrichment. 
    \begin{proposition}
        There is a linear natural transformation whose natural transformation component, and monoidal natural transformation component, are given by $\comp{\lambda}$ in the following two diagrams.
        \[
            \begin{tikzcd}[column sep=3.2em, row sep=3.2em]
            \blank
                & \BB(B,B)\times \BB(A,B)
                \arrow[rd, "\circ", 
                ]\\
            *\times \BB(A,B)
            \arrow[rr, "{\sim}"{swap}, ""{name=bottom}]
            \arrow[ru, "{\id_B\times \BB(A,B)}", 
            ]
                &\blank
                    &\BB(A,B)
                    \arrow[Rightarrow, from=\tikzcdmatrixname-1-2, to=bottom, "\comp{\lambda}"]
            \end{tikzcd}
        \] 
        \[
            \begin{tikzcd}[column sep=3.2em, row sep=3.2em]
            \blank
                & \BB(I,I)\times \BB(I,I)
                \arrow[rd, "\circ"
                ]\\
            *\times \BB(I,I)
            \arrow[rr, "{\sim}"{swap}, ""{name=bottom}]
            \arrow[ru, "{\id_I\times \BB(I,I)}"]
                &\blank
                    &\BB(I,I)
                    \arrow[Rightarrow, from=\tikzcdmatrixname-1-2, to=bottom, "\comp{\lambda}"]
            \end{tikzcd}
        \]
    \end{proposition}
    \begin{proof}
    By appealing to strictification we may treat $\comp{\lambda}$ as the identity 2-cell. Thus, the restriction of $\comp{\lambda}$ to scalars is, trivially, monoidal and the diagram of \autoref{def:linearNT} commutes trivially. 
    \end{proof}

    \begin{proposition}
        There is a linear natural transformation whose natural transformation component, and monoidal natural transformation component, are given by $\comp{\rho}$ in the following two diagrams.
        \[
            \begin{tikzcd}[column sep=3.2em, row sep=3.2em]
            \blank
                & \BB(A,B)\times \BB(A,A)
                \arrow[rd, "\circ"]\\
            \BB(A,B)\times *
            \arrow[rr, "{\sim}"{swap}, ""{name=bottom}]
            \arrow[ru, "{\BB(A,B)\times \id_A}"]
                &\blank
                    &\BB(A,B)
                    \arrow[Rightarrow, from=\tikzcdmatrixname-1-2, to=bottom, "\comp{\rho}"]
            \end{tikzcd}
        \] 
        \[
            \begin{tikzcd}[column sep=3.2em, row sep=3.2em]
            \blank
                & \BB(I,I)\times \BB(I,I)
                \arrow[rd, "\circ"]\\
            \BB(I,I)\times *
            \arrow[rr, "{\sim}"{swap}, ""{name=bottom}]
            \arrow[ru, "{\BB(I,I)\times\id_I}"]
                &\blank
                    &\BB(I,I)
                    \arrow[Rightarrow, from=\tikzcdmatrixname-1-2, to=bottom, "\comp{\rho}"]
            \end{tikzcd}
        \]
    \end{proposition}
    \begin{proof}
    This follows for the same reasons as above.
    \end{proof}

    \begin{proposition}
        There is a linear natural transformation whose natural transformation component, and monoidal natural transformation component, are given by $\comp{\alpha}$ in the following two diagrams.
        \[
            \begin{tikzcd}[column sep=3.0em, row sep=3.2em, cramped]
            \BB(C,D)\times (\BB(B,C)\times \BB(A,B))
                \arrow[rr, "{\BB(C,D)\tensor \circ}", ""{swap, name=top}]
                        && \BB(C,D)\times \BB(A,C)
                        \arrow[d, "\circ"]\\
            (\BB(C,D)\times \BB(B,C))\times \BB(A,B)
            \arrow[u, "\alpha"]
            \arrow[rd, "{\circ \times \BB(A,B)}"{swap}]
                        &
                        \blank &\BB(A,C)\\
            \blank
                &\BB(B,D)\times \BB(A,B)
                \arrow[ru, "\circ"{swap}]
            \arrow[from=top, to=\tikzcdmatrixname-3-2, Rightarrow, "\comp{\alpha}", vertical]
            \end{tikzcd}
        \]
        \[
            \begin{tikzcd}[column sep=3.0em, row sep=3.2em, cramped]
            \BB(I,I)\times (\BB(I,I)\times \BB(I,I))
                \arrow[rr, "{\BB(I,I)\tensor \circ}", ""{swap, name=top}]
                        && \BB(I,I)\times \BB(I,I)
                        \arrow[d, "\circ"]\\
            (\BB(I,I)\times \BB(I,I))\times \BB(I,I)
            \arrow[u, "\alpha"]
            \arrow[rd, "{\circ \times \BB(I,I)}"{swap}]
                        &
                        \blank &\BB(I,I)\\
            \blank
                &\BB(I,I)\times \BB(I,I)
                \arrow[ru, "\circ"{swap}]
            \arrow[from=top, to=\tikzcdmatrixname-3-2, Rightarrow, "\comp{\alpha}", vertical]
            \end{tikzcd}
        \]
    \end{proposition}
    \begin{proof}
        Again, this follows from strictification.
    \end{proof}
    We now have all of the data necessary to enrich over the bicategory of closed representations, and we will use the constructed 2-functors in the last section to change the base of our enrichment
    \scalarEnrichment
    \label{section:cotrace}
    \begin{proof}
        The five above propositions give all of the necessary data to enrich over the bicategory of closed $\BB(I,I)$-iterated representations. To be explicit, we let $\underline{\BB}^\times$ be the $\Repp{\VV}_{\cll}^\times$-bicategory where $\underline{\BB}^{\times}(A,B)$ is given by the $\VV$-representation defined in \autoref{lemma:closedStructureDef}. We define the enriched identity 1-cell to be
        \[
          (\id,\id)\colon (*,*)\to(\BB(A,A),\BB(I,I));  
        \]
        we define the enriched composition 1-cell to be
        \[
        (\circ,\circ) \colon (\BB(B,C)\times \BB(A,B), \BB(I,I)\times \BB(I,I))\to (\BB(A,C), \BB(I,I));  
        \]
        and the enriched unitors and enriched associator are given by the pointwise unitors and associators. Clearly the underlying bicategory of $\underline{\BB}^\times$ is $\BB$, and note that, since the axioms of enriched bicategories hold component-wise in $\Cat$, they hold for $\underline{\BB}^\times$.

        Now using \autoref{enrichingFuncMon} and \autoref{collapsingFuncMon} and Garner and Shulman's~\cite[13.2]{garner2016enriched} base-change trifunctor associated to a monoidal 2-functor, we can replace $\Repp{\VV}_{\cll}^\times$ with $\cat{\VV}$ in our enrichment to get the $\cat{\VV}$-bicategory $\underline{\BB}$. To see that the underlying bicategory of $\underline{\BB}$ is $\BB$, firstly note that, as pointed out in \autoref{theorem:enrichingfunctor}, the enriching functor commutes with the underlying category 2-functor. The collapsing 2-functor does not necessarily commute with the underlying category 2-functor. However, for all $A$ and $B$, $\underline{\BB}^\times(A,B)$ is a $\VV$-representation, and the collapsing functor is the identity on $\VV$-categories.
        \end{proof}
    \begin{example}
        The enriched bicategory $\underline{\Rel}$ has  enrichment given by
        \[
        \underline{\Rel}(A,B)(R,S)=\begin{cases}
            * &\text{ if } R\subseteq S\\
            \varnothing &\text{ otherwise.}
        \end{cases}  
        \]
        In other words, single 2-cells are replaced by truth values.
    \end{example}
    \begin{example}
        The enriched bicategory $\underline{\Bim_R}$ has enrichment given by equipping the set of bimodule homomorphisms with the usual $R$-module structure. 
    \end{example}
    \begin{example}
        The enriched bicategory $\underline{\DGBim}$ has enrichment given by replacing taking the 2-cells between bimodules $M$ and $N$ to be the elements of the chain complex of $R$-modules $\RHom(M,N)$.
    \end{example}
    \begin{example}
        The enriched bicategory \underline{$\VV$-Prof} has enrichment given by replacing the sets of $\VV$-natural transformations between profunctors with $\VV$-natural transformation objects. Let  $P,Q\colon A\profto B$ be profunctors. By the Fubini theorem for ends, the scalar between them is given by
        \[
        \endint_{A}\endint_{B} \VV(P(A,B), Q(A,B))\cong \endint_{(A,B)} \VV(P(A,B), Q(A,B))\cong\Nat(P,Q)
        \]
        which can be found, for example, in Loregian's~\cite[thm.~1.3.1]{Fosco} book. 
    \end{example}
    \begin{example}
        The scalars for the bicategory $\Span(\Set)$ are given by $\Set$ and therefore the enrichment must return the same bicategory. If $\CC$ is some other category that is locally cartesian closed with finite limits, then the enrichment can be recovered by taking the natural strictification, where the objects of $\Span(\CC)$ are thought of as slice categories, and taking the natural transformation objects in $\CC$.
    \end{example}
    \begin{example}
        The enriched bicategory $\underline{\Path(G)}$, for a topological group $G$, has enrichment given by replacing homotopy classes of homotopies with loops. Explicitly, if $p,q\colon x\to y$ are paths in the space, then the scalar between them is the loop that they create, with basepoint shifted to be the unit. Explicitly, the hom-object between $p$ and $q$ is $(\overline{p}\mathbin{;}q)\cdot x^{-1}$, where the semicolon denotes concatenation of paths. 
    \end{example}
    Now that we have the enrichment, we can define the enriched cotrace functor.
    \begin{definition}
        The \textdef{enriched cotrace} functor at an object $A$ in a left-composition-closed right-monoidal-closed bicategory is the functor
        \begin{align*}
            \dTr\coloneqq \underline{\BB}(A,A)(\id_A,-) \colon \underline{\BB}(A,A)&\to \underline{\BB}(I,I) 
        \end{align*}
    \end{definition}
    Recall that the 2-trace of an endomorphism $E\colon A\to A$, as studied by Bartlett~\cite[def.~7.8]{brucethesis}, and Ganter and Kapranov~\cite[def.~3.1]{ganter}, was given by $\BB(A,A)(\id_A,E)$. The cotrace functor above then gives an enriched version of the 2-trace. It is worth noting that the 2-trace as defined by Bartlett was given specifically for the 2-category of 2-Hilbert spaces, and has a Hilbert space structure. We conjecture that the 2-category of 2-Hilbert spaces has a left-composition-closed right-monoidal-closed structure, and that its scalars are given by Hilbert spaces, but we do not prove this.
    
    Recall also that the enriched cotrace functor gives the cotrace in the sense of Day and Street~\cite[def.~8]{street1997monoidal}. In this way we unify the 2-trace.

    Before we move on to compact-closed categories and traces, it is worth noting that the enriched bicategory retains many of the properties of the original bicategory: for example, local completeness and cocompleteness.
    \begin{proposition}
        If $\BB$ is locally complete then so is $\underline{\BB}$. If $\BB$ is locally cocomplete then so is $\underline{\BB}$.
    \end{proposition}
    \begin{proof}
        Note that for every $A,B\in \BB$, $\underline{\BB}(A,B)$ is powered and copowered, with power and copower given by $\Spr(s)\circ (-)$ and $\Spr(s)\lift (-)$ respectively. Since weighted limits are composed of conical limits and powers -- see for example, Kelly's~\cite[thm.~3.73]{kelly1982basic} monograph -- and weighted colimits are composed of conical limits and copowers, if all conical limits exist in $\BB(A,B)$ then all weighted limits exist in $\underline{\BB}(A,B)$ and similarly for weighted colimits.
    \end{proof}
    We don't give a definition here of what it means to be a monoidal-closed \emph{enriched} bicategory. As such, whilst it seems likely that the tensor and closed structures exist for $\underline{\BB}$ we don't give a proof that they do. But it is still useful to know that the tensor gives rise to $\BB(I,I)$-functors between the hom-categories of $\BB$, and that the name and realisation functors have enriched functor analogues.
    \begin{proposition}
        For every pair of objects $A,B\in \BBB$ there is an enriched functor
        \[
            -\tensor B\colon \BBB(A,A)\to \BBB(A\tensor B, A\tensor B)  
        \]
        that we call the enriched right tensor functor, whose underlying functor is the usual tensor functor.
    \end{proposition}
    \begin{proof}
        It suffices to show that there is a strong linear functor \[\BB(A,A)\to\BB(A\tensor B, A\tensor B)\] whose functor component is given by the usual $-\tensor B$. 

        In order to do this we need to give a natural transformation
        \[
            \Spr_{A\tensor B}(s)\circ (f\tensor B)\xrightarrow{m}(\Spr_A(s)\circ f)\tensor B
        \]
        in $\BB(A,A)$ that satisfies the axioms of a linear functor.

        To simplify things let's treat $\BB$ as a Gray monoid. Then we need a natural transformation
        \[
        (s\tensor A\tensor B) \circ (I\tensor f\tensor B) \xrightarrow{m} ((s\tensor A)\circ (I\tensor f))\tensor B    
        \]
        but in a Gray monoid the left and right-hand sides above are equal, so we can choose $m$ to be the identity 1-cell and the axioms hold trivially.
    \end{proof}
    \begin{proposition}
        \label{prop:eName}
        For every pair of objects $A,B\in \BBB$, there is an equivalence defined by a pair of enriched functors
        \[
        \name{(-)}\colon \BBB(A,B)\rightleftarrows \BBB(I,\hom{A,B})\cocolon \unname{(-)}
        \]
        that we call the enriched name and realisation functors, whose underlying functors are the usual name and realisation functors.
    \end{proposition}
    \begin{proof}
        In the proof of \autoref{prop:altAction}, we defined a $\BB(I,I)$-linear structure for the name and realisation functors and showed they were equivalent as linear functors. Thus, by the enriching 2-functor, they are the underlying functors of $\BB(I,I)$-functors.  Since 2-functors preserve equivalences, the enriched name and realisation functors give an equivalence.
    \end{proof}
    \begin{proposition}
        \label{prop:ePostComp}
        Let $f\colon B\to C$ be a 1-cell in $\BB$. Then for every $A\in \BB$ there is an enriched postcomposition functor
        \[
            f\ecirc (-)\colon \BB(A,B)\to \BB(A,C)
        \]
        whose underlying functor is $f\circ (-)$,
        and for every $D\in \BB$ there is an enriched precomposition functor
        \[
            (-)\ecirc f\colon \BB(C,D)\to \BB(B,D)  
        \]
        whose underlying functor is $(-)\circ f$.
    \end{proposition}
    \begin{proof}
        Since $\BB$ is the underlying category of $\BBB$ we know that there is a $\BB(I,I)$-functor
        \[
            f\colon *\to \BBB(B,C)
        \]
        \[
            f\colon *\to \BB(B,C)    
        \]
        that picks out the 1-cell $f$. Then we construct $(-)\ecirc f$ as the composite
        \[
        \BBB(A,B)\xrightarrow{\sim}*\tensor \BBB(A,B)  \xrightarrow{f\tensor \BBB(A,B)}\BBB(B,C) \tensor \BBB(A,B)\xrightarrow{\ecirc} \BBB(A,C)
        \]
        whose underlying functor is $f\circ (-)$ since the underlying functor of enriched composition is usual composition. A similar construction gives the enriched precomposition functor.
    \end{proof}
    The enriched bicategory $\BBB$ is left-composition-closed, like $\BB$, in the sense that the postcomposition functor defined above has a right adjoint.
    \begin{lemma}
    Let $\BB$ be a left-composition-closed, right-monoidal-closed bicategory. The $\cat{\BB(I,I)}$-bicategory $\underline{\BB}$ is left-composition-closed. That is to say, for every 1-cell $f\colon A\to B$ there is a $\BB(I,I)$-adjunction
    \[
        f\ecirc (-) \colon \BBB(A,B) \rightleftarrows \BBB(A,C)\cocolon f\elift (-)
    \]
    where $\circ$ denotes the enriched composition. Furthermore, $f\lift (-)$ is the underlying functor of the $\BB(I,I)$-functor $f\elift(-)$.
    \end{lemma}
    \begin{proof}
    We could give an indirect proof by constructing a linear functor from the usual lift functor. However, in this case a direct approach is just as simple. Firstly, we define a candidate for an enriched functor $f\elift(-)$ and show that it is, in fact, an enriched functor. We then show that this enriched functor is adjoint to enriched post-composition functor. Finally, by uniqueness of adjoints we have that the underlying functor must be the usual lift functor.
    
    We begin by defining a functor $f\elift (-)$. On objects, we define $f\elift(-)$ to be the same as $f\lift (-)$. To define the functor morphism we need to show that for any $g\colon A\to C$ and $h\colon A \to C$ there is a morphism in $\cat{\BB(I,I)}$,
    \[
        \dTr(g\lift h) \to \dTr((f\lift g)\lift (f\lift h)).
    \]
    Since $\dTr$ is a functor, it suffices to know that there is a morphism
    \[
        (g\lift h)\xrightarrow{\zeta} (f\lift g)\lift (f\lift h).
    \]
    We define $\zeta$ is given as the adjunct to the canonical morphism
    \[
    (f\lift g)\circ (g\lift h)\to (f\lift h),   
    \]
    which is itself the adjunct of the map
    \[
        f\circ (f\lift g)\circ (g\lift h)\xrightarrow{\epsilon\circ (g\lift h)} g\circ (g\lift h) \xrightarrow{\epsilon} h.
    \]
    Here $\epsilon$ is the counit for the composition-lift adjunction.
    Thus, $\zeta$ is a candidate for the functor morphism. We now need to show that $\zeta$ adheres to the axioms of enriched functors.

    We firstly need to show that the functor preserves the enriched vertical composition. This means we need to show that the following diagram commutes.
    \begin{figure}[H]
    \[        
    \begin{tikzcd}[ampersand replacement=\&, column sep=6.5em, row sep=3.2em, every matrix/.append style={nodes={font=\scriptsize}}]
        \begin{cdaligned}
            \BBB(A&,C)(h,k)\\
            &\circ \BBB(A,C)(g,h)
        \end{cdaligned}
            \arrow[r, "{f\elift(-)\circ f\elift(-)}"]
            \arrow[d, "\bullet"{swap}]
                \& \begin{cdaligned}
                    \BBB(A&,B)(f\elift h, f\elift k)\\ 
                    &\circ \BBB(A,B)(f\elift g, f\elift h)
                \end{cdaligned}
                \arrow[d, "\bullet"]\\
            \BBB(A,C)(g,k)
            \arrow[r, "{f\elift(-)}", swap]
                \&\BBB(A,B)(f\elift g, f\elift k)
        \end{tikzcd}  
        \]
        \caption{}
        \label{diag:eliftComp}
    \end{figure}
    Firstly note how the enriched vertical composition is defined in the $\BB(I,I)$-category $\BBB(A,C)$. By construction, given $f,g,h\colon A\to C$ we have that the functor morphism for vertical composition
    \[
        \underline{\BB}(A,C)\tensor \BB(A,C)\to \BB(A,C)
    \]
    is given by a morphism
    \[
    \dTr(g\lift h)\circ \dTr(f\lift g)\to \dTr(f\lift h)    
    \]
    which is the composite of three canonical morphisms in $\BB(I,I)$. Firstly we have the braiding in $\BB(I,I)$ given by
    \[
        \dTr(g\lift h)\circ \dTr(f\lift g) \xrightarrow{\altchi} \dTr(f\lift g)\circ \dTr(g\lift h).
    \]
    Secondly, note that the spread functor is strong monoidal and so the enriched cotrace functor, being right adjoint to the spread functor, has a canonical lax monoidal structure. This means there is a canonical natural transformation
    \[
       \dTr(f\lift g)\circ \dTr(g\lift h) \xrightarrow{\phi} \dTr((f\lift g)\circ (g\lift h)).
    \]
    The third and final morphism is the canonical morphism
    \[
        \dTr((f\lift g)\circ (g\lift h)) \xrightarrow{\dTr(\xi)} \dTr(f\lift h),
    \]
    and note that $\xi$ is the adjunct of $\zeta$.

    Thus, to show that \autoref{diag:eliftComp} commutes we need to show that the following diagram commutes.
    \[
        \begin{tikzcd}[ampersand replacement=\&, column sep=6.5em, row sep=3.2em, every matrix/.append style={nodes={font=\scriptsize}}]
            \dTr(h\lift k)\circ \dTr(g\lift h)
            \arrow[d, "\altchi", swap]
            \arrow[r, "f\elift(-)\circ f\elift(-)"]
                \&
                \begin{cdaligned}
                    \dTr((f\lift &h)\lift(f\lift k))\\
                        &\circ 
                        \dTr((f\lift g)\lift (f\lift h))
                \end{cdaligned}
                \arrow[d, "\altchi"]\\
            \dTr(g\lift h)\circ \dTr(h\lift k)
            \arrow[d, "\phi", swap]
                \&
                \begin{cdaligned}
                    \dTr((f\lift &g)\lift (f\lift h)) \\
                    &\circ \dTr((f\lift h)\lift(f\lift k))
                \end{cdaligned}
                \arrow[d, "\phi"]\\
            \dTr((g\lift h)\circ (h\lift k))
            \arrow[d, "\xi", swap]
                \&
                \begin{cdaligned}
                    \dTr((f\lift &g)\lift (f\lift h)\\
                    &\circ (f\lift h)\lift(f\lift k))
                \end{cdaligned}
                \arrow[d, "\xi"]\\
            \dTr(g\lift k)
            \arrow[r, "f\elift(-)"{swap}]
                \&\dTr((f\lift g)\lift (f\lift k))
            \end{tikzcd}
    \]
    However, due to the naturality of $\altchi$ and $\phi$, and the fact that $\dTr$ is a functor, it suffices to show that the following diagram commutes.
    \[
        \begin{tikzcd}
            (g\lift h)\circ (h\lift k)
            \arrow[r, "\zeta"]
            \arrow[d, "\xi"{swap}]
                &((f\lift g)\lift (f\lift h))\circ ((f\lift h)\lift (f\lift k)) 
                \arrow[d, "\xi"]\\
            (g\lift k)
            \arrow[r, "\zeta", swap]
                &(f\lift g)\lift(f\lift k)
        \end{tikzcd}
    \]
    But, since $\xi$ is the adjunct of $\zeta$, this diagram is the adjunct to the following diagram
    \[
        \begin{tikzcd}
        (f\lift g)\circ(g\lift h)\circ (h\lift k)\arrow[r, "\xi\circ (h\lift k)"]
        \arrow[d, "(f\lift g)\circ \xi", swap]
            &(f\lift h)\circ (h\lift k)
            \arrow[d, "\xi"]\\
        (f\lift g)\circ (g\lift k)
        \arrow[r, "\xi", swap]
            &f\lift k
        \end{tikzcd}
    \]
    which commutes because of the naturality of $\xi$. Thus, \autoref{diag:eliftComp} commutes.

    Next we need to show that the morphism preserves identities. In other words we need to show that the following diagram commutes.
    \begin{figure}[H]
    \[
        \begin{tikzcd}
        *
        \arrow[r, "\iota"]
        \arrow[rd, "\iota", swap]
            &\BBB(A,C)(g,g)
            \arrow[d, "\zeta"]\\
        &\BBB(A,B)(f\lift g, f\lift g)
        \end{tikzcd}    
    \]
    \caption{}
    \label{diag:eliftId}
    \end{figure}
        Firstly note that the enriched identity morphism, by construction, is given by the following composite
    \[
        \id\xrightarrow{\psi} \dTr(\id) \xrightarrow{\dTr(\eta)} g\lift g,
    \]
        where $\psi$ is the unit map that gives $\dTr$ a lax monoidal structure, and $\eta$ is the unit of the composition-lift adjunction. Then we need to show that the diagram below commutes.
        \[
        \begin{tikzcd}
                \id
                \arrow[r, "\psi"]
                    & \dTr(\id)
                    \arrow[r, "\dTr(\eta)"]                
                    \arrow[rd, "\dTr(\eta)", swap]
                        &\dTr(g\lift g)
                        \arrow[d, "\dTr(\zeta)"]\\
                        &&\dTr((f\lift g)\lift (f\lift g))
            \end{tikzcd}
        \] 
        Thus, it suffices to show that the following diagram commutes.
        \[
        \begin{tikzcd}
            \id
            \arrow[r, "\eta"]
            \arrow[rd, "\eta"{swap}]
                &g\lift g
                \arrow[d, "\zeta"]\\
            &
                (f\lift g)\lift (f \lift g)
        \end{tikzcd}    
        \]
        But this diagram commutes because the top right path is, by definition equal to the top right path in the following diagram.
        \[
        \begin{tikzcd}
            \id
            \arrow[d, "\eta"]
            \arrow[r, "\eta"]
                & g\lift g
                \arrow[d, "\eta"]\\
            (f\lift g)\lift (f\lift g)
            \arrow[r, "\eta"]
            \arrow[rd, equal]
                &f\lift g\lift ((f\lift g)\circ (g\lift g))
                \arrow[d, "\zeta"]\\
            \blank
                & (f\lift g)\lift (f\lift g)
        \end{tikzcd}    
        \]
        This diagram commutes by naturality and the zigzag identities for the composition-lift adjunction. Thus, we know that \autoref{diag:eliftId} commutes, and so our enriched lift functor $f\elift (-)$ is, in fact, a functor.
        
        Now all that remains is to prove that the enriched lift functor $f\elift (-)$ is adjoint to the enriched post-composition functor $f\ecirc (-)$. This means giving a natural isomorphism
        \[
            \dTr((g\circ f)\lift h)\cong \dTr(f\lift (g\lift h))  
        \]
        and so it suffices to show that there is a natural isomorphism
        \[
            (g\circ f)\lift h\cong f\lift (g\lift h).
        \]
        This follows from the Yoneda lemma and the fact that for any $k\colon A\to A$
        \begin{align*}
            \BB(A,A)(k, (g\circ f)\lift h)&\cong \BB(A,C)(g\circ f\circ k, h)
            \\ &\cong \BB(A,B)(f\circ k, g\lift h)
            \\ &\cong \BB(A,A)(k, f\lift (g\lift h)).
        \end{align*}
        Thus, $f\ecirc(-)$ is left adjoint to $f\elift(-)$ in the category of $\BB(I,I)$-categories. Finally, since the underlying functor of $f\ecirc(-)$ is the usual postcomposition functor $f\circ(-)$, the forgetful 2-functor
        \[
            (-)_0\colon \cat{\BB(I,I)}\to \Cat
        \]
        preserves adjoints, and adjoints are unique up to isomorphism, we know that the underlying functor of $f\elift(-)$ is the usual lift functor $f\lift (-)$.
    \end{proof}
    \begin{proposition}
        \label{prop:epostcompEquiv}
        If $f\colon B \to C$ has a pseudoinverse $f^\bullet$ then the enriched postcomposition functor $f\ecirc (-)$ has a pseudoinverse $f^\bullet \ecirc (-)$, and the enriched precomposition functor $(-)\ecirc f$ has a pseudoinverse $(-)\ecirc f^\bullet$.
    \end{proposition}
    \begin{proof}
        Let $f^\bullet$ be an adjoint pseudoinverse for $f$. We show that the functor $f^\bullet \ecirc (f\ecirc (-))$ is isomorphic to the identity functor. This is sufficient since $f^\bullet$ being adjoint pseudoinverse to $f$ implies that $f$ is adjoint pseudoinverse to $f^\bullet$.
        
        Firstly note that since $\BB$ is the underlying bicategory of $\BBB$ we know that $
            *\xrightarrow{f^\bullet \circ f} \BBB(B,B)
        $
        is isomorphic to the functor $\eid$ which picks out the identity. Since the underlying functor of the enriched composition functor $(-)\ecirc (-)$ is the usual composition functor $(-)\circ(-)$ there is an invertible 2-cell that gives an isomorphism in the diagram below.
        \[
            \begin{tikzcd}
                *
                \arrow[r, "\sim", ""{swap, name=top}]
                \arrow[d, "(f^\bullet \circ f)"{swap}]
                    &*\tensor *
                    \arrow[d, "f^\bullet \tensor f"]\\
        \BBB(B,B)
            &\BBB(C,B)\tensor \BBB(B,C)
            \arrow[l, "\ecirc", ""{swap, name=bottom}]
        \arrow[phantom, "\cong", from=top, to=bottom]
            \end{tikzcd}
        \]
        Since $\BBB$ is an enriched bicategory, we also have an isomorphism of functors
        \[
        \ecomp{\alpha}^{-1}\colon (-)\ecirc ((-)\ecirc (-) )\xRightarrow{\sim}     ((-)\ecirc (-))\ecirc (-).
        \]
        Composing these isomorphisms we have the following isomorphism
        \[
            f^\bullet\ecirc (f\ecirc (-))\cong (f^\bullet\ecirc f)\ecirc (-) \cong (f^\bullet \circ f)\ecirc (-) \cong \id\ecirc (-)\cong \id.
        \]
        Hence, $f\ecirc (-)$ and $f^\bullet \ecirc(-)$ give an equivalence of $\BB(I,I)$-categories. The case for precomposition is analogous.
    \end{proof}

\chapter{Compact-Closed Bicategories and the Two Traces}
Up until now our discussions have focused on the concept of a cotrace. The prefix `co-' here indicates that the cotrace ought to be, in some sense, opposite to some other construction that we can call the trace. This intuition is, of course, backed up by our examples. In different contexts the cotrace gives the centre, the homology and the end, and all of these have dual constructions -- namely the abelianisation, the cohomology or the coend.

It turns out that this dual construction follows from the usual trace in category theory, first defined by Dold and Puppe~\cite[def.~4.1]{DoldPuppe}, that can be defined in any compact-closed category. A compact-closed category is a symmetric monoidal category with duals. That is to say, $\AA$ is compact-closed if every object $A\in \AA$ comes equipped with an object $A^*$ and morphisms
\[
\ev\colon A^*\tensor A\to I\text{ and } \coev\colon I\to A\tensor A^* 
\]
that we draw as cups and caps
such that the usual yanking identities for adjunctions hold.
\[
\fbox{\tikzsetnextfilename{Diagrams/DualInvariance/Cup1}
\begin{tikzpicture}[baseline={(current bounding box.center)}, xscale=1
]

\drawStartNodes{
Ad/$A^*$/1,
A/$A$/0
}

\strCupX{Ad}{A}

\end{tikzpicture}
} \qquad \fbox{\tikzsetnextfilename{DualInvariance/Cap1}
\begin{tikzpicture}[baseline={(current bounding box.center)}, xscale=1
]

\drawStartNodes{
c//0
}
\place{A}{-0.5}{c}
\place{Ad}{0.5}{c}
\strCapX{A}{Ad}

\drawEndNodes{
    Ad/$A^*$,
    A/$A$
}

\end{tikzpicture}
}    
\]
In other words, $\AA$ is a one-object bicategory, equipped with a symmetry, such that every 1-cell comes equipped with a choice of right adjoint.

The motivating example of a compact-closed category is the category of finite dimensional vector spaces. It turns out the trace of any endomorphism in this category can be defined purely in terms of  the compact-closed structure. Given any linear endomorphism $E$, the trace is simply given by the following diagram:
\[
    \myinput{tikz_Introduction_rountrace}
\]
where the braiding in this diagram represents the symmetry. In this chapter we will see how the analogous definition for composition-closed compact-closed bicategories gives rise to a functor called the trace functor. This functor is, in some sense, dual to the cotrace functor.

The first section begins with a review of two different types of commutativity on a monoidal bicategory: braided and symmetric. We give two semi-strictification theorems by Gurski, and prove, using our string diagram language, that the braid for scalars is symmetric whenever the monoidal bicategory is braided. 

The second section focuses on compact-closed bicategories. It is in this section that we define the trace functor and its right adjoint, the cospread functor, and give some examples of the trace functor in certain bicategories. We also show that these functors have scalar-enriched counterparts that form an enriched adjunction. 

The third section is dedicated to a comparison between dagger compact categories -- a category structure that captures the structure of finite dimensional Hilbert spaces -- and composition-closed compact-closed bicategories. We give a heuristic account of why the trace for composition-closed compact-closed bicategories is `split' into the trace and cotrace.

The fourth section gives an account of the properties that the trace and cotrace share with the linear trace. This includes linearity, dual invariance, cyclicity, and preservation of tensors. This also gives us some insight into how the trace and cotrace might be formally dual.

The final section looks at the properties of the codimension and dimension of a given object in the bicategory. The dimension of a finite dimensional Hilbert space $A$ can be given as the trace of the identity morphism at $A$. Defining dimension and codimension analogously, it turns out that codimension and dimension form a monoid-module pair in the monoidal category of scalars. We give an account of the structures that this gives rise to for each of our motivating examples.
\section{Braided and Symmetric Monoidal Bicategories}
\begin{figure}[H]
            \centering
            \begin{minipage}{.5\textwidth}
                \[\myinput{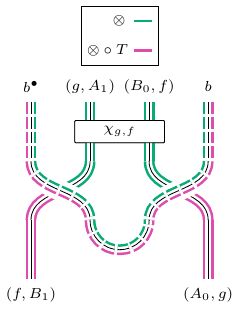}\]
                \caption{}
                \label{scBraid1}
            \end{minipage}%
            \begin{minipage}{.5\textwidth}
                \[\myinput{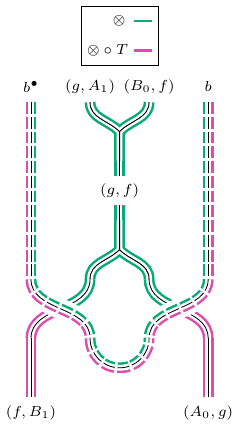}\]
            \caption{}
            \label{scBraid2}
            \end{minipage}
            \end{figure}

\begin{figure}[p]
                \centering
                \begin{minipage}{.5\textwidth}
                    \[\myinput{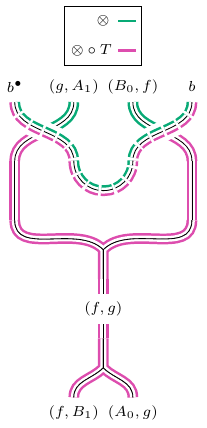}\]                    \caption{}
                    \label{scBraid3}
                \end{minipage}%
                \begin{minipage}{.5\textwidth}
                    \[\myinput{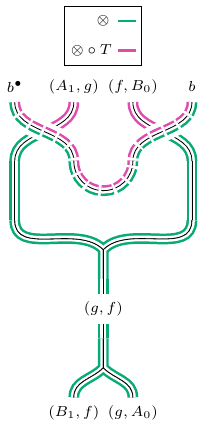}\]       
                    \caption{}
                    \label{scBraid4}
                \end{minipage}
            \end{figure}

\begin{figure}[p]
             \[
              \myinput{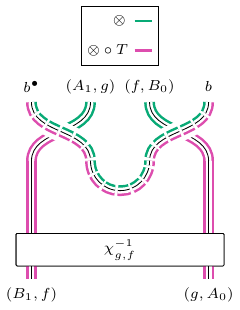}  
            \]
            \caption{}
            \label{scBraid5}
            \end{figure}

In this section we give an account of braided and symmetric monoidal bicategories and the semi-strictification theorems of Gurski. In \autoref{section:Scalars} we saw that the category of scalars is always braided monoidal. In this section we will see that if a monoidal bicategory is braided, then its monoidal category of scalars is symmetric.
\begin{definition}
    A \textdef{braided monoidal bicategory} is a monoidal bicategory $\BB$ equipped with a pseudonatural adjoint equivalence, called the braid, \[b_{A,B}\colon A\tensor B\to B\tensor A\]
    as well as two invertible modifications, called the left and right 2-braids, whose components are given by the diagrams below.
    \[
        \myinput{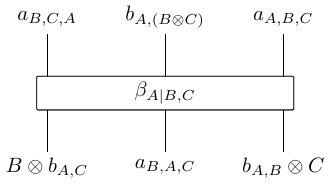}
    \]
    \[
        \myinput{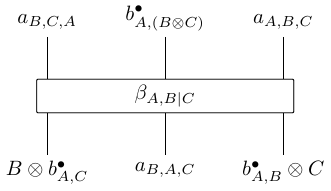}.  
    \]
    These are subject to coherence axioms that can be found, for example, in Stay's~\cite[def.~4.5]{stay2016} paper on compact-closed bicategories. We choose to omit these since we will not appeal to them directly.
    \end{definition}
    \begin{remark}
    Note that if we think of the underlying bicategory as a Gray monoid, the above invertible modifications reduce to the isomorphisms given below.
    \[
        \myinput{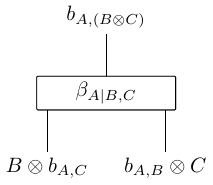}\qquad \myinput{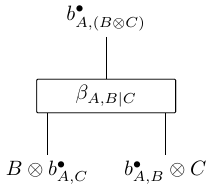}
    \]
    \end{remark}
    In a braided monoidal category $\BB$ we can show that for any $A\in \BB$ the following diagram commutes
    \[
        \begin{tikzcd}
        A\tensor I
        \arrow[rd, "r"{swap}]
        \arrow["b", r]
            & I\tensor A
            \arrow[d, "l"]\\
            & A
    \end{tikzcd}
    \]
    which corresponds to the identity axiom for pseudonatural transformations. An analogous result for braided monoidal bicategories holds, but proving it is far from straightforward.
    \begin{proposition}
        In a braided monoidal bicategory $\BB$ there is an invertible modification that relates the unitors and the braid, given by the diagram below.
            \[
                \fbox{\tikzsetnextfilename{BraidMonBiCat/SymScalars/rWrapAround/2}
\begin{tikzpicture}[baseline={(current bounding box.center)}, xscale=1
]

\drawStartNodes{
    rI/$r$/0
}

\strIdX{rI}[//c1]

\strModXX{rI///c1}{$\gamma$}{l///c3;b//c3/c1}[1]

\strIdX{l}[//c3]
\strIdX{b}[/c3/c1]

\drawEndNodes{
    l/$l$,
    b/$b_{-,I}$
}

\strKeyXX{
    $-\tensor I$/c1;
    $I\tensor -$/c3
}
\end{tikzpicture}
}
            \]
    \end{proposition}
    \begin{proof}
        Firstly notice that there is a modification $\gamma_0$ given by \autoref{gammaFig}.This can be used to construct a modification $\gamma_1$, given by \autoref{gammaFig1}. Finally, we can wrap $\gamma_1$ in naturalisors, as in \autoref{gammaFig2}, to construct the desired 2-cell.
        \begin{figure}[p]
            \[
                \fbox{\tikzsetnextfilename{BraidMonBiCat/SymScalars/rWrapAround/3}
\begin{tikzpicture}[baseline={(current bounding box.center)}, xscale=1
]

\drawStartNodes{
    b/$b_{A,-}$/2,
    l/$l$/0
}

\strBraidXX{b}{l}[/c7/c5/white][c5//c1][2]

\drawEndNodes{
    b/$b_{A,-}$,
    l/$l$
}
\strKeyXX{
    $-\tensor I$/c1;
    $-\tensor (I\tensor I)$/c2;
    $A\tensor -$/c5;
    $-\tensor A$/c7
}
\end{tikzpicture}}
            \]
        \caption{$\gamma_0$}
        \label{gammaFig}
        \end{figure}
        \begin{figure}[p]
            \[\fbox{\tikzsetnextfilename{BraidMonBiCat/SymScalars/rWrapAround/1}
\begin{tikzpicture}[baseline={(current bounding box.center)}, xscale=1
]

\drawStartNodes{
    rI/$r$/0
}
\place{r}{-1.5}[-2]{rI}
\place{rB}{1.5}[-2]{rI}
\strCrossCap{r}{rI}{rB}[//c1/white][//c1/]

\strLabelX{r}{$r$}[//c1/white]
\strLabelX{rB}{$r^\bullet$}[/c1//c1]

\strTwoCellX{rI/c1//c1}{-}{rI//c1/c1,c1}

\strModXX{rI//c1/c1,c1}{$\gamma_1$}{l//c1/c3,c1;b//c1,c3/c1,c1}[1]

\strTwoCellX{l//c1/c3,c1}{-}{l/c1//c3}
\strTwoCellX{b//c1,c3/c1,c1}{-}{b/c1/c3/c1}

\strIdToX{r}{l}[//c1/white]
\strIdToX{rB}{b}[/c1//c1]

\strBraidXX{r}{l}[//c1/white][c1//c3]
\strBraidXX{rB}{b}[/c1//c1][c1/c3/c1]

\strCupX{r}{rB}[//c1/c3]

\strIdX{l}[//c3]
\strIdX{b}[/c3/c1]

\drawEndNodes{
    l/$l$,
    b/$b_{-,I}$
}

\strKeyXX{
    $-\tensor I$/c1;
    $I\tensor -$/c3
}
\end{tikzpicture}}
            \]
        \caption{$\gamma_2$}
        \label{gammaFig2}
        \end{figure}
        \begin{figure}[p]
        \[\myinput{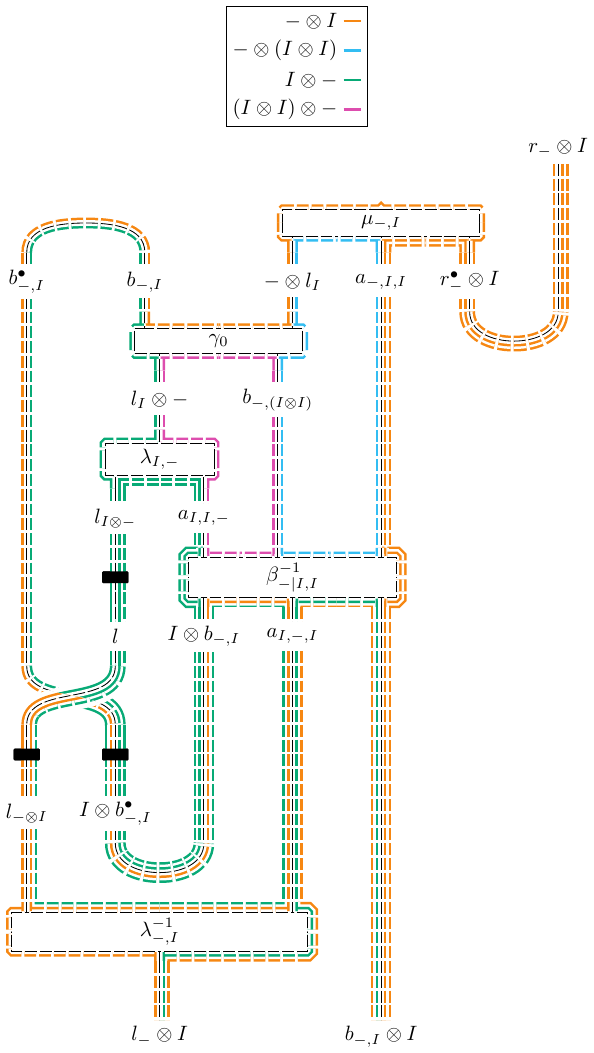}\]
        \caption{$\gamma_1$}
        \label{gammaFig1}
        \end{figure}
    \end{proof}
    Clearly, then, reasoning about braided monoidal bicategories is rather unwieldy. Luckily Gurski gives a strictification result that vastly simplifies proofs.
    \begin{definition}
        A \textdef{strict braided monoidal bicategory} is a braided monoidal bicategory such that the underlying monoidal bicategory is a Gray monoid, and any braid or 2-braid with a unit in its index is an identity map. To be clear:
        \begin{itemize}
            \item $A= A\tensor I \xrightarrow{b_{A,I}} I\tensor A= A$ and $A=I\tensor A \xrightarrow{b_{I,A}} A\tensor I=A$ are both the identity;
            \item for all $A$, $B$, and $C$, the 2-cells in the following diagrams are identities.
            \[
                \myinput{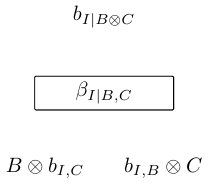}\qquad\myinput{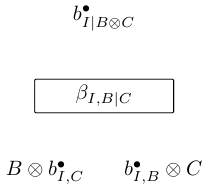}
            \]
            \[
                \myinput{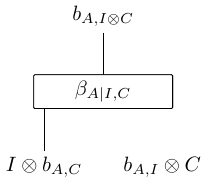}\qquad\myinput{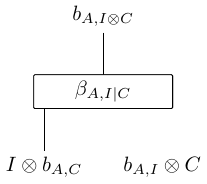}
            \]
            \[
                \myinput{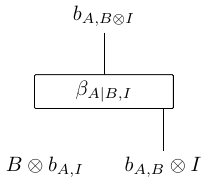}\qquad\myinput{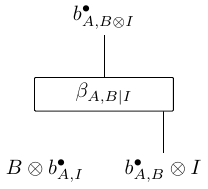}
            \]
        \end{itemize}
    \end{definition}
    \begin{theorem}
        Every braided monoidal category is braided monoidally biequivalent to a strict braided monoidal category.
    \end{theorem}
    \begin{proof}
        See Gurski's~\cite[thm.~27]{GurskiBraided} paper: \emph{Loop spaces, and coherence for monoidal and braided monoidal bicategories}. 
    \end{proof}
    Baez and Dolan~\cite[sec.~5]{baez1995tqft}, in their periodic table of $n$-categories, conjectured that a one-object braided monoidal bicategory ought to be a symmetric monoidal category. It was then pointed out by Baez and Neuchl~\cite[sec.~5]{BaezNeuchl} that if a braided monoidal bicategory could be strictified in the sense of Gurski above, then the conjecture would hold. In light of Gurski's strictification theorem we give an explicit proof of that conjecture using our string diagram language.
    \begin{lemma}
        If $\BB$ is a braided monoidal bicategory then the braiding for its category of scalars is a symmetry.
    \end{lemma}
    \begin{proof}
        Recall that if we strictify our monoidal bicategory then the braiding for the scalars is given by the interchangerator. Let $T\colon \BB\times \BB\to \BB\times \BB$ be the twist 2-functor that reorders pairs of objects, 1-cells and 2-cells. Given $f\colon A_0\to A_1$ and $g\colon B_0\to B_1$ consider the 2-cell in \autoref{scBraid1}. By definition this 2-cell is equal to the 2-cell given in \autoref{scBraid2}. \autoref{scBraid2} is equal to \autoref{scBraid3} by naturality, which is equal to \autoref{scBraid4} and \autoref{scBraid5} by definition. But if $f$ and $g$ are scalars, and we think of $\BB$ as being strict braided monoidal, then we know that $b$ and $b^\bullet$ are identity 1-cells. Thus, we know that $\chi_{g,f}=\chi_{f,g}^{-1}$ and so the braid for scalars is a symmetry.
    \end{proof}

    \newpage
    \begin{definition}
    A \textdef{symmetric monoidal bicategory} is a braided monoidal bicategory equipped with an invertible modification called the syllepsis,
        \[\myinput{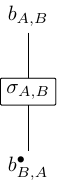}\]
    which is subject to coherence axioms that can be found, for example, in Stay's~\cite[defs.~4.7,~4.8]{stay2016} paper on compact-closed bicategories. We choose to omit these, once again, since we will not appeal to them directly.
    \end{definition}    
    \begin{remark}
        Depending on which coherence axioms we include there is a weaker notion of a monoidal bicategory with a commutativity structure. This is called a \textdef{sylleptic monoidal category}, hence the name syllepsis. This is part of the story of the periodic table of $n$-categories, but we will not make use of this definition here.
    \end{remark}
    \begin{remark}
        If a symmetric monoidal category is right-monoidal-closed then it is also left-monoidal-closed, and so we simply refer to the category as monoidal-closed.
    \end{remark}
    In a symmetric monoidal category the following two morphisms are equal. 
    \[
        \myinput{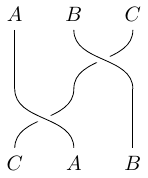} \quad =\quad         \myinput{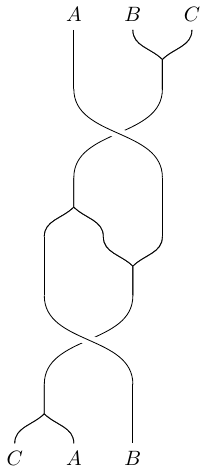}
    \]
    The next proposition gives an analogue of this fact for symmetric monoidal bicategories, and will be useful in our account of compact-closed bicategories.
    \begin{proposition}
        \label{prop:symmetry}
    In a symmetric monoidal bicategory, $\BB$, there is a modification with components given by the diagram below.
    \[
        \myinput{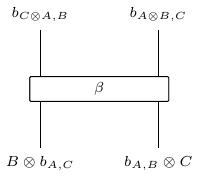}\]
    \end{proposition}
    \begin{proof}
        The 2-cell is given by the composite in \autoref{ch6.1}, which is a modification since each of its constituent 2-cells are modifications.
        \begin{figure}[h]
            \centering
            \[\myinput{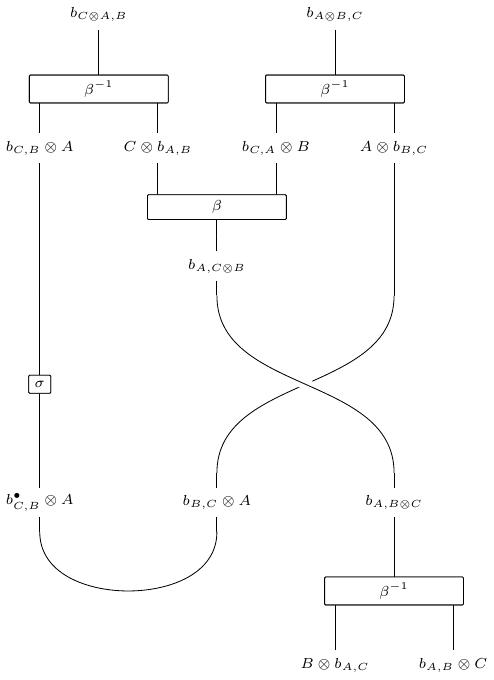}\]
            \caption{}
            \label{ch6.1}
        \end{figure}
    \end{proof}
    \begin{definition}
        A symmetric monoidal bicategory is \textdef{strict} if it is strict as a braided monoidal bicategory and the syllepsis is the identity.
    \end{definition}
    Gurski and Osorno have also given a strictification theorem for symmetric monoidal bicategories.
    \begin{theorem}
        Every symmetric monoidal category is symmetric monoidally biequivalent to a strict braided monoidal category.
    \end{theorem}
    \begin{proof}
        See Gurksi and Osorno's~\cite[thm.~1.13]{GurskiSymmetric} paper: \textbook{Infinite loop spaces, and coherence for symmetric monoidal bicategories}.
    \end{proof}
\section{Compact-Closed Bicategories and the Trace}
\label{section:CompactClosed}
In this section we explore the properties of compact-closed bicategories with the ultimate aim of defining the trace functor. Note here that, as with `monoidal-closed', we hyphenate `compact-closed' to distinguish this closed structure from composition-closed structure. Throughout this section we will apply semi-strictification to improve the readability of certain definitions and proofs.
\begin{definition}
Let $\BB$ be a symmetric monoidal bicategory. A \textdef{dual} for an object $A\in \BB$ is an object $A^*$ equipped with a 1-cell
\[
\coev_A\colon I\to A\tensor A^*    
\]
that we call \textdef{coevaluation}, and a 1-cell
\[
\ev_A\colon A^*\tensor A\to I    
\]
that we call \textdef{evaluation}, such that $-\tensor A$ is left pseudoadjoint to $-\tensor A^*$ with unit given by $-\tensor \coev_A$ and counit given by $-\tensor \ev_A$.
\end{definition}
To be explicit about what this means, if $A$ is dual to $A^*$ then there are adjoint equivalences given by the following two pairs of 1-cells that we call the \textdef{zigzag 1-cells}.
\[
\myinput{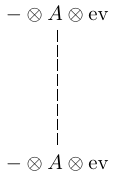} \qquad \myinput{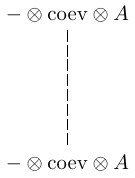}
\]
\[
\myinput{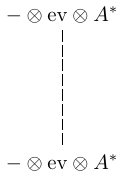} \qquad\myinput{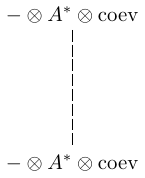}
\]
Recall that the definition of pseudoadjunction we are working with here does not stipulate any additional coherence conditions. Thus, this definition does not include the swallowtail diagrams that can be found in, for example, Stay's~\cite[def.~4.11]{stay2016} account of compact-closed bicategories. Again, as shown by Verity~\cite[lem.~1.3.9]{verity} this is essentially immaterial since we can always replace the 2-cell counit or unit of the hom-equivalence
\[
\BB(B\tensor A, C)\cong \BB(B,C\tensor A^*)    
\]
so that it becomes an adjoint equivalence. This corresponds to changing either the cups or the caps that show that evaluation and coevaluation are adjointly equivalent. If the cups and caps give an adjoint equivalence between the hom-categories, then they adhere to the swallowtail identities. In short, we can always choose our cups or caps in such a way as the swallowtail identities hold. Thus, if $A$ and $A^*$ are dual in the sense of the definition above, then they are dual in the sense of Stay~\cite[def.~4.11]{stay2016}.
\begin{proposition}
    In every symmetric monoidal bicategory the unit object is self-dual.
\end{proposition}
\begin{proof}
    This follows by taking coevaluation to be $l^{\bullet}$ and evaluation to be $l$.
\end{proof}

\begin{proposition}
    Let $A$ be an object in a symmetric monoidal bicategory. If $A$ is dual to both $A'$ and $A^*$ then there is a canonical equivalence between $A'$ and $A^*$.
\end{proposition}
\begin{proof}
    This follows from the uniqueness of pseudoadjoints. Since $-\tensor A'$ and $-\tensor A^*$ are both right pseudoadjoint to $-\tensor A$, we know that there is a canonical equivalence $B\tensor A'\cong B\tensor A^*$ for all $B\in \BB$. In particular, there is a canonical equivalence $I\tensor A'\cong I\tensor A^*$.
\end{proof}
\begin{proposition}
    \label{prop:dualBothWays}
    Let $A$ be an object in a symmetric monoidal bicategory. If $A$ is dual to $A^*$ then $A^*$ is dual to $A$.
\end{proposition}
\begin{proof}
    We will show that $b_{A,A^*}\circ \coev$ and $\ev \mathbin{\circ} b_{A,A^*}$ give the coevaluation and evaluation maps for $A^*$.

    In order to prove this we need to show that our chosen coevaluation and evaluation maps form the unit and counit for a pseudoadjunction between $-\tensor A^*$ and $-\tensor A$. In other words, we need to check that the following functors
    \[
    ((-)\tensor A)\circ (B\tensor b_{A^*,A})\circ (B\tensor  \coev)\colon \BB(B\tensor A^*,C)\to \BB(B,C\tensor A)
    \]
    \[
    (C\tensor \ev)\circ (C\tensor b_{A^*,A})\circ ((-)\tensor A^*)\colon \BB(A,C\tensor A)\to \BB(B\tensor A^*, C)
    \]
    form an equivalence. Note that if we compose the above functors we get the following functors:
    \begin{equation}
    \label{functor:Duals}
        (C\tensor \ev)\circ (C\tensor {b_{A,A^*}})\circ ((-)\tensor A\tensor A^*)\circ (B\tensor b_{A,A^*}\tensor A^*)\circ (B\tensor \coev \tensor A^*);
    \end{equation}
    \begin{equation}
    \label{functor:Duals2}
        (C\tensor \ev \tensor A)\circ (C\tensor {b_{A,A^*}\tensor A})\circ ((-)\tensor A^*\tensor A)\circ (B\tensor b_{A,A^*})\circ (B\tensor \coev).
    \end{equation}
    Thus, it suffices to show natural isomorphisms from (\ref{functor:Duals}) to the identity, and from (\ref{functor:Duals2}) to the identity. The first natural isomorphism is defined by 2-cell in \autoref{ch6.2}  where $\beta$ is the modification given by \autoref{prop:symmetry}. The two bottom modifications follow from strictification for braided monoidal bicategories.
    \begin{figure}[h]
        \centering
        \[\myinput{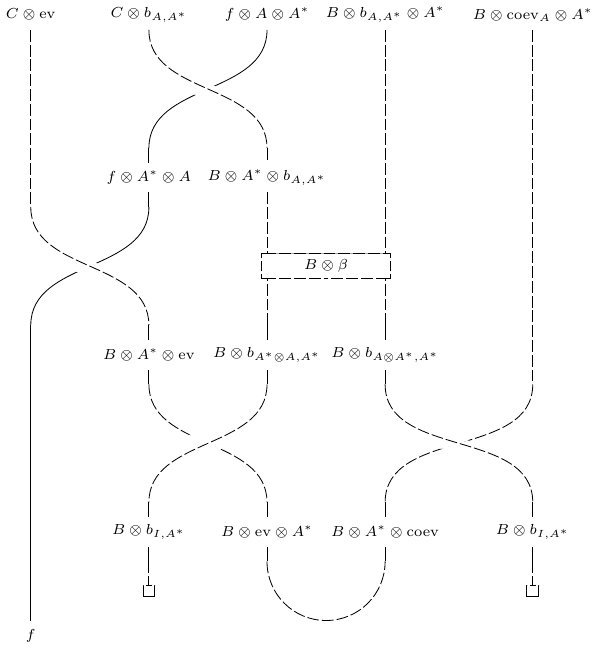}\]
        \caption{}
        \label{ch6.2}
    \end{figure}

     The fact that this constitutes a natural transformation follows from the fact that the 2-cell is a modification. 
     
     The second natural isomorphism is given analogously. Thus, $B\tensor (b\circ \coev)$ and $(\ev\circ b)\tensor C$ give a pseudoadjunction $-\tensor A^*$ $-\tensor A$, and so $A^*$ is dual to $A$ whenever $A$ is dual to $A^*$. 
\end{proof}
\begin{corollary}
    Let $A$ be an object in a symmetric monoidal bicategory. If $A^*$ is dual to $A$, and $A^{**}$ is dual to $A^*$ then there is a natural equivalence from $A$ to $A^{**}$.
\end{corollary}
\begin{proof}
    In this case $A$ and $A^{**}$ are both dual to $A$ and so, by uniqueness of duals, they must be equivalent.
\end{proof}
\begin{definition}
A \textdef{compact-closed bicategory} is a symmetric monoidal bicategory in which every object comes equipped with a choice of dual.
\end{definition}
Note that this definition is given in terms of a structure, but we might as well define a compact-closed bicategory as a symmetric monoidal bicategory in which every object \emph{has} a dual. Since duals are unique up to equivalence, any choice of duals would be equivalent. For this reason, since $I$ is self-dual, we will assume that the choice of dual for $I$ is always given by $I^*=I$. This is a sort of partial strictification that will be convenient for discussing the impact that duals have on scalars.
\begin{proposition}
    In a compact-closed bicategory $\BB$ there is a pseudofunctor, called the \textdef{dualising pseudofunctor},
    \[
        (-)^*\colon \BB^{\op}\to \BB
    \]
    which takes every object to its chosen dual and has hom-functors given by
    \[
    \BB(B,A)\xrightarrow{A^*\tensor -\tensor B^*} \BB(A^*\tensor B\tensor B^*,A^*\tensor A\tensor B^*)\xrightarrow{(\ev\tensor B^*)\circ - \circ (A^*\tensor\coev)} \BB(A^*, B^*).
    \]
\end{proposition}
\begin{proof}
    We need to define the compositor and the identitor for the pseudofunctor. By definition of a dual, $\ev\tensor A^*$ is adjoint equivalent to $\coev\tensor A^*$, and $A\tensor \ev$ is adjoint equivalent to $\coev \tensor A$. Then the canonical 2-cell, defined by the diagram below, gives the identitor.
    \[
        \myinput{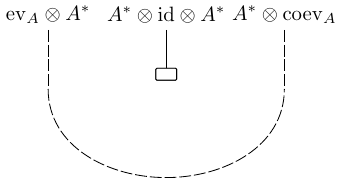}
    \]
    Now, given 1-cells $f\colon B\to A$ and $g\colon C\to B$, the compositor for the dualising pseudofunctor is defined by the following diagram.
    \[
        \myinput{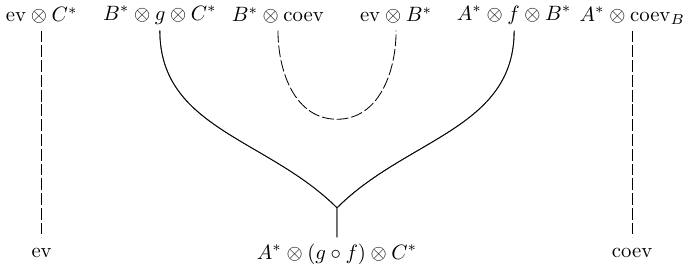}  
    \]
    The fact that these 2-cells adhere to the axioms of pseudofunctors follows from the fact that they are constructed from the tensor identitor and compositor and properties of adjoint equivalences. 
\end{proof}
\begin{remark}
    Note that the source of this pseudofunctor is $\BB^{\op}$ \emph{not} $\BB^{\operatorname{coop}}$, so the dualising pseudofunctor only flips 1-cells, not 2-cells. 
\end{remark}
The above definition is useful for showing that the dualising pseudofunctor is, in fact, a pseudofunctor. But the associated hom-functor can be reinterpreted as composite of name and realisation functors, and this interpretation will simplify working with duals.
\begin{proposition}
    Let $(-)^*\colon \BB(B,A)\to \BB(A^*,B^*)$ denote the hom-functor for the dualising pseudofunctor. There is an isomorphism, natural in $A$ and $B$, between $(-)^*$ and the following functor:
    \[
    \BB(B,A)\xrightarrow{\name{(-)}} \BB(I, A\tensor B^*) \xrightarrow{\unname{(-)}} \BB(A^*, B^*).  
    \]
\end{proposition}
\begin{proof}
    This follows because the naming functor above is given by
    \[
    \BB(B,A)\xrightarrow{B^*\tensor (-)} \BB(B^*\tensor B, B^*\tensor A) \xrightarrow{\coev\circ (-)} \BB(I, B^*\tensor A),
    \]
    the realisation functor above is given by
    \[
    \BB(I, A\tensor B^*)\xrightarrow{(-)\tensor B} \BB(B,A\tensor B^*\tensor B)\xrightarrow{(A\tensor \ev)\circ (-) } \BB(B,A),
    \]
    and the rest follows by pseudonaturality. 
\end{proof}
\begin{corollary}
    The dualising pseudofunctor is a biequivalence.
\end{corollary}
\begin{proof}
    This follows from the fact that both the name and realisation functors give equivalences.
\end{proof}
Note that, in contrast to the various symmetry structures applied to monoidal bicategories, a compact-closed structure on a bicategory will not guarantee any kind of duality structure on the monoidal category of scalars.
\begin{proposition}
    At the unit, the hom-functor for the dualising pseudofunctor
    \[(-)^*\colon \BB(I,I)\to \BB(I,I)\] 
    is naturally isomorphic to the identity functor. 
\end{proposition}
\begin{proof}
This follows from the definition of $(-)^*$ and the fact that evaluation and coevaluation for the unit are given by the left unitor and its adjoint pseudoinverse.
\end{proof}

\begin{proposition}
A compact-closed bicategory is monoidal-closed with the monoidal-closed structure given by $-\tensor A^*$.
\end{proposition}
\begin{proof}
This follows from the definition of a dual.
\end{proof}
In \autoref{section:ClosedMonoidal} we gave a number of examples of monoidal-closed bicategories. In fact, all of our examples were compact-closed. The unit and counit maps given there can be used to reverse-engineer the coevaluation and evaluation maps. Below we give an explicit account of what the dualising pseudofunctor does to objects and 1-cells in each of our examples.
\begin{example}
    The bicategory $\Rel$ is compact-closed. The dual of a set $A$ is the same set $A$. The dual of a relation $R\subset A\times B$ is the transpose relation $R^T$.
\end{example}
\begin{example}
    The bicategory $\Bim_R$ is compact-closed. The dual of an algebra $A$ is $A^{\op}$, the opposite algebra. The dual of an $A$-$B$-bimodule is the same bimodule, thought of as a $B^{\op}$-$A^{\op}$-bimodule.
\end{example}
\begin{example}
    The bicategory $\DGBim_R$ is compact-closed. The dual of an algebra $A$ is $A^{\op}$, the opposite algebra. The dual of a differential graded bimodule works similarly to the above.
\end{example}
\begin{example}
    The bicategory $\VV$-Prof is compact-closed. The dual of a category is $\AA$ is $\AA^{\op}$, the opposite category. A profunctor $P\colon \AA\profto \BB$ is a functor
    \[
    P\colon \BB^{\op}\tensor \AA\to \VV    
    \]
    and by symmetry in $\VV$ this can be identified with a functor
    \[
    P\colon \AA\tensor \BB^{\op}\to \VV    
    \]
    which gives a profunctor $P\colon \BB^{\op}\profto \AA^{\op}$.
\end{example}
\begin{example}
    Given a category $\CC$ with finite limits, the bicategory $\Span(\CC)$ is compact-closed. The dual of an object $C$ is the same object $C$. The dual of a span $A\xleftarrow{f} S\xrightarrow{g} B$ is the span $B\xleftarrow{g} S\xrightarrow{f} A$.
\end{example}
\begin{example}
    Given a topological abelian group $G$, the bicategory $\Path(G)$ is compact-closed. The dual of a point $x$ is its inverse $x^{-1}$. The dual of a path $p\colon x\to y$ is path $(y^{-1}\cdot p\cdot x^{-1})\colon y^{-1}\to x^{-1}$.
\end{example}
\begin{definition}
A symmetric monoidal bicategory $\BB$ is called \textdef{coclosed} if its opposite bicategory is closed. That is to say, $\BB$  comes equipped with a pseudofunctor 
\[
\hom{-,-}\colon \BB^{\op}\times \BB\to \BB    
\]
such that for every $A\in \BB$, $\hom{A,-}$ is left pseudoadjoint to $-\tensor A$.
\end{definition}
\begin{proposition}
    A compact-closed bicategory is coclosed with coclosed structure given by $-\tensor A^*$.
\end{proposition}
\begin{proof}
    Since we know that $A^{**}$ is naturally isomorphic to $A$ we have a sequence of isomorphisms $
    \BB(B\tensor A^*,C)\cong \BB(B, C\tensor A^{**})\cong \BB(B,C\tensor A)
    $ that give coclosedness.
\end{proof}
It is worth noting here that we can use the name and realisation functors in the `opposite direction'. Suppose that we have a 1-cell $f\colon A\tensor B^*\to I$. Then the name of such a 1-cell is given by
\[
\name{f}\colon A\to I\tensor B^{**}\xrightarrow{\sim} B    
\]
and so given a 1-cell $f\colon A\to B$ we can take the `realisation' of $f$ to be a map \[\unname{f}\colon A\tensor B^*\to I\] which is more of a `coname' than a realisation, since it is the name functor for the coclosed structure.

Now that we have given some basic properties of a compact-closed bicategory, we will investigate how the compact-closed structure interacts with composition-closed structure.
\begin{proposition}
    \label{prop:weakdaggercompact}
    A left-composition-closed compact-closed bicategory is also right-composition-closed with the right extension of $g$ along $f$ given by
        \[(g\extend f)\cong (g^*\lift f^*)^*.\]
\end{proposition}
\begin{proof}
This follows from the following chain of adjoint natural isomorphisms
\begin{align*}
    \BB(A,C)(h\circ g,f)&\cong \BB(C,A)((h\circ g)^*, f^*)\\
                        &\cong \BB(C,A)(g^*\circ h^*, f^*)\\
                        &\cong \BB(C,B)(h^*, g^*\lift f^*)\\
                        &\cong \BB(B,C)(h, (g^*\lift f^*)^*),
\end{align*}            
since, recall, the dualising functor $(-)^*$ only reverses 1-cells.
\end{proof}
\begin{remark}
    This means that if we have a left-composition-closed compact-closed bicategory we can simply refer to it as composition-closed compact-closed.

    It is worth noting here that, whilst a definition for compact-closed bicategory was yet to be given formally, Betti and Walters~\cite[def.~5.6]{betti} suggested that the name \textdef{closed bicategory} should refer to composition-closed compact-closed bicategory in which all hom-categories are finitely complete and cocomplete.
\end{remark}
\begin{remark}
    In the next section we will see how this particular property relates to the notion of dagger compactness.
\end{remark}
A priori, right-composition-closedness has the potential to give another scalar enrichment. We have shown that every left-composition-closed, right-monoidal-closed bicategory can be enriched over the bicategory of scalar enriched categories. But by considering this result on $\BB^{\op}$ we could just have easily shown the result for right-composition-closed, right-monoidal-\emph{coclosed} bicategories. Thus, a left-composition-closed compact-closed bicategory might have two distinct scalar enrichments. Given 1-cells $f,g\colon A\to B$ the scalar between them would then be the extension in the following diagram.
\[
\begin{tikzcd}[column sep=3.2em]
        I
        & I
        \arrow[l, "\unname{g}\extend \unname{f}"{swap}, dashed]
            & A\tensor B^*
            \arrow[l, "\unname{f}"{swap}]
            \arrow[ll, "\unname{g}", ""{swap,name=bottom}, bend left=65]
    \arrow[Leftarrow, to=\tikzcdmatrixname-1-2, from=bottom, verticalup]
\end{tikzcd}    
\]
It turns out that this just gives us the enrichment defined in the previous chapter.
\begin{proposition}
    If $\BB$ is a composition-closed compact-closed bicategory then for 1-cells $f,g\colon A\to B$ there are natural isomorphisms
    \[
        \unname{g}\extend \unname{f} \cong \BBB(A,B)(f,g).
    \]
\end{proposition}
\begin{proof}
    Recall that $(-)^*\colon \BB(I,I)\to \BB(I,I)$ is naturally isomorphic to the identity and so for any scalar $s$ we have
    \begin{align*}
        \BB(I,I)(s, \unname{g}\extend \unname{f})&\cong \BB(I,I)(s, ((\unname{f})^*\lift (\unname{g})^*)^*)\\
        &\cong \BB(I,I)(s^*, (\unname{f})^*\lift (\unname{g})^*)\\
        &\cong \BB(A\tensor B^*,I)((\unname{f})^*\circ s^*, (\unname{g})^*)\\
        &\cong \BB(A\tensor B^*,I)((s\circ \unname{f})^*, (\unname{g})^*)\\
        &\cong \BB(A\tensor B^*,I)(s\circ \unname{f}, \unname{g}).
    \end{align*}
    But now note that, since $\unname{(-)}$ is the name functor for $\BB^{\op}$, by \autoref{prop:altAction}, we have natural isomorphisms 
    \begin{align*}
        \BB(A\tensor B^*,I)(s\circ \unname{f}, \unname{g})&\cong \BB(A\tensor B^*,I)(\unname{f\circ \Spr(s)}, \unname{g})\\
        &\cong \BB(A,B)(f\circ \Spr(s),g)\\ &\cong \BB(A,B)(\Spr(s), f\lift g)\\&\cong \BB(I,I)(s, \ctr(f\lift g)).
    \end{align*}
    Thus, there is a natural isomorphism 
    \[
        \BB(I,I)(s, \unname{g}\extend \unname{f})\cong \BB(I,I)(s,\ctr(f\lift g))  
    \]
    which, by the Yoneda lemma, implies that there is a natural isomorphism \[\unname{g}\extend \unname{f}\cong \ctr(f\lift g).\]
\end{proof}

\begin{definition}
The \textdef{trace} functor for an object $A$ in a compact-closed bicategory $\BB$ is the functor
\[
    \rTr\colon \BB(A,A)\to\BB(I,I)
\]
given by taking the composite
\[
\BB(A,A)\xrightarrow{\name{(-)}}\BB(I,A\tensor A^*)\xrightarrow{((\ev\circ b)\circ (-))} \BB(I,I).
\]
\end{definition}
If we want to be explicit about how this functor acts on endomorphisms, the functor takes in an endo-1-cell $f\colon A\to A$ and returns the scalar given by the composite
\[
    I\xrightarrow{coev} A\tensor A^*\xrightarrow{f\tensor A^*}A\tensor A^*\xrightarrow{b} A^*\tensor A\xrightarrow{\ev} I.
\]
Thus, if $\BB$ is a compact-closed bicategory without any non-identity 2-cells, the trace functor is just the usual trace function for compact-closed categories.
Of course, by \autoref{prop:eName} and by \autoref{prop:ePostComp} there is an enriched version of this functor.
\begin{definition}
    The \textdef{enriched trace} functor for an object $A$ in a compact-closed bicategory $\BB$ is the functor
    \[
        \rTr\colon \BBB(A,A)\to\BBB(I,I)
    \]
    given by taking the composite
    \[
    \BBB(A,A)\xrightarrow{\name{(-)}}\BBB(I,A\tensor A^*)\xrightarrow{(\ev\circ b)\ecirc (-)} \BBB(I,I).
    \]
\end{definition}
Most of the following examples of traces in compact-closed bicategories are due to Willerton~\cite{willertonTalk}.
\begin{example}
    In the compact-closed bicategory $\Rel$ the trace of a relation $R\colon A\to A$ is given by the relation
    \[
    \Tr(R)=\begin{cases}
        * & \text{if }\exists a\in A\; aRa\\
        \varnothing&\text{otherwise.} 
    \end{cases}    
    \]
    In other words, $\Tr(R)=*$ if and only if $A^c$ is not reflexive.
\end{example}
\begin{example}
    In the compact-closed bicategory $\Bim_R$ the trace of an $A$-$A$-bimodule $M$ is given by the $R$-module of coinvariants, or the abelianisation, of the algebra. This is the quotient $R$-module
    \[
    M/\langle ma-am \mid m\in M, a\in A \rangle.
    \]
    This follows from the fact that the trace can be thought of as the composite
    \[
    R \xrightarrow{\name{A}} A\tensor A^{\op} \xrightarrow{\unname{M}} R  
    \]
    which given as a tensor is
    \[
    A\tensor_{A^{\op}\tensor A} M. 
    \]
    By definition of the tensor product we get the quotient described above.
\end{example}
\begin{example}
    In the compact-closed bicategory $\DGBim_R$, when $R$ is a field, the trace of a differential graded $A$-$A$-bimodule $M$ is given by the Hochschild homology of $A$ with coefficients in $M$. This follows from the fact that the trace can be thought of as a composite
    \[
    R \xrightarrow{\name{A}} A\tensor A^{\op} \xrightarrow{\unname{M}} R  
    \]
    which is given by the derived tensor
    \[
     R\tensor^L_{A\tensor A^{\op}} M.
    \]
    This derived tensor has $n$'th component given by $\Tor_n^{A\tensor A^{\op}}(A,M)$. As shown by Cartan and Eilenberg~\cite[ch.~IX]{cartan1956homological} this is equivalent to the definition of the Hochschild homology of $A$ with coefficients in $M$.
\end{example}
\begin{example}
    In the compact-closed bicategory $\VV$-Prof the trace of a profunctor $P\colon \AA\profto \AA$ is the coend of $P$. This follows from the Fubini theorem for coends -- see for example Loregian's~\cite[thm.~1.3.1]{Fosco} book --  and the co-Yoneda lemma:
    \begin{align*}
    \Tr(P)\cong\unname{P}\circ \name{\id}\cong \endint^{(A,A')\in \AA^{\op}\tensor \AA} \Hom(A,A')\tensor P(A,A')&\cong\endint^A\endint^{A'}\Hom(A,A')\tensor P(A,A')\\&\cong \endint^{A\in \AA} P(A,A).
    \end{align*}
\end{example}
\begin{example}
    Given a category $\CC$ with finite limits, in the compact-closed bicategory of spans $\Span(\CC)$, the trace of a span $A \xleftarrow{f} S\xrightarrow{g} A$ is given by the pullback in the following diagram.
    \[
        \begin{tikzcd}
            S\times_A S
            \arrow[r]
            \arrow[d]
            \arrow[rd, phantom, very near start, "\ulcorner"]
                &S
                \arrow[d, "f"]\\
            S
            \arrow[r, "g"{swap}]
                &A.
        \end{tikzcd}
    \]
    If $\CC$ is the category of sets then this trace is the set $\{(s,s')\in S\times S\mid f(s)=g(s')\}$.
\end{example}
\begin{example}
    If $G$ is a topological group, in the compact-closed bicategory $\Path(G)$ the trace of a loop $p$ at $x$ is given by the loop $(p\cdot x^{-1})\colon e\to e$.
\end{example}
\begin{remark}
Comparing the above examples of trace to our examples of cotrace it seems that there is some sort of duality between the two: cohomology vs homology, centre vs abelianisation, ends vs coends. Formalising this isn't exactly straightforward. For a start we need some kind of context for this duality -- some kind of category structure where the two functors can be compared. Assuming that $\BBB(I,I)$ is complete and cocomplete we can construct the functor category,
\[
 \hom{\BBB(A,A), \BBB(I,I)}
\]
and since both $\BBB(A,A)$ and $\BBB(I,I)$ are monoidal, this has a monoidal structure given by Day~\cite[p.~19]{dayconvolution} convolution. This would seem to be the obvious choice of category in which to compare the two functors. But since $\BBB(A,A)$ is enriched over $\BBB(I,I)$, the unit for Day convolution is $\BBB(\id,-)$. This is, of course, the cotrace\footnote{Thank you to Ciaran Cassidy for pointing out that the Day convolution of the trace and cotrace functors is `probably not interesting'.}. Thus, the trace and the cotrace cannot be dual in this category, otherwise they would have to be isomorphic as functors. In the final sections of this chapter we investigate how the trace and cotrace are related, in order to find some formal explanation for this putative duality.
\end{remark}

Clearly by \autoref{prop:eName} and by \autoref{prop:ePostComp} the underlying functor of the enriched trace functor is the trace functor. If our bicategory is composition-closed as well as compact-closed then the trace, like the cotrace, has an adjoint.
\begin{proposition}
    For every $A\in \BB$ the enriched trace functor has a right adjoint
    \[
        \rSpr_A\colon \BB(I,I)\to \BB(A,A)
    \]
    that we call the enriched cospread functor.
\end{proposition}
\begin{proof}
    The enriched cospread functor is given by the composite
    \[
    \BBB(I,I)\xrightarrow{(\ev\circ b)\elift (-)}\BBB(I, A\tensor A^*) \xrightarrow{\unname{(-)}} \BB(A,A) 
    \]
    and the fact that this is right adjoint to the trace follows from the fact that the enriched name and realisation functors give an equivalence, whilst the enriched postcomposition and lift functors give an adjunction.
\end{proof}
\begin{remark}
    Note that we now have two pairs of enriched functors that give adjoints. Firstly, there is trace and cospread
    \[
        \rTr\colon \BB(A,A)\rightleftarrows\cocolon \BB(I,I) \cocolon \rSpr.
    \]
    Secondly there is spread and cotrace
    \[
        \Spr\colon \BB(I,I)\rightleftarrows \BB(A,A)\cocolon \dTr.
    \]
    In order to remember names and symbols note that in both cases the functor whose name starts with ``co-" is the right adjoint, has a $\lift$ over the top and is constructed using the lift functor.
\end{remark}
This immediately implies an underlying functor to the cospread functor, which is right adjoint to the unenriched trace functor.
\begin{corollary}
    The trace functor has a right adjoint that we call the cospread functor.
\end{corollary}
\section{Dagger Compact Categories}
\label{section:daggerComp}
In a composition-closed compact-closed bicategory we have a spread and a cospread, a trace and a cotrace. It is here that our analogy with linear algebra breaks down somewhat. Linear endomorphisms have a trace but no cotrace. This is because, in some sense, the trace and the cotrace seem to coincide. Heuristically this seems to follow from the idea that weighted limits and colimits both correspond to sums, and left and right adjoints both correspond to linear adjoints.

To try and give some formal backing to this idea we will compare composition-closed compact-closed bicategories to dagger compact categories. Dagger compact categories, in some sense, completely describe the category of finite dimensional Hilbert spaces. It turns out that the combination of a composition-closed structure and a compact-closed structure gives something very similar to the structure of a dagger compact category, and under certain conditions the trace and cotrace will coincide.
\begin{definition}
    \label{def:daggerCompCat}
    A \textdef{dagger compact category} is a compact-closed category $\CC$ equipped with a strict symmetric monoidal involutive functor
    \[
        (-)^\dagger\colon \CC^{\op}\to \CC
    \]
    that is the identity on objects, such that for every object $A\in \CC$ the following diagram commutes.
    \begin{equation}
        \label{diag:dagcompCondition}
        \begin{tikzpicture}[commutative diagrams/every diagram]
            \node (P0) at (90:1.8cm) {$I$};
            \node (P1) at (90+120:1.8cm) {$A\tensor A^*$} ;
            \node (P2) at (90+240:1.8cm) {$A^*\tensor A$};
            \path[commutative diagrams/.cd, every arrow, every label]
            (P0) edge node[swap] {$\ev^\dagger$} (P1)
            (P0) edge node {$\coev$} (P2)
            (P1) edge node[swap] {$b$} (P2);
            \end{tikzpicture}
        \end{equation}
\end{definition}
\begin{definition}
    If $\CC$ is a dagger category -- that is a category equipped with an identity-on-objects involution $(-)^\dagger\colon \CC^{\op}\to \CC$ -- then we call a morphism $f\colon A\to B$ in $\CC$ \textdef{unitary} if $f^\dagger$ is the inverse of $f$.
\end{definition}
Then we can spell out the definition of a dagger compact category in terms of unitary morphisms.
\begin{proposition}
    A dagger compact category is a compact-closed category equipped with a dagger structure such that
    \[
    (f\tensor g)^\dagger=f^\dagger\tensor g^\dagger,
    \]
    the unitors, associator and braid are unitary, and (\ref{diag:dagcompCondition}) holds.
\end{proposition}
\begin{proof}
    Suppose $\CC$ is a dagger compact category as defined in \autoref{def:daggerCompCat}, then since $(-)^\dagger$ is strict symmetric monoidal, $(f\tensor g)^\dagger=f^\dagger\tensor g^\dagger$. Now note that the coherence maps for $\CC^{\op}$ are the inverses of those for $\CC$. But also note that since $(-)^{\dagger}$ is strict symmetric monoidal the associated coherence maps are identities. That means that the commutativity of the coherence diagrams reduce to the commutativity of the following diagrams.\[
    \begin{minipage}{0.45\textwidth}
        \centering
        \[\begin{tikzcd}
            {(A\tensor B)\tensor C} & {A\tensor(B\tensor C)}
            \arrow["a", curve={height=-18pt}, from=1-1, to=1-2]
            \arrow["{(a^{-1})^\dagger}"', curve={height=18pt}, from=1-1, to=1-2]
        \end{tikzcd}\]
    \end{minipage}
    \begin{minipage}{0.45\textwidth}
        \centering
        \[\begin{tikzcd}
	A\tensor B & B\tensor A
	\arrow["b", curve={height=-18pt}, from=1-1, to=1-2]
	\arrow["(b^{-1})^\dagger"', curve={height=18pt}, from=1-1, to=1-2]
\end{tikzcd}\]
    \end{minipage}
\]
\[
    \begin{minipage}{0.45\textwidth}
        \centering
\[
    \begin{tikzcd}
        I\tensor A & A
        \arrow["l", curve={height=-18pt}, from=1-1, to=1-2]
        \arrow["(l^{-1})^\dagger"', curve={height=18pt}, from=1-1, to=1-2]
    \end{tikzcd}
\]
\end{minipage}
\begin{minipage}{0.45\textwidth}
    \centering
\[
    \begin{tikzcd}
        A\tensor I& A
        \arrow["r", curve={height=-18pt}, from=1-1, to=1-2]
        \arrow["(r^{-1})^\dagger"', curve={height=18pt}, from=1-1, to=1-2]
    \end{tikzcd}
\]
\end{minipage}
\]
    Thus, any dagger compact category as defined in \autoref{def:daggerCompCat} is a dagger compact category as defined in the proposition above.

    The converse follows straightforwardly, since the stipulation that $(f\tensor g)^\dagger=f^\dagger\tensor g^\dagger$, as well as the unitarity of the coherence maps, forces $(-)^\dagger$ to be strict symmetric monoidal.
\end{proof}
Categories with this precise structure seem to have first been defined by Abramsky and Coecke~\cite[def.~7.1]{abramsky2004quantum}, under the name `strongly compact categories', in the context of quantum protocols. Prior to this Doplicher and Roberts~\cite[sec.~1]{doplicher1989new} studied $C^*$-categories with `conjugates', which are dagger compact categories enriched over Banach spaces. Baez and Dolan~\cite[p.~25]{baez1995tqft} also previously referred to monoidal categories `with duals', which are like dagger compact categories, but where the operation $(-)^*$ is only defined on objects. Selinger~\cite[thm.~2.2]{selinger2012finite} proved that all dagger compact categories are, in some sense, described by the category of finite dimensional Hilbert spaces. Formally speaking, Selinger developed a string diagram language for dagger compact categories and showed that if a given statement in that language holds in the category of finite dimensional Hilbert spaces, then it holds for all dagger compact categories. 

In a dagger compact category we refer to $f^\dagger$ as the adjoint of $f$, since for finite dimensional Hilbert spaces the dagger functor takes every linear map to its linear adjoint. Note that adjoints and duals interact well.
\begin{proposition}
    If $\CC$ is a dagger compact category then the following diagram commutes 
    \[
    \begin{tikzcd}
        \CC
        \arrow[r, "(-)^\dagger"]
        \arrow[d, "(-)^*"{swap}]
            &\CC^\op
            \arrow[d, "(-)^*"]\\
        \CC^\op
        \arrow[r, "(-)^\dagger"{swap}]
            &\CC
    \end{tikzcd}
    \]
\end{proposition}
\begin{proof}
    On objects these two functors clearly agree since they both take an object $A$ to its dual $A^*$. Thus, it suffices to show that $(f^{\dagger})^*$ is equal to $(f^*)^\dagger$ for all $f\colon A\to B$. Note that $f^*$ is equal to
    \[
    B^*
    \xrightarrow{B^*\tensor \coev_A } 
    B^* \tensor A\tensor A^* 
    \xrightarrow{B^*\tensor f\tensor A^*} 
    B^*\tensor B\tensor A^*
    \xrightarrow{\ev_B\tensor A^*} A^*  
    \]
    and so $(f^*)^{\dagger}$ is equal to the top path in \autoref{ch6.3}. But the bottom path gives $(f^\dagger)^*$, and we know that each of the internal polygons in \autoref{ch6.3} commutes, either by the commutativity of (\ref{diag:dagcompCondition}) or properties of the braid natural transformation.
    \begin{figure}[p]
    \[
    \rotatebox{90}{
    \begin{tikzcd}[column sep=5.4em, row sep=5.4em, ampersand replacement=\&]
        B^*
            \&B^*\tensor A\tensor A^*
            \arrow[l, "B^*\tensor \coev_A^\dagger"{swap}]
                \&B^*\tensor B\tensor A
                \arrow[l, "B^*\tensor f^\dagger \tensor A^*"{swap}]
                    \&A^*
                    \arrow[l, "\ev_B^\dagger \tensor A^*"{swap}]
                    \arrow[ld, "\coev_B\tensor A^*"{description}]
                    \arrow[dd, "\coev_B\tensor A^*"]\\
        B^*\tensor A^*\tensor A
        \arrow[u, "B^*\tensor \ev_A"{description}]
        \arrow[ur, "B^*\tensor b"{description}]
            \&A\tensor B^*\tensor A^*
            \arrow[u, "b\tensor A^*"{swap}]
                \&B\tensor B^*\tensor A^*
                \arrow[u, "b\tensor A^*"]
                \arrow[l, "f^{\dagger}\tensor B^*\tensor A"]
                    \&\blank\\
        A^*\tensor A\tensor B^*
        \arrow[uu, "\ev_A\tensor B^*", bend left=80]
        \arrow[ur, "b"]
        \arrow[u, "b"{description}]
            \&\blank
                \&\blank
                    \&A^*\tensor B\tensor B^*
                    \arrow[lu, "b"{swap}]
                    \arrow[lll, "A^*\tensor f^\dagger\tensor B^*"]
    \end{tikzcd}
    }
    \]
    \caption{}
    \label{ch6.3}
\end{figure}
\end{proof}
There is a similarity here with \autoref{prop:weakdaggercompact} given above. Note that, in a bicategory, if the right adjoint of $f$ exists then it is given by $f\lift \id$. Similarly, if the left adjoint of $f$ exists, then it is given by $\id\extend f$. So $f\lift \id$ and $\id \extend f$ behave like `weak' adjoints for $f$. 
\begin{corollary}
    \label{prop:weakdaggercompact2}
    If $\BB$ is a composition-closed compact-closed bicategory and $f$ has a weak right adjoint $f^\dagger$ then $(f^\dagger)^*$ is a weak left adjoint to $f^*$.
\end{corollary}
There are, however, two key differences between a dagger compact category and a composition-closed compact-closed bicategory that seem to explain the fact that dagger compact categories have a unique trace. Firstly, in a dagger compact category we know that the adjoint of evaluation is equal to coevaluation modulo braiding: $b\circ \ev^\dagger\cong \coev$, or equivalently, $(\ev\circ b)^\dagger\cong \coev$. Secondly, in a dagger compact category there is a single notion of adjoint. If $f$ is `left' adjoint to $g$ then it is also `right' adjoint to $g$. There is no distinction.
\begin{proposition}
    Let $\BB$ be a composition-closed compact-closed bicategory $\BB$ such that, for every $A$, the composite
    \[
        A\tensor A^*\xrightarrow{b} A^*\tensor A\xrightarrow{\ev_A}I
    \]
    is left adjoint to the coevaluation map. Then the spread functor is naturally isomorphic to the cospread functor.
\end{proposition} 
\begin{proof}
    Note that by \autoref{prop:altAction} it suffices to show that $\unname{\coev \circ s}$ is naturally isomorphic to $\unname{(\ev\circ b)\lift s}$ but $\coev\circ s$ is naturally isomorphic to $(\ev\circ b)\lift s$ by the fact that $\ev\circ b$ is the left adjoint to the coevaluation map.
\end{proof}
In cases where this holds, we have the following chain of adjunctions.
\[
\begin{tikzcd}[column sep=3.2em]
    \BB(I,I)
    \arrow[r, "\operatorname{Spr}"{description}, ""{swap, name=botmid}, ""{name=topmid}]
        &\BB(A,A)
        \arrow[l, bend right=75, ""{name=top}, "\tr"{swap}]
        \arrow[l, bend left=75, ""{swap, name=bottom}, "\ctr"]
    \arrow[phantom, "\dashv"{sloped}, from=top, to=topmid, vertical]
    \arrow[phantom, "\dashv"{sloped}, from=botmid, to=bottom, vertical]
\end{tikzcd}    
\]
\begin{proposition}
    Let $\BB$ be a composition-closed compact-closed bicategory $\BB$ such that, for every $A$, the composite
    \[
        A\tensor A^*\xrightarrow{b} A^*\tensor A\xrightarrow{\ev_A}I
    \]
    is right adjoint to the coevaluation map. Then the cotrace functor is naturally isomorphic to the trace functor.
\end{proposition} 
\begin{proof}
    This follows from the sequence of natural isomorphisms
    \[
    \tr(f)\cong \ev\circ b \circ (f\tensor A^*)\circ \coev \cong \coev\lift ((f\tensor A^*)\circ \coev)\cong \name{\id}\lift \name{f} \cong \ctr(f)
    \]
    which are due to the definition of each of the traces, the assumption above and the definition of the name functor.
\end{proof}
In cases where this holds, we have the following chain of adjunctions.
\[
\begin{tikzcd}[column sep=3.2em]
    \BB(A,A)
    \arrow[r, "\Tr"{description}, ""{swap, name=botmid}, ""{name=topmid}]
        &\BB(I,I)
        \arrow[l, bend right=75, ""{name=top}, "\spr"{swap}]
        \arrow[l, bend left=75, ""{swap, name=bottom}, "\cspr"]
    \arrow[phantom, "\dashv"{sloped}, from=top, to=topmid, vertical]
    \arrow[phantom, "\dashv"{sloped}, from=botmid, to=bottom, vertical]
\end{tikzcd}    
\]
Of course, if $(\ev\circ b)$ is both left and right adjoint to the coevaluation map -- as is formally the case for dagger compact categories -- then there is a single trace functor and a single spread functor. On its own $\ev\circ b$ being adjoint to coevaluation isn't enough to collapse a composition-closed compact-closed bicategory to a dagger compact category. We can reverse 1-cells by taking $f\lift \id$, but this isn't functorial: unless $f$ and $g$ have right adjoints, there is no isomorphism 
\[(g\circ f)\lift \id\xRightarrow{\sim} (f\lift \id)\circ (g\lift \id).\] 
\begin{proposition}
    Let $\BB$ be a composition-closed compact-closed bicategory, such that every $f\colon A\to B$ has unique ambi-adjoint $f^\vee$ and the adjoint to $(\ev\circ b)$ is coevaluation. Then the 2-skeleton of $\BB$ has an induced dagger compact structure.
\end{proposition}
\begin{proof}
    We firstly show that the 2-skeleton has a symmetric monoidal structure. We then show that this monoidal structure is compact-closed. After that we show that the 2-skeleton has a dagger structure, and finally we show that the dagger compact condition holds.

    Let $(\AA)_0$ denote the 2-skeleton of a bicategory $\AA$. Firstly note that given a pseudofunctor $F\colon \CC\to \DD$ there is a functor $F_0\colon \CC_0\to \DD_0$ given by $F_0(A)=F(A)$ on objects and $F(\hom{f})=\hom{F(f)}$ on morphisms. This is well-defined since $F$ preserves isomorphisms between 1-cells, and the compositor and identitor are isomorphisms. Secondly note that for any pair of bicategories $\CC$ and $\DD$, $(\CC\times \DD)_0=\CC_0\times \DD_0$ since any invertible 2-cell in $\CC\times \DD$ is given by a pair of invertible 2-cells, one in $\CC$ and one in $\DD$. Thus, we can define the tensor product on $\BB_0$ to be $\tensor_0$.

    Note that if $f$ and $f^\bullet$ give an equivalence then $\hom{f\circ f^\bullet}=\hom{\id}$ and $\hom{f^\bullet\circ f}=\hom{\id}$. Thus, the unitors and associator for the tensor can be given by $\hom{l}$, $\hom{r}$ and $\hom{a}$. The fact that the 2-associator and 2-unitors are isomorphisms means that the axioms of a monoidal category hold. Hence, the 2-skeleton has a monoidal structure. Similarly, we have braid morphisms $\hom{b}\colon A\tensor B\to B\tensor A$ with an inverse. The syllepsis and the 2-braids are both isomorphisms and so $\hom{b}$ is self inverse and adheres to the braid conditions.
    
    If $A$ and $A^*$ are dual in $\BB$ with chosen evaluation and coevaluation map. The isomorphism classes of these adhere to the yanking conditions and so $A$ and $A^*$ are dual in $\CC$. This gives the compact-closed structure.

    To prove that there is a dagger structure, we define the functor
    \[
    (-)^\dagger\colon \BB_0^{\op}\to \BB_0    
    \]
    to be the functor that is the identity on objects and sends $\hom{f}$ to $\hom{f^\vee}$. This functor preserves identities since $\id^\vee$ is $\id$ and preserves composition since if $g$ is adjoint to $g$, $f^v$ is adjoint to $f$, then $g^\vee\circ f^\vee$ is adjoint to $f\circ g$.

    We now show that the dagger functor is strict symmetric monoidal. In order to do this, we first characterise the unitary morphisms. A morphism $\hom{f}$ is unitary in $\BB_0$ if and only if
    \[
        \hom{f\circ f^\vee}=\hom{\id}\text{ and }\hom{f^\vee\circ f}=\hom{\id},
    \]
    but this holds if and only if $f$ and $f^\vee$ give an equivalence. So $f^\vee$ and $f$ form an adjoint equivalence. Similarly, if $f$ and $f^\bullet$ give an adjoint equivalence, then $\hom{f}^\dagger=\hom{f^\bullet}$ and we have that
    \[
        \hom{f\circ f^\vee}=\hom{\id}\text{ and }\hom{f^\vee\circ f}=\hom{\id}     
    \]
    and so $f$ and $f^\bullet$ are unitary. In particular the unitors and the associator are all unitary. Since pseudofunctors preserve adjoints
    \[
    (\hom{f}\tensor \hom{g})^\dagger=\hom{(f\tensor g)^\vee}=\hom{f^\vee\tensor g^\vee}=\hom{f}^\dagger\tensor \hom{g}^\dagger,
    \]
    and thus the dagger functor is a strict symmetric monoidal functor.

    Finally, we want to show that the diagram (\ref{diag:dagcompCondition}) commutes. In other words we need to know that $\hom{b}\circ\hom{\ev}^\dagger=\hom{\coev}$ holds. But this follows from the fact that $b\circ \ev$ is ambi-adjoint to $\coev$.

    Thus, we have shown that $\BB_0$ has the structure of a dagger compact category.
\end{proof}
\begin{example}
    The 2-skeleton of $\Path(G)$, for $G$ a topological group, is the groupoid given by points in $G$ and homotopy classes of paths between them. The dual of a point is the inverse of that point. The dagger is given by reversing the paths.
\end{example}
So it is the case that a small class of composition-closed compact-closed bicategories can be reduced to dagger compact categories. Interestingly, sets and relations form a dagger compact category where the dagger is given by taking the transpose relation. This can be used to describe the composition-closed structure of the bicategory $\Rel$, but the converse does not seem to be the case: there is no relation $Q\colon A\to A$ such that $R\lift Q$ is the transpose of $R$.

Of course not all dagger compact categories arise in this way. Take, for example, the dagger compact category of finite dimensional Hilbert spaces.

Let $z\colon \mathbb{C}\to \mathbb{C}$ in $\FDHilb$. We know that $z$ corresponds to a complex number, and we also know that $z^\dagger$ is given by the complex conjugate. Now suppose that $\FDHilb\cong \BB_0$ for some composition-closed bicategory $\BB$ with ambi-adjoints. We know that any non-zero $z\colon \mathbb{C}\to \mathbb{C}$ has an ambi-adjoint given by $\overline{z}$, but we also know that $z$ has an inverse. If $z^{-1}$ is the inverse of $z$ then $z^{-1}$ is ambi-adjoint to $z$ and thus, by uniqueness of adjoints $z^{-1}\cong \overline{z}$ which is a contradiction, since complex conjugates do not give inverses.

Despite this, it is possible that there is a mutual generalisation of composition-closed compact-closed bicategories and dagger compact categories which is enough to give trace and cotrace structures. Just as dagger compact categories have `formal adjoints', it may well be possible to introduce the concept of compact-closed bicategories equipped with a `formal' lift and extension structure, that is compatible with the dualising pseudofunctor. A further discussion can be found in the afterword.
\begin{remark}
    In \autoref{section:CompClosed} there was some discussion of how in the bicategories $\Rel$ and --- if $\VV$ is star-autonomous -- $\VV$-Prof, the composition-closed structure behaves in the same way as a star-autonomous category. To be explicit, let us define a composition-autonomous bicategory $\BB$ to be a composition-closed bicategory such that for every $A,B\in \BB$ there is a functor
    \[
    (-)^\vee\colon \BB(A,B)\to \BB(B,A)    
    \]
    which is its own inverse, and for all $f\colon B\to C$, $g\colon A\to C$ and $h\colon B\to D$, the two closed structures are given by
    \[f\lift g\cong (g^\vee\circ f)^\vee
    \text{ and } h\extend f\cong (f\circ h^\vee)^\vee.
    \]
    Examples of such bicategories include star-autonomous categories, thought of as one-object bicategories, as well as $\Rel$ and $\VV$-Prof when $\VV$ is star-autonomous. 

    Note that $f\lift \id^\vee \cong f^\vee\cong \id^\vee \extend f$. For $\Rel$ and $\VV$-Prof it's also the case that $f^{*\vee}\cong f^{\vee *}$ for all $f$. This suggests something resembling a dagger compact structure. But note that $(-)^\vee$ does not extend to a pseudofunctor: $\id^\vee$ is not the identity in either of our examples. Additionally, $(f\tensor g)^\vee$ is not necessarily isomorphic to $f^\vee \tensor g^\vee$ -- although in the case of $\VV$-Prof, it is if $\VV$ is compact-closed. It seems then that whilst $\Rel$ and $\VV$-Prof don't collapse to dagger compact categories, they are somehow closer to being dagger compact than some of our other examples.
\end{remark}
\section{Properties of the Trace and Cotrace}
\begin{figure}[H]
    \[
    \myinput{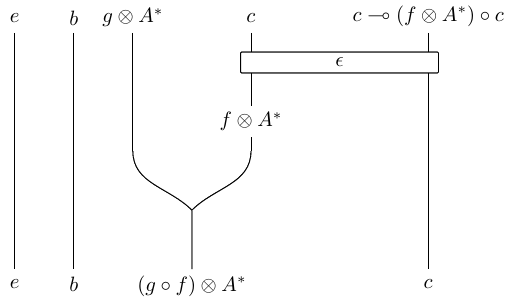}
    \]
    \caption{}
    \label{ch6.4}
    \end{figure}

\begin{figure}
\[
\myinput{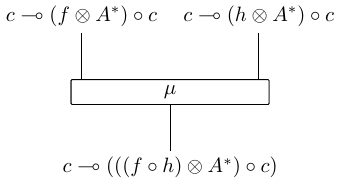}    
\]
\caption{}
\label{ch6.5}
\end{figure}

\begin{figure}
\[
    \myinput{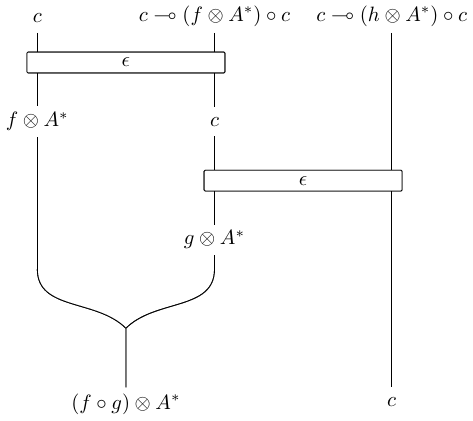}    
\]
\caption{}
\label{ch6.6}
\end{figure}

\begin{figure}
\[
\myinput{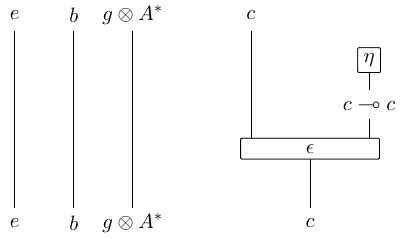}    
\]
\caption{}
\label{ch6.7}
\end{figure}

\begin{figure}
\[  
\myinput{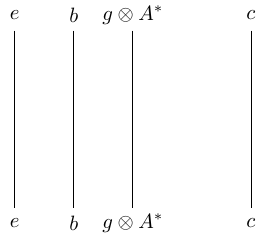}
\]
\caption{}
\label{ch6.8}
\end{figure}

\begin{figure}
    \[\myinput{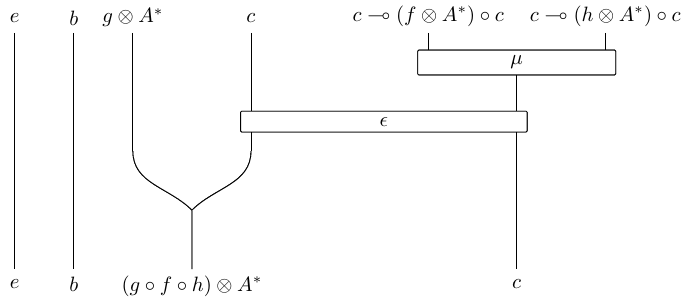}\]
    \caption{}
    \label{ch6.9}
    \end{figure}

\begin{figure}
    \[\myinput{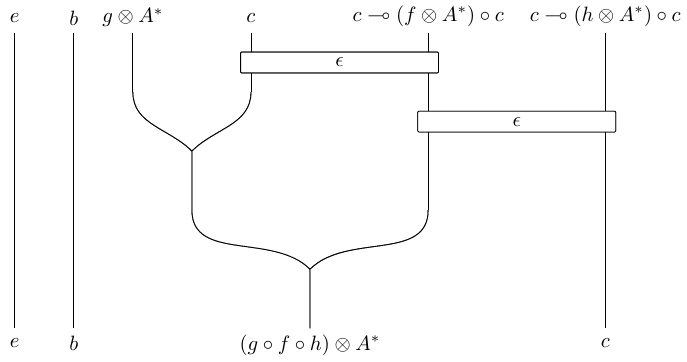}\]
    \caption{}
    \label{ch6.10}
    \end{figure}

In the previous section we shed some light on why it is the case that bicategories have traces and cotraces, whereas in linear algebra there is a single trace that works analogously to both. In this section we see what trace-like properties the cotrace and trace possess. For compact-closed categories the trace is known to have certain properties already -- such as dual invariance, cyclicity, and tensor preservation -- but with the inclusion of a composition-closed structure we see further trace-like properties begin to emerge. To be specific the trace has the following properties,
\begin{itemize}
    \item the trace of a scalar is the scalar itself (this is a special case of the following);
    \item the trace is linear, $
    \Tr(\sum_{i=1}^n \lambda_i f_i)= \sum_{i=1}^n \lambda_i \Tr(f_i);$
    \item the trace is dual invariant,
    \(\Tr(E^*)=\Tr(E)\);
    \item the trace preserves adjoints,
    \(
        \Tr(E^\dagger)=\Tr(E)^{\dagger}
    \), since the adjoint of a scalar is its complex conjugate;
    \item the trace is cyclic, \(\Tr(K\circ H)=\Tr(H\circ K)\);
    \item the trace preserves tensors, $\Tr(F\tensor E)= \Tr(F)\cdot \Tr(E).$
\end{itemize}
We will show that the trace and cotrace share many analogous properties. We also try to provide some intuition at to why the two traces always seem to give dual results. For the remainder of this section we will be working in the context of a composition-closed compact-closed bicategory $\BB$.
\begin{proposition}
    \label{prop:TraceOfScalar}
    Let $s\colon I\to I$ be a scalar in $\BB$. Then we have natural isomorphisms
    \[ 
    \rTr(s)\cong s\cong \dTr(s).
    \]
\end{proposition}
\begin{proof}
    To see that the trace preserves scalars, recall that $I$ is self-dual with the left unitor giving evaluation and coevaluation. Then both the trace of $s$ and the scalar promotion of $s$ are given by
    \[
        l\circ (s\tensor I)\circ l^\bullet.  
    \]
    But note that there is a natural isomorphism given by the diagram below.
    \[
        \myinput{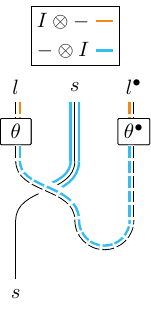}
    \]
    So both $\Spr_I$ and $\rTr_I$ are naturally isomorphic to the identity functor. Thus, by the uniqueness of adjoints both $\dTr_I$ and $\rTr_I$ are naturally isomorphic to the identity functor.
\end{proof}
We next prove that the trace and cotrace are both, in some sense, linear. There are two ways to view this. In our analogy between category theory and linear algebra we let weighted limits and colimits play the role of weighted sums.
\begin{lemma}
    The enriched trace preserves weighted colimits.
\end{lemma}
\begin{proof}
    This follows from the fact that the enriched cotrace functor is a right adjoint. See Kelly's~\cite[sec.~3.2]{kelly1982basic} monograph.
\end{proof}
\begin{lemma}
    The enriched cotrace preserves weighted limits.
\end{lemma}
\begin{proof}
    This follows from the fact that the enriched trace functor is a left adjoint.  See Kelly's~\cite[sec.~3.2]{kelly1982basic} monograph.
\end{proof}
In particular this means that the cotrace and trace functors preserve copowers and powers respectively. Thus, the trace and cotrace preserve scalar multiplication in the following sense.
\begin{corollary}
    For every scalar $s$ and every endo-1-cell $f\colon A\to A$ there is a natural isomorphism
    \[
        \tr(\Spr(s)\circ f)\cong s\circ \tr(f).
    \]
\end{corollary}
\begin{proof}
    This follows from the fact that the copower in $\BB(A,A)$ is given by $\Spr(s)\circ (-)$.
\end{proof}
Note that this makes the trace functor linear in the sense of $\BB(I,I)$-representations.
\begin{corollary}
    For every scalar $s$ and every endo-1-cell $f\colon A\to A$ there is a natural isomorphism
    \[
        \ctr(\Spr(s)\lift f)\cong s\lift \ctr(f).
    \]
\end{corollary}
\begin{proof}
    This follows from the fact that the power in $\BB(A,A)$ is given by $\Spr(s)\lift (-)$.
\end{proof}
Thus, both the cotrace and the trace preserve scalar multiplication of a kind. There is also a different notion of linear functor, \autoref{def:linearFunctor}. Note that $\BB(A,A)$ is a $\BB(A,A)$-representation by precomposition, so the domain and codomain of the trace functor can both be viewed as $\BB(I,I)$-representations.
\begin{lemma}
    \label{lemma:TracesLinear}
    For every $f,g\colon A\to A$ there is a 2-cell \[
        m\colon \tr(g)\circ \ctr(f)\Rightarrow \tr(g\circ f),
    \] natural in $f$ and $g$, such that $
    (\tr,\ctr,m)\colon \BB(A,A)\to \BB(I,I)$ 
    defines a linear functor.
\end{lemma}
\begin{proof}
    Recall that $\ctr(f)\cong (\name{\id}\lift \name{f})\cong(\coev\lift (f\tensor A^*)\circ \coev)$ and so if we expand the composite $\tr(g)\circ \ctr(f)$ we have
    \[
        \ev\circ b \circ (g\tensor A^*)\circ \coev \circ (\coev\lift (f\tensor A^*)\circ \coev).
    \]
    In the interest of saving space, let $c\coloneqq \coev$ and $e\coloneqq \ev$. Then we define $m$ to be the 2-cell in \autoref{ch6.4}, where $\epsilon$ is the counit for the composition-lift adjunction.

    Before we show that this adheres to the axioms of a linear functor, it will be useful to have an explicit description of the natural transformations that make the cotrace into a lax monoidal functor. They are defined as follows. Firstly, we have the unit natural transformation which is given by the unit for the composition-lift adjunction:
    \[
        \id \xrightarrow{\eta} (c\lift c).
    \]
    Secondly, we have the multiplication transformation of \autoref{ch6.5} which is defined as the adjunct of the \autoref{ch6.6}.
    
    In order to show that $m$ adheres to the unit axiom for linear functors we must show that \autoref{ch6.7} is equal to \autoref{ch6.8}, but this follows from the zigzag identities for $\eta$ and $\epsilon$.

    Now we must show that the associativity axiom holds. This means we must show that, for a given $h\colon A\to A$, \autoref{ch6.9} is equal to \autoref{ch6.10}. But this follows from the definition of $\mu$ and the associativity of the compositor for the $-\tensor A^*$ functor. Thus, $(\tr,\ctr,m)$ defines a linear functor.
\end{proof}
This is a form of linearity, and it has some interesting consequences later in the chapter, but it is somewhat at odds with the case for the trace of linear endomorphisms. Since the trace for linear endomorphisms plays the role of both trace and cotrace, the analogous result in $\FDHilb$ would be that there was a homomorphism of representations \[(\Tr,\Tr)\colon \FDHilb(A,A)\to \FDHilb(\mathbb{C},\mathbb{C}),\]
where $\FDHilb(A,A)$ acts on itself via composition. In other words, given linear endomorphisms $g$ and $f$ the analogous result would be that $\Tr(g\circ f)=\Tr(g)\circ \Tr(f)$, which is not the case.

This goes to show that we cannot capture the trace of finite dimensional Hilbert spaces using composition-closed compact-closed bicategories. What is true of finite dimensional Hilbert spaces, however, is that $\FDHilb(A,A)$ is acted on by $\FDHilb(A,A)$ via \emph{addition}. We also know that $\Tr(F+G)=\Tr(F)+\Tr(G)$ which means that the trace map
\[
(\Tr,\Tr)\colon \FDHilb(A,A)\to \FDHilb(\mathbb{C},\mathbb{C})
\]
is a homomorphism of representations, but where the action is given by addition and not composition.

The next property that we want to explore is the extent to which traces preserve duals and adjoints. The case for duals is quite straightforward: both the trace and the cotrace are dual invariant. In the case of compact-closed \emph{categories}, the fact that the trace is dual invariant follows straightforwardly from the string diagram language for compact-closed categories, where strings denote \emph{objects} and beads denote 2-cells. Firstly, writing the evaluation map and coevaluation map as cups and caps, the dual of a map $f\colon A\to A$ is given by the following 2-cell,
\[
    \myinput{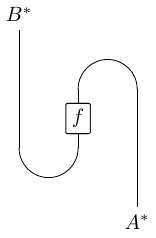}
\]
and so, by the zigzag identities have the following equalities of 2-cells.
\[
 \myinput{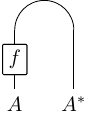}\qquad=\qquad\myinput{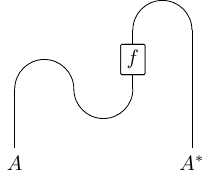}\qquad=\qquad\myinput{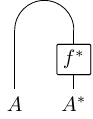}
\]
Now, by, symmetry, we can also see that the following 2-cells are equal,
\[
    \myinput{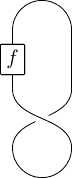}\qquad=\qquad\myinput{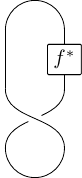}\qquad=\qquad\myinput{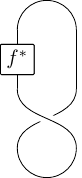}
\]
and so we know that the trace of $f$ is equal to the trace of its dual. Proving that this is the case for the trace in a compact-closed bicategory is similar, but with natural isomorphisms instead of equalities. In what follows we return to using string diagrams for bicategories, so strings are 1-cells and beads are 2-cells.
\begin{proposition}
    \label{prop:coevSwap}
In a compact-closed bicategory, given any $f\colon A\to A$ there is a natural isomorphism, $\zeta$, in the diagram below.
\[
    \myinput{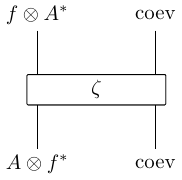}
\]
\end{proposition}
\begin{proof}
    The natural isomorphism is given by \autoref{ch6.11}. Naturality follows from the fact that caps and braids define modifications.
    \begin{figure}[p]
    \[
    \myinput{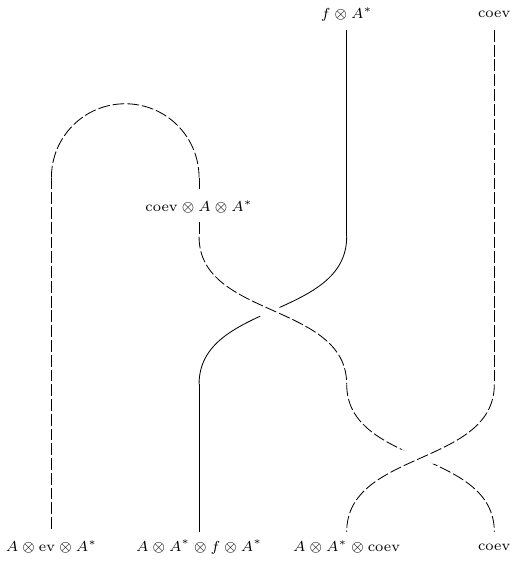}    
    \]
    \caption{}
    \label{ch6.11}
    \end{figure}
\end{proof}
\begin{lemma}
    For every 1-cell $f\colon A\to A$ there is a natural isomorphism
    \[
    \tr(f)\cong \tr(f^*).
    \]
\end{lemma}
\begin{proof}
    The transformation is given by \autoref{ch6.12}, where the natural isomorphisms $\xi$ are given by the fact that the coevaluation map for $A^*$ is $b\circ \coev_{A}$ -- see the proof of \autoref{prop:dualBothWays} -- and similarly for the evaluation map.
    \begin{figure}[p]
    \[
        \myinput{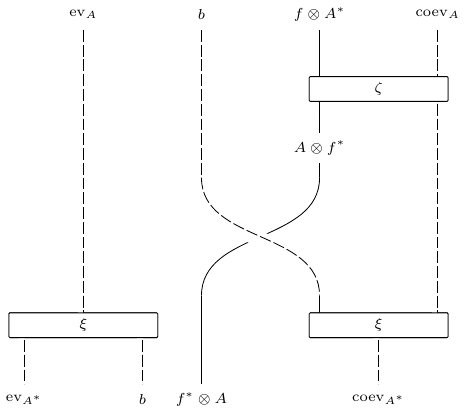}    
    \]
    \caption{}
    \label{ch6.12}
    \end{figure}
\end{proof}
The case for the cotrace is comparatively straightforward.
\begin{lemma}
    For every $f\colon A\to A$ there is a natural transformation
    \[
    \ctr(f)\cong \ctr(f^*).
    \]
\end{lemma} 
\begin{proof}
We know that there is a string of natural isomorphisms
    \begin{align*}
        \ctr(f)\cong \BBB(\id, f)\cong \BBB(\id^*, f^*)\cong \BBB(\id, f^*)\cong \ctr(f^*)
    \end{align*}
and so their composite gives the natural transformation above.
\end{proof}
We now move on to adjoints. For maps between finite dimensional Hilbert spaces, duals and linear adjoints are essentially the same thing. If we choose a basis for $A$ and $B$ then the adjoint $f^\dagger\colon B\to A$ of $f\colon A\to B$ is given by the composite
\[
        B\xrightarrow{\sim} B^*\xrightarrow{f^*} A^* \xrightarrow{\sim} A
\]
and as a result the dual $f^*$ of $f\colon A\to B$ can also be written in terms of its adjoint.

Recall that if a 1-cell $f$ has a right adjoint then it is given by $f\lift \id$. The converse isn't necessarily true -- it is not the case that $f\lift \id$ is always adjoint to $f$ -- but nonetheless lifting the identity through $f$ gives an adjoint in some slightly more general sense.

Neither the cotrace nor the trace will, in general, preserve adjoints. But as a pair we can measure to what extent they fail to preserve adjunctions. 
\begin{proposition}
    In any composition-closed compact-closed bicategory, for all 1-cells $f\colon A\to B$ and all scalars $s\colon I \to I$ there is a natural isomorphism
    \[
        \ctr(f\lift \cspr(s)) \cong \tr(f)\lift s.
    \]
\end{proposition}
\begin{proof}
By using the fact that the cotrace at $I$ is isomorphic to the identity, and the fact that the trace and cospread are $\BB(I,I)$-adjoint we have a string of isomorphisms
\begin{align*}
    \tr(f)\lift s &\cong \ctr(\tr(f)\lift s)\\
    &\cong \BBB(I,I)(\tr(f), s)\\
    &\cong \BBB(A,A)(f, \cspr(s))\\
    &\cong \ctr(f\lift \cspr(s))
\end{align*}
that proves the above claim.
\end{proof}
If we wanted to make this look little more symmetric, since the spread of a scalar at $I$ is the scalar itself we can express the above proposition as
\[
    \ctr(f\lift \cspr(s)) \cong \tr(f)\lift \spr(s).
\] 
This gives us a direct relationship between the cotrace and the trace and gives some sort of idea of how they are dual. Note that in the particular case where we take $s$ to be $\id$ we have the following corollary.
\begin{corollary}
    In any composition-closed compact-closed bicategory, for all $f\colon A\to B$ there is a natural isomorphism
    \[
        \ctr(f\lift \cspr(\id)) \cong \tr(f)\lift \id.
    \]
\end{corollary}
Thus, the extent to which the cotrace-trace pairing fails to preserve adjoints is, in some sense, the extent to which the cospread fails to preserve the identity. If it is the case that $\ev\circ b$ is left adjoint to the coevaluation map, then the spread and cospread are isomorphic, but the spread is strongly monoidal. Thus, $\cspr(\id_I)\cong \id_A$ and so in that particular case the trace and cotrace have a duality property in the sense that
\[
\ctr(f^\dagger)\cong \tr(f)^\dagger 
\]
where $(-)^\dagger$ denotes the right adjoint. If, in addition, $\ev\circ b$ is right adjoint to the coevaluation map, then the cotrace and trace are isomorphic and so the unique trace preserves right adjoints. 

Recall that in \autoref{section:daggerComp} we pointed out that two of our examples are `composition-autonomous' in the sense that the closed structure is given in terms of functors
\[
  (-)^\vee\colon \BB(A,B)\to \BB(B,A)  
\]
and composition. If our bicategory is composition-autonomous then a different kind of duality is preserved.
\begin{proposition}
    Suppose that $\BB$ is a composition-autonomous compact-closed category such that there is a natural isomorphism in the following diagram.
    \[
    \begin{tikzcd}
        \BB(A,B)
        \arrow[r, "(-)^\vee", ""{swap, name=top}]
        \arrow[d, swap, "\name{(-)}"]
            &\BB(B,A)
            \arrow[d, "\unname{(-)}"]\\
        \BB(I,B\tensor A^*)
        \arrow[r, "(-)^\vee"{swap}, ""{name=bottom}]
            &\BB(B\tensor A^*, I)  
        \arrow[phantom, "\cong", from=top, to=bottom]
    \end{tikzcd}    
    \]
    Then, for all $f\colon A\to A$, there is a natural isomorphism
    \[
    \ctr(f^\vee)\cong \tr(f)^\vee.
    \]
\end{proposition}
\begin{proof}
This follows from the string of isomorphisms
        \[\name{\id}\lift \name{f^\vee}\cong (\name{(f^\vee)}^\vee\circ \name{\id})^\vee\
    \cong \unname{f} \circ \name{\id} \cong \Tr(f)^\vee\]
which follow from the supposition above and the definition of a composition-autonomous category.
\end{proof}
\begin{example}
    In the bicategory $\Rel$ this corresponds to the fact that, given a relation $R\colon A\to A$, there is an $a\in A$ such that $aR^ca$ if and only if $R$ is not reflexive.
\end{example}
\begin{example}
    In the bicategory $\VV$-Prof, for star-autonomous $\VV$, this corresponds to the fact that, given a profunctor $P\colon \AA\profto \AA$ we have the following natural isomorphism:
    \[
    \endint_{A\in \AA} \left(P(A,A)\right)^* \cong \left(\endint^{A\in \AA} P(A,A)\right)^*.  
    \]
\end{example}
\begin{remark}
    Note that in these particular cases this means we can define our scalar enrichment in terms of the trace, rather than the cotrace.
\end{remark}
Next, we want to investigate the extent to which the trace and cotrace have the cyclicity property. Like with dual invariance, the trace for compact-closed categories is cyclic, and the argument for compact-closed bicategories is analogous. In the case of compact-closed categories, for every $f\colon A\to A$ the morphisms given by the diagrams below are equal,
\[
    \myinput{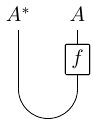}  \quad=\quad  
    \myinput{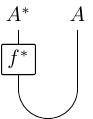}    
\]
similar to the case for coevaluation. This means that we have a sequence of equalities given below.
\[
            \myinput{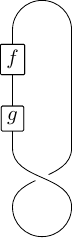} \quad=\quad \myinput{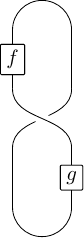}\quad =\quad \myinput{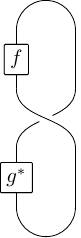}\quad =\quad \myinput{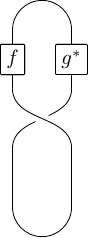}\quad = \quad \myinput{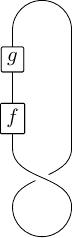}
\]
\begin{proposition}
    In a compact-closed bicategory, given any $f\colon A\to A$ there is a natural isomorphism, $\zeta$, in the diagram below.
\[
    \myinput{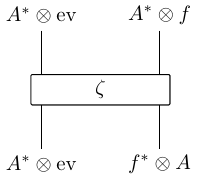}
\]
\end{proposition}
\begin{proof}
    The construction of this isomorphism is analogous to \autoref{prop:coevSwap}.
\end{proof}

\begin{lemma}
    For every pair of 1-cells $f\colon A\to B$ and $g\colon B\to A$ there is a natural isomorphism
    \[
    \tr(f\circ g)    \cong \tr(g\circ f).
    \]
\end{lemma}
\begin{proof}
    The isomorphism is given by \autoref{ch6.13}. Invertibility and naturality follow from the fact that every 2-cell in the composite is invertible and natural.
    \begin{figure}[p]
    \[
    \myinput{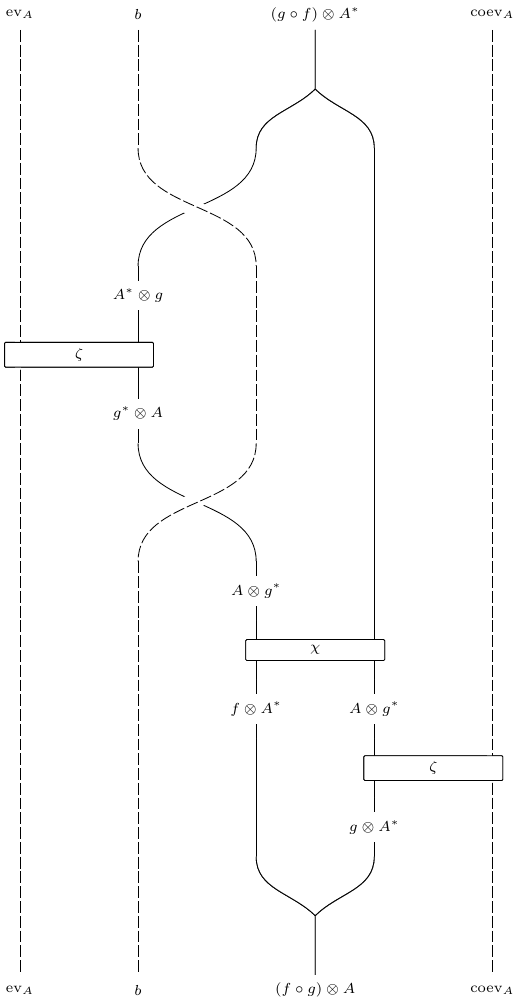}
    \]
    \caption{}
    \label{ch6.13}
    \end{figure}
\end{proof}
Cyclicity does not necessarily hold for the cotrace. Consider, for example the bicategory $\Span(\Set)$. Let $A$ be a set and let spans $S\colon A\to *$ and $T\colon *\to A$ be given by the following diagrams,
\[
\begin{tikzcd}
&
    A
    \arrow[rd,"!"]
    \arrow[ld, "\id"{swap}]\\
A
    && *
\end{tikzcd}
\text{ and }
\begin{tikzcd}
    &
        A
        \arrow[rd,"\id"]
        \arrow[ld, "!"{swap}]\\
    *
        && A
    \end{tikzcd}
\] 
where $*$ denotes the one-object set. Then $T \circ S$ is isomorphic to the span in the following diagram,
\[
\begin{tikzcd}
    &
    A\times A
    \arrow[rd,"p_2"]
    \arrow[ld, "p_1"{swap}]\\
A
    && A
\end{tikzcd}    
\]
where $p_1$ and $p_2$ denote the projection functions. Thus, we know that \[\ctr(T\circ S)\cong\Span(\Set)(\id_A, T\circ S)=\{f\colon A\to A\times A\mid p_1\circ f=\id= p_2\circ f\},\]
which is the singleton $\{\Delta\}$ containing the diagonal map. On the other hand, $S\circ T$ is isomorphic to the span in the following diagram.
\[
\begin{tikzcd}
    &
    A
    \arrow[rd,"!"]
    \arrow[ld, "!"{swap}]\\
*
    && *
\end{tikzcd}    
\]
Thus, we know that 
\[\ctr(T\circ S)\cong\Span(\Set)(\id_A, T\circ S)=\{f\colon *\to A\mid \mathrm{!}\circ f=\mathrm{!}\},\]
which is in bijection with the set $A$. And so, as long as $A$ has cardinality of more than $1$, the two cotraces do not agree.

\begin{figure}[p]
    \[\myinput{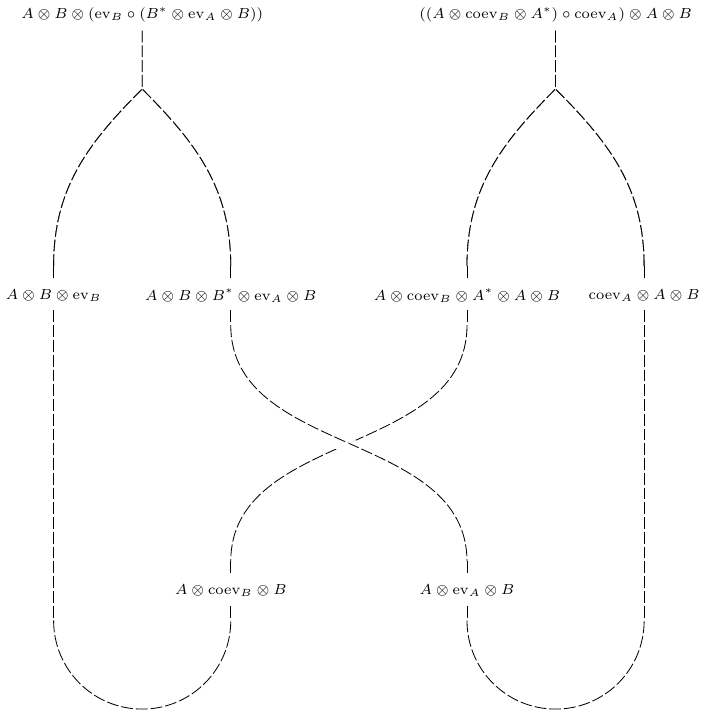}\]
    \caption{}
    \label{diag:dualProp}
\end{figure}

\begin{figure}[p]
    \[
        \myinput{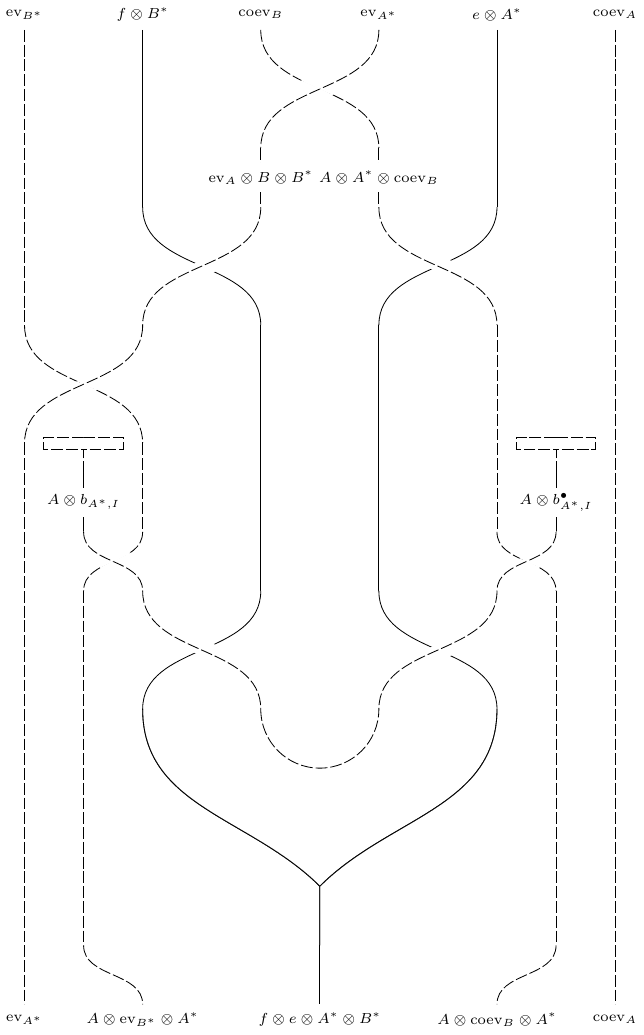}
    \]
    \caption{}
    \label{diag:tensorProp}
    \end{figure}

However, the cotrace has a property that resembles cyclicity.
\begin{lemma}
    For every pair of 1-cells $f\colon A\to B$ and $g\colon B\to A$ there is a natural isomorphism
    \[
    \ctr(f\lift g)    \cong \ctr(g\extend f).
    \]\end{lemma}
\begin{proof}
    This follows from the isomorphisms \[
        \BBB(\id, f\lift g)\cong \BBB(f,g)\cong \BBB(g\extend f),\]
    which follow from the definition of pseudoinverse.
\end{proof}
In particular, if $f\colon A\to B$ has left adjoint $f_\dagger$ and right adjoint $f^\dagger$ then for every $g\colon B\to A$, there is a natural isomorphism \[\ctr(f^\dagger\circ g)\cong \ctr(g\circ f_\dagger).\]
\begin{remark}
Interestingly, the fact that the trace has a cyclicity property gives it the structure of a shadow functor. This means it can be used to define the trace of a 2-cell in the sense of Ponto~\cite{pontoOriginal}. The cotrace does not have the structure of a shadow functor, but it does have enough of a cyclicity property that it can still be used to define a Ponto trace.
\end{remark}
Note that one consequence of cyclicity for the linear trace is that it is conjugate independent. By this we mean that, for any $g\colon B\to A$ that has an inverse
\[
    \Tr(g^{-1}fg)\cong \Tr(f).
\]
The same holds true for the cotrace.
\begin{proposition}
    In any composition-closed compact-closed bicategory, for all $f\colon A\to B$ and all $g\colon B\to A$ with pseudoinverse $g^\bullet$, there is a natural isomorphism
    \[
        \ctr(g^\bullet \circ f \circ g)\cong \ctr(f).
    \]
\end{proposition}
\begin{proof}
    This follows from the natural isomorphisms \[\BBB(\id, g^\bullet\circ f \circ g)\cong \BBB(g\circ g^\bullet, f)\cong \BBB(\id, f).\]
\end{proof}
The final trace-like property we will investigate is the preservation of the tensor product. This line of inquiry was originally suggested to the author by Cranch~\cite{cranch}. Like the above properties, tensor preservation holds strongly for the trace, but somewhat laxly for the cotrace.
\begin{lemma}
    If $A^*$ is dual to $A$, and $B^*$ is dual to $B$, then $B^*\tensor A^*$ is dual to $A\tensor B$ with coevaluation given by 
    \[
        I \xrightarrow{\coev_A} A\tensor A^* \xrightarrow{A\tensor \coev_B\tensor A^*} A\tensor B\tensor B^*\tensor A^*;
    \]
    and evaluation given by
    \[
    B^*\tensor A^* \tensor A \tensor B \xrightarrow{B^*\tensor \ev_A \tensor B} B^*\tensor B \xrightarrow{\ev_B} I.   
    \]
\end{lemma}
\begin{proof}
    To prove this we must prove that there are cups and caps giving an adjoint equivalence between the zigzag 1-cells. The first cup is given by \autoref{diag:dualProp}. The second cup and the caps are defined analogously. The fact that the cups and caps adhere to the zigzag identities follows from \autoref{pseudoinvertible}.
\end{proof}
\begin{lemma}
For every pair of 1-cells $e\colon A\to A$ and $f\colon B\to B$ there is a natural isomorphism
\[
\tr(f)\circ \tr(e) \cong \tr(f\tensor e).   
\]
\end{lemma}
\begin{proof}
    The natural isomorphism is given by \autoref{diag:tensorProp}.
\end{proof}
\begin{lemma}
    For every pair of 1-cells $e\colon A\to A$ and $f\colon B\to B$ there is a natural transformation
\[
\ctr(f)\circ \ctr(e) \rightarrow \ctr(f\tensor e).   
\]
\end{lemma}
\begin{proof}
    This follows from the fact that $\tensor$ gives an enriched pseudofunctor, and so there is a natural map
    \begin{align*}
        \underline{\BB}(B,B)(\id_B,f)\circ \underline{\BB}(A,A)(\id_A, e)&\to \underline{\BB}(A\tensor B,A\tensor B)(\id_A\tensor \id_B, e\tensor f) \\&\cong \underline{\BB}(A\tensor B,A\tensor B)(\id_{A\tensor B}, e\tensor f).
    \end{align*}
\end{proof}\newpage
\section{Dimension and Codimension}
In this final section we use the trace and cotrace to give the definition of the dimension and codimension of an object in a composition-closed compact-closed bicategory. In linear algebra we can define the dimension of a Hilbert space $A$ as the trace of the identity endomorphism on $A$: \[\Dim(A)=\Tr(\id_A).\] 
Taking a similar approach with the trace and cotrace we discover that the dimension and codimension have a particularly structured relationship. It was Willerton~\cite{willertonTalk} who first observed that the 2-dimension of $A$ has a monoid structure and that it acts on the dimension of $A$. Here we give an enriched version of that observation. In all of the following results, $\BB$ is a compact-closed composition-closed bicategory.
\begin{definition}
    Let $A$ be an object in a composition-closed compact-closed bicategory $\BB$. Then we define the \textdef{dimension} of $A$ and the \textdef{codimension} of $A$ to be given by
    \[
        \rDim(A)\coloneqq \rTr(\id_A) \text{ and }\dDim(A)\coloneqq \dTr(\id_A).
    \]
\end{definition}
\begin{proposition}
    For every object $A$, the codimension of $A$, $\BBB(A,A)(\id_A,\id_A)$, has a monoid structure in the symmetric monoidal category of scalars, with multiplication given by enriched composition and unit given by the underlying identity 2-cell
    \[
        \iota\colon \id_I \to \underline{\BB}(A,A)(\id_A,\id_A).
    \]
\end{proposition}
\begin{proof}
    Associativity and unitality hold since they hold for the enriched composition.
\end{proof}
\begin{proposition}
    \label{cotraceMonoid}
    For every object $A$, the 1-cell $\id_A$ has a monoid structure in the category $\BB(A,A)$ where the unit is given by the identity, and composition is given by either of the left or right unitor. This induces a monoid structure on the codimension of $A$, since the cotrace is a lax monoidal functor.
\end{proposition}
\begin{proof}
Unitality follows from the properties of identity 2-cells, associativity follows from the triangle equation in the definition of a bicategory.
\end{proof}
 A classical result in category theory is the fact that any pair of adjunct functors
 \[
F\colon \AA\leftrightarrows \BB \cocolon G    
\]
induces a monad structure for the functor $G\circ F$. In other words, $G\circ F$ comes equipped with the structure of a monoid object in the category of endofunctors on $\AA$. Less well known is the concept of a codensity monad, first introduced by Kock~\cite[p.~3]{Kock}. The codensity monad is the monad structure that comes from taking the right Kan extension of $F\colon \AA\to \BB$ along itself. Taking the right extension, or right lift, in a bicategory also gives a monoid object in the category of endomorphisms.
\begin{proposition}
    \label{monoidDiag2}
    For every object $A$, the codimension of $A$, $\name{\id_A}\lift \name{\id_A}$, has a monoid structure in the symmetric monoidal category of scalars, where the unit is given by the adjunct of the identity 2-cell
    \[
    \id_A \rightarrow{\id_A}
    \]
    and composition is given by the adjunct of the evaluation map
    \[
    \name{\id}\circ (\name{\id}\lift \name{\id}) \circ (\name{\id}\lift \name{\id}) \xrightarrow{\epsilon} \name{\id}\circ (\name{\id}\lift \name{\id})\xrightarrow{\epsilon}\name{\id}. 
    \]
\end{proposition}
\begin{proof}
    See, for example, Leinster's~\cite[sec.~5]{LeinsterCodense} exposition, for an overview of how taking a right Kan extension of a functor along itself gives rise to a monad structure. The same argument holds for right extensions in a bicategory and, by duality, for right lifts as well.
\end{proof}
\begin{lemma}
    The three monoids given in the propositions above are all isomorphic.
\end{lemma}
\begin{proof}
To see that the first and second monoids are isomorphic note that the isomorphism given in \autoref{prop:altEnrichment} gives an isomorphism between $\underline{\BB}(\id,\id)$ and $\ctr(\id_A)$. The fact that this isomorphism is an isomorphism of monoids simply follows from the definition of the enriched composition.

To see that the first and third monoids are isomorphic, firstly note that, by definition, \[\ctr(\id_A)=\name{\id_A}\lift \name{\id_A}.\] The lax multiplication for the cotrace is given by the adjunct of the following composite:
\begin{align*}
    \name{\id} \circ (\name{\id} \lift \name{g}) \circ (\name{\id} \lift \name{f})   & \xrightarrow{\epsilon} \name{g}\circ (\name{\id}\lift \name{f})\\
    &\xrightarrow{\sim}\hom{A,g}\circ \name{\id}\circ (\name{\id}\lift \name{f})\\
    &\xrightarrow{\epsilon}\hom{A,g}\circ \name{f}\\
    &\xrightarrow{\sim}\name{g\circ f}.
\end{align*}
Then the multiplication for the monoid defined in \autoref{cotraceMonoid} is given by adjunct of the top right map in \autoref{monoidDiag}. But the diagram commutes by naturality and the bottom left map defines the multiplication for the monoid in \autoref{monoidDiag2}. Thus, the two monoids are isomorphic. 
\end{proof}
\begin{figure}[hb]
\[\begin{tikzcd}
    \name{\id} \circ (\name{\id} \lift \name{\id}) \circ (\name{\id}\circ \name{\id})
    \arrow[r,"\epsilon"]
        & \name{\id}\circ (\name{\id}\lift \name{\id})
        \arrow[r, "\comp{\rho}^{-1}"]
        \arrow[rdd, "\eta"{swap}]
            &\hom{A,\id}\circ \name{\id}\circ \name{\id}\lift \name{\id}
            \arrow[d, "\epsilon"]\\
        &   
            &\hom{A,\id}\circ\name{\id}
            \arrow[d, "\comp{\rho}"]\\
        &
            & \name{\id}
\end{tikzcd}
\]
\caption{}
\label{monoidDiag}
\end{figure}
\begin{lemma}
    For every $f\colon A\to A$, the trace of $f$, $\tr(f)$, is a right $\dDim(A)$-module in the symmetric monoidal category of scalars.
\end{lemma}
\begin{proof}
    This follows from the fact that the identity is a monoid object in $\BB(A,A)$, every $f\colon A\to A$ is an $\id$-module, \autoref{lemma:linearPreservesModules} which states that linear functors preserve modules, and \autoref{lemma:TracesLinear} which states that the trace and cotrace form a linear functor.
\end{proof}
\begin{corollary}
    The dimension of an object $A$ is a right $\dDim(A)$-module in the symmetric monoidal category of scalars.
\end{corollary}
Once again, the majority of the following examples are due to Willerton~\cite{willertonTalk}.
\begin{example}
    In $\Rel$ the dimension and codimension are both just the point and so the monoid and module structure are both trivial. If we take the cotrace of the empty relation then this has a $*$-module structure where the action map is the unique map
    \[
    \varnothing \times *\cong \varnothing \to \varnothing    
    \]
\end{example}
\begin{example}
    In $\Bim_R$ the dimension and codimension of an algebra are given by the abelianisation of the algebra and the centre of the algebra. The monoid structure on the centre is given by the algebra multiplication. The action of the centre on the abelianisation is also given by the algebra multiplication.
\end{example}
\begin{example}
    In $\DGBim_R$ the dimension and codimension of an algebra are the Hochschild homology and the Hochschild cohomology of the algebra. The monoid structure on cohomology is the cup product. The action of cohomology on homology is given by the cap product.
\end{example}
\begin{example}
    In $\VV$-Prof the cotrace of a profunctor $P\colon \AA\profto \AA$ is given by taking the end. The codimension of $\AA$, then, is given by the natural transformation object from the identity functor to itself
    \[
    \endint_{A\in \AA} \AA(A,A) =\Nat(\Id,\Id).
    \]
    The monoid structure is given by composition. To understand the dimension, let us focus on the case where $\VV$ is $\Set$. The coend of the Hom profunctor is given by the quotient of all endomorphisms
    \[
        \endint^{A\in \AA} \AA(A,A)=\left( \bigsqcup_{A\in A} \AA(A,A)\right)/\sim
    \]
    where $(f\colon A\to A)\sim (g\colon A'\to A')$ if there exists some $\phi\colon A\to A'$ or some $\psi\colon A'\to A$ such that one of the following squares commutes.
    \[
        \begin{tikzcd}
            A
            \arrow[r, "f"]
            \arrow[d, "\phi"{swap}]
                &A
                \arrow[d, "\phi"]\\
            A'
            \arrow[r, "g"{swap}]
                &A'
        \end{tikzcd}
        \begin{tikzcd}
            A'
            \arrow[r, "g"]
            \arrow[d, "\psi"{swap}]
                &A'
                \arrow[d, "\psi"]\\
            A
            \arrow[r, "f"{swap}]
                &A
        \end{tikzcd}  
    \]
    In the case the $\AA$ is a monoid this is the centre of the monoid, so it is reasonable to think of this as the cocentre of the category. The action of natural transformations on the cocentre of the category is just given by pre- or post-composition: the two are equal by naturality.
    \end{example}
\begin{example}
    In $\Span(\Set)$ the codimension of a set $A$ is the set $\{\id_A\}$. The monoid structure is trivial. The dimension of a set $A$ is given by $A$ itself. The action is also trivial.
\end{example}
\begin{example}
    In $\Path(G)$, for a topological group $G$, the dimension and codimension are both given by the constant loop at $e$, the unit. The monoid and module structures are both given by the concatenation of loops.
\end{example}

\chapter*{Further Work}
\addcontentsline{toc}{chapter}{Further Work}
In this afterword we consider some of the further directions that this work may take. The first and most obvious line of work would be further examples. As mentioned earlier, one candidate would be the bicategory of $2$-Hilbert spaces. Another would be the bicategory of `spaces' and Fourier--Mukai kernels, as described by C\u{a}ld\u{a}raru and Willerton~\cite{willertonMukai}, which Willerton~\cite{willertonTalk} has pointed out should also yield Hochschild homology and cohomology as trace and cotrace.

There are also two questions that this thesis fails to answer in full. Firstly, can we truly formalise the analogy between composition-closed compact-closed bicategories, and dagger compact categories? Secondly, is there some sense in which the cotrace is formally dual to the trace?

We have already given a partial answer to the first question -- we have shown that not every dagger compact category arises as the 2-skeleton of a composition-closed compact-closed bicategory, and not every composition-closed compact-closed bicategory gives rise to a dagger compact category. As such, it seems likely that, were a formal analogy to exist, it would take the form of a mutual generalisation.

There are two potential hints as to what a mutual generalisation might look like. Firstly, $\FDHilb$ failed to be the 2-skeleton of a composition-closed compact-closed bicategory on account of the fact that, in category theory, inverses are a type of adjoint. This is not the case in linear algebra, and it seems likely that any mutual generalisation would be forced to replace lifts with some kind of formal structure like the dagger.

The second hint comes from the fact that there is a linear functor given by the pairing of the trace with the cotrace, and that the trace is linear but with respect to a different action. Composition-closed, compact-closed categories are defined in terms of a closed action of $\BB(I,I)$ on $\BB(A,A)$ using composition. But perhaps we could define some structure involving actions without necessarily enforcing that the actions are given by composition. This might then allow for the formal lifts to be defined as adjoint, or at least related, to the action.

To answer the second question we need a context in which the cotrace and trace can be compared. We have already shown that the trace and cotrace are not dual in the context of the functor category $\hom{\BB(A,A), \BB(I,I)}$, equipped with Day convolution. If they were, they would have to be isomorphic since the cotrace is the monoidal unit.

In some sense the trace and cotrace are orthogonal, rather than dual. The trace arises from the monoidal-closed and monoidal-coclosed structures, which are `horizontal' in the sense of being defined in terms of the tensor product. On the other hand the cotrace is defined in terms of monoidal-closed and composition-closed structures, which are `vertical' in the sense of being defined in terms of composition. Capturing this formally, or at least more formally than we already have done, might prove difficult. That being said, the majority of our examples have stemmed from $\VV$-Prof. Note that if we have a profunctor $P\colon \AA\profto \AA$ then the trace of $P$ is the cotrace of
\[
P^{\op}\colon \AA \times \AA^{\op}\to \VV^{\op}    
\]
which is a dualisation of $P$, but where $P$ is thought of as a functor and not a profunctor. This might suggest that we consider our examples not as composition-closed compact-closed bicategories, but instead as compact-closed proarrow equipments in the sense of Wood~\cite[]{woodproarrows}. This may give another way to compare the trace and cotrace. In the context of bicategories the functor $\underline{\VV\text{-}\Prof}(*,-)$ takes categories to their presheaf categories, and profunctor to cocontinuous functors. But note that $\cat{\VV}(*,-)$ is isomorphic to the identity functor. If we worked instead with enriched proarrow equipments, we might be able to replace this with a double functor $\underline{\VV\text{-}\Prof}(*,-)$ which is the identity double functor. In other words, it would send the 1-cells of $\VV$-Prof to profunctors, rather than functors. 

If this were the case, it might then mean that we would have a way to turn 1-cells in an arbitrary $\BB$ into profunctors, rather than functors. Assuming this double functor preserved the composition-closed and compact-closed structures, this might give us a way to view the trace and cotrace as an end and a coend in all cases. Such a perspective would firstly give an obvious duality between the trace and cotrace, but might also open the door to other categorical analogues. For example, the profunctor associated to an endo-1-cell might be seen to encode the `eigenvalues' of that endo-1-cell. 

\printbibliography
\addcontentsline{toc}{chapter}{Bibliography}
\end{document}